\documentclass[fleqn,leqno]{book}%
\usepackage{amsmath}
\usepackage{graphicx}%
\usepackage{amsfonts}%
\usepackage{amssymb}
\usepackage{hyperref}
\usepackage{hyperref}

\usepackage[noload]{qtree}
\usepackage{tikz}
\usetikzlibrary{shapes,backgrounds}

\newtheorem{theorem}{Theorem}[chapter]

\newtheorem{corollary}[theorem]{Corollary}

\newtheorem{definition}[theorem]{Definition}
\newtheorem{example}[theorem]{Example}

\newtheorem{lemma}[theorem]{Lemma}
\newtheorem{notation}[theorem]{Notation}

\newtheorem{proposition}[theorem]{Proposition}
\newtheorem{remark}[theorem]{Remark}

\newenvironment{proof}[1][Proof]{\textbf{#1.} }{\ \rule{0.5em}{0.5em}}
\newcommand{\bi}{\bigskip}
\newcommand{\sm}{\smallskip}

\newcommand{\wh}{\widehat}

\newcommand{\ee}{\end{equation}}
\newcommand{\eea}{\end{eqnarray}}
\newcommand{\bean}{\begin{eqnarray*}}
\newcommand{\eean}{\end{eqnarray*}}
\newif\ifpctex
\newcommand{\noi}{\noindent}

\newcommand{\vp}{\varphi}
\newcommand{\ve}{\varepsilon}
\newcommand{\wt}{\widetilde}

 \newcommand{\rand}[1]{}

     \newcommand{\Rand}[1]{}

\newcommand{\be}[1]{\Rand{\vspace{0,6cm}\tt #1}\begin{equation}\label{#1}}

\newcommand{\bea}[1]{\Rand{\vspace{0,7cm}\tt
#1\vspace{-0,7cm}}\begin{eqnarray}\label{#1}} \marginparwidth2.5cm

\newcommand{\beL}[1]{\Rand{\vspace{0,6cm}\tt #1}\begin{lemma}\label{#1}}
\newcommand{\beD}[1]{\Rand{\vspace{0,6cm}\tt #1}\begin{definition}\label{#1}}
\newcommand{\beT}[1]{\Rand{\vspace{0,6cm}\tt #1}\begin{theorem}\label{#1}}
\newcommand{\beP}[1]{\Rand{\vspace{0,6cm}\tt #1}\begin{proposition}\label{#1}}
\newcommand{\beC}[1]{\Rand{\vspace{0,6cm}\tt #1}\begin{corollary}\label{#1}}

\newcommand{\vr}{\varrho}
\newcommand{\suml}{\sum\limits}
\newcommand{\intl}{\int\limits}

\newcommand{\supl}{\sup\limits}

\newcommand{\X}{{\bf X}}

\newcommand{\F}{\mathbb{F}}
\newcommand{\K}{\mathbb{K}}
\newcommand{\R}{\mathbb{R}}
\newcommand{\N}{\mathbb{N}}

\newcommand{\I}{\mathbb{I}}
\newcommand{\Z}{\mathbb{Z}}

\newcommand{\Q}{\mathbb{Q}}
\newcommand{\mb}{\mbox}

\newcommand{\CF}{{\mathcal F}}

\newcommand{\CL}{{\mathcal L}}

\newcommand{\CM}{{\mathcal M}}

\newcommand{\CR}{{\mathcal R}}
\newcommand{\CP}{{\mathcal P}}

\newcommand{\jjto}{\;\mathop{\Longrightarrow}_{j \to \infty}\;}

\newcommand{\tto}{{_{\D \Longrightarrow \atop t \to \infty}}}
\newcommand{\ttoo}{{_{\D \longrightarrow \atop t \to \infty}}}

\newcommand{\Nto}{{_{\D \Longrightarrow \atop N \to \infty}}}
\newcommand{\nto}{{_{\D \Longrightarrow \atop n \to \infty}}}

\newcommand{\la}{\longrightarrow}

\newcommand{\D}{\displaystyle}

\usepackage{amsmath}
\usepackage{amsfonts}
\usepackage{amssymb}
\usepackage{graphicx}%

 \definecolor{1}{rgb}{1,0,0}
 \definecolor{2}{rgb}{0,1,0}
 \definecolor{3}{rgb}{0,0,1}
\definecolor{4}{cmyk}{0,.3,1,0}

 \definecolor{1red}   {rgb}{0.80,0.00,0.00}
 \definecolor{1blue}  {rgb}{0.00,0.00,0.80}
 \definecolor{1green} {rgb}{0.00,.70,.30}
 \definecolor{1brown} {rgb}{0.60,0.30,0.00}
 \definecolor{1pink}  {rgb}{0.40,0.00,0.40}
 \definecolor{1cyan}  {rgb}{0.00,0.30,0.30}
 \definecolor{1orange}{rgb}{0.80,0.1,0.40}
 \definecolor{1gray}  {rgb}{0.30,0.30,0.30}

 \def\red{\color{1red}}
 \def\blue{\color{1blue}}
 \def\green{\color{1green}}

\begin{document}


\title{Introductory   Lectures on Stochastic Population Systems}

\author{Donald A. Dawson\\ \\ {\em School of Mathematics and Statistics}\\{\em Carleton University}\\{\em Ottawa, Canada}}

\date{}

\maketitle

\pagenumbering{arabic}
 \tableofcontents



\chapter{Introduction}

These notes are based on lectures and courses given at the University of Erlangen, the Britton Lectures at McMaster University, the 2009 PIMS Summer School in Probability and Carleton University.

Historically, the modelling of biological populations has been an
important stimulus for the development of stochastic processes.
The revolutionary changes in the biological sciences over the past
50 years have created many new challenges and open problems. At
the same time  probabilists have developed new classes of
stochastic processes such as interacting particle systems and
measure-valued processes and made advances  in stochastic analysis
that make possible the modelling and analysis of populations
having complex structures and dynamics. This course will focus on
these developments.  In particular stochastic processes that model
populations distributed in space as well as their genealogies and
interactions will be considered. This will include branching
particle systems, interacting Wright-Fisher diffusions,
Fleming-Viot processes and superprocesses.   Basic methodologies
including martingale problems, diffusion approximations, dual
representations, coupling methods, random measures and particle
representations will be introduced.

\bigskip

\chapter{Stochastic models in biology: a historical overview}

\section{Classical deterministic population dynamics}

We begin our historical review with some basic models from demography, ecology and epidemiology.

The mathematical formulation of the growth  of an age-structured population
was developed  by Euler (1760) (\cite{E-60}).
The  (female) birth rate at time $t$  $\{B(t)\}_{t\geq 0}$ satisfies the renewal equation

\be{}  B(t)=\int_0^t B(t-s)(1-L(s))m(s)ds. \ee

This leads to exponential growth $\qquad  B(t) \sim e^{\alpha t}$ where the {\em Malthusian parameter} $\alpha$  is given by the {\em characteristic equation of demography} (Euler-Lotka equation)
\be{EL} 1=\int_0^\infty  e^{-\alpha s}(1-L(s))m(s)ds, \ee
where  $m(s)$ is the average birth rate for an individual of age $s$ and $L(s) $ is the cumulative distribution function of the
lifetime of an individual.

The implications of exponential growth of the human population was the subject of
the famous writings of Thomas Malthus - \textit{Essay on the
Principle of Population - } (1798) which had a major impact and which was one of the influences on Darwin.

\subsection{The Logistic Equation}

 Verhulst (1838) (\cite{Ve-38})  introduced the {\em logistic  equation} which describes the more realistic situation in which resources are limited and the death rate increases as the resources are exhausted.

\bea{logistic}
\frac{dx}{dt}  &&  =\alpha x(1-\frac{x}{N}),\quad x(0)\geq 0,\quad \alpha \geq 0\\
\; x(t)  &&  =\frac{Nx(0)e^{\alpha t}}{N+x(0)(e^{\alpha
t}-1)}\rightarrow N\;\;\text{as }t\rightarrow\infty \eea

Here $N<\infty$ is interpreted as the {\em carrying capacity} of
the environment in which the population lived.

\subsection{The Lotka-Volterra Equations for Competing Species}

Equations to model the competition between species in ecology were
proposed by (Lotka (1925) \cite{L-25} and  Volterra (1926) \cite{V-26}):

\bea{L-V}
&&\frac{dx_1}{dt}=r_1x_1\left(1-\frac{x_1}{K_1}-a_{12}\frac{x_2}{K_1}\right)\\&&
\frac{dx_2}{dt}=r_2x_2\left(1-\frac{x_2}{K_2}-a_{21}\frac{x_1}{K_2}\right)
\eea

Coexistence, that is, a {\em stable equilibrium} with both species present
occurs if  $\frac{1}{a_{21}}>\frac{K_1}{K_2}> a_{12}$.

\begin{remark}     Gause (1934) proposed the {\em
competitive-exclusion principle} that states two  species cannot
stably coexist if they occupy the same niche, for example, if
$a_{12}=a_{21}$.
\end{remark}

\subsection{The SIR Epidemic Model}

A classical model for the progress of an epidemic due to  Kermack
and McKendrick (1927) \cite{KMcK-27} is given by the system of
ode:

\bea{SIR} &&\frac{dS}{dt}=-\beta SI,\quad \frac{dI}{dt}=\beta
SI-\gamma I, \quad \frac{dR}{dt}=\gamma I,\\&& S(0)>0 ,\; I(0)>0,\; R(0)=0.\eea Here $S$ denotes
the population of susceptible individuals, $I$ the population of
infectious individuals and $R$ the population of removed
individuals.

The {\em epidemiological  threshold} quantity is defined by \be{}
\emph{R}_0=\frac{\beta S(0) }{\gamma} (\text{  reproductive ratio)}\ee
If $\emph{R}_0<1$, then the infected population never increases whereas
if $\emph{R}_0>1$ the  epidemic ``will spread''.

\subsection{Population models and dynamical systems}

As suggested by these elementary examples, the modeling of interacting multitype populations leads to a rich area of dynamical systems and there exists an immense literature in this field.

For example the  extension of the Lotka-Volterra equations to  N interacting species is given by the system:
\be{}
\frac{x_i(t)}{dt}=r_ix_i\left(1-\sum_{j=1}^N\alpha_{ij}x_j\right),\quad i=1,\dots,N.
\ee

These multispecies dynamical systems can have very complex behavior
including limit cycles or chaotic behaviour.  In fact S. Smale
(1976) \cite{S-76} proved that for $N\geq 5$ these systems can exhibit any
asymptotic behavior.

\section{Small population effects}

Deterministic models provide good approximations to the growth of
large (noncritical) populations but for small populations and
``nearly critical populations'' it is essential to take account of
their inherent discrete nature and randomness.

\subsection{The Bienamy\'e-Galton-Watson Branching Process (BGW)}
\ The importance of the fact that individuals produce a random number of offspring and the possibility exists that the population can become extinct
 led  Bienaym\'{e} (1845) \cite{B-45}, and  Galton-Watson (1874) \cite{WG-74} to introduce this probabilistic model.

The population size at generation $n$ is denotes $X_n$.  Starting with $X_1=1$, at each generation each individual gives rise to a random number of children as follows:

\[ X_{n+1}=\sum_{k=1}^{X_{n}}\xi_{k},\text{  where the }\quad\{ \xi_{k}\}_{k=1,\dots,X_n} \text{ are }\textbf{ independent}.%
\]%

Generating functions provide a basic tool for developing these processes.  The generating function for the offspring distribution is given by:
\begin{align*}
f(s)  &
=E[s^{\xi}]=\sum_{k=0}^{\infty}p_{k}s^{k},\quad 0\leq s\leq 1\\
\;\;f^{\prime}(1)  &  =m \quad = \quad\text{mean offspring size}.\\
\end{align*}
The key relation is
\[ E[s^{X_n}]=f_n(s),\text{   where   } f_{n+1}=f[f_n(s)].\]

The {\textbf{extinction probability}}, is defined by  $q:= P(X_n=0\;\text{for some }n<\infty)$.

 \beT{} (Steffensen (1930,1932)) If $m\leq 1$ (critical, subcritical branching), the $q=1$.  If $m>1$ (supercritical branching), the $q$ is the  unique nonnegative solution in $[0,1)$
\be{}
s=f(s).
\ee
\end{theorem}

If $m<\infty$, then

\[ \frac{X_n}{m^n} \text{   is a martingale }.\]

\subsubsection{Propagation of initial randomness}

\beT{}(Hawkins and Ulam (1944), Yaglom (1947), Harris (1948), \cite{H-63})\\
If $m>1,\quad \sigma^2 =f^{\prime\prime}(1)+f^\prime(1)-(f^\prime(1))^2 <\infty$, then
\[  \frac{X_n}{m^n}\to W,\;\text{  in  } L^2\text{ and  a.s. as }n\to\infty \]
and
\[ E W=1,\quad \text{Var}(W)=\frac{\sigma^2}{(m^2-m)},\;
P(W=0)=q.\]
\end{theorem}

\subsection{Reed-Frost epidemic model}

A probabilistic analogue of the SIR epidemic model known as the { Reed-Frost model } is given as follows.  We consider an initial population of susceptible individuals  $S_0=N$ and one infected individual $I_0=1$.

\[S_{t+1}\sim
\text{Bin}(S_t,(1-p)^{I_t}),\;\; t\in\mathbb{N},\]
that is, in each time unit a susceptible individual has probability $p$ of meeting each infected individual and one such contact results in infection. Individuals are infected during one time period so that $I_{t+1}= S_t-S_{t+1}$.
If  ${ 1-p=e^{-\lambda/N}}$, then Von Bahr and Martin-L\"of (1980) \cite{BML-80} showed that as $N\to \infty$ the critical threshold is $\lambda =1$.

\subsection{Multitype populations and the Wright-Fisher Model}

The celebrated work of  Mendel (1865) \cite{M-65} on the  inheritance of traits and its rediscovery around 1900 led to the development of the field of genetics.
The modern theory of mathematical genetics was initiated in the work of  Wright (1931), (1932)
\cite{W-31},\cite{W-32} and  Fisher (1930) \cite{F-30}. They introduced a probabilistic model of
finite population sampling that serves as a starting point  for modern population genetics. This model deals with a population of individuals of different types.  As a mathematical idealization they assume that the total population is constant in time and they focus on the changes in the relative proportions of the different types of individual.
The key ingredients are:

\begin{itemize}
\item Fixed finite population size $N$

\item Typespace \ (\textit{alleles})%
\[
E_K=\{1,\dots,K\}
\]
\item
$X_{n}(i)$ is the number of individuals of type $i,$ at generation
$n$.

\end{itemize}

Let $\mathcal{N}(E_K)$ denote the counting measures on $E$.
Then dynamics are defined by a Markov chain $X_n=(X_{n}(1),\dots,X_{n}(K))$ with state space
\[ \{(x_1,\dots,x_K)\in\mathcal{N}:\;\sum_{i=1}^K x_i=N\}.\]

 The intuitive idea leading to the transition mechanism for the {\em neutral model} is that first each individual in the nth generation produces a large number of
 potential offspring. Then in a second stage the population is pruned back (culling) so that the total population remains $N$ (this can be thought of as an analogue of carrying capacity).   Based on the {\em neutral assumption}, that is each of the individuals in produced in the first stage has equal probability of being selected,  the $(n+1)$st generation consists of $N$ individuals of types $\{1,\dots,K\}$ obtained by

\begin{itemize}
\item \emph{multinomial sampling} from the empirical distribution
\end{itemize}

\begin{align*}
P  &  (X_{n+1}=(y_{1},\dots,y_{K})|X_{n}=(x_{1},\dots,x_{K}))\\
&  =\frac{N!}{y_{1}!y_{2}!\dots y_{K}!}\left(
\frac{x_{1}}{N}\right)
^{y_{1}}\dots\left(  \frac{x_{K}}{N}\right)  ^{y_{K}}%
\end{align*}
\bigskip

An important feature of this process is the loss of information (diversity) leading to  \textbf{fixation}, that is, the long time survival of exactly one type. To see this note that
\[
p_{n}(i)\text{ is a { martingale} where
}p_n(i)=\frac{X_n(i)}{N}\] and \[p_{n}(i)\rightarrow\text{ }0\text{ \
or
\ }1\text{ as }n\to\infty \text{  for each }i\text{ w.p.1. }%
\]

\underline{The dual perspective}
\medskip

If we choose $k$ individuals at random from generation $n+1$ and look backwards in time to identify the parents in the nth generation, by an elementary conditional probability calculation, we see that   { \textit {each individual in generation $n+1$  picks its
parents ``at random''.}}  This naturally leads to the notion of {\em identity by descent} introduced by   Mal\'ecot (1941) \cite{M-41}, that is,   two individuals are identical by descent if they  have a common ancestor (and no mutations have occurred).
\begin{center}


\end{center}
\bigskip

\bigskip

\section{The Role of Stochastic Analysis}

\textbf{Basic developments in stochastic analysis}:

The development of stochastic population modelling was made possible by the remarkable developments in stochastic analysis.

\begin{itemize}

\item Markov chains and processes (1906) \cite{M-06},  Kolmogorov (1931), \cite{K-31}

\item Brownian motion Wiener (1923) \cite{W-23}, L\'evy (1948) \cite{L-48}

\item  Ito stochastic calculus  (1942), (1946), (1951) \cite{I-51}

\item Markov processes and their semigroup characterization\\ Feller  (1951) \cite{F-51}, It\^o-McKean (1965).\\

Given a Markov process $\{X(t)\}_{t\geq 0}$ with state space $E$ (for example, compact metric space) and $f\in C(E)$ (bounded continuous functions on $E$) and $x\in E$, let
\[
T_{t}f(x)=E_{x}(f(X(t))\]
$\{T_t\}$ is said to be strongly continuous if $\|T_tf -f\|\to 0$ as $t\downarrow 0$ for $f\in C(E)$.

Then for some class $f\in D(G)\subset C(E)$ the generator $G$ acting on $f$ is defined by
\[Gf(x)=\lim_{t\downarrow
0}\frac{T_tf(x)-f(x)}{t} \text{  exists and }\in C(E).
\]
\end{itemize}
Conditions under $(D(G),G)$  defines a strongly continuous semigroup were obtained in the celebrated Hille-Yosida Theorem (1948)
\cite{H-48}, \cite{Y-48}.

\subsection{Diffusion approximations of branching and Wright-Fisher processes}

In two seminal papers Feller (1939) \cite{F-39}, Feller (1951) \cite{F-51} developed diffusion process approximations to the branching process and the Wright-Fisher model.  These are now referred to as the Feller continuous state branching process (CSBP) and Wright-Fisher diffusion process, respectively. These serve as the ``bridge'' between the discrete world of individuals and generations and the world of differential equations and dynamical systems. It simultaneously maintains the power of analysis of the latter world and the random finite population effects of the former world.

We now state the two basic results.

\begin{theorem} \underline{Nearly critical BGW processes to Feller CSBP
Branching}\newline Consider a sequence of  BGW processes with mean offspring sizes $m_{N}=1+\frac{m}{N}$ and constant finite variance. Assume that
$N^{-1}X_{0}^{N}\rightarrow X_{0}$ as $N\rightarrow \infty$. Then
\[
\{
\frac{1}{N}X_{\lfloor Nt\rfloor}^{N},\;t\geq0\}\Longrightarrow\{X_{t}:t\geq0\}
\]
where the convergence is in the sense of weak convergence of c\`adl\`ag processes.
The limiting process is a continuous process with state space $[0,\infty)$ and the generator of the associated semigroup
is given by
\[ Gf(x)= mx\frac{\partial f}{\partial x} +\frac{\gamma}{2}
x\frac{\partial^2 f}{\partial x^2}\quad\text{for some   }\gamma >0.\]
\end{theorem}
See Section 4.3 for details.

\bigskip
\begin{theorem}
{\underline{Wright-Fisher Diffusion}}\newline Consider a sequence of $K$-type  Wright-Fisher Markov chains
$X^N$ with total population $N$ and assume that $N^{-1}X^N_{0}\rightarrow \mathbf{p}_{0}\in\mathcal{P}(E_K)$ (probability measures on $E_K$) as $N\rightarrow\infty$  . Then
\[
\{\mathbf{p}_{N}(t):t\geq0\}\equiv\{\frac{1}{N}X_{\lfloor Nt\rfloor},\;t\geq
0\}\Longrightarrow\{\mathbf{p} (t):t\geq0\}
\]
where $\{\mathbf{p}(t):t\geq0\}$ is a Markov diffusion process
with values in
the simplex%
\[
\Delta_{K-1}=\{(p_{1},\dots,p_{K}):\;p_{i}\geq0,\sum_{i=1}^{K}p_{i}=1\}
\]
and generator acting on functions $f(p)=f(p_{1},\dots,p_{K}):$%
\bigskip

{
\begin{align*}
G^{(K)}f(\mathbf{p})  &  = \frac{1}{2}\sum
_{i,j=1}^{K}p_{i}(\delta_{ij}-p_{j})\frac{\partial^{2}f(\mathbf{p})}{\partial
p_{i}\partial p_{j}}.%
\end{align*}
}
\end{theorem}

\begin{remark}In the case $K=2$ it suffices to keep track of $p_1\in [0,1]$ and then the generator is
\[  Gf({p_1})=\frac{1}{2}p_1(1-p_1)\frac{\partial^{2}f({p_1})}{\partial
p_1^2}.\] In terms of It\^o's SDE, this satisfies
\[  dp_1(t)= \sqrt{p_1(t)(1-p_1(t))}dw(t),\; p_1(0)\in [0,1],\]
where $\{w(t)\}_{t\geq 0}$ is a standard Brownian motion.

\end{remark}

See Section 5.2 for details.


\section{Darwinian selection}

Darwin's theory of evolution ({\em On the Origin of Species}, (1859 \cite{D-59}) was based on the concept on natural selection based on the differential reproductive success of the different types of individual. In a seminal paper in 1924 J.B.S. Haldane \cite{H-24}, \cite{H-32} initiated the modern synthesis of Darwinism evolution and Mendelian genetics (Mendel (1865)) and formulated the notion of {\em fitness}.

\subsection{Fisher's large population approximation}

A deterministic mathematical model incorporating the notion of fitness was developed by Fisher \cite{F-30} as follows. Consider an infinite diploid population (organisms have a type $(i,j)\in E\times E$) which reproduce sexually with random mating, that is, the offspring type is obtained by choosing two individuals at random (parents) and choosing one of the homologous pairs from each parent. \ The  type space is $E\times E$ ( genotype
determined by the gametes $i$ and $j)$.

Let $x_{i}(t)$ be the amount of  gamete $i$ in the population at
time $t$ and $p_{i}$ denote the frequency $p_{i}=\frac{x_{i}}{\sum
x_{i}}.$

Let $V(i,j)=V(j,i)$ $=$ ``diploid fitness'' \ of the genotype $(i,j)$. The
instantaneous fitness,
$V(i,p)$ of the ith gamete is defined by%
\[
V(i,p)=\sum_{j}p_{j}V(i,j)
\]
and the mean fitness is defined by%
\[
\bar{V}=\sum_{i}V(i,p)p_{i}=\sum_{ij}p_{i}p_{j}V(i,j)\text{.}%
\]

(The {\em  haploid case} is similar  to the additive case $V(i,j)=V(i)+V(j)$.)

In Fisher's formulation the population sizes $x_{i}$ satisfy the differential equations%
\[
\frac{dx_{i}}{dt}=x_{i}V(i,p),\;\;\;i=1,\dots,K
\]
and therefore the proportions $\{p_{i}\}$ satisfy the equations:%
{{
\be{}
\frac{d p_{i}}{dt}=p_{i}(V(i,p)-\bar{V}),\;i=1,\dots,K
\ee
}}
\begin{theorem}
(Fisher's Fundamental Theorem) (\cite{F-30})
\end{theorem}
\begin{itemize}
\item (a) Mean fitness $\bar V(t)$ increases on the trajectories
of \ $\mathbf{p}(t)$.

\item (b) The rate of change of the mean $\bar V(t)$ along orbits
is proportional to the variance.
\end{itemize}

See Section 12.1 for details.

\subsection{Selection and Genetic Drift}

In contrast to Fisher, Wright \cite{W-31}, \cite{W-32}
considered genetic drift (due to finite population effects) to play an important role in evolution and developed his
{\em shifting balance theory of evolution}.  The relative importance to evolution of different mechanisms remains a subject of
investigation and debate (see for example  Barton-Turelli (1997) \cite{BT-97}, Ohta-Gillespie (1996) \cite{OG-96}).

To introduce the interplay between selection and {\em genetic drift} we consider the diffusion approximation to the finite population model with  two types $(1,2)$ and {\em haploid selection}. Let  $p_1(t)$ denote the proportion of type 1. Let type 1 have
fitness $s>0$ and type 2 have fitness $0$. Then the diffusion approximation limit $p_1(t)$ satisfies the SDE
\be{wfd}
dp_1(t)=sp_1(t)(1-p_1(t))dt+\sqrt{\gamma p_1(t)(1-p_1(t))}dw(t)
\ee
Here $\gamma$ is proportional to the inverse of the {\em effective population size}.

The relation between the probability laws of $p_1(t)$ with and without selection (i.e. setting $s=0$) follows from the
Cameron-Martin-Girsanov representation (\cite{RY-91},  Chapt. VIII) as follows.

\begin{theorem}   Let $P_{[0,t]}^s$ be the probability
law on $C_{[0,1]}([0,t])$ of the solution of (\ref{wfd}). \\ The Radon-Nikodym derivative on
$\mathcal{P}(C_{[0,1]}([0,t]))$ is given by: { \[
\frac{dP_{[0,t]}^{s}}{dP_{[0,t]}^{0}}=\exp\left(  \frac{s}{\gamma
}(p_1(t)-p_1(0))-\frac{s^{2}}{\gamma}\int_{0}^{t}p_1(s)(1-p_1(s))ds\right).
\]}
\end{theorem}

The deviation  of very large but finite populations from the deterministic (infinite population) limit can be analysed
by considering the asymptotics as $\gamma\to 0$ and using
Freidlin-Wentzell (see \cite{FW-84}) large deviation methods.

 \section{Spatially structured population systems} {

The above models assume that any new individual can be chosen from
any member of the population and the members of the population
interact in an exchangeable manner at any time.  However real populations are distributed in
space and reproduce and compete locally.  Spatial models play an
essential role in the study of population systems. We begin with a
basic formulation with discrete ``geographic space''.

To begin we consider a population with subpopulations located at sites on the a finite or countable set $S$, for example the lattice  $S=\mathbb{Z}^d$, in which individual migrate
between sites with  migration rates given by a symmetric random walk kernel $\{ p_{\xi-\xi^{\prime}}\}_{\xi,\xi^\prime\in S}$.

 \subsection {Super Random Walk}
The branching random walk (see Athreya-Ney \cite{AN-72}) extends
the basic Bienamy\'e-Galton-Watson process to the situation in which
individuals are located in a countable space $S$.

 The diffusion
limit of branching random walks on $S$ leads to a system of
stochastic differential equations, now called {\em super random walks} (SRW),
as follows:

\begin{align*}
dx_{\xi}(t)  &  =\sum_{\xi^{\prime}\in S}p_{\xi-\xi^{\prime}}(x_{\xi^{\prime}%
}(t)-x_{\xi}(t))dt+\sqrt{\gamma x_{\xi}(t)}dw_{\xi}(t)\\
&  x_{\xi}(0)\geq 0
\end{align*}
where \ $\{w_{\xi}(\cdot)\}_{\xi\in S}$ is a system of independent
Brownian motions. We note that if $S$ is infinite, then the space
of configurations $[0,\infty)^S$ is infinite dimensional.  We can
also identify this with the set of locally finite  measures on $S$,
$\mathcal{M}(S)$. The state at time $t$ is then given by a random
measure on $S$ (see Moyal (1962) \cite{M-62} for an early
formulation of spatial population processes). A key tool in the
study of an important class of random measures is the Laplace
functional.
\bigskip

\bigskip

\begin{theorem}
The transition Laplace Functional of the SRW is given by
 \[
E_{\mathbf{x}(0)}\exp\left(  -\sum_{\xi\in
S}\varphi(\xi)x_{\xi}(t)\right) =\exp\left(  -\sum _{\xi\in
S}v_{t}(\xi)x_{\xi}(0)\right) \]%
for $\varphi\in C_{+}(S)$, where $v_t$ is the unique solution of \[
v_{t}(\xi)=S_{t}\varphi(\xi)-\int_{0}^{t}S_{t-s}\left(
v_{s}^{2}(\xi)\right) ds
\]
where $\{S_{t}\}$ is semigroup on $C_b(S)$ with generator
\be{}  Gf(\xi)= \sum_\xi p_{\xi-\xi^{\prime}}(f(\xi')-f(\xi)).\ee
\end{theorem}

 The measure-valued analogue of this system on $\mathbb{R}^d$,
now called super-Brownian motion,
 was introduced by S. Watanabe (1968) \cite{W-68} in the context of branching processes and by Dawson (1975)\cite{D-75} in the context of stochastic evolution equations.

\subsection{Stepping Stone Models}

The introduction of spatial models in population genetics goes
back to Wright's island model  \cite{W-31} and the work of
Mal\'{e}cot (1941) \cite{M-41}, (1948) \cite{M-48}, (1949) \cite{M-49}, Kimura
(1953) \cite{K-53}, and  Sawyer (1976) \cite{S-76a}.

We will consider below the Wright-Fisher two-type diffusion
stepping stone model with selection {
\begin{align*}
dx_{\xi}(t)  &  =c\sum_{\xi^{\prime}\in S}p_{\xi-\xi^{\prime}}(x_{\xi^{\prime}%
}(t)-x_{\xi}(t))dt\\&+sx_\xi(t)(1-x_\xi(t))dt+\sqrt{2x_{\xi}(t)(1-x_{\xi}(t))}dw_{\xi}(t)\\
x_{\xi}(0)  &  = \theta\in [0,1] \quad\forall \xi
\end{align*}
}
This arises from a collection of Wright-Fisher populations at {\em demes} $\xi\in S$
in which there is probability  $\frac{cp_{\xi-\xi^\prime}}{N}$ that an individual in generation $n+1$
is the offspring of an individual at  deme $\xi^\prime \ne \xi$ in generation $n$.

The stepping stone model is closely related to the voter model
which was introduced by  Clifford, Sudbury (1973) \cite{CS-73},
Holley and Liggett (1975) \cite{HL-75}. The voter model  is a
${\{0,1\}}^S$-valued Markov jump process and is one of the
principal  examples of the class of  interacting particle systems
introduced by Spitzer (1970) \cite{S-70} and  Dobrushin
\cite{DOB-71} (1971) and extensively developed over the past 35
years.

\subsection{Spatial spread of advantageous genes and epidemics}%

An important feature of spatial models is spatial spread, for
example,  of a mutant type or epidemic.  A classic example is the
wave of advance of an advantageous gene modelled by the celebrated
Fisher-KPP equation (1937) (Fisher  (\cite{F-37},  Kolmogorov,
Petrovsky and Piscounov \cite{KPP-37}):

\[
u(t,x)=\text{ proportion of advantageous type }1\text{ at }x\text{
at time }t
\]%
\begin{align*}
 u_{t}  &  =\frac{1}{2}u_{xx}+u(1-u) \\
u(0,x)  &  =1_{(-\infty,0]}(x)
\end{align*}
The fundamental results of Kolmogorov, Petrovsky and Piscounov
establish the existence of  travelling wave solutions at speed
$\sqrt{2}$.  Relations between this travelling wave and the
maximal displacement of branching Brownian motion were used by
Bramson (1983) \cite{B-83} to obtain fine properties of this
phenomenon.

\section{Complex population dynamics}

A generalization of the N species Lotka-Volterra equations is give by the nonlinear dynamical system:
\be{}
\frac{x_i(t)}{dt}=F_i(x_1(t)),\; i=1,\dots,N.
\ee
As mentioned above these multispecies dynamical systems can have very complex behavior
including limit cycles or chaotic behaviour.   Questions of the  complexity, robustness, diversity and spatial structure have been the subject
of much research and debate much of it stimulated by the work of Hutchinson (1957) \cite{H-57} and MacArthur (1955) \cite{Mac-55}, \cite{MW-67}, \cite{Mac-72}.  Some information can be gained by analyzing the behaviour near stationary points, that is $(x_1,\dots,x_N)$ where $F_i(x_1,\dots,x_N)=0\;\forall i$ by considering the linearized system around these points and in particular the spectra of the resulting matrices.  One approach  to  these questions was developed in the 1973 paper of May \cite{M-73} in which he assumed that the matrices are random and used results on random matrix theory to look at the relation between species number $N$ and complexity (measured by the proportion of non-zero matrix elements). Limitations of this analysis have been pointed out \cite{CN-85} and the assumption of randomness does not reflect the dynamical mechanisms involved.  A reasonably robust ecological systems is constantly being tested by the emergence of new mutant types.  For a robust system most such mutations and deleterious and are eliminated.  However rare mutant types can cause the collapse of the system or move the system to a new attractor and a higher or lower order of organization. {\em This means that the ecosystem dynamics is itself subject to evolutionary forces.}  A theoretical framework to classify the types of dynamical system to develop and their stability remains elusive. However stochastic effects are important here in that the generation and survival at low population sizes of mutant types clearly plays an important role.

\subsection{Basic questions on modelling complex and evolving populations}

We next take a quick look at some of the basic questions on the
modelling of complex evolving populations.

\begin{itemize}

\item How do we model a complex multilevel population?

\item  Individuals - type, internal state and dynamics, geographical location,
  fitness, interaction matrix.

\item  reproduction and modification mechanisms: birth,  death and mutation rates,

\item Family structure and genealogy

\item migration dynamics

\item Spatial distribution of types: role of spatial and
hierarchical network structures

\item Networks of interactions between individuals and groups of individuals, for example, competition and or cooperation between
types, for example, in ecology and economics. Stability and collapse.

\item How do we relate the different levels of description: microscopic, mesoscopic and macroscopic?

\item How do we describe the development of  population composition and
structure in different space and time scales going from the
microscopic to evolutionary scales.

\item emergence of new levels of organization

\end{itemize}

\subsection{The role and analysis of stochasticity}

We have seen some of the classical stochastic population models in in our historical survey.  Today, in view of the revolutionary developments in biology over the past 50 years there is an endless richness of biological phenomena and models. There is a huge literature on both deterministic and stochastic models. Once can ask:
what is the role of stochasticity in the development of
complex populations? Where do deterministic models fail and the
role of randomness is essential? There is no simple answer but we make the following four
observations:

\begin{itemize}

\item Evolution is an interplay of nonlinear dynamics, finite population fluctuations and rare random events  at which  fundamental transitions occur in the
population composition and dynamics.

\item Extinction and the loss of information due to finite population resampling.

\item The creation of diversity and information is a result of selective forces acting on randomly produced
mutant types. The latter can be viewed as a  random search through a potentially infinite search space.

\item Demographic and environmental stochasticity is  ubiquitous
and plays a role analogous to molecular motion in statistical
physics. The levels of these sources of randomness influences the
nature of quasiequilibria and non-equilibrium phase transitions.

\end{itemize}

Fortunately over the past fifty year, stochastic analysis has undergone major  developments in many directions, partly in response to these challenges. This has produced a range of tools that have proved effective in addressing some significance issues in the sciences including the biological sciences and hold potential to address even more of these questions in the future.  The objective of these notes is to introduce some of the basic ideas and tools and to provide some pointers to the growing literature in this field.



\bigskip

\chapter[Branching Processes I]{Branching Processes I:\\Supercritical growth and population structure}

The fundamental characteristic of biological populations is that
individuals undergo birth and death and that individuals carry
information passed on from their parents at birth. Furthermore
there is a randomness in this process in that the number of births
that an individual gives rise to is in general not deterministic
but random. Branching processes model this process under
simplifying assumptions but nevertheless provide the starting
point for the modelling and analysis of such populations.   In
this chapter we present some of the central ideas and  key results
in the theory of branching processes.

\section{Basic Concepts and Results on Branching Processes}


\subsection{Bienamy\'e-Galton-Watson processes}

The Bienamy\'e-Galton-Watson branching process (BGW process) is a Markov chain
on $\N_0:=\{0,1,2,\dots\}$. The discrete time parameter is interpreted
as the generation number and $X_n$ denotes the number of
individuals alive in the n'th generation. Generation $(n+1)$
consists of the offspring of the nth generation as follows:
\begin{itemize}
\item  each individual $i$  in the nth generation produces a
random number $\xi_i$ with distribution
\[  p_k=P[\xi_i=k],\;k\in \N_0\]
\item $\xi_1,\xi_2,\dots,\xi_{X_n}$ are independent.
\end{itemize}

Let $X_0=1$. Then for $n\geq 0$
\[
X_{n+1}=\sum_{i=1}^{X_{n}}\xi_{i},\quad\{ \xi_{i}\} \text{ independent}%
\]%
We assume that the mean number of offspring \[ m=\sum_{i=1}^\infty
ip_i<\infty.\]

The BGW process is said to be {\em subcritical} if $m<1$, {\em
critical } if $m=1$ and {\em supercritical} if $m>1$.

A basic tool in the study of branching processes is the {\em generating
function}
\be{}
f(s)   =E[s^{\xi}]=\sum_{k=0}^{\infty}p_{k}s^{k},\quad 0\leq
s\leq 1.\ee
Then
\be{}
 f^{\prime}(1)    =m,\quad f^{\prime\prime}(1)=E[\xi(\xi-1)]\geq 0.\\
\ee
Let
\[ f_n(s)=E[s^{X_n}],\; n\in\N.\]
   Then conditioned on $X_n$, and using the independence of the $\{\xi_i\}$,
\[ f_{n+1}(s)=
E[s^{\sum_{i=1}^{X_n}\xi_i}]=E[f(s)^{X_n}]=f_n(f(s))=f(f_n(s)).\]
Note that $f(0)=P[\xi=0]=p_0$ and
\[ P[X_{n+1}=0]=f(f_n(0))=f(P[X_n=0])\]

Then if $m> 1,\; p_0>0$,    $P[X_n=0]=f_n(0)\uparrow q$ where
$q$ is the smallest nonnegative root of
\[ f(s)=s,\]
and if $m\leq 1$, $P[X_n=0]\uparrow  1$. Note that $1$ and $q$ are  the only roots of $f(s)=s$.

Since $E[X_{n+1}|X_n]=mX_n$,
\be{}  W_n:=\frac{X_n}{m^n}\text{   is a martingale and } lim_{n\to\infty} W_n =W\text{ exists } a.s.
\ee

\beP{q0} We have  $P[W=0]= q\text{  or  } 1$, that is,  conditioned on nonextinction either $W=0$
a.s. or $W>0$ a.s.
\end{proposition}
\begin{proof}
 It suffices to show that $P[W=0]$ is a
root of $f(s)=s$.  The ith individual of the first generation has
a descendant family with a martingale limit which we denote by $W^{(i)}$. Then
$\{W^{(i)}\}_{i=1,\dots, X_1}$  are independent and have the same
distribution as $W$. Therefore \be{} W=\frac{1}{m}\sum_{i=1}^{X_1}
W^{(i)}\ee and therefore $W=0$ if and only if for all $i\leq X_1$,
$W^{(i)}=0$. Conditioning on $X_1$ implies that
\be{}P[W=0]= E(P(W^{(i)}=0)^{X_1})=f(P[W=0]). \ee Therefore $P[W=0]$ is a root of $f(s)=s$.
\end{proof}

\begin{remark} In the case $\rm{Var}(X_1)=\sigma^2<\infty$ we can show by induction that
\be{} \rm{Var}(X_n)=\begin{split}& \frac{\sigma^2 m^n(m^n-1)}{m^2-m}, \quad m\ne 1,\\
&n\sigma^2,\quad m=1\end{split}
\ee
\end{remark}

Then if $m>1$ the martingale $\frac{X_n}{m^n}$ is uniformly integrable and $E(W)=1$. Moreover $\frac{X_n}{m^n}\to W$ in $L^2$ and
\be{} \rm{Var}(W)=\frac{\sigma^2}{m^2-m} >0 \qquad \text{(see Harris \cite{H-63} Theorem 8.1)}.\ee

If $m>1$, $\sigma^2=\infty$,  a basic question concerns the nature of the random variable $W$ and the question whether
or not $\frac{X_n}{m^n}\to W$ in $L^1$. The question was settled by a celebrated result of Kesten and Stigum which we present in Theorem \ref{KS} below. We first introduce some further basic notions.

\subsubsection{Bienamy\'e-Galton-Watson process with immigration (BGWI)}
The Bienamy\'e-Galton-Watson process with offspring distribution $\{p_k\}$
and immigration process $\{Y_n\}_{n\in\N_0}$ satisfies

\be{}   X_{n+1}=\sum_{i=1}^{X_{n}}\xi_{i} +Y_{n+1},\ee where the
$\xi_i$ are iid with distribution $\{p_k\}$.

Let $\mathcal{F}^{Y}$ be
the $\sigma$-field generated by $\{Y_k:k\geq 1\}$ and  $X_{n,k}$
be the number  of descendants at generation $n$ of the individuals who
immigrated in generation $k$.  Then the total number of individuals
in generation $n$ is $X_n=\sum_{k=1}^n X_{n,k}$.

For $k<n$ the random variable $W_{n,k}=X_{n,k}/m^{n-k}$
has the same law as $\wt X_{n-k}/m^{n-k}$
 where $\wt X_n$ is the BGW process with  $Y_k$ initial particles.
Therefore
\be{} E[\frac{X_{n,k}}{m^{n-k}}]=Y_k.\ee

Now consider the subcritical case $m<1$.  If $\{Y_i\}$ are i.i.d. with $E[Y_i]<\infty$, then the Markov chain $X_n$ has a stationary measure
with mean $\frac{ E[Y]}{1-m}$.

Next consider the supercritical case $m>1$.  Then
\be{3.8} E[\frac{X_n}{m^n}|\mathcal{F}^{Y}] =
E[\frac{1}{m^n}\sum_{k=1}^nX_{n,k}|\mathcal{F}^{Y}]=
\sum_{k=1}^n\frac{1}{m^k}
E[\frac{X_{n,k}}{m^{n-k}}|\mathcal{F}^{Y}] =\sum_{k=1}^n\frac{Y_k}{m^k}.\ee
If $\sup_k E[Y_k]<\infty$, then
\be{}  \lim_{n\to\infty}\frac{E[X_n]}{m^n}=\sum_{k=1}^\infty\frac{E[Y_k]}{m^k}<\infty.\ee

A dichotomy in the more subtle case $E[Y_k]=\infty$ is provided by the following theorem of Seneta.

\beT{T.Sen} (Seneta (1970) \cite{Sen-70}) Let $X_n$ denote the BGW process with mean offspring
$m>1$, $X_0=0$  and with i.i.d.
immigration process $Y_n$.\\(a)  If $E[ \log^+
Y_1]<\infty $, then $\lim \frac{X_n}{m^n}$ exists and is finite
a.s.\\(b) If $E[\log ^+ Y_1]=\infty$, then
$\limsup\frac{X_n}{c^n}=\infty $ for every constant $c>0$.

\end{theorem}

\begin{proof} The theorem is a consequence of the following elementary result.

\beL{L1} Let  $Y,Y_1,Y_2,\dots$ be nonnegative iid rv. Then a.s.
\be{}\limsup_{n\to\infty} \frac{1}{n}Y_n=\left\{%
\begin{array}{ll}
    0, & \hbox{ if } E[Y]<\infty\\
    \infty, & \hbox{ if } E[Y]=\infty\\
\end{array}%
\right. \ee

\end{lemma}

\begin{proof}
Recall that $E[Y]=\int_0^\infty P(Y>x)dx$. This gives $\sum_n
P(\frac{Y}{n}> c) <\infty $ for any $c>0$   if $E[Y]<\infty$  and
the result follows by  Borel-Cantelli. If $E[Y]=\infty $, then
$\sum P(\frac{Y}{n}> c) =\infty $ for any $c>0$ and the result
follows by the second Borel-Cantelli Lemma since the $Y_n$ are
independent.
\end{proof}
\medskip

Proof of (a).  By (\ref{3.8})  \be{3.x}
E[\frac{X_n}{m^n}|\mathcal{F}^Y]=\sum_{k=1}^n\frac{Y_k}{m^k}.\ee
Since here we assume $E[\log^+ Y_1]<\infty$,   Lemma
\ref{L1} gives  $\limsup_{k\to\infty}\frac{Y_k}{c^k}<\infty$ for any
$c>0$.   Therefore the series given by the last expression in (\ref{3.x}) converges a.s.  and
therefore $\lim_{n\to\infty} E[\frac{X_n}{m^n}|\mathcal{F}^Y]$ exists and is finite a.s. This implies (a)

Proof of (b). If $E[\log^+ Y_1]=\infty$,
then by  Lemma \ref{L1} $\limsup_{n\to\infty} \frac{\log^+ Y_n}{n}=\infty$
a.s. Therefore for any $c>0$ \be{}\limsup_{n\to\infty} \frac{Y_n}{c^n}=\infty\ee
a.s. Since $X_n\geq Y_n$, (b) follows.

\end{proof}

\subsection{Bienamy\'e-Galton-Watson trees}%

In addition to the keeping track of the  total population of
generation $n+1$ in a BGW process it is useful to incorporate genealogical
information, for example, which individuals in generation $n+1$
have the same parent in generation $n$.  This leads to a natural
family tree structure which  was introduced in the papers of Joffe and  Waugh (1982), (1985),
\cite{JW-82}, \cite{JW-85} in their determination of the
distribution of kin numbers and developed in the papers of Chauvin (1986) \cite{C-86} and Neveu (1986) \cite{N-86}.

A convenient representation of the BGW random tree is as follows.
Let $u=(i_1,\dots,i_n)$ denote an individual in generation $n$ who
is the $i_n$th child of the $i_{n-1}$-th child of $\dots$ of the
$i_1$-th child of the ancestor, denoted by $\emptyset$.  The
{\em space of individuals} (vertices) is given by
\be{ind} \mathcal{I}=\{\emptyset\}\cup \cup_{n=1}^{\infty}\mathbb{N}^{n}.\ee

Given  $u=(u_{1},\dots,u_{m}),\;v=(v_{1}%
,\dots,v_{n})\;\in \mathcal{I},$ we denote the composition by
$uv:=(u_{1},\dots,u_{m},v_{1},\dots,v_{n})$
\medskip

A {\em plane rooted tree} $\mathcal{T}$ with root $\emptyset$ is a
subset of $\mathcal{I}$ such that

\begin{enumerate} \label{PRT}
\item $\emptyset\in\mathcal{T},$

\item If $v\in\mathcal{T}$ and $v=uj$ for some $u\in\mathcal{I}$
and $j\in \mathbb{N}$, then $u\in\mathcal{T}$.

\item For every \ $u\in\mathcal{T}$, there exists a number \
$k_{u}(\mathcal{T})\geq0,$ \ such that $uj\in\mathcal{T}$ if and
only if \ $1\leq j\leq k_{u}(\mathcal{T})$.
\end{enumerate}
A plane tree can be given the structure of a graph in which a parent is connected by an edge to each of its offspring.

Let \ $\mathbf{T}$ be the set of all plane trees.  If $t\in
\mathbf{T}$ let $[t]_n$ be the set of rooted trees whose first $n$
levels agree with those of $t$. Let $\mathbf{V}$ denote the set of connected sequences in $\mathcal{I}$, $\varnothing, v_1,v_2,\dots$, which do not backtrack. Given $t\in\mathbf{T}$, let $V(t)$ denote the set of paths in $t$. If $v_n$ is a vertex at the nth
level, let $[t;v]_n$ denote the set of trees with distinguished
paths such that the tree is in $[t]_n$  $v\in{V}(t)$ and the path  goes through $v_n$.

Given a finite plane tree $\mathcal{T}$ the
{\em height} $h(\mathcal{T})$ is the maximal generation of a vertex in
$\mathcal{T}$ and $\#(\mathcal{T})$ denotes the number of vertices
in $\mathcal{T}$. Let $\mathbf{T}_n$ be the set of trees of height $n$.

A {\em random tree} is given by a probability measure on $\mathbf{T}$.
Given an offspring distribution $\mathcal{L}(\xi)=\{p_k\}_{k\in\N}$, the
corresponding  BGW tree is constructed  as follows:

Let the initial individual be labelled $\emptyset$. Give it a
random number of children denoted $1,2,\dots, \xi_{\emptyset}$.

Then each of these has a random number of children, for example
$i$ has children denoted $ (i,1),\dots,(i,\xi_i)$ etc. Each of
these has children, for example $(i,j)$ has  $\xi_{i,j}$ children
labelled $(i,j,1),\dots, (i,j,\xi_{i,j})$, etc.  Then considering the first n generations
in this way we obtain a probability measure $P^{BGW}_n$ on $\mathbf{T}_n$.

The  probability measures, $P^{BGW}_n$ form a consistent family and induce a
 probability measure $P^{BGW}$ on $\mathbf{T}$, the law of the BGW random tree.

Let \be{}Z_n= \text{
number of vertices in the tree at level }n.\ee Then by the construction it follows that $Z_n$ is a version of the
BGW process and we can think of  the BGW tree as an enriched version
of the BGW process.

\begin{figure}[h]
\begin{center}
\setlength{\unitlength}{0.01cm}
\begin{picture}(1500,600)
\put(700,550){$\varnothing$}
\put(1020,550){$Z_0=1$}
\put(700,550){\line(-3,-4){110}}
\put(700,550){\line(-2,-5){180}}
\put(700,550){{\line(0,-1){150}}}
\put(700,550){\line(2,-5){180}}

\put(580,400){{{{\small $ 1$}}}}
\put(640,400){{\small {$-2$}}}
\put(760,400){{\small{$4$}}}
\put(700,400){{\small{$3$}}}
\put(1020,400){$Z_1=4$}
\put(700,400){{\line(-2,-5){60}}}
\put(700,400){\line(2,-5){60}}
\put(760,400){\line(4,-5){110}}

 \put(580,250){{\small 21}}
\put(640,250){{\small 31}}
\put(760,250){{\small 32}}
\put(820,250){{\small 41}}
\put(880,250){{\small 42}}
\put(1020,250){$Z_2=5$}
\put(640,250){{\line(-2,-5){60}}}
 \put(640,250){\line(0,-1){150}}
\put(640,250){{\line(2,-5){60}}}
 \put(760,250){\line(0,-1){150}}
\put(820,250){\line(0,-1){150}}
\put(520,100){\small{$_{211}$}}
\put(580,100){$_{311}$}
\put(640,100){$_{312}$}
\put(700,100){$_{313}$}
\put(760,100){$_{321}$}
\put(820,100){$_{411}$}
\put(880,100){$_{412}$}
\put(1020,100){$Z_3=7$}
\end{picture}
\end{center}
\caption{BGW Tree}
\end{figure}



\subsubsection{The size-biased BGW tree}

The fundamental notion of size-biasing has many applications.
It will be used below in the proof of  Lyons,
Pemantle and Peres (1995) \cite{LPP-95}
of some basic results on Bienamy\'e-Galton-Watson processes (see Theorem \ref{KS} below).

 To exploit  this  notion for branching processes we
consider the {\em size-biased offspring distribution}

\be{} \wh p_k=\frac{kp_k}{m},\quad k=1,2,\dots.\ee
We denote by $\wh\xi$ a random variable having the size biased
offspring distribution. The  size-biased BGW tree $\wh T$ is
constructed as follows:
\begin{itemize}
   \item  the initial individual is labelled $\varnothing$; $\varnothing$ has a
random number  $\wh \xi_{\varnothing}$ of children (with the size-biased offspring distribution) $\wh p$,
   \item   one of the children of $\varnothing$ is selected at random and denoted  $v_1$ and given an
             independent size-biased number $\wh \xi_{v_1}$ of children,
   \item  the other children of $\varnothing$ are independently assigned ordinary BGW descendant trees with offspring number $\xi$,
   \item   again  one of the children of $v_1$ is selected at random and denoted $v_2$ and given an
             independent size-biased number $\wh \xi_{v_2}$ of children,
   \item  this process is continued
 and produces the size-biased BGW tree  $\wh T$ which is {\em immortal} and infinite distinguished path $v$ which we call the {\em backbone}.
    \end{itemize}
\strut


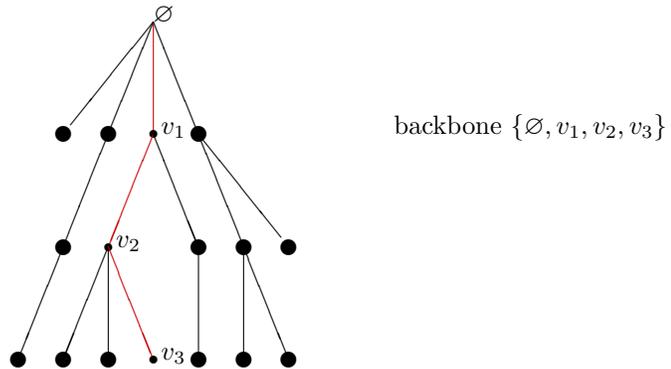
\begin{figure}[h]
\begin{center}
\setlength{\unitlength}{0.01cm}
\begin{picture}(1500,600)
\put(700,550){{$\varnothing$}}
\put(700,550){\line(-4,-5){110}}
\put(700,550){\line(-2,-5){180}}
\put(700,550){{\red\line(0,-1){150}}}
\put(700,550){\line(2,-5){180}}

\put(580,400){{{\circle*{20}}}}
\put(640,400){\circle*{20}}
\put(760,400){{\circle*{20}}}
\put(700,400){{\circle*{10}$v_1$}}
\put(1020,400){backbone $\{\varnothing,v_1,v_2,v_3\}$}
\put(700,400){\red{\line(-2,-5){60}}}
\put(700,400){\line(2,-5){60}}
\put(760,400){\line(4,-5){110}}

 \put(580,250){\circle*{20}}
\put(640,250){\circle*{10}{ $v_2$}}
\put(760,250){\circle*{20}}
\put(820,250){\circle*{20}}
\put(880,250){\circle*{20}}
\put(640,250){{\line(-2,-5){60}}}
 \put(640,250){\line(0,-1){150}}
\put(640,250){\red{\line(2,-5){60}}}
 \put(760,250){\line(0,-1){150}}
\put(820,250){\line(0,-1){150}}
\put(520,100){\circle*{20}}
\put(580,100){\circle*{20}}
\put(640,100){\circle*{20}}
\put(700,100){\circle*{10}{$v_3$}}
\put(760,100){\circle*{20}}
\put(820,100){\circle*{20}}
\put(880,100){\circle*{20}}
\end{picture}
\end{center}
\caption{Size-biased BGW Tree}
\end{figure}

Define the measure $\bar P^{BGW}_*\in\mathcal{P}(\mathbf{T}\times\mathbf{V})$ to be the joint distribution of
the random tree $\wh T$ and backbone
$\{v_0,v_1,v_2,\dots\}$. Let  $\bar P^{BGW}$ denote the marginal distribution of $\wh T$.
We can  view the  vertices off the  backbone
$(v_0,v_1,\dots)$ of the size-biased tree as a branching process
with immigration in which the immigrants are the siblings of the
individuals on the backbone.  The distribution of the number of
immigrants at generation $n$, $Y_n$, is given by the law  $\wh \xi
-1$.

Given a tree $t$ let $[t]_n$ denote the tree restricted to generations $1,\dots,n$.  Let
$Z_n(t)$ denote the number of vertices in the tree  at the nth level (generation) and $\mathcal{F}_n=\sigma([t]_n)$.
Let  \be{} W_n(t):=\frac{Z_n(t)}{m^n}\ee  denote the martingale
associated to a tree $t$ with $Z_n(t)$ vertices at generation  $n$.
\beL{3.RN} (a)  The Radon-Nikodym derivative of the marginal distribution $\bar
P^{BGW}|_{\mathcal{F}_n}$ of $\wh {\mathcal{T}}$   with respect to $P^{PGW}|_{\mathcal{F}_n} $  is given by  \be{} \frac{d\bar
P^{BGW}_n}{dP^{BGW}_n}(t)=W_n(t).\ee
(b)
Under the measure $\bar P^{BGW}_*$, the vertex $v_n$ at the $n$th
level  of the tree $\wh{T}$  in the random path $(v_0,v_1,\dots)$
is uniformly distributed on the vertices at the nth level of
$\wh{T}$.
\end{lemma}
\begin{proof}

We will verify that
\be{RN}\; \bar P_*^{BGW}[t,v]_n= \frac{1}{m^n} P^{BGW}[t]_n\ee
and therefore
\be{RN}\; \bar P^{BGW}[t]_n=W_n(t)P^{BGW}[t]_n.\ee
First observe that the \bea{}\; &&\bar P^{BGW}_*(Z_1=k,v_1=i)= \frac{kp_k}{m}\cdot
\frac{1}{k}\\&&= \frac{p_k}{m}=\frac{1}{m}P(\xi=k),\text{ for } i=1,\dots,k.\nonumber\eea
since $v_1$ is  randomly chosen from the offspring $(1,\dots,\wh {\xi}_{\varnothing})$.

Now consider $[\wh{T},v]_{n+1}$.  We can construct this by first selecting $\wh{\xi}_0$ and $v_1$ and then following the next $n$ generations of the resulting descendant tree and backbone  as well as the BGW descendant trees  of the remaining $\wh{\xi}_0-1$ vertices in the first generation.  If $\wh{\xi}_{\varnothing}(t)= k$  we denote the resulting descendant trees
by $t^{(1)},t^{(2)},\dots,t^{(k)}$.

Let $v_{n+1}(t)$ be a vertex (determined by a position in the lexicographic order) at level
$n+1$. It determines   $v_1(t)$ and    the descendant tree  $t^{(v_1)}$  that it belongs to. If $\wh \xi_{\varnothing}(t)=k,\;v_1(t)=i$,  then we obtain
\be{} \bar
P^{BGW}_*[t;v]_{n+1}=\frac{p_k}{m}\cdot\bar
P^{BGW}_*[t^{(i)};v_{n+1}]_n\cdot \prod_{j=1,j\ne i}^kP^{BGW}[t^{(j)}]_n.\ee Then by
induction for each $n$ \be{} \bar P^{BGW}_*[t;v]_n=\frac{1}{m^n}P^{BGW}[t]_n \ee
for each of the $Z
_n(t)$ positions $v$ in the lexicographic order {  at level }n
and $[t]_n$. Consequently we have obtained the {\em martingale change of measure}
\be{} \bar P^{BGW}_*[t]_n=\frac{Z_n(t)}{m^n}P^{BGW}[t]_n\ee
and
\be{} \bar P^{BGW}_*[v=i|t] =\frac{1}{Z_n(t)}\text{ for }i=1,\dots,Z_n(t).\ee

\end{proof}

For an infinite tree $t$
 we define \be{}  W(t):=\limsup_{n\to\infty} W_n(t).\ee
Note that in the critical and subcritical cases the measures $P^{BGW}$ and $\bar P^{BGW}$ are
singular since the $P^{BGW}$- probability of nonextinction is zero.  The question as to whether
or not they are singular in the supercritical case will  be the focus of the next subsection.



\subsection{Supercritical branching}

As mentioned above if $0< m<\infty$, then under $P^{BGW}$  \be{}W_n =\frac{Z_n}{m^n}\ee
is a martingale and converges to a random variable $W$ a.s. as $n\to\infty$. The
characterization of the limit $W$ in the supercritical case, $m>1$,
under minimal conditions was obtained  in the following theorem of
Kesten and Stigum (1966) \cite{KS-66}. The proof given below
follows the ``conceptual proof'' of Lyons, Pemantle and Peres
(1995) \cite{LPP-95}.

\beT{KS} (Kesten-Stigum (1966)  \cite{KS-66})  Consider the BGW process with offspring $\xi$ and mean offspring size $m$. If $1<m<\infty$, the following are
equivalent

(i) $P^{BGW}[W=0]=q$

(ii) $E^{BGW}[W]=1$

(iii) $E[\xi\log^+ \xi]<\infty$
\end{theorem}

\begin{proof}
By Lemma \ref{3.RN}  \be{} \frac{d\bar
P^{BGW}_n}{dP^{BGW}_n}(t)=W_n(t)\ee where the left side denotes the
Radon-Nikodym derivative wrt
$\mathcal{F}_n=\sigma([t]_n)$.

Note that $P^{BGW}(W=0)\geq q$ where $q=P^{BGW}(E_0)$ where $E_0:=\{Z_n=0\text{  for some }n<\infty\}$ (extinction probability).
Moreover, since $\mathcal{F}_n\uparrow \mathcal{F}=\sigma(t)$,  we have the Radon-Nikodym dichotomy (see  Theorem \ref{RNT})  \be{RN1} W=0,\quad P^{BGW}-a.s. \quad
\Leftrightarrow P^{BGW}\bot \bar P^{BGW}\quad \Leftrightarrow
W=\infty\quad \bar P^{BGW}-a.s.\ee and \be{RN2} \int WdP^{BGW}=1 \quad
\Leftrightarrow \bar P^{BGW} \ll  P^{BGW}\quad \Leftrightarrow
W<\infty\quad \bar P^{BGW}-a.s.\ee

Now recall (\ref{3.x}) that the size-biased tree can be represented as a branching process
with immigration in which the distribution of the number of
immigrants at generation $n$, $Y_n$, is given by the law  $\wh \xi
-1$, that is
\be{}
E[Z_n|\mathcal{Y}]= \sum_{k=1}^n\frac{Y_k}{m^k}.
\ee

If
$E[\log^+ \wh \xi]=\frac{1}{m}E[\xi \log^+ \xi]=\frac{1}{m}\sum_{k=1}^\infty kp_k\log k =\infty$, then \be{}  W=\lim_{n\to\infty} \frac{Z_n}{m^n}  =\infty,\quad \bar{P}^{BGW}a.s.\ee
by Theorem \ref{T.Sen} (b).  Therefore  $P^{BGW}(W=0)=1$
by (\ref{RN1}).

If $E[\xi\log^+ \xi]<\infty$, then $E[\log^+ \wh \xi]=\sum_{k=1}^\infty kp_k\log k <\infty$.
and by  Theorem \ref{T.Sen}(a)
\be{}  \lim_{n\to\infty} E(\frac{Z_n}{m^n}|\mathcal{F}^{\mathcal{Y}}) =\sum_{k=1}^\infty\frac{Y_k}{m^k} <\infty,\quad \bar{P}^{BGW}a.s.\ee
  and therefore
  \be{} W=\lim_{n\to\infty} \frac{Z_n}{m^n} <\infty, \bar P^{BGW}-a.s.\ee
  Then
$E^{BGW}[W]=\int W dP^{BGW} =1$ by (\ref{RN2}).

Finally, since  by Proposition \ref{q0} $ P^{BGW}(W=0)= q\text{  or  }0$, we obtain (i).

\end{proof}
\begin{remark} The supercritical branching model is the basic model for a growing population with unlimited resources.  A more realistic model is a spatial model in which resources are locally limited but the population can grow by spreading spatially. A simple deterministic model of this type is the Fisher-KPP equation.  We will consider the analogous spatial stochastic models in a later chapter.

\end{remark}

\subsection{The general branching model of Crump-Mode-Jagers}

We now consider a far-reaching generalization of the Bienamy\'e-Galton-Watson
process known as a Crump-Mode Jagers (CMJ) process (\cite{CM-68},
\cite{J-75}). This is a process with time parameter set
$[0,\infty)$ consisting of finitely many individuals at each time.

 With each individual $x$ we denote its birth time $\tau_x$, lifetime $\lambda_x$  and
 reproduction process $\xi_x$. The latter is a point process which gives
 the sequence of birth times of individuals.  $\xi_x(t)$ is the number of offspring
 produced (during its lifetime)  by an individual $x$ born at time 0 during $[0,t]$. The intensity of
 $\xi_x$, called the {\em reproduction function}  is defined by
 \be{} \mu(t)=E[\xi(t)].\ee

 The lifetime distribution is defined by
 \be{}  L(u)=P[\lambda \leq u].\ee

We begin with one individual $\emptyset$ which we assume is born
at time $\tau_{\emptyset}=0$.  The reproduction processes $\xi_x$
of different individuals are iid copies of
 $\xi$.

The basic probability space is \be{}
(\Omega_{\mathcal{I}},\mathcal{B}_{\mathcal{I}},P_{\mathcal{I}})=\prod_{x\in\mathcal{I}}
(\Omega_x,\mathcal{B}_x,P_x)\ee where $\mathcal{I}$ is given as in (\ref{ind}) and
$(\xi_x,\lambda_x)$ are random variables defined on $(\Omega_x,\mathcal{B}_x,P_x)$ with
distribution as  above.

We then determine the birth times $\{\tau_x, \;x\in\mathcal{I}\}$
as follows:

 \bea{} && \tau_{\emptyset}=0,
 \\ && \tau_{(x',i)}=\tau_{x'}+\inf{ \{u:\xi_{x'}(u)\geq i\}}.\nonumber \eea
Note that for individuals never born $\tau_x=\infty$.

Let \be{ZI} Z_t= \sum_{x\in\mathcal{I}} 1_{\tau_x\leq
t<\lambda_x},\quad T_t= \sum_{x\in\mathcal{I}} 1_{\tau_x\leq t}\ee
that is, the {\em number of individuals alive at time $t$} and {\em total
number of births before time $t$}, respectively.

For $\lambda >0$ we define
 \be{} \xi^\lambda(t) :=\int_0^t e^{-\lambda u}\xi(du).\ee

 The {\em Malthusian parameter} $\alpha$ is defined by the equation
 \be{}  E[\xi^\alpha(\infty) ]=1\ee
 that is,
\be{C1}  \int_0^\infty e^{-\alpha t}\mu(dt)=1.\ee

The {\em stable average age of child-bearing} is defined as

\be{C2} \beta= \int_0^\infty t \wt{\mu}(dt)  \text{   where   }\wt{\mu}(dt)= e^{-\alpha t}\mu(dt).\ee

\begin{example}
Consider  a population in which individuals have an internal state space, say $\N$.  Assume that the individual starts in state $0$ at its time of birth and and its internal state changes according to a Markov transition mechanism.  Finally assume that when it is in state $i$ it produces offspring at rate $\lambda_i$.
\end{example}

\begin{definition}  {\em Characteristics of an individual:} A {\em characteristic} of an individual is given by a  process
 $\phi:\mathbb{R}\times \Omega\to \mathbb{R}_+$ which is given by a
 $\mathcal{B}(\mathbb{R})\times
\sigma(\xi)$-measurable non-negative function satisfying  $\phi(t)=0$ for $t<0$, let
\be{}  Z^\phi_t=\sum_{x\in\mathcal{I}}\phi_x(t-\tau_x)\ee
denote the process {\em counted with  characteristic} $\phi$.
\end{definition}

\begin{example}\label{age}
If  $\phi^a(t)=1_{[0,\inf(a,\lambda))}(t)$, then $Z^{\phi^a}_t$ counts the number of individuals alive at time $t$ whose ages are less than $a$.

\end{example}

The following fundamental  generalization of the
Kesten-Stigum theorem was developed in  papers of Doney (1972),(1976) \cite{Do-72}, \cite{Do-76},
and Nerman (1981) \cite{N-81}.

\beT{Doney-Nerman} Consider a CMJ process with malthusian parameter $\alpha $ and assume that
$\beta <\infty$.

(a) \cite{Do-72} Then as
$t\to \infty$, $e^{-\alpha t}Z_t$ converges in distribution to  $m W_\infty$
where
\be{} m=\frac{\int_0^\infty e^{-\alpha s}(1-L(s))ds}{\beta}
\ee and  $W_\infty$ is a random variable (see Proposition \ref{N-Mart}) and

(b) The following are equivalent: \bea{}
&&E[\xi^\alpha(\infty)\log^+\xi^\alpha(\infty)]<\infty\\&& E[W]>0\\&& E[e^{-\alpha
t}Z_t]\to E[W]\quad\text{as  } t\to\infty\\&& W>0 \text{  a.s. on
} \{T_t\to \infty\}.\eea

(c)  \cite{N-81} Under the condition that there
exists a non-increasing integrable function $g$ such that

\be{N-5.3}
E[\sup_t \frac{(\xi^\alpha(\infty)-\xi^\alpha(t))}{g(t)}]<\infty,\ee then
$e^{-\alpha t}Z_t$ converges a.s. as $t\to\infty$.

\end{theorem}

\begin{remark} A sufficient condition is the existence of non-increasing integrable function $g$ such that
\be{N-5.4} \int_0^\infty \frac{1}{g(t)}e^{-\alpha t}\mu(dt)<\infty.\ee
(See Nerman \cite{N-81} (5.4)).
\end{remark}

\textbf{Comments on Proofs}

(b) The equivalence statements can be proved in this general case following the same lines as that of Lyons, Pemantle and Peres - see Olofsson (1996) \cite{O-96}.

(a) - convergence in distribution was proved by Doney (1972) \cite{Do-72}. However the almost sure convergence  required some basic new ideas since  we can no longer directly use the martingale
convergence theorem  since $Z_t$ is not a martingale in the general case.  The a.s. convergence was proved by Nerman \cite{N-81}.
 We will not give Nerman's  long
detailed technical proof of this result but will now introduce the
key tool used in its proof and which is of independent interest,
namely, an underlying  intrinsic martingale $W_t$ discovered by
Nerman \cite{N-81} and then give an intuitive idea of the remainder of the proof.

Denote the mother of $x$ by $m(x)$  and let
\be{}\mathcal{I}_t=\{x\in\mathcal{I}:\tau_{m(x)}\leq
t<\tau_x<\infty\},\ee the set of individuals whose mothers are
born before time $t$ but who themselves are born after $t$

Consider the individuals ordered by their times of birth

\be{} 0=\tau_{x_1}\leq \tau_{x_2}\leq \dots\ee
Define  $\mathcal{A}_n=\sigma\text{-algebra generated by }
\{(\tau_{x_i},\xi_{x_i},\lambda_{x_i}):i=1,\dots,n\}$
Recall (\ref{ZI}) and let $\mathcal{F}_t= \mathcal{A}_{T_t}$.

Define \be{} W_t:=\sum_{x\in\mathcal{I}_t}e^{-\alpha \tau_x}.\ee

 \beP{N-Mart} (Nerman (1981) \cite{N-81}) (a) The process
$\{W_t,\, \mathcal{F}_t\}$ is a non-negative martingale with
$E[W_t]=1$.

(b) There exists a random variable $W_\infty <\infty$ such that
$W_t\to W_\infty$ a.s. as $t\to\infty$. \end{proposition}

\begin{proof} Define \bea{}   && R_0=1,\\
&& R_n=1+\sum_{i=1}^n e^{-\alpha
\tau_{x_i}}(\xi^\alpha_{x_i}(\infty)-1),\quad n=1,2,\dots
\nonumber \eea

Equivalently, letting $\tau_{(x_i,k)}$ denote the time of birth of the kth offspring of $x_i$,  \be{}
R_n=1+\sum_{i=1}^n\sum_{k=1}^{\xi_{x_i}(\infty)} e^{-\alpha
\tau_{(x_i,k)}}  -\sum_{i=1}^n e^{-\alpha \tau_{x_i}}\ee so that
$R_n$ is a weighted (weights $e^{-\alpha \tau_x}$) sum of
children of the first $n$ individuals.

We next show that $(R_n, \mathcal{A}_n)$ is a non-negative
martingale. $R_n$ and $\tau_{x_{n+1}}$ are
$\mathcal{A}_n$-measurable and $\xi^\alpha_{x_{n+1}}$ is
independent of $\mathcal{A}_n$ and

\be{} E[\xi^\alpha_{x_{n+1}}(\infty)] =\mu_\alpha(\infty)=1.\ee
Therefore \be{} E[R_{n+1}-R_n]=
e^{-\alpha\tau_{x_{n+1}}}E[\xi^{\alpha}_{x_{n+1}}-1]=0.\ee

Next we observe that since $\mathcal{I}(t)$ consists of exactly
the children of the first $T_t$ individuals to be born after $t$,
it follows that  $W_t= R_{T_t}$.

Note that for fixed $t$, $\{T_t > k\}=\{\tau_{x_n}\leq t\}\in\mathcal{A}_n$ and therefore $T_t$ is an increasing family of integer-valued stopping times with
respect to  $\{\mathcal{A}_n\}$. Therefore  $\{W_t\}$ is a supermartingale with respect to the filtration
$\{\mathcal{A}_{T_t}\}$.

Since $E[T_t]<\infty$ and
\be{}  E[|R_{n+1}-R_n|\;|\mathcal{A}_n]= e^{-\alpha
\tau_{x_{n+1}}}E[|\xi^\alpha(\infty)-1|]\leq 2.\ee
 a standard argument  (e.g.
Breiman \cite{B-68} Prop. 5.33) implies that $E[W_t]=E[R_{T_t}]=1$  and $\{W_t\}$ is actually a martingale.

(b) This follows from (a) and the martingale convergence theorem.

\end{proof}

\begin{remark} We now sketch an intuitive explanation for the proof of the a.s. convergence of $e^{-\alpha t}Z_t$ using Proposition \ref{N-Mart}. This is based on the relation between $W_t$ and $Z_t$ which is somewhat
indirect.
 To give some idea of this, let
 \be{} W_{t,c} =\sum_{x\in\mathcal{I}_{t,c}}e^{-\alpha
\tau_x},\ee where \be{}\mathcal{I}_{t,c}=\{x=(x',i):\tau_{x'}\leq
t,\; t+c< \tau_x<\infty\}.\ee
 Note that if we consider the characteristic $\chi^c$ defined by
 \be{} \chi^c(s)=(\xi^\alpha(\infty)-\xi^\alpha(s+c))e^{\alpha
s}\text{ for
 } s\geq 0,\ee
 then
 \be{ }W_{t,c} =e^{-\alpha t}Z^{\chi^c}_t\ee
where
 \be{}  Z^{\chi^c}_t=\sum_{x\in \mathcal{I}} \chi^c_x(t-\tau_x),\quad \chi^c_x(s)=(\xi_x^\alpha(\infty)-\xi_x^\alpha(s+c))e^{\alpha
s}.\ee

Note that $\lim_{c\to 0} W_{t,c}=W_t$ and $\lim_{c\to 0} Z^{\chi^c}_t= Z^\chi_t$
where
\be{}\chi(s)=\int_s^\infty e^{-\alpha (u-s)}\xi(du).\ee

Then $Z^\chi_t\to m_\chi W_\infty,\; a.s.$ where
\be{}
m_\chi =\frac{\int_0^\infty e^{-\alpha t}(1-L(t))dt}{\beta}.
\ee

 In the special case where   $\xi$ is stationary then the distribution of $\chi(s)$ does not depend on $s$. Then $Z_t^\chi$ is a sum of $Z_t$ i.i.d. random variables and therefore as $t\to\infty$, $Z^\chi_t$ should approach  a constant times $Z_t$ by the law of large numbers.
\end{remark}

\subsubsection{Stable age distribution}

The notion of the stable age distribution of a population is a basic concept in demography going back to Euler.  The stable age distribution in the deterministic setting of the Euler-Lotka equation (\ref{EL}) is

\be{}  U(\infty,ds)= \frac{(1-L(s))e^{-\alpha s}ds}{\int_0^\infty(1-L(s))e^{-\alpha s}ds}.\ee

 It was introduced into the study of branching processes by  Athreya and Kaplan (1976)
\cite{AK-76}.
 Let  $Z^a_t$ denote the number of individuals of age $\leq a$. The normalized age distribution at time $t$ is defined by
\be{}
U(t,[0,a)):= \frac{Z^a_t}{Z_t},\quad a\geq 0.
\ee



\begin{theorem}

(Nerman \cite{N-81} Theorem 6.3 - Convergence to stable age distribution)  Assume that $\xi$ satisfies the conditions of Theorem \ref{Doney-Nerman}.
  Then on the event $T_t\to\infty$,
\be{}  U(t,[0,a)) \to \frac{\int_0^a(1-L(u))e^{-\alpha u}du}{\int_0^\infty(1-L(u))e^{-\alpha u}du} \text{   a.s.  as
}t\to\infty.\ee
\end{theorem}

\subsection{Multitype branching}

A central idea in evolutionary biology is the differential growth
rates of different types of individuals. Multitype branching
processes  provide a starting point for our discussion of this
basic topic.

Consider a multitype BGW process with $K$ types. Let
$\xi^{(i,j)}$ be a random variable representing the number of
particles of type $j$ produced by one type $i$ particle in one
generation.

 Let
$Z^{(j)}$ be the number of particles of type $j$ in generation $n$
and ${\bf Z}_n:=(Z_n^{(1)},\dots,Z_n^{K})$.

For ${\bf k}=(k_1,\dots,k_K)$, let  $p^{(i)}_{\bf k}=P[
\xi^{(i,j)}=k_j,\;j=1,\dots,K]$. Assume that \bea{} &&{\bf
M}=(m_{(i,j)})_{i,j=1,\dots,K},\\&& m_{(i,j)}=E[\xi^{(i,j)}] <\infty \quad \forall\; i,j.\nonumber\eea  Then

\be{}  E(\mathbf{Z}_{m+n}|\mathbf{Z}_m)= \mathbf{Z}_m {\bf
M}^n,\quad m,n\in\mathbb{N}.\ee

The behaviour of $E[\mathbf{Z}_n]$ as $n\to\infty$ is then obtained
from the classical Perron-Frobenius Theorem:

\beT{Per-Frob}(Perron-Frobenius) Let ${\bf M}$ be a nonnegative
$K\times K$ matrix.  Assume that ${\bf M}^n$ is strictly positive
for some $n\in\mathbb{N}$.  Then ${\bf M}$ has a largest positive
eigenvalue $\rho$ which is a simple eigenvalue with
positive right and left normalized eigenvectors ${\bf u}=(u_i)$
($\sum u_i=1$) and ${\bf v}=(v_i)$ which are the only nonnegative
eigenvectors. Moreover

\be{}  {\bf M}^n=\rho^n{\bf M}_1+{\bf M}_2^n\ee where ${\bf
M}_1=(u_iv_j)_{i,j\in\{1,\dots,K\}}$ normalized by $\sum_i,j u_i v_j=1$. Moreover ${\bf
M}_1{\bf M}_2= {\bf M}_2{\bf M}_1=0,\; {\bf M}_1^n={\bf M}_1$.

Finally, \be{} |M^n_2|=O(\alpha^n)\ee for some $0<\alpha <\rho$.

\end{theorem}

The analogue of the Kesten-Stigum theorem stated above is given as
follows.

\begin{theorem} (Kesten-Stigum (1966) \cite{KS-66}), (Kurtz,
Lyons, Pemantle and Peres (1997) \cite{KLPP-97})\\  (a) There is a
scalar random variable $W$ such that \be{} \lim_{n\to\infty}
\frac{{\bf Z}_n}{\rho^n}=W{\bf{u}}\;\; a.s.\ee and $P[W>0] >0$ iff
\be{} E[\sum_{i,j=1}^J \xi^{(i,j)}\log^+ \xi^{(i,j)}]<\infty. \ee

(b) Almost surely, conditioned on nonextinction, \be{}
\lim_{n\to\infty}\frac{{\bf Z}_n}{|{\bf Z}_n|}={\bf{u}}. \ee

\end{theorem}

\section{Multilevel branching}

Consider a {\em host-parasite population} in which the individuals in the host population reproduce by BGW branching and the population of parasite on a given host also develop by an independent BGW branching.  This is an example of a {\em multilevel branching system}.

A {\em multilevel population system} is a hierarchically structured  collection
of objects at different levels as follows:\\
 ${ E_0} $ denotes the set of possible types of level 1 object, \\
for $n\geq 1$ each level $(n+1)$ object is given by a collection of level $n$
object including their their multiplicities.\\

\noindent \underline {Multilevel branching dynamics }   \\
Consider a continuous time branching process such that
\begin{itemize}
\item for $n\geq 1$, when a level $n$ object branches,
all its offspring are copies of it
\item  if $n\geq 2$, then
the offspring contains a copy of the set of level-$n-1$ objects contained in the parent level $n$ object.
\item let $\gamma_n$
the level n branching rate and by $f_n(s)$ the level n
offspring generating function.
\end{itemize}
Then the questions of extinction, classification into critical, subcritical and supercritical case  and growth asymptotics in
the supercritical case are more complex than the single level branching case. See for example, Dawson and Wu (1996) \cite{DW-96}.

\chapter[Branching Processes II]{Branching Processes II: \\Convergence of critical branching to Feller's CSB}


\section{Birth and Death Processes}

\subsection{Linear birth and death processes}
Branching processes can be studied in discrete or
continuous time. We now consider a classical continuous time version. This
is a continuous time Markov chain, $\{X_{t}\}_{t\geq0}$ with
state space $\N_0$ and with linear birth and death
rates, $b$ and $d$ and let $V=b+d\geq 0$.  This corresponds to
a branching system in which (independently) each  particle can die or be replaced by two offspring in the
interval $[t,t+\Delta t)$ with probability $V\Delta t +o(\Delta t)$. This means  that
the time until the first branch (birth-death event) is an exponential random variable
with mean $\frac{1}{V}$. $V$ is called the branching rate. When
the particle ``branches'' it dies with probability $\frac{d}{b+d}$ and
is replaced by two descendants with probability $\frac{b}{b+d}.$ Note
that this process can be built directly on a probability space
containing a sequence of iid exponential (1) rv's and a sequence
of iid Bernoulli ($p=\frac{b}{b+d}$) rv's (or a sequence of iid Uniform
$[0,1]$ rv's) and this description can be used to generate a
simulation of the model.
The special case in which $d =0$ is called the {\em Yule} process.

The birth and death process can also be realized on a probability space
$(\Omega,\mathcal{F},P)$
 on which independent Poisson random measures $N_1,N_2$ on $\mathbb{R}^2_+$ are defined. Then the birth and death process  is defined via a stochastic differential equation
driven by the Poisson noises, namely,

\bea{bdsde}&&\\\;\; X_t=x_0+\int _0^t\int_0^{b X(s_-)} N_1(du,ds) -
\int _0^t \int_0^{d X(s-) }N_2(du,ds).\nonumber\eea This equation has
a pathwise unique c\` adl\`ag solution which is a continuous time Markov chain
with the required transition rates.
 See  Li-Ma \cite{LM-08}.

 Let  $P_{x_0}$ denote the resulting probability law on $D_{\N_0}([0,\infty))$, the space of c\`adl\`ag functions from $[0,\infty)$ to $\N_0$.

\subsection{Semigroups and generating functions}

Given the Markov process $X_t$ we can associate a Markov semigroup $\{T_{t}:t\geq0\}$ of operators on the Banach space $C_{0}%
(\N_0)$ (the space of bounded functions on $\N_0$, with limits at infinity) as follows:
\[
T_{t}f(x_0):=E_{x_0}(f(X_{t}))=\int f(x) P_{x_0}(X_t\in dx).
\]
This semigroup determines the finite dimensional distributions of the Markov
chain. This semigroup satisfies the conditions of the Hille-Yosida theorem with
generator given by%
\begin{align*}
Gf(n)  &  =\frac{dT_{t}f(n)}{dt}|_{t=0}\\
&  = b n(f(n+1)-f(n))+ dn(f(n-1)-f(n))
\end{align*}

Now consider the {\em Laplace function} of $X_t$ starting with one
particle at time $0$:
\[
L(t,\theta):= E_{x_0}(e^{-\theta X_{t}}),\;\text{   with  }x_0=1,\;\theta\geq0
\]
 Noting the outcome  at the first branching time and using the independence of
the particle and its offspring when a birth occurs, we obtain
the nonlinear renewal-type equation
\[
L(t,\theta)=e^{-Vt}e^{-\theta}+\frac{d}{b+d}(1-e^{-Vt})+V\frac{b}{b+d}\int_{0}%
^{t}e^{-Vu}L^{2}(t-u,\theta)du
\]
Alternately, note that we can represent the jump in $X_t$ at a branching time by  the addition of an
independent random variable $\zeta$  with Laplace transform  $E(e^{-\theta \zeta})=\frac{b}{b+d}
e^{-\theta}+\frac{d}{b+d}e^{\theta}$. Since the branching occurs at linear rate $VX_{t}$ at time $t$, we get
\begin{align*}
\frac{\partial L(t,\theta)}{\partial t}  &  =\lim_{\Delta\rightarrow
0}\frac{L(t+\Delta,\theta)-L(t,\theta)}{\Delta}\\
&  =\lim_{\Delta\rightarrow0}\left(\frac{E(E(e^{-\theta X_{t+\Delta}}|X
_{t-}))-E(e^{-\theta X_{t-}})}{\Delta}\right)\\
&  =\lim_{\Delta\rightarrow0}\frac{V\Delta\lbrack
 E(X_{t-}E(e^{-\theta(X_{t-}+\zeta)}|X_{t-}))-E(X_{t-} e^{-\theta X_{t-}})\rbrack +o(\Delta)}{\Delta}\\
&
=-\{V\frac{b}{b+d}(e^{-\theta}-1)+V\frac{d}{b+d}(e^{\theta}-1)\}\frac{\partial
L(t,\theta)}{\partial\theta}.%
\end{align*}
Here we have used $E(Xe^{-\theta X})= -\frac{\partial L(\theta)}{\partial \theta}$.
So we then have the first order PDE
\be{GPDE}
\frac{\partial L(t,\theta)}{\partial t}+ V[\frac{b}{b+d}(e^{-\theta
}-1)+\frac{d}{b+d}(e^{\theta}-1)]\frac{\partial L(t,\theta
)}{\partial\theta}=0,\quad L(0,\theta)=e^{-\theta}.
\ee

We can solve this by finding the  {\em characteristic curves} $(t(s),\theta(s))\;$ in the $(t,\theta)$ plane along which $L(t(s),\theta(s))$ is constant (refer to Garabedian (1964) \cite{Ga-64}, John (1982) \cite{J-82} or Delgado (1997) \cite{De-97}). We write this as%
\be{CE}
\frac{\partial}{\partial s}L(t(s),\theta(s))=L_{1}\frac{\partial
t(s)}{\partial s}+L_{2}\frac{\partial\theta(s)}{\partial s}=0
\ee
where $L_1,L_2$ denote the first partial derivatives with respect to $t,\theta$ respectively.
Comparing (\ref{CE}) with  (\ref{GPDE}) leads to the {\em characteristic equations}%
\be{}
\frac{\partial L(s)}{\partial s}=0,\;\;\;\frac{\partial
t(s)}{\partial s}=1,\;\;\;\frac{\partial\theta(s)}{\partial
s}=h(\theta)=b(e^{-\theta
}-1)+d(e^{\theta}-1)
\ee
For $b\ne d$ we obtain the characteristic curve
\be{} \frac{(e^{-\theta}-1)e^{(b-d)t}}{b e^{-\theta}-d}=\text{ constant}.\ee
and general solution
\be{} L(t,\theta)=\Psi( \frac{(e^{-\theta}-1)e^{(b-d)t}}{b e^{-\theta}-d})\ee
where $\Psi$ is a differentiable function.  From the initial condition we have
\be{}  \Psi(\frac{e^{-\theta}-1}{be^{-\theta}-d})=e^{-\theta X_0}.\ee

Solving for $\Psi$ we obtain  for $b\ne d$ the solution
\be{cfgo} L(t,\theta)= \left(\frac{d(e^{-\theta}-1)e^{(b-d)t}-(be^{-\theta}-d)}{b(e^{-\theta}-1)e^{(b-d)t}-(be^{-\theta}-d)}\right)^{X_0}
\ee
and for $b=d$

\be{cgf}  L(t,\theta)=
\left( \frac{1-(bt-1)(e^{-\theta}-1)}{1-bt(e^{-\theta}-1)}\right)^{X_0}.\ee

\begin{remark} Note that the form of the Laplace transforms (\ref{cfgo}), (\ref{cgf})implies
the {\em branching property}, namely, if $X_0=X_{0,1}+X_{0,2}$, then the probability law of $X_t$ is identical to the distribution of the sum of independent random variables $X_{t,1} +X_{t,2}$ where $X_{t,i}$ are versions of the linear birth and death process with initial conditions $X_{0,1}, X_{0,2}$.

\end{remark}

\subsubsection{Distribution function, moments, extinction probability}  Setting $b,\;d$ as the birth and death rates.  Then replacing $\theta$ by $-\ln z$ in $L_t(\theta)$ we obtain the probability generating function
\be{} G_t(z) =L(t,-\ln z)=\sum_{k=0}^\infty z^k p_k(t).\ee
Then expanding in a power series in $z$ we can obtain the standard formula
\be{bdp1}  p_0(t)= f(t),\ee
\be{bdp2} p_n(t)= (1-f(t))(1-g(t))g(t)^{n-1},\;n\geq 1\ee
where
\be{}  f(t)=\frac{d(e^{(b-d)t}-1)}{be^{(b-d)t}-d},\quad g(t)= \frac{b(e^{(b-d)t}-1)}{be^{(b-d)t}-d}.\ee
Similarly if $b=d=\frac{V}{2}$, then
\be{} p_n(t)=\frac{(bt)^{n-1}}{(1+bt)^{n+1}},\;n\geq 1,\ee
\be{criticalp0} p_0(t)=\frac{bt}{1+bt}.\ee

Then the extinction probability is
\be{extprob} \lim_{t\to\infty} p_0(t)=\lim_{t\to\infty} \frac{d(e^{(b-d)t}-1)}{be^{(b-d)t}-d}=\left\{\begin{split} &1\quad\text{if  } b\leq d,\\ &\frac{d}{b}\quad\text{if  }b>d.\end{split}\right.
\ee

Recalling that
\be{} E(X_t)= -\left. \frac{\partial L_t(\theta)}{d\theta}\right|_{\theta=0},\ee

\be{} E((X_t)^2)= \left. \frac{\partial^2 L_t(\theta)}{d\theta^2}\right|_{\theta=0}.\ee
we can obtain

\be{} E(X_t)= X_0e^{(b-d)t},\ee
\be{BD2M}   E((X_t)^2)= (X_0)^2e^{2(b-d)t}+ \frac{X_0(b+d)}{b-d}e^{(b-d)t}(e^{(b-d)t}-1),\; b\ne d\ee
\be{}E((X_t)^2)= (X_0)^2+2bt,\; \text{  if  }b=d.\ee




\section{Critical branching}

Exponential growth of a population is unrealistic and therefore
supercritical branching models describe only the growth of a
population as long as the resources are unlimited. Otherwise logistic competition comes into play.
  We will return
to this circle of questions throughout this course.

Only critical branching processes have the property that the mean
population size is stable but as shown above the critical
branching process actually suffers extinction with probability
one.  Nevertheless critical branching processes have played a key
role in the development of stochastic population models. We will later see that a key feature of critical branching
is the limiting behavior of the process conditioned on non-extinction up to time $t$ and letting  $t\to\infty$.  We now give two formulations
of the resulting behavior.

\begin{theorem}  Consider the BGW process $Z_n$ with mean offspring size $m=1$.  Suppose that
$\sigma^2:=\rm{Var}(\xi)=E[\xi^2]-1\leq \infty$.  Then

(i) Kolmogorov \be{} \lim_{n\to\infty}
nP[Z_n>0]=\frac{2}{\sigma^2} \ee

(ii) Yaglom:  If $\sigma <\infty$, then the conditional distribution
of $\frac{Z_n}{n}$ given $Z_n>0$ converges as $n\to\infty$ to an
exponential law with mean $\frac{\sigma^2}{2}$.
\end{theorem}

We refer the proof of this to the literature \cite{AN-72},
\cite{LPP-95}.

\beT{} Consider the
critical linear birth and death process, $\{X_t\}$ with
$\alpha=\frac{1}{2}$, $b=d=\frac{V}{2}$. Then

(i)  Extinction probability: $\lim_{t\rightarrow\infty}p_{0}(t)
 =1$.
 \medskip

 (ii) Expected extinction time: Let $\tau:= \inf\{t:X_t =0\}$ Then $E[\tau]=\infty$.

\medskip

(iii)  Exponential limit law: conditioned on $X_t\ne 0$,
\[ \frac{X_t}{t} \Rightarrow Y\]
where $Y$ is exponential with mean $b$.

\end{theorem}

\begin{proof}
The proof is based on the explicit form of the generating function (\ref{cgf}).

(i) From (\ref{criticalp0}), $p_0(t)=\frac{bt}{bt+1}\to 1$ as $t\to \infty$.

(ii) The expected extinction time is infinite%
\begin{align*}
E(\tau)  &  =\int_{0}^{\infty}(1-p_{0}(t))dt=\int_{0}^{\infty}\frac{1}%
{bt+1}dt=\infty.\\
\end{align*}

(iii) From (\ref{cgf}), \bea{} &&\\&& E(e^{-\frac{X_t\theta}{t}}|X_t\ne
0)=\frac{L(t,\frac{\theta}{t})-P(X_t=0)}{1-P(X_t=0)}\nonumber\\&&
=
\frac{1-(bt-1)(e^{-\frac{\theta}{t}-1)}{1-bt(e^{-\frac{\theta}{t}-1})}
-\frac{bt}{bt+1}}
{\frac{1}{bt+1}}\nonumber
\\&& \lim_{t\to\infty}E(e^{-\frac{X_t\theta}{t}}|X_t\ne
0)  =\frac{1}{1+ b{\theta}} \nonumber\eea and
which is the Laplace transform of the  exponential distribution with mean $b$.
\end{proof}

\section[Feller's continuous state branching]{Feller's continuous state branching process (CSBP)}

Consider  the It\^o stochastic differential equation (SDE)
\[
dX_{t}=mX_{t}dt+{ \sqrt{\gamma X_{t}}}dW_{t},\quad X_0=x\geq 0%
\]
where $\{W_{t}\}$ is a standard Brownian motion.
This equation has a non-Lipschitz coefficient but its pathwise uniqueness
follows from the  Yamada-Watanabe theorem \cite{YW-71}.

Using It\^o's lemma one can then check that the generator of the resulting diffusion process
 acting on $D(G)= \{f\in C_0^2(\mathbb{R}_+),\; xf_x,\;xf_{xx}\in C_0(\mathbb{R}_+^d)\}$ satisfies
\be{FCSBG} Gf(x)= mx\frac{\partial f}{\partial x} +\frac{1}{2}\gamma
x\frac{\partial^2 f}{\partial x^2}\ee
and therefore $X_t$ is a realization of the Feller CSBP process.

\beP{FCSBLT} (Laplace transform and extinction probability)\\
(a)The  Laplace transform is given by
\be{}
L(\theta,t)=\mathbb{E}_x\exp(-\theta X_t)=\exp(-u(t)x)\ee
where $u(s)$ satisfies the equation:
\begin{equation}
\frac{\partial u}{\partial s}=m u-\frac{\gamma}{2}u^2\quad
u(0)=\theta.
\end{equation}

(b) In the critical case $m=0$
\be{EXT} P_{x}(x_t =0)=\exp\left(-\frac{x}{
\gamma t}\right).\ee

\end{proposition}
\begin{proof} Assume that $\theta(s)\geq 0$ is differentiable. Then applying  It\^o's lemma (\cite{RY-91}, Theorem 3.3, Remark 1) to $F(\theta,x)=e^{-\theta x}$ we have
\bea{} &&F(\theta(t),X_t)-F(\theta(0),X_0)\\
&& = m\int_0^t X_s F_2(\theta(s),X_s)ds + \int_0^t F_2(\theta(s),X_s)\sqrt{\gamma X_s}dW_s\nonumber\\&&
+\int_0^t F_1(\theta(s),X_s)d\theta(s) +\frac{\gamma}{2}\int_0^t X_s F_{22}(\theta(s),X_s)ds\nonumber
\eea
Noting that $E(X_se^{-\theta X_s})=-L_1(\theta,s)$ we obtain
\bea{}&&
\frac{\partial L(\theta(s),s)}{\partial s}=L_1(\theta(s),s)\frac{d\theta(s)}{ds}-m\theta(s)
L_1(\theta(s),s)+\frac{\gamma}{2}\theta(s)^2
L_1(\theta(s),s)\\&& \quad\text{with}\quad L(\theta,0)=e^{-\theta x}\nonumber
\eea
If $ u$ is a solution of
\be{lle}
\frac{\partial u(\theta,s)}{\partial s}= m
u(\theta,s)-\frac{\gamma}{2}u^2(\theta,s),\qquad u(\theta,0)=\theta
\ee
then the derivative with respect to $s$   \begin{equation*}
\frac{\partial}{\partial s} L(u(\theta,t-s),s)=0,\quad 0\leq s\leq t
\end{equation*}
and therefore
\be{LT}
\mathbb{E}_x(e^{-\theta
X_t})=L(\theta,t)=L(u(\theta,t),0)=e^{-u(\theta,t)x}.
\ee

(b) Solving (\ref{lle})  we get
\be{}
u(\theta,t)=\frac{\theta}{(1+t\gamma\theta)},\text{  if  }m=0\ee
\be{LTI}
u(\theta,t)=\frac{\theta me^{mt}}{m+\gamma \theta(e^{mt}-1)},\text{  if  }m\ne 0.
\ee
If $m=0$
\be{} P_{X_0}(x_t =0)= \lim_{\theta\to\infty} e^{-x_0 u(\theta,t)}=e^{-\frac{x_0}{
\gamma t}}.\ee

\end{proof}

\begin{remark}
An immediate consequence of (\ref{LT}) is that for each $t$,  $X_t$ is an infinitely divisible random variable. In fact the law of $X_t$ corresponds to the law of the sum of a Poisson distributed number of independent exponential random variables.  These facts will provide an important tool for the study of these processes and their infinite dimensional generalizations.
\end{remark}
\bigskip
\subsubsection{Feller CSBP with immigration}

 Adding an immigration term $c t$ to $X_t$, one obtains
the continuous state branching with immigration process (CBI), and
can verify (see e.g. Li (2006) \cite{Li-06})  the following:
\medskip
\medskip
\beP{} Consider the continuous subcritical branching process with
immigration (CBI), given by the SDE:
\begin{equation}
 dY_t=cdt - bY_tdt+\sqrt{\gamma Y_t}dW_t, \quad
 Y_0=y_0,\quad b,c>0.
\end{equation}

(a) The Laplace transform of the distribution of $Y_t$ is given by:
\be{FLTI}
\mathbb{E}_{y_0}\exp(-\theta Y(t))=e^{-y_0u(t) -\int_0^t c u(s)\,ds};\\
\frac{\partial u}{\partial t}=-b u- \frac{\gamma}{2}u^2\quad
u(0)=\theta >0.
\ee

(b) In the subcritical case  $Y_t$ converges to equilibrium, $Y_t\Rightarrow Y_\infty$ as $t\to\infty$, where $Y_\infty$ has the {\em gamma distribution}
with Laplace transform
\be{} L(\theta)= \frac{c}{[(b+\gamma\theta)}.\ee

\end{proposition}
\begin{proof}
(a)  This can be proved using the method of Theorem \ref{FCSBLT}. Alternately, we can prove this by consider the process with immigrants coming according to  $\frac{1}{K}\sum \delta_{y_i}$ where $\{y_i\}$ are the points of a Poisson process with rate $K$  and letting $K\to\infty$.

(b) We obtain

\be{} L_t(\theta)= \frac{c}{[(b+\gamma\theta)-\gamma\theta e^{-bt}]^{1/\gamma}}.\ee
   from (a) by  simple integration  of (\ref{FLTI}).  (b) then follows by taking $t\to\infty$.
\end{proof}

\medskip

\begin{remark}
The critical Feller CSBP with immigration

\be{CSBI}
\begin{split}
 dY_t&=\beta dt+2 \sqrt{Y_t}dW_t\\
 Y_0&=y_0
 \end{split}
\ee
is the square of a $\beta$-dimensional Bessel process. (See Revuz
Yor \cite{RY-91} where this is called a BESQ$^\beta$ process).  For $\beta \geq
2$, $\{0\}$ is polar.  For $0<\beta<2$, $\{0\}$ is instantaneously
reflecting.  For $0<\beta <1$ the set $\{t:X_t=0\}$ is a perfect
set. (See Revuz Yor \cite{RY-91} Chap. XI.)
\end{remark}

\section[Diffusion process limits]{Diffusion limits of critical and nearly critical branching processes}

\subsection{Convergence to Feller's continuous state branching process}

In a celebrated paper Feller (1951) \cite{F-51} developed the
diffusion approximation to branching processes using semigroup
methods.

\begin{theorem} ({Convergence of B+D and BGW processes to Feller CSBP})\newline

(a) Consider the sequence of birth and death process, $\{X^K_t\}, \;K\in\N$, with linear birth and death rates $b_K=1+\frac{m}{2K}, \; d_K=1-\frac{m}{2K}$ with
$X^K_0 =\lfloor Kz\rfloor$.  Assume that $\frac{\lfloor Kz\rfloor}{K}\to x$  and let
\be{} Z^K_t:= \frac{1}{K}X^K_{Kt}.\ee
Then as $K\to\infty$
\be{}  \{Z^K_t\}_{t\geq 0} \Longrightarrow \{Z_t\}_{t\geq 0},\ee
where $\{Z_t\}_{t\geq 0}$ is a CSBP with generator $G$ given by (\ref{FCSBG}) with $\gamma =1$  and $Z_0=x$. The convergence is in the sense of weak convergence of probability measures on $D_{[0,\infty)}([0,\infty))$ and the limiting process is a.s. continuous.

(b) Consider a sequence of BGW processes $\{X_k^N\}$ with mean
offspring sizes \be{ekc1}E(\xi^N)= m_{N}= 1+\frac{m}{N}\ee and offspring variances
 \be{ekc2}\rm{Var}(\xi^N)=\gamma >0.\ee  Let
 \be{}  Z^N_t:=\frac{1}{N}X^N_{\lfloor Nt\rfloor}.
 \ee
 Assume that $Z^N_0\rightarrow Z_{0}$ as $N\rightarrow \infty$. Then
\be{} \{Z^N_t\}_{t\geq 0} \Longrightarrow\{Z_t\}_{t\geq 0},
\ee
that is, $Z^N_t$ converges in distribution on $D_{[0,\infty)}([0,\infty))$ to a
Markov diffusion process, $\{Z_t\}_{t\geq 0}$,  called the Feller continuous state branching
process (CSBP).  The generator of the CSBP $\{Z_t\}$ acting on functions
$f\in C^2_0([0,\infty))$  is given by
\be{CSBPG} Gf(x)= mx\frac{\partial f}{\partial x} +\frac{1}{2}\gamma
x\frac{\partial^2 f}{\partial x^2}.\ee

\end{theorem}

\begin{proof}(a) The proof follows a standard program for weak convergence of processes, namely,
\begin{itemize}
\item the convergence of the finite dimensional distributions,
\item proof that the laws of the processes $P^K\in\mathcal{P}(D_{[0,\infty)}([0,\infty)))$ are tight.
\end{itemize}

To show that the finite dimensional distributions converge, first substitute birth and death rates $b=1+\frac{m}{2K}$,
$d=1-\frac{m}{2K}$,  in (\ref{cfgo}) to obtain the Laplace transform of $Z^K_t$ with  $Z^K_0=\lfloor Kz\rfloor$ as follows:

\bea{LTZ}&&\\ \qquad &&E(e^{-\theta Z^K_t})=L^K(t,\theta)\nonumber\\&&= \left(-\frac{K (e^{-\frac{\theta}{K}}-1)(e^{m t}-1) -\frac{m}{2}(e^{-\frac{\theta}{K}}-1)e^{m t}-\frac{m}{2}(e^{-\frac{\theta}{K}}+1)}{ K(e^{-\frac{\theta}{K}}-1)(e^{m t}-1) +\frac{m}{2}(e^{-\frac{\theta}{K}}-1)e^{m t}-\frac{m}{2}(e^{-\frac{\theta}{K}}+1)}   \right)^{\lfloor Kz\rfloor}\nonumber
\\&&= \left(-\frac{ (-\theta) (e^{m t}-1)+\frac{\theta^2}{2K} -\frac{m}{2}(e^{-\frac{\theta}{K}}-1)e^{m t}-\frac{m}{2}(e^{-\frac{\theta}{K}}+1)+O(K^{-2})}{ -\theta(e^{m t}-1)+\frac{\theta^2}{2K} +\frac{m}{2}(e^{-\frac{\theta}{K}}-1)e^{m t}-\frac{m}{2}(e^{-\frac{\theta}{K}}+1)+O(K^{-2})}   \right)^{\lfloor Kz\rfloor}\nonumber
\\& & \longrightarrow \exp\left(-\frac{m \theta z e^{m t}}{m+\theta (e^{m t} -1)}  \right).\nonumber\eea
This coincides (see Proposition \ref{FCSBLT}) with the Laplace transform at time $t$ of the diffusion process with $Z_0=x$ and with generator
\be{} Gf(x)= m x\frac{\partial f}{\partial x} +\frac{1}{2}
x\frac{\partial^2 f}{\partial x^2}.\ee

Using the Markov property and the continuity of the transition probability in $x$ we can then obtain convergence of the finite dimensional distributions.

To complete the proof we must verify that the probability laws of $\{Z^K_t\}_{t\geq 0}$ denoted by $P^K\in\mathcal{P}(D_{[0,\infty)}([0,\infty))$ are tight. We will use the Aldous condition.  We first verify that given $\delta>0$ there exists $0<L<\infty$ such
\be{max} \sup_K P^K(\sup_{0\leq t\leq T} X^K(t) > L)\leq \delta.\ee
Note that the generator of $Z^K_t=\frac{X^K_{Kt}}{K}$ is
\be{} G^Kf(\frac{n}{K})=\frac{ n}{K}\cdot K^2[f(\frac{n+1}{K})+f(\frac{n-1}{K})-2f(\frac{n}{K})]+\frac{m n}{2K}\cdot K[f(\frac{n+1}{K})-f(\frac{n-1}{K})].\ee
Then
\be{} M^K_t:=Z^K_t-m \int_0^t Z^K_s ds\quad\text{is a martingale}.\ee
By Gronwall's inequality
\be{} \sup_{0\leq t\leq T} Z^K_t\leq \sup_{0\leq t\leq T} |M^K_t|e^{mt}.\ee
Applying Doob's maximal inequality to $M^K_t$
\be{} P(\sup_{0\leq t\leq T} |M^K_t|\geq R)\leq \frac{E((M^K_T)^2)}{R^2}.\ee
It remains  to compute $E((M^K_T)^2)$. We have
\be{} E((M^K_T)^2)\leq   E(Z^K_T)^2)+ 2|m| \int_0^T E(Z^K_sZ^K_T)ds+\int_0^T\int_0^TE(Z^K_sZ^K_t)dsdt.\ee
Using (\ref{BD2M}) we can check that
\be{}  E[(Z^K_t)^2]\leq Z^K_0 e^{mt}\frac{(e^{mt}-1)}{m} + (Z^K_0)^2 e^{2mt}.\ee

 A simple calculation then yields
\be{} E((M^K_T)^2)\leq  C(T,z)\ee
where $C(T,z)$ does not depend on $K$ which proves (\ref{max}).

We can then apply the Aldous sufficient condition for tightness, namely, given stopping times $\tau_K\leq T$ and $\delta_K\downarrow 0$ as $K\to\infty$
\be{} \lim_{K\to\infty} P^K(|Z^K_{\tau_K+\delta_K}-Z^K_{\tau_K}|>\ve)=0.\ee
First note that $X^K_{\tau_K}$ is tight so we can take a convergent subsequence. Then by Skorohod's representation we can put these on a common probability space so that there is a.s. convergence.  In this setting assume that $X^{K_n}_{\tau_{K_n}}\to x$.  It now suffices to prove that $X^{K_n}_{\tau_{K_n}+\delta_{K_n}}$ converges in distribution to $x$.  Then
by the strong Markov property we have
\bea{}&&\\ &&E(e^{-\theta (Z^K_{\tau_K+\delta_K})}-e^{-\theta Z^K_{\tau_K}}|Z^K_{\tau_K})\nonumber\\&& =  \left(-\frac{K (e^{-\frac{\theta}{K}}-1)(e^{m \delta_K}-1) -\frac{m}{2}(e^{-\frac{\theta}{K}}-1)e^{m \delta_K}-\frac{m}{2}(e^{-\frac{\theta}{K}}+1)}{K (e^{-\frac{\theta}{K}}-1)(e^{m
\delta_K}-1) +\frac{m}{2}(e^{-\frac{\theta}{K}}-1)e^{m \delta_K}-\frac{m}{2}(e^{-\frac{\theta}{K}}+1)}   \right)^{\lfloor KZ^K_{\tau_K}\rfloor }-e^{-\theta Z^K_{\tau_K}}\nonumber\\&& \longrightarrow 0 \text{   on   } \{\sup_{0\leq t\leq T} Z^K(t) \leq L\}\text{  as }K\to\infty.\nonumber\eea
Therefore  $Z^K_{\tau_K+\delta_K}-Z^K_{\tau_K} \to 0$ in distribution and for $\ve, \eta>0$ we can find $K_0$ such that
 \be{} P^K(|Z^K_{\tau_K+\delta_K}-Z^K_{\tau_K}|>\ve) <2\eta,\quad \forall \;K\geq K_0.\ee
 This completes the proof of tightness.

(b) See Ethier and Kurtz (\cite{EK-86} Chapter 9, Theorem 1.3) for a
proof based on a  semigroup convergence theorem (e.g. \cite{EK-86}, Chap. 1, Theorem 6.5). This involves showing that
\be{} \lim_{N\to\infty}\sup_{x=\frac{\ell}{N},\;\ell\in\N}|N(T_Nf(x)-f(x))-Gf(x)|=0\quad\forall\; f\in C^\infty_c([0,\infty)),\ee
where
\be{}
T_Nf(x)=E[f(\frac{1}{N}\sum_{k=1}^{Nx}\xi^N_k)],\quad x\in\{\frac{\ell}{N},\;\ell\in \N\}\ee
and where $\{\xi^N_k\}$ are i.i.d. satisfy (\ref{ekc1}), (\ref{ekc2}).

\end{proof}

\begin{remark}
These results can also be proved using the martingale problem formulation in the same way as is carried out below for the Wright-Fisher model.

\end{remark}

\section{The critical BGW tree}

\subsection{The rooted BGW tree as a metric space}

We begin by recalling that  a BGW tree $\mathcal{T}\in\mathbf{T}$ with root $\varnothing$ is a graph in which the vertices are  a subset of
\be{}\mathcal{I}= \varnothing \cup \cup_{n=1}^\infty \mathbb{N}_0\ee
 satisfying conditions (\ref{PRT}).
  Recall that if $x=(i_1,\dots, i_n)\in\mathcal{T}$ is said to be in generation $n$, denoted by $H_N(x)=n$ where $N=\#(\mathcal{T})$.
The edges are given by the set of pairs of the form $((i_1,\dots,i_n),(i_1,\dots,i_n,j) $.

The {\em lexicographic order} is  an order relation on the vertices of $\mathcal{T}$ defined as follows.
We say that  $x=(i_1,\dots,i_n)$ and $y=(j_1,\dots,j_m)$  have a {\em last common ancestor} at generation $\ell \geq 1$ if
\be{} (i_1,\dots,i_\ell)=(j_1,\dots,j_\ell) \text{ and } i_{\ell+1}  \ne j_{\ell+1} \text{( or is empty)}.  \ee

Given $\mathcal{T}$ with $\#(\mathcal{T})=N$ we can order the vertices in lexicographic order $\varnothing,x_1,x_2,\dots,x_{N-1}$.
We can then embed it  in the plane so $x_i$ appears to the left of $x_j$ if $i<j$.

The corresponding {\em height function} $H_N(k) $ of a tree of size $\#(\mathcal{T})=N$ is defined by
\be{} H_N(k):=|x_k|,\quad 0\leq k\leq N-1\ee
where $|x|$ denotes the generation of $x$.

Note that  the number of visits of $H_N(k)$  to $n$ gives the population size at generation $n$, that is, \be{}Z_n =\sum_{k=0}^{N-1}1_{\{n\}}(H_N(k))\ee
where  $1_{\{n\}}$ denotes the indicator function.

 We now define a distance between the individuals in $\mathcal{T}$. If we assign length $1$ to each edge then a metric $d_{\mathcal{T}}(x,y)$ can be defined on $\mathcal{T}$ by
\be{} d_{\mathcal{T}}(x,y):= \text{the length of the shortest path in  }
\mathcal{T} \text{  from } x \text{  to } y.\ee
Since the critical BGW tree is a.s. finite this produces a compact metric space and is an example of random compact rooted real tree which we define below.

\begin{remark}  Note that a reordering  of the offspring (in the lexicographic order) defines a root preserving isometry.  We can then associate to  $\mathcal{T}$ the corresponding equivalence class of plane trees (modulo the family of root preserving isometries). This equivalence class is characterized by $(\#(\mathcal{T}),\varnothing,d_{\mathcal{T}}(.,.))$.
\end{remark}

We now briefly introduce the { reduced tree} at generation $n$.  We denote the set of nth generation individuals
\be{}  X_n =\mathcal{T}\cap \mathbb{N}_0^n.\ee
The {\em reduced tree} \be{}\mathcal{T}^R_n:= \{x\in\mathcal{T}: x=(i_1,\dots,i_r),\;r=1,\dots,n, \text{ such that }\exists (i_1,\dots, i_n)\in X_n\}.\ee
We also define a metric on $X_n$ by $d_n(x,y) :=n-\ell$ if the last common ancestor of $x,y$ is in generation $\ell <n$.
It is easy to verify that $d_n$ is an {\em ultrametric}, that is,
\be{} d_n(x,y)\leq \max(d_n(x,z),d_n(z,y))\quad \text{ for any } z\in X_n.\ee

 \subsection{The contour functions}

Given a tree $\mathcal{T}$ with $\#(\mathcal{T}) <\infty$ we define the {\em contour function} \be{}C^{\mathcal{T}}=C^{\mathcal{T}}(t):0\leq t\leq 2(\#(\mathcal{T})-1)\ee
which is obtained by taking a particle that starts from the root
of $\mathcal{T}$ and visits continuously all edges at speed one,
moving away from the root if possible otherwise going backwards along the edge leading to the root and respecting the
lexicographical order of vertices. The domain of $C^{\mathcal{T}}$ can be extended to $[0,\infty)$ by setting $C^{\mathcal{T}}(t)=0$ for $t>2(\#(\mathcal{T})-1))$. In other words,  $C^{\mathcal{T}}$  is a piecewise linear process given by the distance from the root as we move through the tree.

We have considered above the Yaglom conditioned limit theorem (Theorem 4.1) for a critical BGW process.  Similarly it is if interest to consider the conditioned BGW process conditioned on $\#(\mathcal{T})$.  In order to formulate results for this we need to introduce two additional notions, real trees and the Gromov-Hausdorf metric.



\begin{figure}[h!]
\caption{BGW Tree and contour function, $N=$ 10}
\begin{center}
\setlength{\unitlength}{0.1cm}
\begin{picture}(40,80)(20,0)
\setlength{\unitlength}{2cm}
\put(0,0){$\text{root}$}
\put(0,0){\line(-1,1){1}}
\put(0,0){\line(1,1){1}}
\put(1,1){\line(1,2){.5}}
\put(1,1){\line(-1,2){.5}}
\put(-1,1){\line(1,2){.5}}
\put(-1,1){\line(-1,2){.5}}
\put(-1,1){\line(0,2){1}}
\put(-1,2){\line(1,2){.5}}
\put(-1,2){\line(-1,2){.5}}
\put(1.1,1){$2$}
\put(-0.9,1){$1$}
\put(1.5,2){$22$}
\put(0.5,2){$21$}
\put(-0.5,2){$13$}
\put(-1,2){$12$}
\put(-1.5,2){$11$}
\put(-0.5,3){$122$}
\put(-1.5,3){$121$}
\end{picture}
\end{center}
\end{figure}
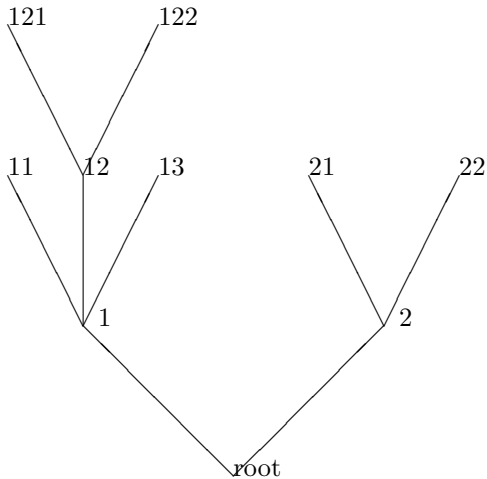

\begin{figure}[h!]
\begin{center}
\setlength{\unitlength}{0.1cm}
\begin{picture}(90,60)
\setlength{\unitlength}{1cm}
\put(3,6){$\large{\text{Contour function}}$}
\put(0,0){\vector(1,0){10}}
\put(0,0){\vector(0,1){5}}
\put(0,0){\line(1,3){0.5}}
\put(1,3){\circle*{0.1}}
\put(0.5,1.5){\line(1,3){0.5}}
\put(0.5,1.5){\circle*{0.1}}
\put(0.6,1.5){$1$}
\put(1,3){\line(1,-3){0.5}}
\put(1.5,1.5){\line(1,3){0.5}}
\put(2,3){\line(1,3){0.5}}
\put(2.5,4.5){\line(1,-3){0.5}}
\put(3,3){\line(1,3){0.5}}
\put(3.5,4.5){\line(1,-3){1}}
\put(4.5,1.5){\line(1,3){0.5}}
\put(5,3){\line(1,-3){1}}
\put(6,0){\line(1,3){1}}
\put(7,3){\line(1,-3){0.5}}
\put(7.5,1.5){\line(1,3){0.5}}
\put(8,3){\line(1,-3){1}}
\put(1.1,3){$11$}
\put(2,3){\circle*{0.1}}
\put(2.1,3){$12$}
\put(2.6,4.5){$121$}
\put(2.5,4.5){\circle*{0.1}}
\put(3.6,4.5){$122$}
\put(3.5,4.5){\circle*{0.1}}
\put(5,3){\circle*{0.1}}
\put(5.1,3){$13$}
\put(6.5,1.5){\circle*{0.1}}
\put(6.6,1.5){$2$}
\put(7,3){\circle*{0.1}}
\put(7.1,3){$21$}
\put(8,3){\circle*{0.1}}
\put(8.1,3){$22$}
\end{picture}
\end{center}
\end{figure}

\subsection{Real trees}
Following Evans \cite{E-07} and Le Gall \cite{LG-06} we now introduce the notion of real
trees and their coding. See Dress and Terhalle \cite{DT-96}, \cite{DMT-96} for general background on ``tree theory''.
 \beD{} A metric space
$(\mathcal{T},d)$ is a {\em real tree} if the following two properties
hold for every $(x,y)\in \mathcal{T}$.
\begin{itemize}
\item there is a unique isometric map $f_{x,y}$ from $[0,d(x,y)]$
into $\mathcal{T}$ such that $f_{x,y}(0)=x$ and
$f_{x,y}(d(x,y))=y$ \item If $q$ is a continuous injective map
from $[0,1]$ into $\mathcal{T}$, such that $q(0)=x,\; q(1)=y$,
then \be{}q([0,1])=f_{x,y}([0,d(x,y)]).\ee A rooted real tree is a
real tree $(\mathcal{T},d)$ with a distinguished vertex
$\emptyset$ called the root.
\end{itemize}

\end{definition}

As explained above it is natural to consider the equivalence class ${\mathbb{T}}$ of real trees $(\mathcal{T},d)$  modulo the family of root preserving isometries. Since this results in a collection of compact metric spaces, it can be furnished with the Gromov-Hausdorff metric $d_{GH}$ (see Appendix II, section \ref{s.AGH}).

(Recall that $d_{GH}((E_1,d_1),(E_2,d_2))$ is given by the infimum of the Hausdorff distances of  the images of $(E_1,d_1),(E_2,d_2)$
under the set of isometric embeddings of $(E_1,d_1),(E_2,d_2)$, respectively,  into a common compact metric space $(E_0,d_0)$.)

\beP{} (Evans, Pitman, Winter (2003) \cite{EPW-06}). The space of real trees furnished with the Gromov-Hausdorff topology, $(\mathbb{T},d_{GH})$,  is Polish.
\end{proposition}

\begin{remark} A metric space $(E,d)$ can be embedded isometrically into a real tree iff the {\em four point condition}
\be{} d(x,y)+d(u,v) \leq \max(d(x,u)+d(y,v),d(x,v)+d(y,u))\ee
is satisfied { for all 4-tuples } $u,v,x,y$   (Dress (1984) ,\cite{Dr-84})

\end{remark}

\subsection{Excursions from zero and real trees}

Consider a continuous function $g:[0,\infty)\to [0,\infty)$ with
non-empty compact support such that $g(0)=0$ and $g(s)=0\;\forall\; s>\inf\{t>:g(t)=0\}$ (we call this an positive excursion from 0). For $s,t\geq 0$, let \be{}
m_g(s,t)=\inf_{r\in[s\wedge t, s\vee t]}g(r),\ee

\be{} d_g(s,t)=g(s)+g(t)-2m_g(s,t).\ee

It is easy to check that $d_g$ is symmetric and satisfies the
triangle inequality.  Let  $\mathcal{T}_g$ denote the quotient
space  $[0,\infty)/\equiv $ where $ s \equiv t$ if $d_g(s,t)=0$.
 Then it can be verified
that the metric space $(\mathcal{T}_g,d_g)$ is a real tree (Le Gall (2006) \cite{LG-06}, Theorem 2.1).

Given $g$ the ancestral relationships can be reconstructed by noting that
$s$ is an ancestor of $t,$ $s\prec t$ iff $g(s)=\inf_{[s,t]}g(r)$

Let $(\mathbb{C},\|\cdot\|):=(\{(g,d_g):g\text{  a positive  excursion from } 0, d_g=  \text{ sup norm metric}\})$.

\bigskip

It can be verified that (e.g. Le Gall (2006) \cite{LG-06}, Lemma 2.3))
  the mapping from $(\mathbb{C},\|\cdot\|)$ to $(\mathbb{T},d_{GH})$ is continuous, that is,
 for two continuous functions $g,g'$ such
that $g(0)=g'(0)=0$:

\be{excon}  d_{GH}(\mathcal{T}_g,\mathcal{T}_{g'})\leq 2\parallel
g-g'\parallel.\ee

\subsection{The Aldous Continuum Random Tree}

Let  $\{B_t\}_{t\geq 0}$ be a standard Brownian motion and
\be{}  \tau_1:= \sup \{t\in [0,1]:B_t=0\},\quad \tau_2:= \inf \{t\geq 1:B_t=0\}.\ee
Then the {\em Brownian excursion} is a nonhomogeneous Markov process defined as follows:
\be{} B^e_t:= \frac{1}{\sqrt{(\tau_2 -\tau_1)}}B(\tau_1 +t(\tau_2 -\tau_1)),\quad 0\leq t\leq 1.\ee
It can be shown (see It\^o-McKean) \cite{IM-65} that the marginal PDF is given by
\be{} f(t,x)= \frac{2x^2}{\sqrt{2\pi t^3(1-t)^3}}e^{-\frac{x^2}{2t(1-t)}}.\ee

\beD{} ({\em Aldous continuum random tree}) Let $(B^e_t)_{0\leq t\leq 1}$ be a normalized Brownian
excursion (extended  to $[0,\infty)$ by setting $B^e_t=0$ for
$t>1$).  The corresponding  random real tree $(\mathcal{T}^e,d^e)$ is called the
{\em continuum random tree} (CRT). We denote by $P^{CRT}\in\mathcal{P}((\mathbb{T},d_{GH}))$ the probability law of
$\mathcal{T}^e$.

 The CRT  was introduced  by Aldous
(1991-1993) in a series of papers \cite{A-91a}, \cite{A-91b} and
\cite{A-93}.

\end{definition}

\subsection{Conditioned limit theorem for the critical BGW tree}

 Consider the special case of a BGW process with geometric offspring distribution, that is, \be{geom} p_k=P(\xi=k)=\frac{1}{2^{k+1 }},\quad k=0,1,2,\dots,\ee
 \beL{} For the offspring distribution (\ref{geom}) the contour process  $C^{\mathcal{T}}$ is given by a
 simple random walk $\{S_k\}$ with \be{}P(S_{k+1}-S_k=\pm 1)=\frac{1}{2}.\ee
 \end{lemma}
 \begin{proof}
 This can be verified by first noting that in this case  the number of jumps  from $0$ to $1$ corresponds to the number of offspring of the initial vertex. Now let $\tau^k_1,\tau^k_2,\dots$ denote the times of visits to height $k$. Consider the $m$th such visit to height $k$, $k\geq 1$.  The corresponding vertex is the offspring of a vertex at height $k-1$, say, the $\ell$th offspring. Then
 \be{}  C(\tau^k_m+1) = \left\{\begin{split}& k+1, \text{ with probability }P(\xi\geq \ell+1|\xi\geq \ell)=\frac{1}{2}\\& k-1\text{ with probability }P(\xi= \ell|\xi\geq \ell)=\frac{1}{2}\end{split}\right.
 \ee
 \end{proof}

 \beP{} Let $P^{BGW}(\cdot|\#(\mathcal{T})=n) \in \mathcal{P}((\mathbb{T},d_{GH}))$ denote the probability law of the BGW tree
 with offspring distribution (\ref{geom}) conditioned to have $n$ vertices. Then
 \be{} P^{BGW}(\frac{\mathcal{T}}{2\sqrt{2n}}|\#(\mathcal{T})=n)\Rightarrow P^{CRT}\ee
 in the sense of weak convergence in $\mathcal{P}((\mathbb{T},d_{GH}))$.
 \end{proposition}
\begin{proof}
  Letting $S_0=0$, and $N=\min\{k>0:S_k=0\}$ and conditioning on $N=2n$ we have  the contour process for this BGW process to have total population $n$.  Note that  this is simply an excursion of the simple random walk conditioned to first return to the origin at time $N=2n$, $S^N_k$.  But it is known that  which rescaled converges as $n\to\infty$ to a Brownian excursion from $0$ (see Durrett, Iglehart and Miller (1977) \cite{DIM-77}).
 \be{} \left(\frac{S^N(\lfloor 2nu\rfloor )}{2\sqrt{2n}}\right)_{0\leq u\leq 1}\Rightarrow (B^{e}_u)_{0\leq u\leq 1}.\ee
 where $B^{e}$ is the standard Brownian excursion.
 Using (\ref{excon}) and the continuous mapping theorem (\cite{Bill-99}, Theorem 2.7) this implies  that the laws of the corresponding BGW trees converge to the CRT as $n\to\infty$.
\end{proof}

\subsection{Aldous Invariance Principle for BGW trees}

A remarkable result of Aldous is the {\em invariance principle} for scaling limit of critical BGW tree, that is, the CRT arises as the limit for the entire class of critical BGW processes with aperiodic offspring distributions having finite second moments.

\begin{theorem} (Invariance principle for BGW trees - Aldous (1993) \cite{A-93},
Theorem 23.) \\ Consider the critical BGW tree with offspring
distribution $\mu.$ Assume that $\mu$ is aperiodic with variance  $\sigma^{2}<\infty$. \ Then
the distribution of the rescaled tree
\[
\frac{\sigma}{2\sqrt{n}}\mathcal{T}%
\]
under the probability measure $P^{BGW}(\cdot|\#(\mathcal{T})=n)$ (i.e.
conditioning that total population up to extinction is $n$)
converges as $n\rightarrow\infty$ to the law of CRT.
\end{theorem}

\begin{proof} The proof is given in \cite{A-93}.  It is too long and complex to include here.
\bigskip

However some of the ideas behind the proof are as follows.
 Using (\ref{excon}) we see that the result would  follow if  the rescaled contour process
\be{}(\frac{\sigma}{2\sqrt{n}}C^{\mathcal{T}}(2nt):\;0\leq t\leq 1)\ee
under the probability measure $P^{BGW}(\cdot|\#(\mathcal{T})=n)$
converges in distribution to the normalized Brownian excursion.  In the general case can no longer be represented by the excursion of a simple random walk.
 Aldous (1993) \cite{A-93}  proof of the invariance result  is based on a characterization of the distribution of the CRT.
Marckert and Mokkadem (2003) \cite{MM-03} gave an alternate proof (assuming the offspring distribution has exponential moments) involving only the contour and height functions. In particular they proved that for any critical offspring distribution with variance $\sigma^2$ the weak convergence of the rescaled contour function and  height processes. Their key idea is to couple the height process to the random walk (``depth-first queue process'')
\be{} S_n(j)=\sum_{i=0}^{j-1} (\xi_i-1),\; 1\leq j\leq n,\ee
that is, with jump distribution is given by $q_i=p_{i+1},\; i=-1,0,1,2,\dots$,
conditioned by $S_n(0)=0,\; S_n(i)\geq 0,\; 1\leq i\leq n-1,\; S_n(n)=-1$.
Then
\be{}  H_n(\ell)= \rm{Card}\{j:0\leq j\leq \ell -1,\;\min_{0\leq k\leq \ell -j}S_n(j+k)=S_n(j) \} , \quad  0\leq \ell < n-1.\ee
They then obtain exponential bounds on deviations between the height process and the conditioned random walk $S_n$ to prove that

\be{} (\frac{H_n(nt)}{\sqrt{n}})_{0\leq t\leq 1}\Rightarrow (\frac{2}{\sigma}B^e(t))_{0\leq t\leq 1}.\ee

\end{proof}

\begin{remark}

 It has also been proved that starting the BGW process with $n$ individuals then the rescaled height function   \be{} \{\frac{1}{n}H_n({\lfloor n^2t\rfloor})\}_{t\geq 0} \to (H_t)_{t\geq 0}\quad \text{with }H_0=1 \ee
where
\be{}  H_t =(B_t-\inf_{0\leq s\leq t} B_s)\ee
where $B_t$ is a Brownian motion, that is, $H_t$ is reflecting Brownian motion. (See \cite{LG-99}).

Recall that  the Ray-Knight Theorem (\cite{RW-87}, 52.1) states that if $B_t$ is a Brownian motion with local time $\{\ell^a_t\}$
\be{} T:= \inf\{u:\ell^0_u >1\},\ee  then the Brownian local time $\{\ell_T^a:a\geq 0\}$ has the same law as the
Feller CSB satisfying
\be{} d Z_t = 2\sqrt{Z_t}dW_t,\quad  Z_0=1.\ee
In other words the local time of the height process is a version of the Feller CSB starting at 1. More precisely,  the initial mass $Z_0=1$ corresponds to the local time at $0$ of a reflecting Brownian motion on $[0,T]$  and for $t\geq 0$ $Z_t =\ell^a_{T}$, that is the occupation density of the reflecting Brownian motion.

\end{remark}

\section{Remark on general continuous state branching}

By the  basic result of Silverstein \cite{S-69} the general continuous state branching process has log-Laplace equation

\be{}  u_t(\lambda)+\int_0^t\psi(u_s(\lambda))=\lambda,\ee
with
\be{} \psi(u) =\alpha u+\beta u^2 +\int_0^\infty (e^{-ru}-1+ru)
\nu(dr)\ee
where $\alpha,\beta\geq 0$ and $\nu$ is a $\sigma$-finite measure on $(0,\infty)$ such that $\int (r\wedge r^2)\nu(dr)<\infty$. This include the class of $(1+\beta)$ CSB which arise as limits of BGW processes in which the offspring distribution has infinite second moments and are related to stable processes and other L\'evy processes. The genealogical structure, stable continuum trees and convergence of the contour process in this general setting have been developed by Duquesne and LeGall \cite{DL-02} but we do not consider this  major topic here.

\chapter{Wright-Fisher Processes}


\section{Introductory remarks}

The BGW processes and birth and death processes we have studied in the previous chapters have the property that
\be{} X_n\to 0\text{  or  } \infty,\quad a.s.\ee

A more realistic model is one in which the population grows at low population densities and tends to a steady state
near some constant value.  The Wright-Fisher model that we consider in this chapter (and the corresponding Moran continuous time model) assume that the total population remains at a constant level $N$ and focusses on the changes in the relative proportions of the different types.  Fluctuations of the total population, provided that they do not become too small, result in time-varying resampling rates in the Wright-Fisher model but do not change the main qualitative features of the conclusions.

The branching model and the Wright-Fisher idealized models are complementary.  The branching process model provides an important approximation in two cases:
\begin{itemize}
\item  If the total population density becomes small then the critical and near critical branching process provides an useful approximation to compute extinction probabilities.
\item If a new type emerges which has a competitive advantage, then the supercritical branching model provides a good approximation to the growth of this type as long as its contribution to the total population is small.
\end{itemize}
Models which incorporate multiple types, supercritical growth at low densities and have non-trivial steady states will be discussed in a later chapter.  The advantage of the idealized models we discuss here is the possibility of explicit solutions.

\section{Wright-Fisher Markov Chain Model}

The classical neutral Wright-Fisher (1931) model is \ a discrete
time model of a population with constant size $N$ and types
$E=\{1,2\}$. \ Let $X_{n}$ be the number of type $1$ individuals
at time $n$. Then $X_{n}$ is a Markov chain with state space
$\{0,\dots,N\}$ and transition probabilities:

\[
P(X_{n+1}=j|X_{n}=i)=\left(
\begin{array}
[c]{c}%
N\\
j
\end{array}
\right)  \left(  \frac{i}{N}\right)  ^{j}\left(
1-\frac{i}{N}\right)  ^{N-j},\quad j=0,\dots,N.%
\]

In other words at generation $n+1$ this involves binomial sampling
with probability $p=\frac{X_n}{N}$, that is,  the current
empirical probability of type $1$.  Looking backwards from the
viewpoint of generation $n+1$ this can be interpreted as having
each of the $N$ individuals of the $(n+1)$st generation ``pick
their parents at random'' from the population at time $n $.

Similarly, the neutral $K$-allele Wright Fisher model with types $E_K=\{e_1,\dots,e_K\}$ is given by a
 Markov chain $X_n$ with state space $\mathcal{n}(E_K)$ (counting measures)
and
\bea{}&&
P(X_{n+1}=(\beta_{1},\dots\beta_{K})|X_{n}=(\alpha_{1},\dots,\alpha
_{K}))\\&& =\frac{N!}{\beta_{1}!\beta_{2}!\dots\beta_{K}!}\left(  \frac{\alpha_{1}%
}{N}\right)  ^{\beta_{1}}\dots\left(  \frac{\alpha_{K}}{N}\right)  ^{\beta
_{K}}\nonumber
\eea
In this case the binomial sampling is simply replaced by
multinomial sampling.

Consider the multinomial distribution with parameters $(N,p_1,\dots,p_K)$. Then the moment generating function is given by

\be{} M(\theta_1,\dots,\theta_K)=E(\exp(\sum_{i=1}^K \theta_iX_i))=\left(\sum_{i=1}^K p_ie^{\theta_i}\right)^N\ee

Then
\be{} E(X_i)=Np_i,\quad \rm{Var}(X_i)= Np_i(1-p_i),\ee
and
\be{} \rm{Cov}(X_i,X_j)=-Np_ip_j,\;i\ne j.\ee

\begin{remark} We can relax the assumptions of the Wright-Fisher model in two ways. First,  if we relax the assumption of the total population constant, equal to $N$, we obtain a Wright-Fisher model with variable resampling rate (e.g.  Donnelly and Kurtz \cite{DK-99a} and Kaj and Krone \cite{KK-03}).

To introduce the second way to relax the assumptions note that we can obtain  the Wright-Fisher model as follows.  Consider a population of $N$ individuals in generation $n$ with possible types in $E_K$, $Y^n_1,\dots,Y^n_N$. Assume  each individual has a Poisson number of offspring with mean m, $(
Z_1,\dots,Z_N)$ and the offspring is of the same type as the parent.  Then \[\text{  conditioned on  }\sum_{i=1}^N Z_i=N,\] the resulting population  $(Y^{(n+1)}_1,\dots,Y^{(n+1)}_N)$  is multinomial $(N;\frac{1}{N};\dots,\frac{1}{N})$,   that is, we have a a multitype (Poisson) branching process conditioned to have constant total population $N$.  If we then define
 \be{} p_{n+1}(i)=\frac{1}{N}\sum_{j=1}^N 1(Y^{(n+1)}_j=i),\; i=1,\dots,K,\ee then $(p_{n+1}(1),\dots,p_{n+1}(K))$ is multinomial $(N;p_n(1),\dots,p_n(K))$ where
 \be{} p_n(i)= \frac{1}{N}\sum_{j=1}^N 1(Y^n_j=i),\; i=1,\dots,K.\ee
  We can generalize this by assuming that the offspring distribution of the individuals is given by a common distribution on $\mathbb{N}_0$.  Then again conditioned the total population to have  constant size $N$ the vector $(Y^{n+1}_1,\dots,Y^{n+1}_N)$   is {\em exchangeable} but not necessarily   multinomial.
    This exchangeability assumption is the basis of   the Cannings Model (see e.g. Ewens \cite{E-04}).

 \end{remark}

 A basic phenomenon of neutral Wright-Fisher without mutation is {\em fixation}, that is, the elimination of all but one
type at a finite random time.  To see this note that for each $j=1,\dots,K$,
$\delta_j\in\mathcal{P}(E_K)$ are absorbing states and
$X_n(j)$ is a martingale.  Therefore $X_n\to X_\infty,\;a.s.$ Since $\rm{Var}(X_{n+1})=NX_n(1-X_n)$, this means that $X_\infty =0\text{  or  }1,\; a.s.$ and  $X_n$ must be $0$ or $1$ after a finite number of
generations (since only the values $\frac{k}{N}$ are possible).

\subsection{Types in population genetics}

The notion of type in population biology is based on the {\em genotype}. The genotype of an individual is specified by the {\em genome} and this codes {\em genetic information} that passes, possibly modified, from parent to offspring (parents in sexual reproduction).  The genome consists of a set of  {\em chromosomes} (23 in humans).  A chromosome is a single molecule
of DNA that contains  many genes, regulatory elements and other nucleotide sequences. A given position on a chromosome is called a {\em locus}  (loci) and may be occupied by one or more {\em genes}. Genes code for the production of a protein. The different variations of the gene at a particular locus are called {\em alleles.}  The ordered list of loci  for a particular genome is called a {\em genetic map}.
The {\em phenotype} of an organism describes its structure and behaviour, that is, how it interacts with its environment. The relationship between genotype and phenotype is not necessarily 1-1.  The field of {\em epigenetics} studies this relationship and in particular the mechanisms during cellular development that produce different outcomes from the same  genetic information.

{\em  Diploid individuals}  have two homologous copies of each chromosome, usually one from the mother and one from the father in the case of sexual reproduction. Homologous chromosomes contain the same genes at the same loci but possibly different alleles at those genes.

\subsection{Finite population resampling in a diploid population}
For a diploid population with $K$-alleles $e_{1},\dots,e_{K}$
at a particular gene we can focus on  the set of types denoted by $E_K^{2\circ }$ consisting of the set of $\frac{K(K+1)}{2}$ unordered pairs
$(e_i,e_j)$.  The genotype $(e_i,e_j)$ is said to be homozygous (at the locus in question) if
$e_i=e_j$, otherwise heterozygous.

Consider a finite population of $N$ individuals. Let
\[
P_{ij}=\text{ proportion of type }(e_{i},e_{j})%
\]

Then, $p_{i}$, the proportion of allele $e_{i}$ is%
\[
p_{i}=P_{ii}+\frac{1}{2}\sum_{j\ne i}P_{ij}. %
\]

The probability $\{P_{ij}\}$ on   $E_K^{2\circ }$ is said to be a
Hardy-Weinberg equilibrium  if \be{}P_{ij} =
(2-\delta_{ij})p_{i}p_{j}.\ee This  is what is obtained if one picks
independently the parent types $e_{i}$ and $e_{j}$ from a
population having proportions $\{p_{i}\}$ ( in the case of sexual
reproduction this corresponds to ``random mating'').

Consider a diploid Wright-Fisher model  with $N$ individuals  therefore $2N$ genes  with random mating. This means that an individual at generation $(n+1)$  has two genes randomly chosen from the $2N$ genes in generation $n$.

In order to introduce the notions of identity by descent and genealogy we assume that in generation $0$ each of the $2N$ genes correspond to different alleles.   Now consider generation $n$.  What is the probability, $F_n$,  that an individual is homozygous, that is, two genes selected at random are of the same type (homozygous)? This will occur only if they are both descendants of the same gene in generation $0$.

First note that in generation 1,  this means that an individual is homozygous only if the same allele must be selected twice and this has probability $\frac{1}{2N}$.  In generation $n+1$ this happens if the same gene is selected twice or if different genes are selected from generation $n$ but they are identical alleles. Therefore,
\be{} F_1=\frac{1}{2N},\qquad F_n=\frac{1}{2N}+(1-\frac{1}{2N})F_{n-1}.\ee

Let  $H_n:=1-F_n$ (heterozygous). Then
\be{} H_1=1-\frac{1}{2N},\qquad  H_n=(1-\frac{1}{2N})H_{n-1},\qquad H_n=(1-\frac{1}{2N})^n\ee

Two randomly selected genes are said to be {\em identical by descent} if they are the same allele. This will happen if they have a {\em common ancestor}.  Therefore if $T_{2,1}$ denotes the time in generations back to the common ancestor we have
\be{} P(T_{2,1}>n)=H_n =(1-\frac{1}{2N})^n, \quad n=0,1,2,\dots,\ee

\be{} P(T_{2,1}=n)=\frac{1}{2N}(1-\frac{1}{2N})^{n-1},\qquad n=1,2\dots.\ee

Similarly, for $k$ randomly selected genes they are identical by descent if they all have a common ancestor.
We can consider the time $T_{k,1}$ in generations back to the most recent common ancestor of $k$ individuals randomly sampled from the population. We will return to discuss the distribution of $T_{k,1}$ in the limit as $N\to\infty$ in Chapter 9.

\subsection{Diploid population with mutation and selection}

In the previous section we considered only the mechanism of resampling (genetic drift). In addition to genetic drift the basic genetic mechanisms include mutation, selection and recombination. In this subsection we consider the Wright-Fisher model incorporating mutation and selection.

For a diploid population of size $N$ with mutation, selection and resampling the {\em reproduction cycle} can be modelled as follows (cf \cite{EK-86}, Chap. 10).
We assume that  in generation $0$ individuals have genotypic proportions
$\{P_{ij}\}$ and therefore the proportion of type $i$ (in the population of $2N$ genes) is
\[
p_{i}=P_{ii}+\frac{1}{2}\sum_{j\ne i}P_{ij}. %
\]

Stage I:\\
In the first stage  diploid  cells undergo meiotic
division producing haploid gametes (single chromosomes), that is, meiosis reduces
 the number of sets of chromosomes from two to one.
 The resulting  {\em gametes}  are haploid cells; that is, they contain one half a complete set of chromosomes.
 When two gametes fuse (in animals
typically involving a sperm and an egg), they form a {\em zygote} that
has two complete sets of chromosomes and therefore is diploid. The
zygote receives one set of chromosomes from each of the two
gametes through the {\em fusion} of the two gametes.
By the assumption of random mating, then in generation 1  this produces zygotes in {\em Hardy-Weinberg
proportions}  $(2-\delta_{ij})p_{i}p_{j}$.

Stage II: Selection and Mutation.

{\em Selection.} The resulting zygotes can have different viabilities for survival. The viability of $(e_i,e_j)$ has viability $V_{ij}$. Then the proportions of surviving zygotes  are proportional to the product of the
viabilities and the Hardy-Weinberg proportions, that is,

\be{}
P^{\rm{ sel}}_{k,\ell}=\frac{V_{k\ell}\cdot(2-\delta_{k\ell})p_{k}p_{\ell}}{\sum_{k'\leq\ell'}%
(2-\delta_{k'\ell'})V_{k'\ell'}p_{k'}p_{\ell'}}
\ee

{\em Mutation}.  We assume that each of the 2 gametes forming zygote can (independently) mutate with  probability $p_m$ and that if a gamete of type $e_i$ mutates then it produces a gamete
of type $e_j$ with probability $m_{ij}$.

\bea{re0}&&\\
P^{\rm{sel,mut}}_{ij}&&=(1-\frac{1}{2}\delta_{ij})\sum_{k\leq\ell}(m_{ki}m_{\ell
j}+m_{kj}m_{\ell i})P_{k\ell}^{\rm{sel}}\nonumber\\&&
=(1-\frac{1}{2}\delta_{ij})\sum_{k\leq\ell}(m_{ki}m_{\ell
j}+m_{kj}m_{\ell i})\frac{V_{k\ell}\cdot(2-\delta_{k\ell})p_{k}p_{\ell}}{\sum_{k'\leq\ell'}%
(2-\delta_{k'\ell'})V_{k'\ell'}p_{k'}p_{\ell'}}\nonumber
\eea




Stage III:  Resampling. Finally random sampling reduces the
population to $N$
adults with proportions $P_{ij}^{\rm{next}}$ where%
\be{re1}
(P_{ij}^{\rm{next}})_{i\leq j}\sim\frac{1}{N}\text{multinomial }(N,(P_{ij}%
^{\rm{sel,mut}})_{i\leq j}).
\ee
We then obtain a population of $2N$ gametes with proportions

\be{re2}
p^{\rm{next}}_{i}=P^{\rm{next}}_{ii}+\frac{1}{2}\sum_{j\ne
i}P^{\rm{next}}_{ij}.
\ee

 Therefore we have defined the process $\{X^N_n\}_{n\in\N}$ with state space
$\mathcal{P^N}(E_K)$.  If $X^N_n$ is a Markov chain we defined the transition function
\[
P(X^N_{n+1}=(p_1^{\rm{next}},\dots,p_K^{\rm{next}})|X^N_{n}=(p_1^{},\dots,p_K^{}))
=\pi_{p_1^{},\dots,p_K^{}}(p_1^{\rm{next}},\dots,p_K^{\rm{next}})\]
where the function $\pi$ is obtained from (\ref{re0}), (\ref{re1}), (\ref{re2}).
See Remark \ref{RCMP}.



\section{Diffusion Approximation of Wright-Fisher}

\subsection{Neutral 2-allele Wright-Fisher model}

As a warm-up to the use of diffusion approximations we consider
the case of 2 alleles $A_1,A_2$, $(k=2)$.  Let  $X^{N}_n$ denote the number of individuals of type $A_1$ at the nth generation. Then as above $\{X^{N}_n\}_{n\in\N}$ is a Markov chain.

\begin{theorem} (Neutral case without mutation)
Assume that $N^{-1}X^N_{0}\rightarrow p_0$ as $N\rightarrow\infty$. Then
\[
\{p_{N}(t):t\geq0\}\equiv\{N^{-1}X^N_{\lfloor Nt\rfloor},\;t\geq0\}\Longrightarrow
\{p(t):t\geq0\}
\]
where $\{p(t):t\geq0\}$ is a Markov diffusion process with state space $[0,1]$ and with  generator
\be{NWF}Gf(p)=\frac{1}{2}p(1-p)\frac{d^2}{d p^2}f(p)\ee if $f\in C^{2}([0,1])$.  This is equivalent
to the pathwise unique solution of the SDE%
\begin{align*}
dp(t)  &  =\sqrt{p(t)(1-p(t))}dB(t)\\
p(0)  &  =p_0.
\end{align*}
\end{theorem}

\begin{proof} Note that in this case $X^N_{n+1}$ is Binomial$(N,p_{n})$ where $p_n=\frac{X^N_n}{N}$.
Then from the Binomial formula, \
\begin{align*}
E_{X^N_n}(\frac{X^N_{n+1}}{N})  &  =\frac{X^N_n }{N}\\
 E_{X^N_n}[\left(  \frac{X^N_{n+1}}{N}-\frac{X^N_n}{N}\right)  ^{2}\left. \right|\frac{X^N_n}{N} ]  &
=\frac{1}{N}\left(\frac{X^N_n}{N}\left(  1-\frac{X^N_n}{N}\right) \right).
\end{align*}
We can then verify that
\be{2mp1n}\{p_{N}(t):=N^{-1}X^N_{\lfloor Nt\rfloor }:t\geq0\}\text{    is a martingale }\ee with
\bea{2moment}
&&E(p_{N}(t_{2})-p_{N}(t_{1}))^{2}    =E\sum_{k= {\lfloor Nt_{1}\rfloor}}%
^{{\lfloor Nt_{2}\rfloor }}(p_N(\frac{k+1}{N})-p_N(\frac{k}{N}))^2 \\&&=
\frac{1}{N} E\sum_{k= {\lfloor Nt_{1}\rfloor}}%
^{{\lfloor Nt_{2}\rfloor }}
p_{N}(\frac{k}{N})(1-p_{N}(\frac{k}{N}))\nonumber
\eea
and then that%
\be{2mp2n}
M_{N}(t)=p_{N}^{2}(t)-\frac{1}{N} \sum_{k= {0}}%
^{{\lfloor Nt\rfloor }}
p_{N}(\frac{k}{N})(1-p_{N}(\frac{k}{N}))
\ee
is a martingale.

Let $P^N_{p_N}\in\mathcal{P}(D_{[0,1]}([0,\infty))$ denote the probability law of  $\{p_N(t)\}_{t\geq 0}$ with $p_N(0)=p_N$.
From this we can prove that the sequence $\{P^N_{p_N(0)}\}_{N\in\N}$  is tight on $\mathcal{P}(D_{[0,\infty)}([0,1]))$. To verify this  as in the previous chapter we use  Aldous  criterion $P^N_{p_N(0)}(p_N(\tau_N +\delta_N)-p_N(\tau_N)>\ve)\to 0$ as $N\to\infty$ for any stopping times $\tau_N\leq T$ and $\delta_N\downarrow 0$.  This follows easily from the strong Markov property, (\ref{2moment}) and Chebyshev's inequality.
Since the processes $p_N(\cdot)$ are bounded it then  follows that for  any limit point $P_{p_0}$ of $P^N_{p_N(0)}$ we have
\be{2mp} \begin{split} &\{p(t)\}_{t\geq 0} \text{  is a bounded martingale with } p(0)=p_0\text{ and with increasing
process} \\& \langle p\rangle_t= \int_{0}^{t}p(s)(1-p(s))ds. \end{split}
\ee
Since the largest jump of $p_N(\cdot)$ goes to $0$ as $N\to\infty$ the limiting process is continuous (see Theorem 17.14 in the Appendix). Also,  by the Burkholder-Davis-Gundy inequality we have
\be{}  E((p(t_2)-p(t_1))^4)\leq \rm{const}\cdot (t_2-t_1)^2,\ee
so that $p(t)$ satisfies  Kolmogorov's criterion for a.s. continuous.

 We can then prove that there is a unique solution to this {\em martingale problem}, that is, for each $p$ there exists a unique probability measure
on $C_{[0,\infty)}([0,\infty))$  satisfying (\ref{2mp}) and therefore this defines a Markov diffusion process with generator (\ref{NWF}).

The uniqueness can proved by determining all joint moments of the form
\be{} E_p((p(t_1)^{k_1}\dots (p(t_\ell))^{k_\ell}),\quad 0\leq t_1< t_2<\dots <t_\ell,\quad k_i\in\N
\ee
by solving a closed system of differential equation.
It can also be proved using duality and this will be done in detail below (Chapter 7) in a more general case.)
\end{proof}

We now give an illustrative application of the diffusion
approximation, namely the calculation of expected fixation times.

\begin{corollary} (Expected fixation time.)
Let $\tau:= \inf\{t: p(t)\in \{0,1\}\}$ denote the fixation time of the diffusion process. Then%
\[
E_{p}[\tau]=g(p)=-[p\log p+(1-p)\log(1-p)].
\]
\end{corollary}

\begin{proof}
Let $f\in C^2([0,1]),\;f(0)=f(1)=0$. Let $g_f(p):=\int_0^\infty
T_sf(p)ds$, and note that as $f\uparrow 1_{(0,1)}$ this converges to
the expected time spent in $(0,1)$.   Since $p(t)\to \{0,1\}\text{
as } t\to\infty,$ a.s., we can show that
\[ G\left(\int_0^t T_sf(p)ds\right) =\int_0^t GT_sf(p)ds =T_tf(p)-f(p)\to 0-f(p) \text{  as  } t\to\infty,\]
that is,
\[ Gg(p)=-f(p)\]
where $G$ is given by (\ref{NWF}).

Applying this to a sequence of $C^2$ functions increasing to
$1_{(0,1)}$
we get%
\begin{align*}
E_{p}(\tau)  &  =\int_{0}^{\infty}P_{p}(\tau>t)dt\\
&  =\int_{0}^{\infty}T_{t}1_{(0,1)}(p)dt\\
&  =g(p)
\end{align*}
We then obtain $g(p)$ by solving the differential equation
$Gg(p)=-1$ with boundary conditions
$g(0)=g(1)=0$ to obtain%
\[
E_{p}[\tau]=g(p)=-[p\log p+(1-p)\log(1-p)].
\]
\end{proof}

Let  $\tau_{N}$ denote the fixation time for $N^{-1}X_{[Nt]}$. We want to show that
\be{}
E_{\frac{X^N_0}{N}}[\tau_{N}]\rightarrow E_{p_0}[\tau]\text{ if }\frac{X^N_0}{N}\rightarrow
p \text{  and  } N\to\infty.
\ee
However note that  $\tau$ is not a continuous function on
$D([0,\infty),[0,1])$.  The weak convergence can be proved  for $\tau^{\ve}=\inf\{t: p(t)\not\in (\ve,1-\ve)\}$ (because there is no ``slowing down'' here). To complete the proof   it can be verified that for $\delta >0$
\be{}
\lim_{\ve \to 0}\limsup_{N\to\infty} P(|\tau^\ve_N-\tau_N|>\delta)=0
\ee
 (see Ethier and Kurtz,
\cite{EK-86}, Chapt. 10, Theorem 2.4).

\subsubsection{2-allele Wright-Fisher with mutation}

For each $N$ consider a Wright-Fisher population of size $M_N$ and
with mutation rates $m_{12}=\frac{u}{N}$  $A_1\to A_2$ and $m_{21}=\frac{v}{N}$ $A_2\to A_1$.

In this case $X^{N}_{n+1}$ is Binomial$(M_N,p_n)$ with
\be{}  p_n=(1-\frac{u}{N})\frac{X^{N}_n}{M_N}+\frac{v}{N}(1-\frac{X^{N}_n}{M_N}).\ee
We now consider
\be{} p_N(t)=\frac{1}{M_N} X_{\lfloor Nt\rfloor}.\ee

If we assume that \be{}\gamma=\lim_{N\to\infty}\frac{N}{M_N},\ee then the diffusion approximation is given by the diffusion process $p_t$ with generator

\be{2WFM} Gf(p)=\frac{\gamma}{2} p(1-p)\frac{\partial^2}{\partial p^2}+[-up+v(1-p)]\frac{\partial}{\partial p}.\ee

In this case the domain of the generator involves boundary conditions at $0$ and $1$ (see \cite{EK-86}, Chap. 8, Theorem 1.1)  but we will not need this.

\begin{remark}
Note that the diffusion coefficient is proportional to the inverse population size. Below for more complex models we frequently think of the diffusion coefficient in terms of  {\em inverse effective population size}.
\end{remark}

\subsubsection{Error estimates}

Consider a haploid Wright-Fisher population of size $M$ with mutation rates $m_{12}=u,\;m_{21}=v$.

Let  $p^{(M,u,v)}_t$ denote the diffusion process with generator (\ref{2WFM}) with $\gamma=\frac{1}{M}$. Then  if $\alpha,\beta\geq 0$,  the law of
\be{} \{Z_t^{(\alpha,\beta)}\}_{t\geq 0}:= p^{(M,\frac{\alpha}{M},\frac{\beta}{M})}_t\ee
is independent of $M$ and is a Wright-Fisher diffusion with generator

 \be{}Gf(p)=\frac{\gamma}{2} p(1-p)\frac{\partial^2}{\partial p^2}+[-\alpha p+\beta(1-p)]\frac{\partial}{\partial p}.\ee

The assumption of mutation rates of order $O(\frac{1}{N})$
corresponds to the case in which both mutation and genetic drift are of the same order and appear in the limit as population sizes goes to $\infty$. Other only one of the two mechanisms appears in the limit as $N\to\infty$.

On the other hand one can consider the diffusion process as an approximation to the finite population model.  Ethier and Norman (\cite{EN-77})  obtained an estimate of the error due to the diffusion approximation in the calculation of the expected value of a smooth function of the nth generation allelic frequency.

To formulate their result consider the Wright-Fisher Markov chain model $\{X^{(M,u,v)}_n\}$ with population size $M$ and one-step mutation probabilities $m_{12}=u,\;\;m_{21}= v$ and $p^{(M,u,v)}_t$ the Wright-Fisher diffusion  with generator (\ref{2WFM}) with $\gamma=\frac{1}{M}$.

\beT{} (Ethier and Norman \cite{EN-77}) Assume that $f\in C^6([0,1])$.
 Then for $n\in\N_0$,
 \be{} \begin{split} &|E_x(f(X_n^{(M,u,v)})- E_x(f(p_n^{(M,u,v)})|\\&
 \leq \frac{\max(u,v)}{2}\cdot \|f^{(1)}\| +\frac{1}{M}\left(\frac{1}{8}\|f^{(2)}\|+\frac{1}{216\sqrt{3}}\|f^{(3)}\|\right)\\&
 +\frac{9\max(u^2,v^2)}{2}\left(\sum_{j=1}^6\|f^{(j)}\|\right)+\frac{7}{16M^2} \sum_{j=2}^6\|f^{(j)}\|
 \end{split}
 \ee
where $\|f^{(j)}\|$ is the sup of the $j$th derivative of $f$.
\end{theorem}

We do not include a proof but sketch the main idea.  Let
\be{} (S_nf)(x):= E_x[f(X_n^{(M,u,v)})],\ee
\be{} (T_tf)(x):= E_x[f(p^{(M,u,v)}_t)].\ee
If $g\in C_b^6([0,\infty))$, then we have the Taylor expansions
\be{} (T_1g)(x)=g(x)+(Gg)(x)+\frac{G^2g(x)}{2}+\omega_2\frac{\| G^3g\|}{6},\quad |\omega_2|\leq 1\ee
and
\be{} (S_1g)(x)= g(x)+\sum_{j=1}^5 E_x[(X^{(M,u,v)}_1-x)^j]\frac{g^{(j)}(x)}{j!}+\omega_1 E_x[(X^{(M,u,v)}_1-x)^6]\frac{\|g^{(6)}\|}{6!},\quad |\omega_1|\leq 1.\ee
We then obtain
\be{}  \|S_1g-T_1g\|_M\leq \sum_{j=1}^6\gamma_j\|g^{(j)}\|\ee
for some constants $\gamma_j$.

The proof is then completed using the inequality
\be{} \|S_nf-T_nf\|_M\leq \sum_{k=0}^{n-1}\|(S_1-T_1)T_k\|_M\ee
where $\|\cdot\|_M$ is the sup norm on $\{\frac{j}{M}:j=\,1,\dots,M\}$.

\subsection{K-allele Wright-Fisher Diffusion}

Now consider the $K$-allele Wright-Fisher Markov chain
$\{X_k^{2N}\}_{k\in\N}$ with $2N$ gametes present in each generation
and assume that the mutation rates and fitnesses satisfy%
\be{}
m_{ij}    =\frac{q_{ij}}{2N} ,\;i\neq j,\quad m_{ii}=1-\frac{m}{N},\; m=\sum_j q_{ij} \ee
\be{VC} V_{ij}    =1+\frac{\sigma_{ij}}{2N}+O(\frac{1}{N^2}).\ee

We now consider the Markov process with state space
\be{}
\Delta_{K-1}:=\{(p_{1},\dots,p_{K}):\;p_{i}\geq0,\sum_{i=1}^{K}p_{i}=1\}.
\ee
defined by
\be{} \{p^{2N}(t):t\geq0\}\equiv\{\frac{1}{2N}X^{2N}_{[2Nt]},\;t\geq0\}.
\ee

\begin{theorem}
Assume that ${2N}^{-1}X^{2N}_{0}\rightarrow p$ as $N\rightarrow\infty$ in  $\Delta_{K-1}$.

Then the laws of the c\`adl\`ag processes $\{p_N(t):= \frac{1}{2N}X^{2N}_{t}\}_{t\geq 0}$ are tight and for any limit point
and function $f(p)=f(p_{1},\dots,p_{K-1})\in C^2(\Delta_{K-1})$,

\be{wfmp} M_f(t):= f(p(t))-\int_0^t G^Kf(p(s))ds \qquad\text{   is a martingale}\ee
 where
\bea{KWFG}
&&G^{K}f(p) \\ &&  =\frac{1}{2}\sum_{i,j=1}^{K-1}p_{i}(\delta_{ij}-p_{j}%
)\frac{\partial^{2}f(p)}{\partial p_{i}\partial p_{j}}\nonumber\\
&&  +\sum_{i=1}^{K-1}\left[ m \left(  (\sum_{j=1,j\ne i}^{K} q_{ji} p_j -
p_{i}\right)  +p_{i}\left(  \sum_{j=1}^{K}\sigma_{ij}p_j-\sum_{k,\ell}%
^{K}\sigma_{k\ell}p_{k}p_\ell\right)  \right]  \frac{\partial
f(p)}{\partial
p_{i}}.\nonumber%
\eea

The  martingale problem (\ref{wfmp}) has a unique solution which determines  a Markov diffusion process
 $\{p(t):t\geq0\}$ called the {\em K-allele Wright-Fisher diffusion}.

\end{theorem}
\begin{proof}  Following the pattern of the 2-allele neutral case the proof involves three steps which we now sketch.

Step 1.  The tightness of the probability laws $P^{N}$ of $\{p^{2N}(\cdot)\}$ on $D_{\Delta_{K-1}}([0,\infty))$  can be proved using Aldous criterion.

Step 2. Proof that for any limit point of $P^N$ and $i=1,\dots,K$

\bea{mpa} &&M_i(t):=p_i(t)-p_i(0)-\int_0^t\left[  m\left(  \sum_{j=1}^{K} q_{ji} p_j(s) -
p_{i}(s)\right)\right. \\&&\left. +p_{i}(s)\left(  \sum_{j=1}^{K}\sigma_{ij}p_j(s)-\sum_{k,\ell}%
^{K}\sigma_{k\ell}p_{k}(s)p_\ell(s)\right)\  \right] ds\nonumber
\eea
is a martingale  with quadratic covariation process
\be{mpb}
\langle M_i,M_j\rangle_t =\frac{1}{2}\int_0^t p_j(s)(\delta_{ij}-p_i(s))ds
\ee
To verify this, let $\mathcal{F}_{\frac{k}{2N}}=\sigma\{ p^{2N}_i(\frac{\ell}{2N}):\ell\leq k,i=1,\dots,K\}$.
Then we have for $k\in\N$

\be{}
\begin{split} & E[p^{2N}_i(\frac{k+1}{2N})-p^{2N}_i(\frac{k}{2N})\left. \right| \mathcal{F}_{\frac{k}{2N}}]\\ &= \frac{1}{2N} \left[m\left(\sum_{j=1,j\ne i}^{K} \frac{q_{ji}}{m} p^{2N}_j(\frac{k}{2N})
 -p^{2N}_{i}(\frac{k}{2N})\right)\right. \\&\left. \qquad\quad +\left(  \sum_{j=1}^{K}\sigma_{ij}p^{2N}_j(\frac{k}{2N})-\sum_{k,\ell=1}
^{K}\sigma_{k\ell}p^{2N}_{k}(\frac{k}{2N})p^{2N}_\ell(\frac{k}{2N})\right)  \right]\\&\qquad +o(\frac{1}{2N}) \end{split}\ee

\be{CVC}
Cov(p^{2N}_{i}(\frac{k+1}{2N}),p^{2N}_{j}(\frac{k+1}{2N})|\mathcal{F}_{\frac{k}{2N}})
=\frac{p^{2N}_{i}}{2N}(\frac{k}{2N})(\delta_{ij}-p^{2N}_{j}(\frac{k}{2N})) +o(\frac{1}{N})
\ee

\begin{remark}\label{RCMP} The Markov property for $X^N_n$ follows if in the resampling step the $\{P^{\rm{sel,mut}}_{ij}\}$ are in Hardy-Weinberg proportions which implies  that the
$\{p_i^{\rm{next}}\}$ are \be{}\text{multinomial}(2N,(p_1^{\rm{sel,mut}},\dots,p_K^{\rm{sel,mut}})).\ee This is true without selection or with multiplicative selection $V_{ij}=V_iV_j$ (which leads to haploid selection in the diffusion limit) but not in general. In the diffusion limit this can sometimes be dealt with by the $O(\frac{1}{N^2})$ term in (\ref{VC}). In general  the diffusion limit result remains true but the argument is more subtle.  The idea is that the selection-mutation changes the allele frequencies more slowly than the mechanism of Stages I and III which rapidly bring the frequencies to Hardy-Weinberg equilibrium  - see \cite{EK-86}, Chap. 10, section 3.

\end{remark}

Then for each $N$ and $i$
\bea{} &&M^{2N}_i(t):=p_i^{2N}(t)-p_i^{2N}(0)-\int_0^t\left[  m\left(  \sum_{j=1,j\ne i}^{K} q_{ji} p^{2N}_j(s) -
p^{2N}_{i}(s)\right)\right. \\&&\left. +p^{2N}_{i}(s)\left(  \sum_{j=1}^{K}\sigma_{ij}p^{2N}_j(s)-\sum_{k,\ell}%
^{K}\sigma_{k\ell}p^{2N}_{k}(s)p^{2N}_\ell(s)\right)\  \right] ds^N +o(\frac{1}{N})\nonumber
\eea
is a martingale and for $i,j=1,\dots,K$
\bea{} &&E[(M^{2N}_i(t_2)-M^{2N}_i(t_1))(M^{2N}_j(t_2)-M^{2N}_j(t_1))]\\&&=\frac{1}{2N} E\sum_{k={\lfloor 2Nt_1\rfloor}}^{{\lfloor 2Nt_2\rfloor}}
p^{2N}_i(\frac{k}{2N})(\delta_{ij}-p^{2N}_j(\frac{k}{2N})) +o(\frac{1}{N}).\nonumber
\eea


Step 3.  Proof that there exists a unique probability measure on $C_{\Delta_{K-1}}([0,\infty))$ such that
(\ref{mpa}) and (\ref{mpb}) are satisfied.

Uniqueness can be proved in the neutral case, $\sigma\equiv0$, by showing that
moments are obtained as unique solutions of a closed system of differential equations.

\end{proof}
\bigskip
\begin{remark} The uniqueness when $\sigma$ is not zero follows from the dual representation developed in the next chapter.
\end{remark}

\section{Stationary measures}

A special case of a theorem in Section 8.3 implies that if the matrix $(q_{ij})$ is irreducible, then the Wright-Fisher  diffusion is
ergodic with unique stationary distribution.

\subsection{The Invariant Measure for the neutral K-alleles WF Diffusion}
Consider the neutral $K$-type Wright-Fisher diffusion with type-independent mutation (Kingman's
``house-of-cards'' mutation model) with generator

 \bea{IMAMP}
&&G^{K}f(p)   =\frac{1}{2}\sum_{i,j=1}^{K-1}p_{i}(\delta_{ij}-p_{j}%
)\frac{\partial^{2}f(p)}{\partial p_{i}\partial p_{j}}+\frac{\theta}{2}\sum_{i=1}^{K-1}\left(\nu_i -
p_{i}   \right)  \frac{\partial
f(p)}{\partial
p_{i}}.\nonumber%
\eea
where the {\em type-independent mutation} kernel is given by $\nu\in \Delta_{k-1}$.

\begin{theorem}\label{IMAEQ} (Wright \cite{W-49}, Griffiths \cite{G-79})
The Dirichlet distribution $D(p_{1},\dots,p_{n})$ on $\Delta_{K-1}$ with
density
\begin{align*}
\Pi_{K}(dp)  &  =\frac{\Gamma(\theta_{1}+\dots+\theta_{K})}{\Gamma(\theta
_{1})\dots\Gamma(\theta_{K})}p_{1}^{\theta_{1}-1}\dots p_{K}^{\theta_{K}%
-1}dp_{1}\dots dp_{K-1}\\
\theta_{j}  &  =\theta\nu_{j},\;\nu\in\mathcal{P}(1,\dots,K)%
\end{align*}
is a reversible stationary measure for the neutral K-alleles WF diffusion with $\gamma=1$.

In the case $K=2$ this is the Beta distribution  \be{}
\frac{\Gamma(\theta)}{\Gamma(\theta_1)\Gamma(\theta_2)}x_1^{\theta_1-1}(1-x_1)^{\theta_2-1}dx_1.\ee
\end{theorem}

\begin{proof} (cf. \cite{EK-81})
Reversibility and stationarity means that when $\Pi_K$ is the initial distribution,
then $\{p(t):0\leq t\leq t_{0}\}$ has the same distribution as
$\{p(t_{0}-t):0\leq t\leq t_{0}\}$.
In terms of the strongly
continuous semigroup $\{T(t)\}$ on $C(\Delta_{K-1})$ generated by
$G$ a necessary and sufficient condition (see Fukushima and
Stroock (1986) \cite{FS-86})  for
reversibility with respect to $\Pi_K$ is that%
\[
\int g\,T(t)\,fd\Pi_K=\int f\,T(t)\,gd\Pi_K\;\forall\;f,g\in
C(\Delta_{K-1}),\;t\geq0
\]
or equivalently that%
\[
\int gGfd\Pi_K=\int fGgd\Pi_K\;\forall\;f,g\in D(G)
\]
or for $f,g$ in a core for $G$ (see Appendix I).

Since the space of polynomials in $p_{1}%
,\dots,p_{K}$ is a core for $G$ it suffices by linearity to show that%
\[
\int gGfd\Pi=\int fGgd\Pi\quad \forall\;f=f_{\alpha},g=f_{\beta}
\]
where $f_{\alpha}=p_{1}^{\alpha_{1}}\dots p_{K}^{\alpha_{K}}$. \ Let
 $|\alpha|=\sum\alpha_{i}.$

Then%
\begin{align*}
&  \int f_{\beta}Gf_{\alpha}d\Pi_K\\
&=\frac{1}{2}\int [\sum_{i=1}^K\alpha_i(\alpha_i+\theta_i-1)f_{\alpha+\beta-e^i}
-|\alpha|(|\alpha|+\sum_{i=1}^K\theta_i-1)f_{\alpha+\beta}]d\Pi_K\\
&  =\frac{1}{2}\left\{  \sum_{i=1}^{K}\frac{\alpha_{i}(\alpha_{i}+\theta_{i}%
-1)}{\alpha_{i}+\beta_{i}+\theta_{i}-1}-\frac{|\alpha|(|\alpha|+\sum\theta
_{i}-1)}{|\alpha|+|\beta|+\sum\theta_{i}-1}\right\} \\
&  \cdot\frac{\Gamma(\alpha_{1}+\beta_{1}+\theta_{1})\dots\Gamma(\alpha
_{K}+\beta_{K}+\theta_{K})}{\Gamma(|\alpha|+|\beta|+\sum\theta_{i}%
-1)}\frac{\Gamma(\sum\theta_{i})}{\Gamma(\theta_{1})\dots\Gamma(\theta_{K})}.
\end{align*}
To show that this is symmetric in $\alpha,\beta$, let $h(\alpha,\beta)$ denote
the expression within $\{...\}$ above. Then%
\begin{align*}
&  h(\alpha,\beta)-h(\beta,\alpha)\\
&  =\sum\frac{\alpha_{i}^{2}-\beta_{i}^{2}+(\alpha_{i}-\beta_{i})(\theta
_{i}-1)}{\alpha_{i}+\beta_{i}+\theta_{i}-1}-\frac{|\alpha|^{2}-|\beta
|^{2}+(|\alpha|-|\beta|)(\sum\theta_i-1)}{|\alpha|+|\beta|+\sum\theta_{i}-1}\\
&  =\sum(\alpha_{i}-\beta_{i})-(|\alpha|-|\beta|)\\
&  =0
\end{align*}
\end{proof}

\begin{corollary}
Consider the  mixed moments:%
\[
m_{k_{1},\dots k_{K}}=\int\dots\int_{\Delta_{K-1}}p_{1}^{k_{1}}\dots
p_{K}^{k_{K}}\Pi_{K}(dp)
\]
Then%
\[
m_{k_{1},\dots k_{K}}=\frac{\Gamma(\theta_{1})\dots\Gamma(\theta_{K})}%
{\Gamma(\theta_{1}+\dots+\theta_{K})}\frac{\Gamma(\theta_{1}+\dots+\theta
_{K}+k_{1}+\dots+k_{K}))}{\Gamma(\theta_{1}+k_{1})\dots\Gamma(\theta_{K}%
+k_{K})}.
\]
\end{corollary}

\subsubsection{Stationary measure with selection}
If selection (as in (\ref{KWFG}) is added then
the stationary distribution is given by the ``Gibbs-like'' distribution%
\be{sm}
\Pi_\sigma(dp)=C\exp\left(
\sum_{i,j=1}^{K}\sigma_{ij}p_{i}p_{j}\right) \Pi_K( dp_{1}\dots dp_{K-1})
\ee
and this is reversible. (This is a special case of a result that will be
proved in a later section.)

\subsection{Convergence of stationary measures of $\{p^N\}_{N\in\N}$}

It is of interest to consider the convergence  of the stationary measures of the Wright-Fisher Markov chains to (\ref{sm}). A standard
argument applied to the Wright-Fisher model is as follows.

\begin{theorem}
{\em Convergence of Stationary Measures.}  Assume that the diffusion limit, $p(t),$ has a unique
invariant measure, $\nu$ and that $\nu_{N}$ is an invariant
measure for $p^{N}(t).$ Then%
\be{}
\nu_{N}\Longrightarrow\nu\text{   as   }N\to\infty.
\ee
\end{theorem}

\begin{proof} Denote by $\{T_t\}_{t\geq 0}$ the semigroup of the Wright-Fisher diffusion.
Since the state space is compact, the space of probability measure is compact.
and therefore the sequence  $\nu_{N}$ is tight $M_{1}(\Delta_{K-1})$. Given a limit point $\tilde{\nu}$ and a subsequence $\nu_{N^{\prime}}$
that converges weakly to $\tilde{\nu}\in M_{1}(\Delta_{K-1}) $ it follows that for $f\in
C(\Delta_{K-1}),$%
\begin{align*}
\int T(t)fd\tilde{\nu}  &  =\lim_{N^{\prime}\rightarrow\infty}\int
T(t)fd\nu_{N^{\prime}}\;(\text{by }\nu_{N^{\prime}}\Longrightarrow\nu)\\
&  =\lim_{N^{\prime}\rightarrow\infty}\int T_{N^{\prime}}(2N^\prime t)fd\nu_{N^{\prime}%
}\;\;(\text{by }p_{N}\Longrightarrow p)\\
&  =\lim_{N^{\prime}\rightarrow\infty}\int fd\nu_{N^{\prime}}\text{ \ (by inv.
of }\nu_{N^{\prime}}\text{)}\\
&  =\int fd\tilde{\nu}\text{ (by }\nu_{N^{\prime}}\text{ }\Longrightarrow
\tilde{\nu}).
\end{align*}
Therefore $\tilde{\nu}$ is invariant for $\{T(t)\}$ and hence
$\tilde{\nu}=\nu$ by assumption of the uniqueness of the invariant measure for
$p(t)$. That is, any limit points of $\{\nu_{N}\}$ coincides with $\nu$ and
therefore $\nu_{N}\Longrightarrow\nu$.
\end{proof}

\subsubsection{Properties of the Dirichlet Distribution}

1. Consistency under merging  of types.

Under $D(\theta_{1},\dots,\theta_{n}),$ the distribution of $(X_{1},\dots,X_{k}%
,1-\Sigma_{i=1}^{k}X_{i})$ is
\[
D(\theta_{1},\dots,\theta_{k},\theta_{k+1}+\dots+\theta_{n})
\]
and the distribution of $\frac{X_{k}}{1-\Sigma_{i=1}^{k-1}X_{i}}=$
$\frac{X_{k}}{\Sigma_{i=k}^{K}X_{i}}$is $Beta(\theta_{k},\Sigma_{i=k+1}^{K}%
\theta_{i}). $

2. {Bayes posterior under random sampling}

Consider the n-dimensional Dirichlet distribution, $\mathcal{D}$($\alpha)$
with parameters $(\alpha_{1},\dots,\alpha_{n}).$ Assume that some phenomena is described
by a random probability vector $p=(p_{1},\dots,p_{n})$. Let \ $\mathcal{D}%
(\alpha)$ be the ``prior distribution of the vector $p$. Now let us assume
that we take a sample and observe that $N_{i}$ of the outcome are $i$. Now
compute the posterior distribution of $p$ given the observations $N=(N_{1}%
,\dots,N_{n})$ as follows: Using properties of the Dirichlet distribution we
can show that it is \ \
\begin{align*}
P(p  &  \in dx|N)=\frac{1}{Z}\frac{x_{1}^{\alpha_{1}}\dots x_{n}^{\alpha_{n}%
}x_{1}^{N_{1}}\dots x_{n}^{N_{n}}}{\int x_{1}^{\alpha_{1}}\dots x_{n}%
^{\alpha_{n}}x_{1}^{N_{1}}\dots x_{n}^{N_{n}}dx_{1}\dots dx_{n}}\\
&  =\frac{1}{Z^{\prime}}x_{1}^{(\alpha_{1}+N_{1})}\dots x_{n}^{(\alpha
_{n}+N_{n})}.
\end{align*}
That is, \be{bayes}P(p\in\cdot|N)\text{ is  }\mathcal{D(}\alpha_{1}+N_{1},\dots,\alpha
_{n}+N_{n}).\ee


\chapter{Infinitely many types models}


\section{Introduction}

\subsection{Motivation}
We have considered particle systems \ above with finitely many
types. In the 1970's with the advent of electrophoresis and
molecular biology, new models were needed in which the number of
types were not fixed.\ In many cases the number of types can be
random and new types can be introduced at random times. Several
models began to appear at that time involving infinitely many
types, for example the {\em ladder} or {\em stepwise mutation model} of Ohta and Kimura (1973) \cite{OK-73} (which could model for
example continuous characteristics). Another model was one in
which no attempt to model the structure of types was made but in
which new types can be introduced (leading to the
\textit{infinitely many alleles} model) (Kimura and Crow (1964) \cite{KC-64}). In this model we take
$[0,1]$ as the type space. Then when a new type is needed we can
choose a type in $[0,1]$ by sampling from the uniform distribution
on $[0,1].$    The {\em infinitely many sites model} introduced by  Kimura  in  1969  provides an idealization of the genome viewed as a sequence of  nucleotides (A,T,C,G). These processes now form the basis for molecular population genetics.

More generally, such infinitely many type models provide the
possibility of coding information at a number of levels and
provide a powerful tool for the study of complex systems. \ For
example we can code historical information, genealogical
information, and information about the random environment that has
been visited. \ In addition it allows for individuals with
internal structure described by an internal state space and state
transition dynamics.

\subsection{Plan}
The objective of this chapter is to construct two basic infinitely-many-type processes, formulated as measure-valued processes,
 by taking the projective limit of the finite type Feller CSB and Wright-Fisher diffusions.  These are the {\em Jirina measure-valued branching process} and {\em infinitely many alleles model} of Crow and Kimura.  The latter has played a central role in population genetics.  We will establish a relation between the invariant measures of these two processes that allows us to obtain the basic properties of the {\em Poisson-Dirichlet} distribution and the {\em Griffiths, Engen and McCloskey (GEM) representation}.

 We begin by considering the diffusion limit of a measure-valued generalization of the Wright-Fisher Markov chain in the setting of semigroup theory. This process, the {\em Fleming-Viot process},   includes the infinitely many alleles model as a special case.
In the next Chapter we reformulate these processes in terms of measure-valued martingale problems and develop techniques for working with  more general classes of measure-valued processes including the class of superprocesses and the class of Fleming-Viot processes with selection, mutation  and recombination.

\section{Measure-valued Wright-Fisher Markov chain}

We now consider a Wright-Fisher model of a population of $N$ individuals in which the space of types is a separable metric space  $E$.
The process is then a Markov chain $\{p^N_n\}_{n\in \N}$   with state space%
\[
\mathcal{P}^{N}(E)=\left\{  \mu=\frac{1}{N}\sum_{i=1}^{N}\delta_{x_{i}}%
,\;\{x_{1},\dots,x_{N}\}\in E\right\}\subset \mathcal{P}(E).
\]

In this case the mutation process is a Markov chain on $E$ with   probability transition function $P(x,dy)$  giving
the type distribution of the offspring of a type $x$ parent if
mutation occurs. Let $V>0$ be a measurable function on $E$ with
$V(x)$ interpreted as the (haploid) fitness of a type $x$ individual.

Then as in the finitely many type case $X^N_n$  is a Markov chain in
$\mathcal{P}^{N}(E)$ with one step transition function
$P(\mu,d\mu^{\rm{next}})$.  This is obtained by noting that
$X^N_{n+1}$ is a random probability measure on
$\mathcal{P}^N(E)$ given by:
\be{WFMC} X^N_{n+1}=\frac{1}{N}\sum_{i=1}^N\delta_{y_i}\ee where
$y_1,\dots,y_N$ are i.i.d. $\mu_n^*(dy)$ where
\[
\mu^{\ast}_n(dy)= \int_E\left(  \frac{V(x)X^N_n(dx)}{\int V(x)X^N_n(dx)
}\right)  P(x,dy),
\]
that is, as before selection first and then mutation and sampling.

\begin{example}  Infinitely many alleles model \cite{KC-64}. $E=[0,1]$ and
\be{}P(x,dz)=(1-m)\delta_{x}(dz)+m\int_0^1\delta_{y}(dz)\lambda(dy)\ee
where $\lambda$ is Lebesgue measure on $[0,1]$.
\end{example}

\begin{example}  The infinitely many sites model was introduced by Kimura \cite{K-69}, \cite{K-71}. (See also Ethier and Griffiths (1987) \cite{EG-87})

Infinitely many sites model. $E=[0,1]^{\mathbb{Z}_{+}}$%
\[
P(\mathbf{x},dy)=(1-m)\delta_{\mathbf{x}}(dy)+m\int_{0}^{1}\delta
_{\{\xi,\mathbf{x)}}\lambda(d\xi)
\]
Here we interpret $\xi$ as the locus on the genome where the last mutation occurred.\bigskip

The number of segregating sites is the number of homologous DNA positions that differ in a sample of $m$ sequences.
They are used  to investigate phylogenetic relationships.  The location of polymorphisms within humans  are also used to determine  the
potential differences in reactions of individuals to medical treatments.

See Section 8.3.3  for the analysis of segregating sites.

\end{example}
\begin{example} Ladder model of Ohta and Kimura (1973) \cite{OK-73}. This stepwise-mutation model was introduced to describe the distribution of allelic types distinguishable as signed electrical charges in gel electrophoresis experiments. \\
Here $E=\mathbb{Z}$ and

\be{}P(x,dy)=(1-m)\delta_{x}(\cdot)+\frac{m}{2}\delta_{x+1}(\cdot)+\frac{m}{2}\delta_{x-1}(\cdot).\ee
\end{example}

\section{The neutral Fleming-Viot process with mutation generator $A$}

The infinitely many alleles diffusion of Crow and Kimura can be studied as an infinite dimensional diffusion (i.e. countably many types) (see Ethier and Kurtz (1981) \cite{EK-81}) in which a mutation always leads to a new type. However it is advantageous to reformulate it as a measure-valued process.
This process was introduced by Fleming and Viot in 1979
\cite{FV-79}. We will show that it arises as the diffusion limit of the
measure-valued Wright-Fisher model.

We will now  we derive the Fleming-Viot process under some simplifying assumptions using semigroup methods. The general case will be dealt with below in the martingale problem setting.

\textbf{Assumptions}
\begin{itemize}
\item Let $E$ be a compact metric space.
\item $V\equiv 1$, that is, we omit the selection effect.
\item  We consider a mutation process
given by a Feller process on $E$ with generator $(D(A),A)$ and semigroup $\{S_t:t\geq 0\}$ on $C(E)$.
$A$ will be called the mutation operator for the
Fleming-Viot process.

 We assume
that $D(A)$ contains an algebra $D_0(A)$ that separates points and $S_t:D_0(A)\to D_0(A)$.
 Then linear
combinations of functions in $D_0(A)$ form an algebra of functions
separating points and therefore is dense in $C(E)$
and therefore measure-determining.  We also assume that $A$ arises as the limit of a sequence of mutation Markov chains on $E$ with transition kernels $\{P_N(.,.)\}_{N\in\N}$, that is for $f\in D(A),$
\be{}
N(\left\langle f,P_{N}\right\rangle -f)\rightarrow Af\text{ as }%
N\rightarrow\infty
\ee
uniformly on $E$.

\end{itemize}

The state space for the Fleming-Viot (FV) process is $\mathcal{P}(E)$, the set of Borel probability
measures on $E$ with the topology of weak convergence.
For $f\in C(E)$, $\mu\in \mathcal{P}$ we denote $\left\langle
f,\mu\right\rangle =\int fd\mu$.

We will now obtain the neutral FV process as the limit of neutral (i.e. $V\equiv1$)  Wright-Fisher Markov chains in the diffusion time scale, \be{} p^N(t)= X^N_{\lfloor Nt\rfloor}\;\in \mathcal{P}(E)\ee where $X^N_n$ is defined by (\ref{WFMC}) with
mutation transition functions $P^{N}(x,dy)$.


 \ In order
to identify the limiting generator for a $\mathcal{P}(E)$-valued diffusion we need a measure-determining
family of test functions. Consider the algebra $\mathcal{D}$ of
nice
functions on $\mathcal{P}(E)$ containing the functions:%
\be{tf} F(\mu)=\left\langle f_{1},\mu\right\rangle
\dots\left\langle f_{n},\mu\right\rangle \ee with $n\geq1$ and
$f_{1},\dots,f_{n}\in D(A).$ This algebra of functions is measure-determining in $\mathcal{P}(\mathcal{P}(E))$.

\begin{notation}
For $F\in \mathcal{D},\;x\in E$ we define
\begin{align*}
\frac{\partial F(\mu)}{\partial\mu(x)}  &  =\lim_{\varepsilon
\rightarrow 0}\frac{F(\mu+\varepsilon\delta_{x})-F(\mu
)}{\varepsilon}|_{\varepsilon=0}=\frac{\partial F(\mu+\ve\delta_x)}{\partial \ve}|_{\ve =0}\\
\frac{\partial^{2}F(\mu)}{\partial\mu(x)\partial\mu(y)}  &
=\frac{\partial^2 F(\mu+\ve_1\delta_x +\ve_2 \delta_y)}{\partial \ve_1 \partial \ve_2}|_{\ve_1=\ve_2 =0}
\end{align*}
\end{notation}

\begin{proposition}  Let $p^N_n$ denote the measure-valued Wright-Fisher Markov chain (\ref{WFMC}) under the above assumptions. Then $p^N(t)=X^N_{\lfloor Nt\rfloor}\Rightarrow p_t$ where $\{p_t\}_{t\geq 0}  $ is a $\mathcal{P}(E)$-valued Markov process with generator
\begin{align*}
&  GF(\mu)\\
&  =\sum_{1\leq i<j\leq n}(\left\langle f_{i}f_{j},\mu\right\rangle
-\left\langle f_{i},\mu\right\rangle \left\langle f_{j},\mu\right\rangle
)\prod_{\ell:\ell\neq i,j}\left\langle f_{\ell},\mu\right\rangle +\sum
_{i}\left\langle Af_{i},\mu\right\rangle \prod_{\ell:\ell\neq i}\left\langle
f_{\ell},\mu\right\rangle \\
&  =\frac{1}{2}\left[  \int\frac{\partial^{2}F(\mu)}{\partial
\mu(x)\partial\mu(y)}\delta_{x}(dy)\mu(dx)-\int\frac{\partial^{2}F(\mu
)}{\partial\mu(x)\partial\mu(y)}\mu(dx)\mu(dy)\right] \\
&  +\int A\frac{\partial F(\mu)}{\partial \mu(x)}\mu(dx).
\end{align*}
for all $F\in \mathcal{D}$.
\end{proposition}

\begin{proof} Here we follow the Ethier-Kurtz semigroup approach.  Using the Kurtz semigroup convergence theorem (\cite{EK-86}, Chap. 1, Theorem 6.5, Proposition 3.7 and Chap. 4. Theorem 2.5  -see  Appendix III Theorems \ref{SGCT}, \ref{SGCT2}).
 Using these results it suffices to show that for $F\in\mathcal{D}$,

\be{KCO}\lim_{N\to\infty}
NE_{\mu}[F(p^N_{\frac{1}{N}})-F(\mu))]=
\lim_{N\to\infty}
NE_{\mu}[F(X^N_1)-F(\mu)]=GF(\mu)
\ee
uniformly in $\mu\in \mathcal{P}(E)$.

First note that for $f_{1},\dots,f_{n}$ $\in C(E)$%
\[
E_\mu(F(X^N_1))=E_{\mu}[\left\langle f_{1}, X^N_1\right\rangle \dots\langle f_{n}%
,X^N_1\rangle ],\quad F(\mu)=  \langle f_{1},\mu\rangle
\dots\langle f_{n},\mu\rangle
\]
where $X^{N}_1=\frac{1}{N}\sum_{i=1}^{N}\delta_{Y_{i}}$ and  $Y_{1}%
,\dots,Y_{N}$ are i.i.d.  $\mu P^N.$ Hence
\begin{align*}
E_{\mu}[\left\langle f_{1},X^N_1\right\rangle \dots\left\langle f_{n}%
,X^N_1\right\rangle ]  &  =\frac{1}{N^{n}}E\left[  \sum_{i=1}^{N}f_{1}%
(Y_{i})\dots\sum_{i=1}^{N}f_{n}(Y_{i})\right] \\
&  =\sum_{k=1}^{n}\frac{N^{[k]}}{N^{n}}\sum_{\beta\in\pi(n,k)}\prod_{j=1}%
^{k}\left\langle \prod_{i\in\beta_{j}}f_{i},\mu P^N\right\rangle \\
&  =\sum_{k=1}^{n}\frac{N^{[k]}}{N^{n}}\sum_{\beta\in\pi(n,k)}\prod_{j=1}%
^{k}\left\langle \left\langle \prod_{i\in\beta_{j}}f_{i},P^N\right\rangle
,\mu\right\rangle
\end{align*}
where $N^{[k]}=\frac{N!}{(N-k)!},\;\pi(n,k)$ is the set of partitions $\beta$
of $\{1,\dots,n\}$ into $k$ nonempty subsets $\beta_{1},\dots,\beta_{k}$,
labelled so that $\min\beta_{1}<\dots<\min\beta_{k}$.

Only the terms involving  $k=n,n-1$ contribute in the limit. To see this note that we
 can choose $n$ different $Y_i$'s in $N(N-1)\dots (N-n+1)= N^n-\frac{n(n-1)}{2}N^{n-1}+O(N^{n-2})$ ways
and $n-1$ different $Y_i$'s in $N(N-1)\dots (N-n+2)=N^{n-1}-O(N^{n-2})$ ways.  For $k= n-2$
we can choose $k$ different $Y_i$'s is $O(N^{n-2})$ ways, etc.

\begin{align*}
&  NE_{\mu}[F(X^N_{1})-F(\mu)]\\
&  =N\left\{  \frac{N^{[n]}}{N^{n}}\prod_{j=1}^{n}\left\langle f_{j},\mu
P_{N}\right\rangle +\frac{N^{[n-1]}}{N^{n}}\sum_{1\leq i<j\leq n}\left\langle
f_{i}f_{j},\mu P_{N}\right\rangle \prod_{\ell:\ell\neq i,j}\left\langle
f_{\ell},\mu P_{N}\right\rangle \right. \\
&  +\left.  O(N^{-2})-\prod_{j=1}^{n}\left\langle f_{j},\mu\right\rangle
\right\} \\
&  =N\left\{  \left(  1-\frac{n(n-1)}{2N}\right)  \prod_{j=1}^{n}\left\langle
f_{j},\mu P_{N}\right\rangle \right. \\
&  +\left.  \frac{1}{N}\sum_{1\leq i<j\leq n}\left\langle f_{i}f_{j},\mu
P_{N}\right\rangle \prod_{\ell\neq i,j}\left\langle f_{\ell},\mu
P_{N}\right\rangle -\prod_{j=1}^{n}\left\langle f_{j},\mu\right\rangle
\right\}  +O(\frac{1}{N})\\
\end{align*}
Note that $\lim_{N\to\infty} \langle f,\mu P_N\rangle =\langle f,\mu\rangle$. Now
let $b_j=\langle f_j,\mu\rangle$ and $a_j=\langle f_j,\mu P_N\rangle$ and recall that
\be{}\lim_{N\to\infty} N(\langle f,\mu P_N\rangle-\langle f,\mu\rangle)=\langle Af,\mu\rangle\ee
Then using this together with the collapsing sum
\begin{align*}
&a_{1}\dots a_{n}+(a_{1}\dots a_{n-1}b_{n}-a_{1}\dots a_{n})    +(a_{1}\dots
a_{n-2}b_{n-1}b_{n}-a_{1}\dots a_{n-1}b_{n})\\
&+(b_{1}\dots b_{n}-a_{1}b_{2}\dots b_{n})-b_{1}\dots b_{n}   =0\\ &
a_{1}\dots a_{n}-b_{1}\dots b_{n}=\sum_k(a_{1}\dots a_{k}b_{k+1}\dots
b_{n}-a_{1}\dots a_{k+1}b_{k+2}\dots b_{n}).
\end{align*}%

or rewriting

\be{}\prod_{j=1}^n \langle f_j,\mu P_N\rangle = \prod_{j=1}^n \left[  \left(\langle f_j,\mu\rangle +(\langle f_j,\mu P_N\rangle-\langle f_j,\mu\rangle) \right)\right]\ee

 we obtain

\begin{align*} NE_{\mu}[F(X^N_{1})-F(\mu)]
&  =\sum_{1\leq i<j\leq n}(\left\langle f_{i}f_{j},\mu P_{N}\right\rangle
-\left\langle f_{i},\mu P_{N}\right\rangle \left\langle f_{j},\mu
P_{N}\right\rangle )\prod_{\ell:\ell\neq i,j}\left\langle f_{\ell},\mu
P_{N}\right\rangle \\
&  +\sum_{i=1}^{n}\left\langle Af_{i},\mu\right\rangle \prod_{j:j<i}%
\left\langle f_{j},\mu\right\rangle \prod_{j:j>i}\left\langle f_{j},\mu
P_{N}\right\rangle +O(N^{-1})\\
&  =GF(\mu)+o(1)
\end{align*}
uniformly in $\mu$.

The completes the verification of condition (\ref{KCO}).
\end{proof}



\section{The Infinitely Many Alleles Model}

This is a special case of the Fleming-Viot process which has
played a crucial role in modern population biology. It has type
space $E=[0,1]$ and
{\em type-independent}  mutation operator with mutation source  $\nu_0\in\mathcal{P}([0,1])$
\begin{align*}
Af(x)  &=  \theta(\int p(x,dy)f(y)-f(x))\\
&  =\theta (\int f(y)\nu_0(dy)-f(x)).
\end{align*}
Since $A$ is a bounded operator we can take indicator functions of intervals in $D(A)$. If we have a partition
$[0,1]=\cup_{j=1}^K B_j$ where the $B_j$ are intervals, consider the set $D(G)$ of functions
\be{tf2} F(\mu)=\left\langle f_{1},\mu\right\rangle
\dots\left\langle f_{n},\mu\right\rangle \ee with $n\geq1$ and where
the functions $f_{1},\dots,f_{n}$ are finite linear combinations of indicator functions of the intervals $\{A_j\}$.  Then  the function
$GF(\mu)$ can be written in the same form and   we can prove that the $\Delta_{K-1}$-valued process $\{p_t(A_1),\dots,p_t(A_K)\}$ is a version of the $K-allele$ process with generator
\bea{KWFG2}
&&G^{K}f(p) \\ &&  =\frac{1}{2}\sum_{i,j=1}^{K-1}p_{i}(\delta_{ij}-p_{j}%
)\frac{\partial^{2}f(p)}{\partial p_{i}\partial p_{j}}  +\theta \sum_{i=1}^{K-1}(\nu_0(A_i) -
p_{i})    \frac{\partial
f(p)}{\partial
p_{i}}.\nonumber%
\eea

We will next give an explicit construction of this process that allows
us to derive a number of interesting properties of this important
model.

\subsection{Projective Limit Construction of the Infinitely Many Alleles Model}

Let  $\mu, \nu_0\in\mathcal{P}(E)$, $\mathcal{C}= C_{[0,\infty)}([0,\infty))$.  Let  $U$ denote the collection of finite partitions $u=(A^u_1,\dots,A^u_{|u|})$ of $E$ into measurable
subsets in $\mathcal{B}(E)$ and $|u|$ denotes the number of sets in the partition $u$. We place a partial ordering on $U$ as follows:%
\[
v\succ u
\]
 if $v$ is a refinement of $u$.  We can also identify partitions with the
finite algebras of subsets of $E$ they generate.  Given a partition we define the probability measure, $P_u$ on $\mathcal{C}^u$ as the law of the Wright-Fisher diffusion with generator

\begin{align*}
G^{(K)}f(p)  &  =\frac{1}{2}\sum_{i,j=1}^{K-1}p_{i}(\delta_{ij}-p_{j}%
)\frac{\partial^{2}f(p)}{\partial p_{i}\partial p_{j}}\\
&  +\frac{1}{2}\sum_{i=1}^{K-1}\theta(\nu_i-p_i)\frac{\partial
f(p)}{\partial p_{i}} \\ \nu_i &  :=\nu_{0}(A_{j})
\end{align*}
and initial measure $\mu$, that is, the law of $(p_t(A^u_1),\dots,p_t(A^u_{|u|})$ (and the additive extension of this to the algebra generated by $u$).

\begin{remark} Recall that the associated Markov transition function is determined by the joint moments as follows.

Since the family of functions  $p_1^{k_1}\dots,p^{k_{K-1}}_{K-1}$
belong to $D(G^{(K)})$ we can apply $G^{(K)}$ and obtain the following
system of equations for the joint moments:
\be{} m_{k_1,\dots,k_{K-1}}(t):= E[p_{1}^{k_{1}}(t)\dots p_{K-1}^{k_{K-1}}(t)], \ee
\begin{align*}
\frac{\partial}{\partial t}m_{k_1,\dots,k_{K-1}}(t)  &
=\frac{1}{2}\sum_{i}k_{i}(k_{i}-1)m_{k_{1},.,k_{i}-1,..k_{K-1}}(t)\\
&  -\frac{1}{2}\sum_{i\neq j}k_{i}k_{j}m_{k_{1},\dots k_{K}}(t)\\
&  +\frac{\theta}{2}\sum_{i=1}^{K-1}\nu_{i}k_{i}m_{k_{1}%
,.k_{j}-1,.,k_{i}+1,..k_{K}}(t)\\
&  -\frac{\theta}{2}\sum_{i=1}^{K-1}k_{i}m_{k_{1},\dots k_{K-1}}(t)
\end{align*}
Since this system of linear equations is closed, there exists a unique solution which characterizes the $K$-allele Wright-Fisher diffusion.
\end{remark}

In a similar way  we can apply this to the function corresponding to the coalescence of two partition elements
\begin{align*}
f(p)  &  =\tilde{f}(\tilde{p})\\
\tilde{p}  &  =(\tilde{p}_{1},\dots,\tilde{p}_{K-1})\\
&  =(p_{1},\dots,p_{\ell-1},p_{\ell+1},\dots,p_{k-1},p_{k+1},\dots
,p_{K-1},(p_{\ell}+p_{k}))
\end{align*}
\bigskip%
\begin{align*}
G^{(K)}f(p)  &  =\frac{1}{2}\sum_{i,j=1}^{K-2}\tilde{p}_{i}(\delta_{ij}%
-\tilde{p}_{j})\frac{\partial^{2}f(\tilde{p})}{\partial\tilde{p}_{i}%
\partial\tilde{p}_{j}}\\
&  +\frac{\theta}{2}\sum_{i=1}^{K-2}(\wt\nu_i-\wt p_i)\frac{\partial\tilde{f}(\tilde
{p})}{\partial\tilde{p}_{i}}\\
&  =G^{(K-1)}\tilde{f}(\tilde{p})
\end{align*}
In other words we have consistency under coalescence of the partition
elements. Because of uniqueness this implies that the process $\tilde
{p}(t)=(\tilde{p}_{1}(t),\dots,\tilde{p}_{K-1}(t))$ coincides with the
$(K-1)-$allele Wright-Fisher diffusion.

We  denote the canonical projections $\pi_u:\mathcal{C}^{\mathcal{B}(E)}\to
\mathcal{C}^{u}$ and $\pi_{uv}:\mathcal{C}^{v}\to
\mathcal{C}^{u}$ if $v\succ u$ such that

$\pi_u=\pi_{uv}\pi_v,\; v\succ u$.

The family $\{P_u\}_{u\in U}$ forms a projective system of probability laws, that is for every pair, $(u,v)$, $v\succ u$,
$\{P_{u}\}$ then satisfies
\be{} \pi_{uv}(P_v)=P_u,\qquad P_u(B)=P_v(\pi^{-1}_{uv}(B)).\ee

Therefore,  by Theorem \ref{PLT} (in Appendix I) there exists a projective limit measure, that is, a probability measure  \ $P_\infty$ on
$\mathcal{C}^{\mathcal{B}([0,1])}$ such that for any $u\in U$, $\pi_u P_\infty  =P_u$.

For fixed $t$ (or any finite set of times)  we can identify the projective limit, \be{}\{\tilde{p}_{t}%
(A):A\in\mathcal{B}([0,1])\}\ee  with an element of $\mathcal{X}([0,1]),$ the space
of all finitely additive, non-negative, mass one measures on
$[0,1],$ equipped with the projective limit topology, i.e., the
weakest topology such that for all Borel subset $B$ of $[0,1]$,
$\mu(B)$ is continuous in $\mu$. Under this topology,
$\mathcal{X}([0,1])$ is Hausdorff.
The $\sigma$-algebra
$\mathcal{B}$ of the space $\mathcal{X}([0,1])$ is the smallest
$\sigma $-algebra such that for all Borel subset $B$ of $[0,1]$,
$\mu(B)$ is a measurable function of $\mu$.

For fixed
$t\in[0,\infty),\; \tilde{p}_{t}(\cdot)$ is a.s. a finitely
additive measure, that is, a member of $\mathcal{X}[0,1]$ and satisfies the conditions of Theorem \ref{PRPOL} in the Appendices (conditions 1,2 follow immediately from the construction, 3 follows since for any $A\in\mathcal{B}([0,1])$ $E(p_t(A))\leq \max(\mu(A),\nu_0(A))$ and (4) is automatic since all measures are bounded by 1).
 Therefore for fixed $t$  this determines a unique countably
additive version $p_{t}(\cdot)$, that is, a random countable additive measure $p_t\in \mathcal{P}([0,1])$ a.s. Similarly, taking two times $t_1,t_2$ we obtain a the joint distribution of a pair $(p_{t_1},p_{t_2})$ of random probability measures. We can then verify that $t\to \int f(x)p_t(dx)$ is a.s. continuous for countable convergence determining class of functions so that there is an a.s. continuous version with respect to the topology of weak convergence.

\begin{remark} We can carry out the same construction assuming that  for each $u\in U$ the Wright-Fisher diffusion starts with the stationary Dirichlet measure and obtain by the projective limit a probability measure on $\mathcal{P}(E)$ which for any partition has the associated Dirichlet distribution.
\end{remark}

\section{The Jirina  Branching Process}

In 1964 Jirina \cite{J-64} gave the first construction of a
measure-valued branching process. The state space is the space of
finite measures on $[0,1],\;$ $M_{f}([0,1]).\;\nu_{0}\in
M_{1}([0,1])$. \ We will construct a version of this process with
immigration by a projective limit construction.

Given a partition $(A_1,\dots,A_K)$ of $[0,1]$ let
 $\{X_{t}(A_i):t\geq0,\;i=1,\dots,K\}$ satisfy the SDE (Feller CSB plus immigration):
\bea{JSDE}
\;\; && dX_{t}(A_i)   =c(\nu_{0}(A_i)-X_{t}(A_i))dt+\sqrt{2\gamma X_{t}(A_i)}dW_{t}^{A_i}\\
&&X_{0}(A_i)    =\mu(A_i)\nonumber
\eea
where $\nu_{0}$ is in $\mathcal{P}([0,1])$ and for each $i$, $W^{A_i}_t$ is a standard Brownian motion and for $i\ne j$ $W^{A_i}_t$ and $W^{A_j}_t$ are independent.

We can then verify that the processes $X_t(A_i):i= 1,\dots, K$ are independent and
as $t\rightarrow\infty$,
$X_{t}(A_i)$ converges in distribution to a stationary measure $X_{\infty}(A_i)$ with
density which satisfies%
\begin{align*}
f_i(x)  &  =\frac{1}{Z}x^{\theta_i-1}e^{-\theta x}%
,\;x>0
\end{align*}
where $\theta =\frac{c}{\gamma},\; \theta_i=\theta \nu_0(A_i)$.

This can be represented by ${X}_{\infty}(A)=\theta^{-1}G(\theta\nu_{0}(A))$
where $\theta=\frac{c}{\gamma}$ and
\bean{}
&&\mathcal{L}\{(X_\infty(A_1),\dots,X_\infty(A_K)) \}=\\&&
\mathcal{L}\{\frac{1}{\theta}[G(\theta_1),G(\theta_1+\theta_2)-G(\theta_1),\dots,G(\theta)-G(\theta-\theta_K)]\}
\eean
where
 $G(s)$
is the {\em Moran subordinator} - see subsection \ref{s.moran} below.

For $u=(A^u_1,\dots, A^u_{|u|}) \in U$ (defined as in the last subsection)  let $\{P_{u}=\mathcal{L}(\{(X_{t}(A_1),\dots,X_t(A_{|u|}):t\geq0,\;A\in u\})\}$.  Then the collection $\{P_u\}_{u\in U}$
 forms a
projective system  and as in the previous section there exists
a projective limit measure \ $P_\infty$ on $(C_{[0,\infty)}([0,\infty)))^{\mathcal{B}([0,1])}$. Moreover for fixed $t\in\lbrack0,\infty),$ $X_{t}(\cdot)$ is a.s. a
finitely additive measure that is regular (on a countable generating subset of
$\mathcal{B}([0,1])$) we obtain a unique countably additive version (recall Theorem \ref{PRPOL}). Thus,
$\{X_{t}(\cdot):t\geq0\}$ is a measure-valued process and again we can obtain an a.s. continuous $M_F([0,1])$-valued version.
This $M_F([0,1])$-valued process is called the {\em Jirina process}.

\begin{corollary} The stationary measure for the Jirina process is given by the random measure
\be{} X_\infty(A)=\frac{1}{\theta} \int_0^1 1_A(x)dG(\theta s),\quad A\in\mathcal{B}([0,1])\ee
where $G(\cdot)$ is the Moran gamma subordinator.
\end{corollary}

\section{Invariant Measures of the IMA and
Jirina Processes}

\subsection{The Moran (Gamma) Subordinator }
\label{s.moran}

We begin by recalling the the {\em Gamma distribution} with parameter $\alpha>0$ given by the density function

\[
g_{\alpha}(u)=u^{\alpha-1}e^{-u}/\Gamma(\alpha)
\]

and Laplace transform of $g_{\alpha}$ is%
\[
\int_{0}^{\infty}g_{\alpha}(y)e^{-\lambda y}dy=\frac{1}{(1+\lambda)^{\alpha}},\;\lambda >-1,
\]

The Moran subordinator $\{G(\alpha): \alpha\geq 0\}$ is an increasing process with stationary independent increments
$G(\alpha_2)-G(\alpha_1),\; \alpha_1<\alpha_2$ given by $g_{\alpha_2-\alpha_1}$.

\subsubsection{L\'evy representation}

\beL{}
\be{MSub}
E\left(  e^{-\lambda G(\alpha)}\right)     =\exp\left(
-\alpha \int_{0}^{\infty }(1-e^{-u\lambda})\frac{e^{-u}}{u}du\right).
\ee
\end{lemma}
\begin{proof}
Note that%
\begin{align*}
\frac{\partial}{\partial \lambda}\int_{0}^{\infty}(1-e^{-\lambda z})z^{-1}e^{-z}dz  &
=\int_{0}^{\infty}(e^{-\lambda z})e^{-z}dz=\frac{1}{1+\lambda}\\
\int_{0}^{\infty}(1-e^{-\lambda z})z^{-1}e^{-z}dz  &  =\log(1+\lambda)
\end{align*}

Hence we have the L\'evy-Khinchin representation with L\'evy measure $\frac{e^{-z}}{z},\; z>0$%
\be{pdlm}
\frac{1}{(1+\lambda)^{\alpha}}=\exp\left\{  -\alpha\int_{0}^{\infty}(1-e^{-\lambda z}%
)z^{-1}e^{-z}dz\right\}.
\ee

\end{proof}

\subsubsection{Poisson representation}

The {\em Poisson random field} with intensity measure $\mu$ is a random counting measure $\Pi$
on a space $S$. \ $\Pi(A_{i}),\Pi(A_{j})$ are independent if $i\neq j $ and
$\Pi(A)$ is Poisson with parameter $\mu(A)$.

\beT{} (Campbell's Theorem.) Let  $\Pi$ be a Poisson random field with intensity $\mu\in M(S)$ and $f:S\rightarrow \mathbb{R},\;\Sigma
=\sum_{x\in\Pi}f(x)$ $=\int f(x)\Pi(dx)$ converges a.s.  if and only if
\[
\int_{S}\min(|f(x)|,1)\mu(dx)<\infty
\]
and then%
\[
E(e^{s\int f(x)\Pi(dx)})=\exp(\int(e^{sf(x)}-1)\mu(dx)),\;\;s\in R
\]
provided the integral on the right exists.
\end{theorem}

Now consider the Poisson random measure on $[0,1]\times
(0,\infty)$
\be{prepn1}
\Xi_\theta=\sum\delta_{\{x,u\}}%
\ee
with intensity measure
\[
\theta \nu_0(dx)\frac{e^{-u}}{u}du.
\]

Let $\wt X_\infty(A):= \int_A\int_0^\infty u \Xi_\theta(dx,du)$.
Then by Campbell's Theorem
\bea{} && E(e^{-\lambda \wt X_\infty(A)})=E(e^{-\lambda \int_A\int_0^\infty u \Xi_\theta(dx,du)})\\&&
= e^{-\theta\nu_0(A)\int (1-e^{-\lambda u})\frac{e^{-u}}{u}du}. \nonumber
\eea
Hence we can  represent equilibrium of the Jirina process by the random measure with Poisson representation $\{{X}_{\infty}(A):A\in
\mathcal{B}([0,1])\}$ by%
\be{prepn2}
{X}_{\infty}(A)=\theta^{-1}\int_{A}\int_{0}^{\infty}u\,\Xi_\theta(dx,du)
\ee
and this can be obtained as the projective limit of the finite systems.

If $\nu_0$ is Lebesgue measure on $[0,1]$ then have that the $\{X_\infty([0,s)\}_{0\leq s\leq 1}=\{G(s)\}_{0\leq s\leq 1}$ where $G(s)$ is the  Moran subordinator with with increments  $G(s_2)-G(s_1)$ having the Gamma $\theta (s_2-s_1)$ distribution  $\theta =\frac{c}{\gamma}$.

\subsection{Representation of the Infinitely Many Alleles
Equilibrium}

Recall (Theorem \ref{IMAEQ}) that the Dirichlet distribution Dirichlet$(\theta_{1},\dots,\theta_{n})$ has the joint density on relative to
$(n-1)$-dimensional Lebesgue measure on $\Delta_{n-1}$ given by%
\[
f(p_{1},\dots,p_{n-1})=\frac{\Gamma(\theta_{1}+\dots +\theta_{n})}{\Gamma
(\theta_{1})\dots\Gamma(\theta_{n})}p_{1}^{\theta_{1}-1}p_{2}^{\theta_{2}%
-1}\dots p_{n}^{\theta_{n}-1}.
\]
Recall that if the $\theta$ are large the measure concentrates away from the
boundary whereas is the $\theta$ are small things concentrate near the
boundary corresponding to highly disparate $p$ with a few large $p_{j}$ and
the others small. \ For example if the $\theta_j$ are small but equal there is a
high probability that at least one of the $p_{j}$ is much greater than
average; and which  value or values of $j$ have large $p_{j}$ is a matter of chance.

\begin{proposition}
Let $X_\infty$ denote the equilibrium random measure for the Jirina process and consider a partition \ $[0,1]=\cup_{i=1}^{n}A_{i}$ and define%
\be{}
Y(A_{i}):=\frac{X_{\infty}(A_{i})}{X_{\infty}([0,1])}=\frac{G(\theta |A_i|)}{G(\theta)}.%
\ee
Then the  family $(Y(A_{1}),\dots,Y(A_{K}))$ is \textbf{independent }of $X_{\infty
}([0,1])$ \ and has as distribution the Dirichlet$(\theta_{1},\dots\theta
_{K})$ where $\theta_{j}=\theta\nu_{0}(A_{j}).$
\end{proposition}

\begin{proof}
Let ${Y}$ be Gamma$(\theta)$ and\\ $(P_{1}%
,\dots,P_{K})$ Dirichlet $(\theta_{1},\dots,\theta_{K})$ with ${Y}$ and
$(P_{1},\dots,P_{K})$ independent, and define $(Y_{1},\dots,Y_{K})$ by

\be{}Y_{i}:={Y}P_{i}.\ee

We will verify that  $(Y_1,\dots,Y_K)$  has the  joint probability density  function
\be{}
g(y_1,\dots,y_K)=\prod_{i=1}^{K}u_{i}^{\theta_{i}-1}e^{-u_{i}}/\Gamma(\theta_i).
\ee

Consider the 1-1 transformation $(Y_{1},Y_{2},\dots,Y_{K})\leftrightarrow
(Y,P_{2},\dots,P_{K})$ with Jacobian \be{}|J|=\left\{|\frac{\partial x_{1},\dots, x_K}{\partial y_{1},\dots, y_K}|,\;x_{1}=y,x_{2}=p_{2}%
,\dots,x_{K}=p_{K}\right\} =\frac{1}{y^{K-1}}.\ee

By independence of $Y$ and $(P_1,\dots,P_K)$, we obtain the  joint density of $(Y_{1},\dots,Y_{K})$ as\\

\bean{}
&&g(y_{1},\dots,y_{K})    =f(p_{1},\dots,p_{K}|{Y})f_{{Y}}(y)|J|\\
&&  =f(p_{1},\dots,p_{K})\frac{1}{\Gamma(\theta)}y^{\theta-1}e^{-y}%
\frac{1}{y^{K-1}}\\
&&  =\frac{\Gamma(\theta)}{\Gamma(\theta_{1})\dots\Gamma(\theta_{K})}\\
&&  .\left(  \frac{y_{1}}{\sum y_{i}}\right)  ^{\theta_{1}-1}\dots\left(
\frac{y_{K}}{\sum y_{i}}\right)  ^{\theta_{K}-1}\frac{1}{\Gamma(\theta)}(\sum
y_{i})^{(\theta-1)}e^{-\sum y_{i}}.(\sum y_{i})^{-(K-1)}\\
&&  =\prod_{i=1}^K\frac{1}{\Gamma(\theta_{i})}y_{i}^{\theta_{i}-1}e^{-y_{i}}%
\eean
Note that this coincides with the Dirichlet$(\theta_1,\dots,\theta_K)$ distribution.
\end{proof}

\beC{} The invariant measure of the infinitely many alleles model can be represented by the
random probability measure
\be{}  Y(A) =\frac {X_\infty(A)}{X_\infty([0,1])},\quad,  A\in \mathcal{B}([0,1]).
\ee
where $X_\infty(\cdot)$ is the equilibrium of the above Jirina process and $Y(\cdot)$ and $X_\infty([0,1])$ are independent.
\end{corollary}

\subsubsection{Reversibility}

Recall that the Dirichlet distribution is a {\em reversible} stationary measure for the $K-type$ Wright-Fisher model with house of cards mutation (Theorem \ref{IMAEQ}).  From this and the projective limit construction it can be verified that $\mathcal{L}(Y(\cdot))$ is a reversible stationary measure for the infinitely many alleles process.  Note that reversibility actually characterizes the IMA model among neutral Fleming-Viot processes with mutation, that is, any mutation mechanism other than the ``type-independent'' or ``house of cards'' mutation leads to a stationary measure that is \underline{not}
reversible (see Li-Shiga-Yau (1999) \cite{LSY-99}).

\subsection{The Poisson-Dirichlet Distribution}

Without loss of generality we can assume that $\nu_0$ is Lebesgue measure on $[0,1]$. This implies that the IMA equilibrium is given by a random probability measure which is pure atomic
 \be{} p_\infty =\sum_{i=1}^\infty  a_i \delta_{x_i},\quad \sum_{i=1}^\infty a_i =1,\; x_i\in [0,1]\ee
  in which the $\{x_i\}$ are i.i.d. $U([0,1])$ and the atom sizes $\{a_i\}$ correspond to the normalized jumps of the Moran subordinator.  Let $(\xi_1,\xi_2,\dots)$ denote the reordering of the atom sizes $\{a_i\}$ in decreasing order.

The Poisson-Dirichlet $PD(\theta)$ distribution is defined to be the  distribution of the
infinite sequence $\xi=(\xi_{1},\xi_{2},\dots)$ which satisfies
\[
\xi_{1}\geq\xi_{2}\geq...,\;\sum_{k}\xi_{k}=1.
\]

This sequence is given by
\be{pdg}
\xi_{k}=\frac{\eta_{k}(\theta)}{G(\theta)}=\frac{\eta_{k}}{\sum \eta_\ell},\quad k=1,2,\dots,
\ee
where $\eta_k=\eta_{k}(\theta)$ is the height of the $kth$ largest jump in $[0,\theta]$
of the Moran process (subordinator), $G$ and $G(\theta) =\sum_{\ell =1}^\infty \eta_\ell$.

\subsubsection{Properties of the Poisson-Dirichlet Distribution}

Recalling (\ref{prepn1}), (\ref{prepn2}) we note that the set of heights of the jumps of $G(\cdot)$ in $[0,\theta]$ form a Poisson random field
$\Pi_{\theta}$ on $(0,\infty)$ with intensity measure
\[
\theta\frac{e^{-u}}{u} du.
\]

We can then  give a direct description of $PD(\theta)$ in terms of
such a Poisson random field. \ If $\eta_{1}\geq\eta_{2}\geq\eta_{3}\geq\dots$ are
the points of such a random field ordered by size then%
\[
\xi_{k}=\frac{\eta_{k}}{\sum_{\ell=1}^\infty \eta_\ell}%
\]
defines a sequence $\xi$ having the distribution $PD(\theta)$.

By the law of large numbers (for the Poisson) we get%
\[
\lim_{t\rightarrow0}\frac{\#\{k:\eta_{k}>t\}}{L(t)}=1
\]
with probability one where%
\[
L(t)=\int_{t}^{\infty}\theta\frac{e^{-u}}{u}du\sim-\theta\log t.
\]
Thus%
\begin{align*}
\#\{k  &  :\eta_{k}>t\}\sim-\theta\log t\\
\eta_{\theta\log\frac{1}{t}}  &  \approx t\quad  \text{as } t\to 0.\\
\eta_{k}  &  \approx e^{-k/\theta}%
\end{align*}

Thus $\xi_{k}$ decays exponentially fast
\[
-\log\xi_{k}\sim\frac{k}{\theta}\quad \text{as  }k\to\infty.%
\]

\subsubsection{The Distribution of Atom Sizes}

We now introduce the random measure on $(0,\infty),$%

\begin{align*}
Z_{\theta}((a,b))  &  =\frac{\Xi_\theta([0,1]\times(a,b))}{G(\theta)}\\
\int_0^\infty uZ_{\theta}(du)  &  =1.
\end{align*}
This is the \textit{distribution of normalized atom sizes} and this just
depends on the normalized ordered atoms and hence is independent of
$X_{\infty}([0,1])$. Intuitively, as $\theta\rightarrow\infty$, $Z_{\theta}(\theta du) $
converges in some sense to
\[
\frac{e^{-u}}{u}du.
\]
To give a precise formulation of this we first note that
\[
\int_{0}^{\infty}u^{k}\left(  \frac{e^{-u}}{u}\right)  du=\Gamma
(k)=(k-1)!
\]
Then one can show (see Griffiths (1979), \cite{G-79}) that
\be{} \lim_{\theta\to\infty} \theta^{k-1} \int_0^\infty u^k Z_\theta( du) =(k-1)!\ee
and there is an associated CLT
\be{} \sqrt{\theta}\; \frac{\theta^{k-1}\int u^kZ_\theta( du)-(k-1)!}{(k-1)!}\Rightarrow N(0,\sigma_k^2),\ee
with $\sigma_k^2 =\frac{(2k-1)!-(k!)^2}{((k-1)!)^2}$
{ Joyce, Krone and Kurtz (2002) \cite{JKK-02}}.  Also see Dawson and Feng (2006) \cite{DF-06} for the related large deviation behaviour.

\subsection{The GEM Representation}

Without loss of generality we can assume $\nu_{0}=U[0,1]$. Consider a
partition of $[0,1]$ into $K$ intervals of equal length. Then the random
probability
\[
\vec{p}_{K}=(p_{1},\dots,p_{K})
\]
has the symmetric Dirichlet distribution $D(\alpha,\dots,\alpha)$ with
$\alpha=\frac{\theta}{K}$.

\subsubsection{Randomized Ordering via Size-biased sampling}

Let $\mathcal{N}$ be a random variable having values in $\{1,2,\dots,K\}$ in
such a way that%
\[
P(\mathcal{N}=k|\vec{p}_{K})=p_{k},\;\;(1\leq k\leq K)
\]
Then a standard calculation shows that the vector
\[
\vec {p^{\prime}}=(p_{\mathcal{N}},p_{1},\dots,p_{\mathcal{N}-1}%
,p_{\mathcal{N}+1},\dots,p_{K})
\]
has distribution (cf. (\ref{bayes}))
\[
D(\alpha+1,\alpha,\dots,\alpha)
\]

It follows that $(p_{\mathcal{N}},1-p_{\mathcal{N}})$ has the Dirichlet
distribution (Beta distribution)%
\[
D(\alpha+1,(K-1)\alpha)
\]
so that it has probability density function%
\[
\frac{\Gamma(K\alpha+1)}{\Gamma(\alpha+1)\Gamma(K\alpha-\alpha)}p^{\alpha
}(1-p)^{(K-1)\alpha-1}.%
\]
Given $v_{1}=p_{\mathcal{N}}$, the conditional distribution of the remaining
components of $\vec p$ is the same as that of $(1-p_{\mathcal{N}})p^{(1)},$ where
the $(K-1)$-vector $\vec p^{(1)}$ has the symmetric distribution $D(\alpha,\dots,\alpha)$.

We say that $p_{\mathcal{N}}$ is obtained from $\vec{p}$ by size-biased
sampling. This process may now be applied to $\vec p^{(1)}$ to produce a component,
$v_{2}$ with distribution
\[
\frac{\Gamma((K-1)\alpha+1)}{\Gamma(\alpha+1)\Gamma((K-1)\alpha-\alpha
)}p^{\alpha}(1-p)^{(K-2)\alpha-1}%
\]
and a $(K-2)$ vector $\vec p^{(2)}$ with distribution $D(\alpha,\dots,\alpha)$. This is an example of Kolmogorov's {\em stick breaking} process.
\


\begin{theorem}
(a) As $K\rightarrow\infty$, with $K\alpha = \theta$ constant,  the distribution of the vector $\vec{q}%
_{K}=(q_{1},q_{2},\dots,q_{K})$ converges weakly to the GEM distribution with
parameter $\theta$, that is the distribution of the random
probability vector $\vec{q}=(q_{1},q_{2},\dots)$ where
\[
q_{1}=v_{1},\;q_{2}=(1-v_{1})v_{2},\;q_{3}=(1-v_{1})(1-v_{2})\nu_{3},\;\dots
\]
with $\{v_{k}\}$ are i.i.d. with Beta density (Beta$(1,\theta)$)%
\[
\theta(1-p)^{\theta-1},\;0\leq p\leq1
\]
\newline (b) If $\vec{q}=(q_{1},q_{2},\dots)$ is reordered (by size) as
$\vec{p}=(p_{1},p_{2},\dots)$, that is i.e. $p_{k}$ is the kth largest of
the $\{q_{j}\},$then\textbf{ }$\vec{p}$ has the Poisson-Dirichlet distribution,
$PD(\theta).$
\end{theorem}

\begin{proof}
(a)  Let $(p^K_1,\dots,p^K_k)\in\Delta_{K-1}$ be a random probability vector obtained by decreasing size reordering of a probability vector sampled  from the  distribution $D(\alpha,\dots,\alpha)$ with $\alpha=\frac{\theta}{K}$. Then let $(q^K_1,\dots,q^K_k)$ be the size-biased reordering of $(p^K_1,\dots,p^K_k)$. Then as shown above we can rewrite this as
\be{} q^K_1=v^K_1,\;q^K_2=(1-v^K_1)v^K_2),\; q^K_3=(1-v^K_1)(1-v^K_2)v^K_3,\dots\ee
where $v^K_1,\dots,v^K_{K-1}$ are independent and  $v^K_r$ has pdf

\be{}
\frac{\Gamma((K-r)\alpha+1)}{\Gamma(\alpha+1)\Gamma((K-r)\alpha-\alpha
)}u^{\alpha}(1-u)^{(K-r-1)\alpha-1},\quad 0\leq u\leq 1.
\ee

Now let $K\rightarrow\infty$ with $K\alpha=\theta$. Then
\[
\frac{\Gamma(K\alpha+1)}{\Gamma(\alpha+1)\Gamma(K\alpha-\alpha)}p^{\alpha
}(1-p)^{(K-1)\alpha-1}\rightarrow\theta(1-p)^{\theta-1}.
\]%
and
\begin{align*}
&  \frac{\Gamma((K-r)\alpha+1)}{\Gamma(\alpha+1)\Gamma((K-r)\alpha-\alpha
)}p^{\alpha}(1-p)^{(K-r-1)\alpha-1}\\
&  =\frac{\Gamma((K-r)\theta/K+1)}{\Gamma(\theta/K+1)\Gamma((K-r)\theta
/K-\theta/K)}p^{\theta/K}(1-p)^{(K-r-1)\theta/K-1}\\
&  \rightarrow\theta(1-p)^{\theta}%
\end{align*}
Thus the distributions of the first $m$ components of the vector
$\vec{q}_{K}$ converge weakly to the distribution of the first $m$
components of the random (infinite)  probability vector $\vec{q}$ defined by
\[
q_{1}=v_{1},\;q_{2}=(1-v_{1})v_{2},\;q_{3}=(1-v_{1})(1-v_{2})\nu_{3},\;\dots
\]
where $\{v_{k}\}$ are i.i.d. with Beta density (Beta$(1,\theta)$)%
\[
\theta(1-p)^{\theta-1},\;0\leq p\leq1.
\]

(b) By the projective limit construction, the PD($\theta$) distribution arises as the limit in distribution of the ordered probability vectors $(p^K_1,\dots,p^K_k)$.  Then the size-biased reorderings converge in distribution to the size-biased reordering $q_1,q_2,\dots$ of the probability vector $p_1,p_2,p_3,\dots$. Clearly the decreasing-size reordering of $q_1,q_2,\dots$ reproduces $p_1,p_2,\dots$.
\end{proof}

\begin{remark} The distribution of sizes of the {\em age ordered alleles} in the infinitely many alleles model is given by the GEM distribution. The intuitive idea is as follows. By exchangeability at the individual level the probability that the $k$th allele at a given time survives the longest (time to extinction) among those present at that time is proportional to $p_k$ the frequency of that allele in the population at that time.  Observing that the ordered survival times correspond to ages under time reversal, the result follows from reversibility. See Ethier \cite{E-90} for justification of this argument and a second proof of the result.

\end{remark}

\begin{remark} There is a two parameter analogue of the Poisson-Dirichlet
introduced by Perman, Pitman and Yor (1992) \cite{PPY-92} that shares some features with the PD distribution.
See Feng (2009) \cite{F-09} for a recent detailed exposition.

\end{remark}

\subsection{Application of the Poisson-Dirichlet distrbution:\\ The Ewens Sampling Formula}

In analyzing population genetics data under the neutral hypothesis it is
important to know the probabilities of the distribution of types obtained in
taking a random sample of size $n$. \ For example, this is used to test for
neutral mutation in a population.

Consider a random sample \ of size $n$ chosen from the random vector $\xi
=(\xi_{1},\xi_{2},\dots)$ chosen from the distribution $PD(\theta).$

We first compute the probability that they are all of the same type.
\ Conditioned on $\xi$ this is $\sum_{k=1}^{\infty}\xi_{k}^{n}$ and hence the
unconditional probability is%
\[
h_{n}=E\left\{  \sum_{k=1}^{\infty}\xi_{k}^{n}\right\}
\]
Using Campbell's formula we get%
\begin{align*}
&  E\left\{  \sum_{k=1}^{\infty}\eta_{k}^{n}\right\} \\
&  =\frac{d}{ds}\left(  E(e^{\int_{0}^{1}\int_{0}^{\infty}(e^{sz^{n}%
}-1)\frac{\theta e^{-z}}{z}dz}\right)  |_{s=0}\\
&  =\int z^{n}\frac{\theta e^{-z}}{z}dz=\theta(n-1)!
\end{align*}
Also%
\[
E((\sum_{k=1}^\infty \eta_k)^{n})=\frac{\Gamma(n+\theta)}{\Gamma(\theta)}.
\]
By the Gamma representation (\ref{pdg}) for the ordered jumps of the Gamma subordinator we
get%
\begin{align*}
E\left\{  \sum_{k=1}^{\infty}\eta_{k}^{n}\right\}   &  =E\left\{  (\sum_{k=1}^\infty \eta_k)
^{n}\sum_{k=1}^{\infty}\xi_{k}^{n}\right\} \\
&  =E((\sum_{k=1}^\infty \eta_k)^{n})E\left\{  \sum_{k=1}^{\infty}\xi_{k}^{n}\right\}
\;\;\text{by independence}%
\end
{align*}

Therefore%
\[
h_{n}=\frac{\theta\Gamma(\theta)(n-1)!}{\Gamma(n+\theta)}=\frac{(n-1)!}%
{(1+\theta)(2+\theta)\dots(n+\theta-1)}%
\]

In general in a sample of size $n$ let%
\begin{align*}
a_{1}  &  =\text{ number of types with 1representative}\\
a_{2}  &  =\text{ number of types with 2 representatives}\\
&  \dots\\
a_{n}  &  =\text{number of types with n representatives}%
\end{align*}
Of course,
\[
a_{i}\geq0,\;\;a_{1}+2a_{2}+\dots+na_{n}=n.
\]
We can also think of this as a partition%
\[
\mathfrak{a}=1^{a_{1}}2^{a_{2}}\dots n^{a_{n}}%
\]
Let \ $P_{n}(\mathfrak{a})$ denote the probability that the sample exhibits
the partition $\mathfrak{a}$. Note that \ $P_{n}(n^{1})=h_{n}.$

\begin{proposition}
(Ewens sampling formula)%
\be{}
P_{n}(\mathfrak{a})=P_{n}(a_{1},\dots,a_{n})=\frac{n!\Gamma(\theta)}%
{\Gamma(n+\theta)}\prod_{j=1}^{n}\left(  \frac{\theta^{a_{j}}}{j^{a_{j}}%
a_{j}!}\right).
\ee
\end{proposition}

\begin{proof}
We can select the partition of $\{1,\dots,n\}$ into subsets of sizes
$(a_{1},\dots,a_{n})$ as follows. Consider $a_{1\text{ }}$boxes of size $1$,
$\dots$ $a_{n}$ boxes of size $n.$ Then the number of ways we can distribute
$\{1,\dots,n\}$ is $n!$ but we can reorder the $a_{i}\;$boxes \ in $a_{i}!$
ways and and there are $(j!)$ permutations of the indices in each of the
$a_{j}$ partition elements with $j$ elements. \ Hence the total number of ways
we can do the partitioning to $\{1,\dots,n\}$ is $\frac{n!}{\prod(j!)^{a_{j}%
}a_{j}!}.$

Now condition on the vector $\xi=(\xi(1),\xi(2),\dots).$ The probability that
we select the types (ordered by their frequencies is then given by)%
\[
P_{n}(\mathfrak{a}|\xi)=\frac{n!}{\prod(j!)^{a_{j}}a_{j}!}\sum_{I_{a_{1},\dots
an}}\xi\left(  k_{11}\right)  \xi(k_{12})\dots\xi(k_{1a_{1}})\xi(k_{21}%
)^{2}\dots\xi(k_{2a_{2}})^{2}\xi(k_{31})^{3}\dots
\]
where the summation is over
\[
I_{a_{1},\dots an}:=\{k_{ij}:i=1,2,\dots;j=1,2,\dots,a_{i}\}
\]
Hence using the Gamma representation we get%
\begin{align*}
&  \frac{\Gamma(n+\theta)}{\Gamma(\theta)}P_{n}(\mathfrak{a})\\
&  =\frac{n!}{\prod(j!)^{a_{j}}a_{j}!}E\left\{  \sum\eta\left(  k_{11}\right)
\eta(k_{12})\dots\eta(k_{1a_{1}})\eta(k_{21})^{2}\dots\eta(k_{2a_{2}})^{2}%
\eta(k_{31})^{3}\dots\right\}
\end{align*}
But%
\begin{align*}
&  E\left\{  \sum_{i,j}\eta\left(  k_{11}\right)  \eta(k_{12})\dots
\eta(k_{1a_{1}})\eta(k_{21})^{2}\dots\eta(k_{2a_{2}})^{2}\eta(k_{31})^{3}%
\dots\right\} \\
&  =\prod_{j=1}^{n}E\left\{  \sum_{k=1}^{\infty}\eta(k)^{j}\right\}  ^{a_{j}%
}=\prod_{j=1}^{n}\left\{  \int_{0}^{\infty}z^{j}\theta\frac{e^{-z}}%
{z}dz\right\}  ^{a_{j}}\text{ \ by Campbell's thm}\\
&  =\prod_{j=1}^{n}\{\theta(j-1)!\}^{a_{j}}%
\end{align*}
Therefore substituting we get%
\[
P_{n}(\mathfrak{a})=\frac{n!\Gamma(\theta)}{\Gamma(n+\theta)}\prod_{j=1}%
^{n}\left(  \frac{\theta^{a_{j}}}{j^{a_{j}}a_{j}!}\right).
\]
\end{proof}

\chapter{Martingale Problems and Dual Representations}

\section{Introduction}
We have introduced above the basic mechanisms of branching, resampling and mutation mainly using generating functions and semigroup methods.  However these methods have limitations and in order to work with a wider class of mechanisms we will introduce some this additional  tools of stochastic analysis in this chapter.  The martingale method which we  use has proved to be a natural framework for studying a wide range of problems including those of population systems. The general framework is as follows:

\begin{itemize}
\item the object is to specify a Markov process on a Polish space $E$ in terms of its probability laws $\{P_x\}_{x\in E}$ where $P_x$ is a probability measure on $C_E([0,\infty))$ or $D_E([0,\infty))$ satisfying $P_x(X(0)=x)=1$.
\item the probabilities $\{P_x\}\in \mathcal{P}(C_E([0,\infty)))$ satisfy a {\em martingale problem} (MP).  One class of martingale problems is defined by the set of  conditions  of the form \be{}F(X(t))-\int_0^t GF(X(s))ds,\quad F\in\mathcal{D}\qquad (\mathcal{D},G)- \text{martingale problem}\ee is a $P_x$ martingale where $G$ is a linear map from $\mathcal{D}$ to $C(E)$,  and  $\mathcal{D}\subset C(E)$ is measure-determining.
\item the martingale problem MP has one and only one solution.

\end{itemize}

Two martingale problems $MP_1$, $MP_2$ are said to be {\em equivalent} if a solution to $MP_1$ problem is a solution to $MP_2$ and vice versa.

In our setting the existence of a solution is often obtained as the limit of a sequence of probability laws of approximating processes.  The question of uniqueness is often the more challenging part.  We introduce the method of dual representation which can be used to establish uniqueness for a number of basic population processes. However the method of duality is applicable only for special classes of models.  We introduce a second method, the Cameron-Martin-Girsanov type change of measure which is applicable to some basic problems of stochastic population systems.  Beyond the domain of applicability of these methods, things are much more challenging.  Some recent progress has been made in a series of papers of Athreya, Barlow, Bass, Perkins \cite{ABBP-02}, Bass-Perkins \cite{BP-03}, \cite{BP-08} but open problems remain.

We begin by reformulating the Jirina and neutral IMA Fleming-Viot  in the martingale problem setting.  We then develop the Girsanov and duality methods in the framework of measure-valued processes and apply them to the Fleming-Viot process with selection.

\section{The Jirina martingale problem}

By our projective limit construction of the Jirina process (with $\nu
_{0}=$ Lebesgue), we have a probability space $(\Omega,\mathcal{F},\{X^{\infty
}:[0,\infty)\times\mathcal{B}([0,1])\rightarrow\lbrack0,\infty)\},P)$ such
that a.s. $t\rightarrow X_{t}^{\infty}(A)$ is continuous and $A\rightarrow
X_{\cdot}^{\infty}(A)$ is finitely additive. We can take a
modification, $X$, of $X^{\infty}$ such that a.s. $X:[0,\infty)\rightarrow
M_{F}([0,1])$ is continuous where $M_{F}([0,1])$ is the space of (countably
additive) finite measures on $[0,1]$ with the weak topology. \ We then define
the filtration
\[
\mathcal{F}_{t}:\sigma\{X_{s}(A):0\leq s\leq t,\;A\in\mathcal{B}([0,1])\}
\]
and $\mathcal{P}$, the $\sigma$-field of predictable sets in $\mathbb{R}%
_{+}\times\Omega$ (ie the $\sigma$-algebra generated by the class of
$\mathcal{F}_{t}$-adapted, left continuous processes).

Recall that for a fixed set $A$ the Feller CSB with immigration satisfies
\be{}X_{t}(A)-X_{0}(A)-\int_{0}^{t}c(\nu_{0}(A)-X_{s}(A))ds
=\int_{0}%
^{t}\sqrt{2\gamma X_{t}(A)}dw_{t}^{A}\ee
which is an $L^2$-martingale.

Moreover, by polarization
\be{jmp3}
\left\langle M(A_{1}),M(A_{2})\right\rangle _{t}=\gamma\int_{0}^{t}X_{s}%
(A_{1}\cap A_{2})ds
\ee
and if $A_{1}\cap A_{2}=\emptyset$, then the martingales $M(A_{1})_{t}$ and
$M(A_{2})_{t}$ are orthogonal. This is an example of an {\em orthogonal martingale measure}.

Therefore for any Borel set $A$
\be{jmp1}
M_{t}(A):=X_{t}(A)-X_{0}(A)-\int_{0}^{t}c(\nu_{0}(A)-X_{s}(A))ds
\ee
is a martingale with increasing process
\be{jmp2}
\left\langle M(A)\right\rangle _{t}=\gamma\int_{0}^{t}X_{s}(A)ds.
\ee

We note that we can define integrals with respect to an orthogonal martingale measure (see next subsection)
and show that (letting $X_t(f)=\int f(x)X_t(dx)$ for $f\in \mathcal{B}([0,1]))$

\be{jmp1}
M_{t}(f):=X_{t}(f)-X_{0}(f)-\int_{0}^{t}c(\nu_{0}(f)-X_{s}(f))ds=\int f(x)M_t(dx)
\ee
which is a martingale with increasing process
\be{} \left\langle M(f)\right\rangle _{t}=\gamma\int_{0}^{t}f^2(x)X_{s}(dx)ds.
\ee

This suggests the martingale problem for the Jirina process which we state in subsection \ref{ss.JMP}.


\subsection{Stochastic Integrals wrt Martingale Measures}

A general approach to martingale measures and stochastic integrals with respect to martingale measures was developed by Walsh \cite{W-86}. We briefly review some basic results.

Consider the collection of $\mathcal{S}$ {\em simple} functions $\psi$ of the form
\be{simf}
\psi(t,\omega,x)=\sum_{i=1}^{K}\psi_{i-1}(\omega)\phi_{i}(x)1_{(t_{i-1}%
,t_{i}]}(t)
\ee
for some $\phi_{i}\in b\mathcal{B}(E)$, $\psi\in b\mathcal{F}_{t_{i-1}}$,
$0=t_{0}<t_{1}\dots<t_{K}\leq\infty$.
The predictable $\sigma$-field $\mathcal{P}r$ on $\mathbb{R}_{+}\times\Omega\times
E$ is the $\sigma$-field generated by the class of simple functions of the form (\ref{simf}).

For  $\psi\in \mathcal{S}$, define%
\[
M_{t}(\psi):=\int_{0}^{t}\int\psi(s,x)dM(s,x)=\sum_{i=1}^{K}\psi
_{i-1}(M_{t\wedge t_{i}}(\phi_{i})-M_{t\wedge t_{i-1}}(\phi_{i}))
\]
Then $M_{t}(\psi)\in\mathcal{M}_{\text{loc}}$ (the space of $\mathcal{F}_{t}$
local martingales) and%
\[
\left\langle M(\psi)_{t}\right\rangle =\int_{0}^{t}X_{s}(\gamma\psi_{s}%
^{2})ds.
\]

Let%
\[
\mathcal{L}_{loc}^{2}=\left\{  \psi:\mathbb{R}_{+}\times\Omega\times
E\rightarrow\mathbb{R}:\psi\text{ is }\mathcal{P}r%
\text{-measurable,\ }\int_{0}^{t}X_{s}(\psi_{s}^{2})ds<\infty,\;\forall
t>0\right\}
\]
\bigskip

\begin{lemma}
For any $\psi\in\mathcal{L}_{\text{loc}}^{2}$ there is a sequence $\{\psi
_{n}\}$ in $\mathcal{S}$ such that%
\[
P\left(  \int_{0}^{n}\int(\psi_{n}-\psi)^{2}(s,\omega,x)\gamma(x)X_{s}%
(dx)ds>2^{-n}\right)  <2^{-n}\text{.}%
\]
\end{lemma}

\begin{proof}
Let $\mathcal{\bar{S}}$ denote the set of bounded $\mathcal{P}r$-measurable functions which can be approximated as above.
\ $\mathcal{\bar{S}}$ is closed under $\rightarrow^{bp}$. Using \newline
$\mathcal{H}_{0}=\{f_{i-1}(\omega)\phi_{i}(x)$,$\phi\in b\mathcal{E}%
,\;f_{i-1}\in b\mathcal{F}_{t_{i-1}},\phi_{i}\in b\mathcal{E}\}$, we see that
$\psi(t,\omega,x)=\sum_{i=1}^{K}\psi_{i-1}(\omega,x)1_{(t_{i-1},t_{i}]}(t)$ is
in $\mathcal{\bar{S}}$ for any $\psi_{i-1}\in b(\mathcal{F}_{t_{i-1}}%
\times\mathcal{E})$. If $\psi\in b(\mathcal{P}r)$, then%
\[
\psi_{n}(s,\omega,x)=2^{n}\int_{(i-1)2^{-n}}^{i2^{-n}}\psi(r,\omega,x)dr\text{
is }s\in(i2^{-n},(i+1)2^{-n}],\;i=1,2,\dots
\]
satisfies $\psi_{n}\in\mathcal{\bar{S}}$ by the above. \ For each
$(\omega,x),\psi_{n}(s,\omega,x)\rightarrow\psi(s,\omega,x)$ for Lebesgue a.a.
$s$ by Lebesgue's differentiation theorem and it follows easily that $\psi
\in\mathcal{\bar{S}}$. Finally if $\psi\in\mathcal{L}_{\text{loc}}^{2}$, the
obvious truncation argument and dominated convergence (set $\psi_{n}%
=(\psi\wedge n)\vee(-n)$ completes the proof.
\end{proof}

\begin{proposition}
There is a unique linear extension of $M:\mathcal{S\rightarrow M}_{\text{loc}%
}$ (the space of local martingales) to a map $M:\mathcal{L}_{\text{loc}}^{2}\rightarrow\mathcal{M}_{\text{loc}%
}$ such that $M_t(\psi)$ is a local martingale with increasing process $\langle M(\psi)\rangle_t$ given by
\[
\left\langle M(\psi)\right\rangle _{t}:=\int_{0}^{t}\gamma X_{s}(\psi_{s}%
^{2})ds\;\forall\;t\geq0\;\; a.s.\forall\;\psi\in\mathcal{L}_{\text{loc}}^{2}.
\]
\end{proposition}

\begin{proof}
We can choose $\psi_{n}\in\mathcal{S}$ \ as in the Lemma. Then%
\begin{align*}
\left\langle M(\psi)-M(\psi_{n})\right\rangle _{n}  &  =\left\langle
M(\psi-\psi_{n})\right\rangle _{n}=\gamma\int_{0}^{n}X_{s}(\gamma(\psi(s)-\psi
_{n}(s))^{2}ds\\
& P\left(\left\langle M(\psi)-M(\psi_{n})\right\rangle _{n} >2^{-n}\right)< 2^{-n}%
\end{align*}
The $\{M_t(\psi_n)\}_{t\geq 0}$ is Cauchy and using Doob's inequality and the Borel-Cantelli Lemma we can define $\{M_t(\psi)\}_{t\geq 0}$
such that
\[
\sup_{t\leq n}|M_{t}(\psi)-M_{t}(\psi_{n})|\rightarrow0\text{ a.s. as
}n\rightarrow\dot{\infty}.%
\]
This yields the required extension and its uniqueness.
\end{proof}

Note that it immediately follows by polarization that if $\psi,\phi
\in\mathcal{L}_{\text{loc}}^{2},$%
\[
\left\langle M(\phi),M(\psi)\right\rangle _{t}=\gamma\int_{0}^{t}X_{s}%
(\phi_{s}\psi_{s})ds
\]
Moreover, in this case $M_{t}(\psi)$ is a $L^{2}$-martingale, that
is,
\[
E(\left\langle M(\psi)\right\rangle _{t})=\gamma\int_{0}^{t}E(X_{s}(\psi
_{s}^{2}))ds<\infty
\]
provided that%
\[
\psi\in\mathcal{L}^{2}=\{\psi\in\mathcal{L}_{\text{loc}}^{2}:E(\int_{0}%
^{t}X_{s}(\psi_{s}^{2})ds)<\infty,\;\forall t>0\}.
\]

\begin{remark}
Walsh (1986) \cite{W-86} defined a more general class of martingale measures on a
measurable space $(E,\mathcal{E)}$ for which the above construction of
stochastic integrals can be extended. $\{M_{t}(A):t\geq0,\;A\in\mathcal{E}\}$
is an $L^{2}$-martingale measure wrt $\mathcal{F}_{t}$ iff

(a) $M_{0}(A)=0\;\;\forall A\in\mathcal{E}$,

(b) $\{M_{t}(A),t\geq0\}$ is an $\mathcal{F}_{t}$-martingale for every
$A\in\mathcal{E}$,

(c) for all $t>0,\;M_{t}$ is an $L^{2}$-valued $\sigma$-finite measure$.$

The martingale measure is \textit{worthy }if there exists a $\sigma$-finite
``dominating measure'' $K(\cdot,\cdot,\cdot,\omega),\,$on $\mathcal{E}%
\times\mathcal{E}\times\mathcal{B}(\mathbb{R}_{+}),\;\omega\in\Omega$ such that

(a) $K$ is symmetric and positive definite, i.e. for any $f\in b\mathcal{E}%
\times\mathcal{B}\mathbb{(R}_{+})$,%
\[
\int\int\int f(x,s)f(y,s)K(dx,dy,ds)\geq0
\]

(b) for fixed $A,B,\;\{K(A\times B\times(0,t]),t\geq0\}$ is $\mathcal{F}_{t}$-predictable

(c) $\exists\;E_{n}\uparrow E$ such that $E\{K(E_{n}\times E_{n}\times
\lbrack0,T]\}<\infty\;\forall n,$

(d) $|\left\langle M(A),M(A)\right\rangle _{t}|\leq K(A\times A\times\lbrack0,t]).$

\end{remark}
\bigskip
\subsection{Uniqueness and stationary measures for the Jirina Martingale Problem}
\label{ss.JMP}

A probability  law, $\mathbb{P_\mu}\in \mathcal{P}(C_{M_{F}([0,1])}([0,\infty)))$, is a solution of the {\em Jirina martingale problem},
if under $\mathbb{P}_\mu$, $X_{0}=\mu$  and
\be{JMPS}\begin{split}
M_{t}(\phi)  &  :=X_{t}(\phi)-X_{0}(\phi)-\int_{0}^{t}c(\nu_{0}(\phi
)-X_{s}(\phi))ds,\;\\
&\text{ is a }L^{2}\text{, }\mathcal{F}_{t}\text{-martingale }\forall\;\phi
\in b\mathcal{B}([0,1])\text{ with increasing process}\\
\left\langle M(\phi)\right\rangle_t  &  =\gamma\int_{0}^{t}X_{s}(\phi^{2})ds, \text{ that is},\\
& M^2_t(\phi)-\langle M(\phi)\rangle _t \text{   is a martingale}.\end{split}
\ee

\begin{remark} This is {\em equivalent} to the martingale problem
\be{}
M_F(t)= F(X_t)-\int_0^t GF(X(s))ds\quad\text{is a martingale}
\ee
for all $F\in\mathcal{D}\subset C(M_F([0,1])$ where
\[ \mathcal{D}=\{F:F(\mu)=\prod_{i=1}^n \mu(f_i),\;\; f_i\in C([0,1]),\;i=1,\dots,,n,\;n\in\N\}\]
and
\bean{} GF(\mu)=&& c\int\left[ \int\frac{\partial F(\mu)}{\partial \mu(x)}\nu_0(dx)-\frac{\partial F(\mu)}{\partial \mu(x)}\right]\mu(dx)\\&& +\frac{\gamma}{2}\int\int \frac{\partial^2 F(\mu)}{\partial\mu(x)\partial\mu(y)}(\delta_x(dy)\mu(dx)-\mu(dx)\mu(dy))
\eean
\end{remark}

\begin{theorem}\label{JMP}
 There exists one and only one solution $\mathbb{P}_\mu \in\mathcal{P}(C_{M_{F}([0,1])}([0,\infty))$ to the martingale problem (\ref{JMPS}).  This defines a continuous $M_F([0,1])$-valued continuous strong Markov process.\newline (b) (Ergodic Theorem) Given any
initial condition, $X_{0},$ the law of $X_{t}$ converges weakly to a limiting
distribution as $t\rightarrow\infty$ with Laplace functional
\be{}
E(e^{-\int_0^1 f(x)X_\infty(dx)})=\exp\left(-\frac{2c}{\gamma}\int_0^1\log(1+\frac{f(x)}{\theta})\nu_0(dx)\right).\ee
This can be represented by
 \be{} X_{\infty}(A)=
\frac{1}{\theta}\int_0^1 f(s)G(\theta ds)\ee where $G$ is the Gamma (Moran)
subordinator (recall (\ref{MSub})).
\end{theorem}

\begin{proof} Outline of method.
As discussed above the projective limit construction produced a solution to this martingale problem.

A fundamental result of Stroock and Varadhan (\cite{SV-79} Theorem 6.2.3) is that in order to prove that
the martingale problem has at most one solution it suffices to show that the
one-dimensional marginal distributions $\mathcal{L}(X_{t}),t\geq0,$ are
uniquely determined. \ Moreover in order to determine the law of a random
measure, $X$, on $[0,1]$ it suffices to determine the Laplace functional.

The main step of the proof is to verify that if $P_\mu$ is a solution to the Jirina martingale problem, $t>0,$ and $f\in C_+([0,1])$, then
\be{}
E_\mu(e^{-X_{t}(f)})=e^{-\mu(\psi(t))-c\nu_{0}(\int_{0}^{t}\psi(t-s)ds)}%
\ee
where%
\begin{align*}
\frac{d\psi(s,x)}{ds}  &  =-c\psi(s,x)-\frac{\gamma}{2}\psi^{2}(s,x),\\
\psi(0,x)  &  =f(x).
\end{align*}

\noindent {\em STEP 1: -discretization}\\
We can choose a sequence of partitions $\{A^n_{1},\dots,A^n_{K_n}\}$ and $\lambda^n_1,\dots,\lambda^n_{K_n}$  such that
\be{}\sum_{i=1}^{K_n}\lambda^n_i 1_{A^n_i}\uparrow f(\cdot).\ee
We next show that for a partition $\{A_{1},\dots,A_{K}\}$ of $[0,1]$ and
$\lambda_{i}\geq0,\;i=1,\dots,K,$
\begin{align*}
\exp(-\sum_{i=1}^{K}\lambda_{i}X_{t}(A_{i}))  &  =\exp\left(  -\sum_{i=1}^K\psi
_{i}(t)X_{0}(A_{i})-\sum c\nu_{0}(A_{i})\int_{0}^{t}\psi_{i}(t-s)ds\right) \\
\frac{d\psi_{i}}{ds}  &  =-c\psi_{i}-\frac{\gamma}{2}\psi_{i}^{2}\\
\psi_{i}(0)  &  =\lambda_{i}.%
\end{align*}

To verify this first note that by It\^o's Lemma, for fixed $t$ and $0\leq s\leq
t$,
\begin{align*}
d\psi_{i}(t-s)X_{s}(A_{i})  &  =X_{s}(A_{i})d\psi_{i}(t-s)+\psi_{i}%
(t-s)dX_{s}(A_{i})\\
&  =X_{s}(A_{i})d\psi_{i}(t-s)+\psi_{i}(t-s)c\nu_{0}(A_{i})\\
&  \;\;\;\;\;\;-c\psi_{i}(t-s)X_{s}(A_{i})+\psi_{i}(t-s)dM_{s}(A_{i})\\
\end{align*}
and
\begin{align*}
&  X_{t}(A_{i})\psi_{i}(0)-X_{0}(A_{i})\psi_{i}(t)\\
&  =-\int_{0}^{t}X_{s}(A_{i})\dot{\psi}_{i}(t-s)ds+c\nu_{0}(A_{i})\int_{0}%
^{t}\psi_{i}(t-s)ds\\
&  \;\;\;\;\;\;\;\;-c\int_{0}^{t}X_{s}(A_{i})\psi_{i}(t-s)ds+N_{t}(A_{i})
\end{align*}
where $\{N\}_{0\leq s\leq t}$ is an orthogonal martingale measure with
\begin{align*}
N_{s}(A_{i})  &  =\int_{0}^{s}\psi_{i}(t-u)M(A_{i},du)\\
\left\langle N(A_{i})\right\rangle _{s}  &  =\frac{\gamma}{2}\int_{0}^{s}%
\psi_{i}^{2}(t-u)X_{u}(A_{i})du\\
\left\langle N(A_{i}),N(A_{j})\right\rangle _{s}  &  =0\;\;\text{if \ }i\neq
j\text{.}%
\end{align*}
Again using It\^o's lemma, for $0\leq s\leq t$

\
\begin{align*}
de^{-X_{s}(A_{i})\psi_{i}(t-s)}  &  =\dot{\psi}_{i}(t-s)e^{-X_{s}(A_i)\psi
_{i}(t-s)}ds-\psi_{i}(t-s)e^{-X_{s}(A_i)\psi_{i}(t-s)}dX_{s}(A_i)\\
&  \;\;\;\;\;\;+\frac{\gamma}{2}e^{-X_{s}(A_i)\psi_{i}(t-s)}\psi_{i}^{2}%
(t-s)X_{s}(A_i)ds\\
&  =\dot{\psi}_{i}(t-s)e^{-
X_{s}(A_i)\psi_{i}(t-s)}ds+c\psi_{i}(t-s)e^{-X_{s}(A_i)%
\psi_{i}(t-s)}X_{s}ds\\
&  \;\;\;\;\;\;+c\nu_{0}(A)\psi_{i}(t-s)e^{-X_{s}(A_i)\psi_{i}(t-s)}ds\\
&  \;\;\;\;\;\;+\frac{\gamma}{2}e^{-X_{s}(A_i)\psi_{i}(t-s)}\psi_{i}^{2}%
(t-s)X_{s}ds+dN_{s}(A_{i})\\
&  =c\nu_{0}(A_{i})\psi_{i}(t-s)e^{-X_{s}(A_i)\psi_{i}(t-s)}ds+dN_{s}(A_{i})
\end{align*}

Then by the method of integrating factors we can get

\[
\wt{N}_{s}(A_i)= e^{\left( - X_{s}(A_{i})\psi_{i}(t-s)+c\nu_{0}(A_{i})\int_{s}^{t}%
\psi_{i}(t-u)du\right)  },\quad 0\leq s\leq t,
\]
is a bounded non-negative martingale that can be represented as

\be{}
\wt{N}_{t}(A_i)-\wt{N}_{0}(A_i)=\int_{0}^{t}e^{-\zeta_i\left(  s\right)  }dN_{s}(A_i).
\ee
where
\be{}
\zeta_i(s)    =\left(  c\nu_{0}(A_{i})\int_{s}^{t}\psi
_{i}(t-u)du\right).
\ee

Noting that the martingales $\wt{N}_t(A_i),\wt{N}_t(A_j)$ are orthogonal if $i\ne j$
we can conclude that

\be{}
e^{-\sum_{i}\left(  X_{s}(A_{i})\psi_{i}(t-s)-c\nu_{0}(A_{i})\int_{s}^{t}%
\psi_{i}(t-u)du\right)  },\quad 0\leq s\leq t,\ee
is a bounded martingale. Therefore for each $n$

\be{}  E\left[e^{-\sum_{i=1}^{K_n}\left(  X_{t}(A^n_{i})\psi^n_{i}(0)\right)}\right]=
e^{-\sum_{i=1}^{K_n}\left(  X_{0}(A^n_{i})\psi^n_{i}(t)-c\nu_{0}(A^n_{i})\int_{s}^{t}%
\psi^n_{i}(t-u)du\right)  }
\ee

\noindent {\em STEP 2: Completion of the proof}\\
Taking limits as $n\to\infty$ and  dominated convergence  we can then show that the Laplace functional%
\[
E(e^{-X_{t}(f)})=e^{-X_{0}(\psi(t))-c\nu_{0}(\int_{0}^{t}\psi(t-s)ds)}%
\]
where%
\begin{align*}
\frac{d\psi(s,x)}{ds}  &  =-c\psi(s,x)-\frac{\gamma}{2}\psi^{2}(s,x),\\
\psi(0,x)  &  =f(x).
\end{align*}

Therefore the distribution of $X_t(f)$ is determined for any non-negative continuous function, $f$, on $[0,1].$ Since the Laplace functional characterizes the law of a random measure, this
proves that the distribution at time $t$ is uniquely determined by the
martingale problem. This completes the proof of uniqueness.





(b)  Recall that  $\psi(\cdot,\cdot)$ satisfies
\begin{align*}
\frac{d\psi(x,s)}{ds}  &  =-c\psi(x,s)-\frac{\gamma}{2}\psi^{2}(x,s),\\
\psi(x,0)  &  = f(x)
\end{align*}
Solving, we get%
\[
\psi(x,t)=\frac{f(x)e^{-ct}}{1+\frac{f(x)}{\theta
}-\frac{f(x)}{\theta}e^{-ct}}\allowbreak,\;\theta=\frac{2c}{\gamma}%
\]

 Next, note that $\psi(x,t)\rightarrow0$ as $t\rightarrow\infty$ and
\begin{align*}
\int_{0}^{\infty}\psi(x,s)ds  &  =\int_{0}^{\infty}\frac{f(x)e^{-ct}%
}{1+\frac{f(x)}{\theta}-\frac{f(x)}{\theta}e^{-ct}}%
dt=\int\frac{(-\frac{f(x)}{c})d(e^{-ct})}{1+\frac{f(x)}{\theta
}-\frac{f(x)}{\theta}(e^{-ct})}\\
&  =\frac{2}{\gamma}\int_{0}^{1}\frac{\frac{f(x)}{\theta}%
du}{1+\frac{f(x)}{\theta}-\frac{f(x)}{\theta}u}=\frac{2}{\gamma
}\int_{0}^{\frac{f(x)}{\theta}}\frac{\frac{f(x)}{\theta}%
du}{1+\frac{f(x)}{\theta}-\frac{f(x)}{\theta}u}\\
&  =\frac{2}{\gamma}\log(1+\frac{f(x)}{\theta})
\end{align*}

Therefore
\begin{align*}
E(e^{-X_{t}(f)})
\rightarrow
e^{-\frac{2c}{\gamma}\int\log(1+\frac{f(x)}{\theta})\nu_0(dx)}
\end{align*}

This coincides with the Laplace functional of
\[
\frac{1}{\theta}\int f(s)G(\theta ds)
\]
where $G(\cdot)$ is the Moran subordinator.
Therefore $X_\infty(f)$ can be represented as \[
X_{\infty}(f)=\frac{1}{\theta}\int f(s)G(\theta ds).
\]
\end{proof}

\begin{remark} A more general class of measure-valued branching processes, known as {\em superprocesses} or {\em Dawson-Watanabe processes} will be discussed in Section 9.4.

\end{remark}

\section{The infinitely many alleles martingale problem}

From the construction above we can show that the probability law of the infinitely many alleles Fleming-Viot process
$\{X_{t}:t\geq0\}$ on $C([0,\infty),M_{1}([0,1]))$ satisfies the
\textit{martingale problem}%

\be{fvmp1}
M_{t}(\phi)    :=X_{t}(\phi)-X_{0}(\phi)-\int_{0}^{t}c(\nu_{0}(\phi
)-X_{s}(\phi))ds,
\ee
 is a $L^{2}$ $\mathcal{F}_{t}$ martingale $\forall\;\phi
\in b\mathcal{E}$  with increasing process
\be{fvmp2} \left\langle M(\phi)\right\rangle    =\gamma\int_{0}^{t}(X_{s}(\phi
^{2})-X_{s}(\phi)^{2})ds.
\ee

We now show that this martingale problem completely characterizes the process.
First note that $M_{t}(\phi)$ extends to a martingale measure with
covariation
\[
\left\langle M(A),M(B)\right\rangle _{t}=\int_{0}^{t}Q(X_{s};A,B)ds
\]
where%
\[
Q(\mu;dx,dy)=\mu(dx)\delta_{x}(dy)-\mu(dx)\mu(dy).
\]
We observe that $M$ is a worthy martingale measure \ with dominating measure
$K(dx,dy,ds)=Q(X_{s};A,B)ds,$ because
\[
|\int_{0}^{t}Q(X_{s};A,A)ds|\leq\int_{0}^{t}|Q(X_{s};A,A)|ds\leq t.
\]

In order to prove that the martingale problem is well-posed we will introduce
moment measures. \

Let $X$ be a random probability measure on the Polish space $(E,\mathcal{E})$.
The nth moment measure is a probability measure on $E^{n}$ defined as follows:%
\[
M_{n}(dx_{1},\dots,dx_{n})=E(X(dx_{1}),\dots,X(dx_{n}))
\]

$M_{n}$ is the probability law of n-exchangeable $E$-valued random variables
$(Z_{1},\dots,Z_{n})$. \

\begin{lemma}
(a) A random probability measure $X$ on $E$ is uniquely determined by its
moment measures of all orders.\newline (b) The sequence $\{X_{n}\}$ of random
probability measures with moment measures $\{M_{n,m},n,m\in\mathbb{N}\}$
converges weakly to a random probability measure $X$ with moment measures
$\{M_{m}\}$ as $n\rightarrow\infty$ if and only if $M_{n,m}\Longrightarrow
M_{m}$ for each $m\in\mathbb{N}$.
\end{lemma}

\begin{theorem}
There exists a {unique} solution, $\mathbb{Q}$, to the martingale
problem (\ref{fvmp1}),(\ref{fvmp2}).
\end{theorem}

\begin{proof}
The existence has been proved above.

By the result of Stroock and Varadhan it suffices to show that the
one-dimensional marginal distributions are uniquely determined. \ But from the
Lemma, to determine the law of the random measure, $Z_{t}$, on $[0,1]$ it
suffices to determine all the moment measures.

But by Ito's Lemma, the moment measures satisfy the following system of
equations%
\begin{align*}
&  \frac{\partial M_{n}(t;dx_{1},\dots,dx_{n})}{\partial t}\\
&  =\sum_{i=1}^{n}c[M_{n-1}(t;dx_{1},..,\not x  _{i},..,dx_{n})\nu_{0}%
(dx_{i})-M_{n}(t;dx_{1},\dots,dx_{n})]\\
&  -\frac{1}{2}\gamma n(n-1)M_{n}(t;dx_{1},\dots,dx_{n})\\
&  +\frac{1}{2}\gamma\sum_{i}\sum_{j\neq i}M_{n-1}(t;dx_{1},\dots
,dx_{i-1},\not d  \not x  _{i},dx_{i+1},\dots,dx_{n})\delta_{x_{j}}(dx_{i})\\
&  M_{n}(0;dx_{1},\dots,dx_{n})=\mu(dx_{1})\dots\mu(dx_{n}).
\end{align*}
Then
\[
M_{1}(t,dx)=e^{-ct}X_{0}(dx)+(1-e^{-ct})\nu_{0}(dx)
\]
and we can then solve the remaining equations recursively. \ This implies that
all the moment measures of $Z_{t}$ are uniquely determined by the martingale
problem. \ Hence the martingale problem has a unique solution.
\end{proof}

\section{Dual martingale problems}

Dual processes play an important role in the study of interacting particle systems (see Liggett \cite{L-85}).
A dual representation for the Fleming-Viot process was introduced in Dawson and Hochberg (1982) \cite{DH-82}.
The following generalization with applications to measure-valued processes was established in (Dawson-Kurtz (1982) \cite{DK-82}).
Here we give the main ideas and refer \cite{D-93}, Sect. 5.5 for the details.

To give the main idea we first present the theorem in a simplified case.

\begin{theorem}\label{DR} (Dual Representation)

Let $E_1,E_2$  be Polish spaces and $F(\cdot,\cdot),\;GF(\cdot,\cdot),\;HF(\cdot,\cdot)\in
\mathcal{B}_b(E_1\times E_2)$, $\beta\in \mathcal{B}_b(E_2)$ and $P_x:E_1\to
\mathcal{P}(D_{E_1}(0,\infty))$ and $Q_y:E_2\to
\mathcal{P}(D_{E_2}(0,\infty))$

Assume that
\bea{}  &&F(X(t),y)-\int_0^t GF(X(s),y)ds\text{   is a
}P_{X(0)}\text{ martingale for each }y\in E_2\\&&
F(x,Y(t))-\int_0^t HF(x,Y(s))ds\text{   is a }Q_{Y(0)}\text{
martingale for each }x\in E_1\nonumber\eea and \be{}
GF(x,y)=HF(x,y)+\beta(y)F(x,y) .\ee

Then
\be{}
E^X_x(F(X(t),y))=E^Y_y\left(F(x,Y(t))\exp(\int_0^t\beta(Y(s))ds\right),\quad
\quad 0<t<T\ee

\end{theorem}

\begin{proof}  Let
\be{}  \Phi(s,t):= E^X_x\otimes
E^Y_y\left(F(X(s),Y(t))\exp(\int_0^t\beta(Y(u))du)\right) \ee
\be{} \Phi(t,0)= E^X_x \left(F(X(t),y)\right)\ee \be{} \Phi(0,t)=
E^Y_y \left(F(x,Y(t))\right) \ee

\be{} \Phi_1(s,t)= E^X_x\otimes
E^Y_y\left(GF(X(s),Y(t))\exp(\int_0^t\beta(Y(u))du)\right) \ee

\be{} \Phi_2(t,s)= E^X_x\otimes E^Y_y\left([HF(X(t),Y(s))+\beta
Y(s)F(X(t),Y(s))]\exp(\int_0^s\beta(Y(u))du)\right) \ee
where $\Phi_1,\Phi_2$ are the first partial derivatives with respect to the first and second variables.
Under the assumptions, $\Phi_1(s,t-s),\Phi_2(s,t-s),\;0\leq s\leq t$ exist and are uniformly bounded .

Therefore

\be{}\Phi(0,t)-\Phi(t,0)=\int_0^t \frac{\partial}{\partial
s}\Phi(s,t-s)=\int_0^t (\Phi_1(s,t-s)-\Phi_2(s,t-s))ds=0\ee

\end{proof}

In applications the assumption that $\beta(\cdot)$ and $GF(\cdot,\cdot)$ are bounded needs to be relaxed. The following extension  (see \cite{D-93}, Cor. 5.5.3) provides the required  conditions.

\beP{DUIS}  Assume that

(i) $F\in C_b(E_1\times E_2),$ and $\{F(\cdot,y):y\in E_2)$ is measure-dtermining on $E_1$

(ii) there exist stopping times $\tau_K\uparrow t$ such that
\bea{DUI}  &&\left\{ (1+\sup_x |GF(x,Y(\tau_K))|)\cdot \exp(\int_0^{\tau_K}|\beta(Y(u))|du)\right\}_K\\
&&\text{ are } Q_{\delta_y}-\text{uniformly integrable for all }y\in E_2\nonumber\eea
and
(iii) $Q_{\delta_y}(Y(s-)\ne Y(s))=0$ for each $s\geq 0$, that is, no fixed discontinuities.

Then the $G$-martingale problem is well-posed and for all $y\in E_2$
\be{} P_\mu(F(X(t),y))=\int_{E_1} \mu(dx)\left(Q_{\delta_y}[F(x,Y(t))\exp\left(\int_0^t\beta(Y(u))du\right)\right).\ee
\end{proposition}

\begin{example} (The Wright-Fisher diffusion with polynomial drift)

Let $\Delta_{d-1}=\{(x_1,\dots,x_d), x_i\geq 0,
i=1,\dots,d,\;\sum_{i=1}^d x_i\leq 1\}$

Then consider the Wright-Fisher diffusion $\{x(t)\}$ with generator
\[  G=\sum_{i,j=1}^da_{ij}(x)\frac{\partial ^2}{\partial
x_i\partial x_j}+\sum_{i=1}^db_i(x)\frac{\partial}{\partial x_i}\]
where $\{a_{i,j}(x)\}$ is the real symmetric non-negative definite matrix,
$\{a_{ij}(x)\}=\{x_i(\delta_{ij}-x_j)\}$  and the drift coefficient $b_i(x)$ is a polynomial satisfying certain natural
boundary conditions on $\Delta_{d-1}$ to ensure that the process remains in $\Delta_{d-1}$.

Shiga (1981)) \cite{S-81} obtained a dual in terms of a family of functions
$\{\phi_\alpha\}_{\alpha\in\Gamma},\quad \phi_\alpha\in D(G)$ defined by

\[ \phi_\alpha(x_1,\dots,x_d) =\prod_{i=1}^d x_i^{\alpha_i},\qquad \alpha=(\alpha_1,\dots,\alpha_d)\in \Gamma \]
and showed that
\[ G\phi_\alpha =\sum_\beta
Q_{\alpha,\beta}(\phi_{\beta}-\phi_{\alpha})+h_\alpha\phi_\alpha\]
where
$Q=\{Q_{\alpha,\beta}\}$ defines a conservative Markov chain
$\alpha_t$ with state space $\Gamma$.  Then  the following identity follows from Proposition \ref{DUIS}:

\be{}  E_x[\phi_{\alpha}(x(t))]=
E_\alpha[\phi_{\alpha_t}(x)\exp(\int_0^{t}h_{\alpha_u}du)],\;
0\leq t\leq t_0\ee {provided that  }
\[ E_\alpha[\exp(\int_0^{t_0}|h_{\alpha_u}|du)]<\infty\quad
\forall\alpha\in \Gamma\]  Therefore  the corresponding Wright-Fisher martingale problem is well-posed.

\end{example}

\begin{example}\label{MCdual} (Markov chains) Consider a continuous time Markov chain with state space $E_K=\{1,\dots,K\}$ and transition rates
\be{}  i\to j \text{ with rate }  m_{ij}, \; i\ne j,\; m_{ii}=-\sum_{j\ne i}m_{ij}.\ee

Let $\mathcal{PK}$ denote the collection of subsets of $E_K$ and define the function $F:\mathcal{PK}\times E_K$ by
\be{} F(A,j)=1_A(j)\ee.

Now consider the Markov $\mathcal{A}_t$ chain with state space $\mathcal{PK}$ and transition rates
\be{} A\to A\cup\{j\} \text{  at rate  }\sum_{\ell\in A} m_{j\ell},\; j\in A^c \ee
\be{} A\to A\backslash\{j\} \text{  at rate  }\sum_{\ell\in A^c} m_{j\ell} \quad \text{if }j\in A\ee

\bea{} &&\\GF(A,j)&&= \sum_\ell m_{j\ell}(1_A(\ell)-1_A(j))=\sum_{\ell\in A}m_{j\ell}(1-F(A,j))-\sum_{\ell\in A^c}m_{j\ell}F(A,j)\nonumber\eea

Then \be{} HF(A)=  \sum_{k\in A^c} [(\sum_{\ell\in A} m_{k\ell})(F(A\cup\{k\})-F(A)]+  \sum_{k\in A}[ (\sum_{\ell\in A^c} m_{k\ell})(F(A\backslash \{k\})-F(A))]\ee
and therefore
\bea{} &&\\HF(A,j)&&= \sum_{k\in A^c} (\sum_{\ell\in A} m_{k\ell}))(1_{(A\cup\{k\})}(j)-1_A(j))+  \sum_{k\in A} (\sum_{\ell\in A^c} m_{k\ell})(1_{(A\backslash \{k\})}(j)-1_A(j)) \nonumber\\&&
= \sum_{\ell\in A}m_{j\ell}(1-F(A,j))-\sum_{\ell\in A^c}m_{j\ell}F(A,j).\nonumber \eea

By duality we have
\be{MCD} E_j(1_{\ell}(X_t))=E_{\{\ell\}}(1_{\mathcal{A}_t}(j)).\ee

\begin{remark} If $\{m_{ij}\}$ is irreducible, then  the Markov chain $\mathcal{A}_t$ has two traps $\emptyset$ and $E_K$.  It is easy to verify that $\mathcal{A}_t$ is absorbed at a trap with probability one. This together with (\ref{MCD}) implies  (the elementary result) that $P_j(x(t)=\ell)$  converges as $t\to\infty$  to a stationary measure $\pi_\ell$ with
\be{} \pi_\ell=P_{\{\ell\}} ({\mathcal{A}}_t \to E_K),\quad \ell\in E_K\ee and that  $\lim_{t\to\infty}P_{j}(X(t)=\ell)$ is independent of $j$.

\end{remark}

\end{example}

\section{Dual representation of the neutral Fleming-Viot process}
\label{sdrfv}

The method of dual representation plays an important role in the study of Fleming-Viot processes and will be frequently used below.
To introduce this we first consider  the special case of a neutral Fleming-Viot process with a nice mutation process.

\subsection{The General Neutral F.V. Process}

Let $E$ be a compact metric space, $A$ be a linear operator defined on
$D(A)\subset C(E)$ and assume that the closure of $A$ generates a
Feller semigroup, $\{S_{t}:t\geq0\}$ on $C(E)$. A probability
measure \ $\mathbb{P}_{\mu}$ on $C([0,\infty),M_{f}(E))$ is said
to be a solution of the \textit{neutral Fleming-Viot martingale
problem} $\mathbb{MP}_{(A,Q,0)}$ with initial
condition $\mu$ if%
\[
\mathbb{P}_{\mu}(X_{0}=\mu)=1
\]
and for each $\phi\in C_{b}^{+}(E)\cap D(A)$%
\[
M_{t}^{0}(\phi):=\left\langle \phi,X_{t}\right\rangle
-\left\langle \phi ,X_{0}\right\rangle -\int_{0}^{t}\left\langle
A\phi,X_{s}\right\rangle ds
\]
where $M_{t}^{0}$ defines a martingale measure $M^{0}(ds,dx)$ with covariance%
\begin{align*}
\left\langle M^{0}(dx),M^{0}(dy)\right\rangle _{t}  &  =\gamma\int_{0}%
^{t}Q(X_{s};dx,dy)ds\\
Q(\mu;dx,dy)  &  =\delta_{x}(dy)\mu(dx)-\mu(dx)\mu(dy).
\end{align*}

\begin{theorem}
There exists a unique solution to the $\mathbb{MP}_{(A,Q,0)}$
martingale problem.
\end{theorem}

\begin{proof}
This will be proved in the following section.
\end{proof}

\subsection{Equivalent Formulation of the martingale problem}

We now turn to an equivalent formulation of the Fleming-Viot process that will be needed for the application of the dual representation in the next chapter.

  Let  $ F\in D(G)\subset \in C^2(\mathcal{P}(E))$

\be{}  GF(\mu)=\int_E \left(A\frac{\delta F(\mu)}{\delta \mu(x)}\right)\mu(dx)+\frac{\gamma}{2}\int_E\int_E \frac{\delta^2 F(\mu)}{\delta\mu(x)\delta\mu(y)} Q(\mu;dx,dy)\ee
where  $Q(\mu,dx,dy):= \mu(dx)\delta_x(dy)-\mu(dx)\mu(dy)$.

Now consider function  $F(\mu,(f,n))=\int\dots\int f(x_1,\dots,x_n)\mu^n(dx) $ with $f\in C(E^n)$, $n\in\N$ and
\be{} \mu^n(dx)=\mu(dx_1)\dots\mu(dx_n).\ee
  Then
\be{} GF(\mu,(f,n))= \langle \mu^n,A^{(n)}f\rangle +\frac{\gamma}{2}\sum_{i\ne j}\left(\langle\mu^{n-1},\wt\Theta_{ij}f\rangle -\langle \mu^n,f\rangle\right)\ee

\be{theta} (\wt \Theta_{ij}f)(y_1,\dots,y_{N-1}):= f(x_1,\dots,x_N)\ee
On the right side of (\ref{theta})
\be{}\begin{split}&
x_k=y_k\text{  for  }k<  i\vee j,\; k\ne i\wedge j\\&
x_{i\vee j}=x_{i\wedge j}=y_{i\wedge j}\\&
x_k= y_{k-1}\text{  for  } k>i\vee j.
\end{split}
\ee

\subsection{The dual representation of the Fleming-Viot process}

The Fleming-Viot process has
state space $\mathcal{P}(E)$. We assume that the mutation process has semigroup $S_t$ with generator  $A$
and there exists an algebra of functions  $D_0(E)$ dense in $C(E)$ and $S_t:D_0(E)\to D_0(E)$.

We can then consider the extension of the mutation process
to $E^n,\;n\geq 1$ corresponding to $n$ i.i.d. copies of the basic mutation process and with generator $A^{(n)}=\sum_{i=1}^n A_i$ where $A_i$ denotes the action of $A$ on the ith variable.

Let  \be{}E_2:= \{(f,n): f\in (D_0(E))^n\cap ,\;n\in\N\}.\ee

Define $F: \mathcal{P}(E)\times E_2\to\mathbb{R}$ by
\be{}
F(\mu,(f,n))=  \int_{E^{n}} f_n(x_1,\dots,x_n)\mu(dx_1)\dots\mu(dx_n.\ee

Now consider the Fleming-Viot process with generator:

\be{} GF(\mu,(f,n))=\int_E \left(A\frac{\partial F(\mu,(f,n))}{\partial \mu(x)}\right)\mu(dx)+\frac{\gamma}{2}\int_E\int_E \frac{\partial^2 F(\mu,(f,n))}{\partial\mu(x)\partial\mu(y)} Q(\mu;dx,dy)\ee

and note that for for each  $\mu\in\mathcal{P}(E)$ this coincides with
\be{} HF(\mu,(f,n))= F(\mu,(A^{(n)}f,n))+ \frac{\gamma}{2} \sum_{j=1}^{n}\sum_{k\ne j}[F(\mu,(\wt\Theta_{jk}f,n)-F(\mu,(f,n)))]\ee
where $\wt\Theta_{jk}:(D_0(E))^n\to (D_0(E))^{n-1}$ is defined by (\ref{theta}).

Then
$H$ is the generator of a c\`adl\`ag process with values in
$E_2$ and law $\{Q_f:f\in E_2\}$ which evolves as follows:

\begin{itemize}
\item $Y(t)$ jumps from $(D_0(E)^n,n)$ to $(D_0(E)^{n-1},n-1)$ at rate
$\frac{1}{2}\gamma n(n-1)$ \item at the time of a jump, $f$ is
replaced by  $\wt\Theta_{jk}f$ \item  between jumps, $Y(t)$ is
deterministic on $D_0(E)^n$ and evolves according to the semigroup
$(S^n_t)$ with generator $A^{(n)}$.
\end{itemize}

\beT{} (a) Let $(\{ X(t)\}_{t\geq 0},\{P_\mu:\mu\in
\mathcal{P}(E)\})$ be a solution to the Fleming-Viot martingale
problem and the process $(\{Y(t)\}_{t\geq 0},\{Q_{(f,n)}:(f,n)\in
E_2\}$ be defined as above. Then

(a) these processes are dual, that is, \be{}
P_\mu(F(X(t),(f,n)))=Q_f(F(\mu,Y(t))),\;
\quad (f,n)\in E_2.\ee

(b) The martingale problem is well-posed and the Fleming-Viot process is a strong Markov process.

\end{theorem}
\begin{proof}
In this case for $(f,n)\in E_2$ $\mu\in\mathcal{P}(E)$, \be{}
GF(\mu,(f,n))= HF(\mu,(f,n))\ee
and
the uniqueness follows from Theorem (\ref{DR}).
(b) follows by the Stroock-Varadhan Theorem.
\end{proof}

\subsection{The Kingman coalescent}

 Consider the special case with no mutation, that is, $A\equiv 0$. Then we can represent the dual process $Y(t)$ with $Y(0)=(f,n)$ as follows.

\be{}Y(t)=(f_t,n_t)\ee
where $n_t\leq n$ and there is a map
\be{} \pi_t:\{1,\dots,n\}\to \{1,\dots,n_t\}\ee
and $f_t\in C(E^{n_t})$ given by
\be{} f_t(y_1,\dots,y_{n_t})=f(x_1,\dots,x_n)\quad\text{with  } x_i=y_{\pi_t(i)},\;i=1,\dots,n.\ee

In other words $\pi_t$ is a process with values in the set of partitions of $\{1,\dots,n\}$ and $n_t$ is a pure death process with deaths rate  $\gamma \left(\begin{array}{c}k \\2\end{array}\right)$ where $n_t=k$. This partition-valued process is the {\em Kingman coalescent} \cite{K-82a} and plays an important role in population genetics.

\section{Interactions via change of measure- the Girsanov
formula}

We have established uniqueness for the Fleming-Viot process with diploid selection using duality in the previous section.  However for more general state dependent fitness functions there is no natural dual process. Instead we can use a change of measure argument based on a Girsanov-type formula. Here  we consider the  Girsanov transformation for measure-valued processes
introduced in \cite{D-78}.  We give here a version  suitable for applications in later chapters on selection and logistic competition.

\subsubsection{Preliminaries}

We begin by reviewing with some general notions of stochastic analysis.

Definition. Let \ $(\Omega,\mathcal{F},\{\mathcal{F}_{t}\}_{t\geq0}%
,\mathbb{P})$ be a complete probability space \ such that
$\mathcal{F}_{0}$ contains all $P$-null sets of $\mathcal{F}$ and
$\mathcal{F}_{t}$ is right continuous.

An $\mathcal{F}_{t}-$adapted c\`{a}dl\`{a}g process, $Y$, is a
(\textit{classical) semimartingale }if there exists processes $N$
and $B$ with
$N_{0}=B_{0}=0$ and%
\[
Y_{t}=Y_{0}+N_{t}+B_{t}%
\]
where $N_{t}$ is a local martingale and $B_{t}$ is a finite
variation process.

A generalization  of the classical Girsanov Theorem to semimartingales due to Meyer is as follows.

\begin{theorem} (\cite{Pr-90}, Chap. III, Theorem 20)\\
Let $(\Omega,\mathcal{F},\{\mathcal{F}_{t}\}_{t\geq0},\mathbb{P})$
be as above. Let $X$ be a classical semimartingale under
$\mathbb{P}$, with decomposition $X=M+A$ where $M$ is a local
martingale. Let $\mathbb{Q}$ be
equivalent to $\mathbb{P}$ with%
\[
Z_{t}=E^{\mathbb{P}}\left[  \frac{d\mathbb{Q}}{d\mathbb{P}}|\mathcal{F}%
_{t}\right]
\]
\ Then $X$ is also a classical semimartingale under $\mathbb{Q}$
and has a
decomposition $X=M^{\mathbb{Q}}+C$ where%
\[
M_{t}^{\mathbb{Q}}=M_{t}-\int_{0}^{t}\frac{1}{Z_{t}}d\left[  Z,M\right]  _{s}%
\]
is a $\mathbb{Q}$-local martingale and $C=X-M^{\mathbb{Q}}$ is a $\mathbb{Q}%
$-finite variation process.
\end{theorem}

\begin{proof}
The idea is to use It\^o's Lemma to show that \
$Z_{t}^{-1}M_{t}^{\mathbb{Q}}$ is a $\mathbb{P}$-martingale.
\end{proof}
\bigskip

Let $\{M_t\}$ be a martingale and $f(\cdot) \in L_{\rm{loc}}^2(M)$. Then by It\^o's Lemma  the {\em stochastic exponential}
\be{SE}
Z_{t}:=\exp\left(  \int_{0}^{t}f_{s}dM_{s}-\frac{1}{2}\int_{0}^{t}f_{s}%
^{2}d\left\langle M\right\rangle _{s}\right)\quad \text{which satisfies  } dZ_t= Z_t f_tdM_t
\ee
is a local martingale.

\begin{proposition}
(Novikov's condition) \
A sufficient condition for the stochastic exponential $Z_t$%
to be a martingale is that%
\[
E\left[ \frac{1}{2} \exp\left(  \int_{0}^{t}f_{s}^{2}d\left\langle
M\right\rangle _{s}\right)  \right]  <\infty.
\]
\end{proposition}
See Ikeda-Watanabe \cite{IW-81} Theorem 5.3.

\subsubsection{Girsanov's Transformation for Measure-Valued Processs}

Let $E$ be a locally compact space,
\begin{itemize}
\item $\mathcal{D}(A)$ be a measure-determining linear subspace of $ C_b(E)$ containing constants and $A$ a linear mapping from $\mathcal{D}(A)\to C_b(E)$
\item $Q:M_F(E)\to M_F(E\times E)$ is continuous, and
\be{} Q(\mu,B,B)\leq K \mu(B),\quad K<\infty,\; B\in\mathcal{B}(E),\ee
\item  $V$ is be a measurable function $V:[0,\infty)\times
M_{F}(E)\times
E\rightarrow R$
\end{itemize}

Then a probability measure \ $\mathbb{P}_{\mu}$ on
$C([0,\infty),M_{f}(E))$ is said to solve the \textit{ martingale
problem} $\mathbb{MP}_{(A,Q,V)}$ with
initial condition $\mu$ if%
\[
\mathbb{P}_{\mu}(X_{0}=\mu)=1
\]
and for each $\phi\in \mathcal{D}(A)$%
\begin{align*}
M_{t}^{V}(\phi)  &  :=\left\langle \phi,X_{t}\right\rangle
-\left\langle \phi,X_{0}\right\rangle -\int_{0}^{t}\left\langle
A\phi,X_{s}\right\rangle
ds\\
&  -\int_{0}^{t}\int\int\phi(x){V(s,X_{s},y)}(
Q(X_{s};dx,dy))ds
\end{align*}
is a $P_\mu$-martingale with increasing process
\be{} \langle M^V(\phi)\rangle_t = \int_0^t\int\int\phi(x)\phi(y)Q(X_s,dx,dy)ds.\ee

Then $M_{t}^{V}$ defines a martingale measure $M^{V}(ds,dx)$ with covariation
\be{}
\left\langle M^{V}(dx),M^{V}(dy)\right\rangle _{t}=\int_{0}^{t}%
Q(X_{s};dx,dy)ds.
\ee

\beT{DGIR}
Assume that $\mathbb{P}_{\mu}$ is the $\ $unique solution of the
martingale
problem $\mathbb{MP}_{(A,Q,0)}$ and that $\mathbb{P}_{\mu}$-a.s.%
\be{GA}
\int_{0}^{t}\int\int
V(s,X_{s},x)V(s,X_{s},y)Q(X_{s};dx,dy)ds<\infty,\;\forall
t>0.
\ee

Define the $\mathbb{P}_{\mu}$ continuous local martingales:
\be{} N^V_t=\frac{1}{\gamma}\int_{0}^{t}\int V(s,X_{s}%
,y)M^{0}(ds,dy)\ee \be{}  \langle N^V\rangle_t= \int_{0}^{t}\int\int V(s,X_{s},x)V(s,X_{s}%
,y)Q(X_{s};dx,dy)ds\ee
and the stochastic exponential
\be{} Z_{t}^{V}    :=\exp\left(N^V_t-\frac{1}{2}\langle N^V\rangle_t\right).\ee


(a) (Existence) Assume that (\ref{GA}) holds $\mathbb{P}_\mu$,-a.s. Then $\mathbb{Q}_\mu:=Z^V_t\mathbb{P}_\mu$ is a solution to the
$(A,Q,V)$ local martingale problem.

(b) (Uniqueness) If $\ \mathbb{Q}_{\mu}$ is any solution of the martingale problem $\mathbb{MP}%
_{(A,Q,V)}$ such that (\ref{GA}) holds $\mathbb{Q}_\mu$ a.s., then
\[
\frac{d\mathbb{Q}_{\mu}}{d\mathbb{P}_{\mu}}|_{\mathcal{F}_{t}}=Z_{t}^{V}%
\]
and therefore there is only one such solution.

\end{theorem}

\begin{proof}
Note that
\be{} M^V_t(\phi)=M^0_t(\phi)-\langle N^V,M^0(\phi)\rangle_t.
\ee

(a)  Assume that $\mathbb{P}_\mu$ is a solution to the $(A,Q,0)$ martingale problem and that (\ref{GA}) holds $\mathbb{P}_\mu$-a.s.
By It\^o's formula we have
\bean{} dZ^V_tM^V_t(\phi)&&=d(Z^V_t(M^0_t(\phi)-\langle N^V_t,M^V_t(\phi))\rangle\\
&&=(M^0_t(\phi)-\langle N^V_t,M^V_t(\phi))dZ^V_t + Z^V_tdM^0_t(\phi)-Z^V_t \langle N^V,M^0(\phi)\rangle_t+dZ_t\cdot dM^0_t\\
&&=(M^0_t(\phi)-\langle N^V_t,M^V_t(\phi))dZ^V_t + Z^V_tdM^0_t(\phi)
\eean
since by (\ref{SE})  $dZ^V_t \cdot dM^0_t=  Z^V_t dN^V_t dM^0_t(\phi) =Z^V_td\langle M^0(\phi),N^V\rangle_t$
\bea{} &&\\ &&d(Z_tM^V_t(\phi))=d(Z_t(M^0_t(\phi)-\frac{1}{2}-\int_{0}^{t}\int\int\phi(x)V(s,X_{s},y)(
Q(X_{s};dx,dy))ds))\nonumber\\&&
= M^0_t(\phi)dZ_t+Z_t(dM^0_t -\int\int\phi(x){V(s,X_{s},y)}(
Q(X_{s};dx,dy))dt)\nonumber\\&& +Z_t  \langle M^0_t(\phi),\int_{0}^{t}\int V(s,X_{s}%
,y)M^{0}(ds,dy) \rangle\nonumber\\&&
=M^0_t(\phi)dZ_t+Z_tdM^0_t(\phi)\nonumber
\eea
so that $M^V_t(\phi)$ is a local martingale under $\mathbb{Q}_\mu:=Z^V_t\mathbb{P}_\mu$.
Moreover the quadratic variation of $M^V_t$ is

\be{} \langle M^V(\phi),M^V(\phi)\rangle_t=\int_0^t\int\int\phi(x)\phi(y)Q(X_s;dx,dy)ds,\; \mathbb{Q}_\mu-a.s\ee
In other words $\mathbb{Q}_\mu$ is a solution to the $(A,Q,V)$ local martingale problem.

(b) Now assume that $\mathbb{Q}_\mu$ is a solution to the $(A,Q,V)$-martingale problem and (\ref{GA}) holds $\mathbb{Q}_\mu$ a.s. The same argument as (a) implies that   $M^{0}_t(\phi)$ is a local martingale under $Z^{-V}_t\mathbb{Q}$ and $M^0_t(\phi)Z^{-V}_t$ is a $\mathbb{Q}_\mu$-local martingale. Let
\be{}
\tau_{n}=\inf\{t:\int_{0}^{t}[\int\int(V(s,X_{s},x)V(s,X_{s},y)+1)Q(X_{s}%
;dx,dy)+1]ds\geq n\}\leq n
\ee
Since $\langle N^V\rangle _{t\wedge\tau_n}$ is bounded,  Novikov's critierion implies that $Z^{-V}_{t\wedge \tau_n}$ is a martingale and
\be{} d\mathbb{P}_{\mu,n}:= Z^{-V}_{t\wedge \tau_n}d\mathbb{Q}_\mu \ee
defines probability and $M^0_t(\phi)$ is a $\mathbb{P}_n$-local martingale
and
\be{} \langle M^0_{t\wedge \tau_n}(\phi),M^0_{t\wedge \tau_n}(\phi)\rangle
=\int_0^{t\wedge\tau_n}\int\int \phi(x)\phi(y)Q(X_s;dx,dy)ds\quad \forall \; t\geq 0,\;\mathbb{P}_n-a.s.\ee
Since this is bounded it is integrable and therefore $M^0_t(\phi)$ is a $\mathbb{P}_n$ martingale.  Let
$\mathbb{P}_n^{ex}|\mathcal{F}_{\tau_n}=  \mathbb{P}_n|\mathcal{F}_{\tau_n}$ and $\mathbb{P}^{ex}_n(X_{\tau_n+\cdot}|\mathcal{F}_{\tau_n})=\mathbb{P}_{X_{\tau_n}}(\cdot)$. Then $\mathbb{P}_n^{ex}$ solves the $(A,Q,0)$ martingale problem and therefore since we assumed that this is well-posed, we have $\mathbb{P}_n^{ex}=\mathbb{P}_\mu$. Therefore (\ref{GA}) implies that
\be{}\mathbb{P}_n^{ex}(\tau_n<t)=\mathbb{P}_\mu(\tau_n<t)\to 0\quad \text{ as  } n\to\infty.\ee

Then since we assume that $\mathbb{Q}_\mu$ solves the $(A,Q,V)$ martingale problem and (\ref{GA}) holds $\mathbb{Q}_\mu$ a.s., $M^{-V}_t$ is a $\mathbb{Q}_\mu$ non-negative local martingale.
Then
\bean{} E^{\mathbb{Q}_\mu}(Z^{-V}_t)&&\geq E^{\mathbb{Q}_\mu}(Z^{-V}_{t\wedge \tau_n}  1_{\tau_n\geq t})
\\&& =E^{\mathbb{Q}_\mu}(Z^{-V}_{t\wedge \tau_n})-E^{\mathbb{Q}_\mu}(Z^{-V}_{t\wedge \tau_n} 1_{\tau_n< t})\\&&
= 1-\mathbb{P}_n(\tau_n<t)\to 1 \text{   as  }n\to\infty.
\eean
Hence $Z^{-V}_t$ is a $\mathbb{Q}_\mu$-martingale and we define \be{}\wt{\mathbb{P}}_\mu|{\mathcal{F}_t}=Z^{-V}_td\mathbb{Q}_\mu|{\mathcal{F}_t}.\ee
But then we can verify that  ${M}^0_t(\phi)$ is a $\wt{\mathbb{P}}_\mu$ local martingale and therefore $\wt{\mathbb{P}}_\mu= {\mathbb{P}}_\mu$.
This implies that
\be{}\mathbb{Q}_\mu|_{\mathcal{F}_t}=\frac{1}{Z^{-V}_t}\mathbb{P}_\mu={Z^{V}_t}\mathbb{P}_\mu\ee
and therefore it is unique.

\end{proof}

\beP{}
 If $\sup|V(s,\mu,x)|\leq V_0\text (constant)$,
then $Z_{t}^{V}$ is a martingale under $\mathbb{P}_{\mu}$ and \
$\mathbb{Q}_{\mu}:= Z^V_{t}\mathbb{P}_\mu$ is the unique law that satisfies
$\mathbb{MP}_{(A,Q,V)}.$
\end{proposition}
\begin{proof}
Define
\be{}
\tau_{n}=\inf\{t:\int_{0}^{t}[\int\int(V(s,X_{s},x)V(s,X_{s},y)+1)Q(X_{s}%
;dx,dy)+1]ds\geq n\}\leq n
\ee

\[
V^{n}(s,X,x)=1(s\leq\tau_{n})V(s,X,x).
\]
Then as in the proof of the Theorem $\mathbb{Q}_{\mu,n}:=Z_{\tau_{n}}^{V}d\mathbb{P}_\mu$
is a solution to the
$\mathbb{MP}_{(A,Q,V^{n})}$-martingale problem. Taking $\phi(x)\equiv 1$ we have
\begin{align*}
E^{\mathbb{Q}_{n}}(X_{t}(1))  &  =\mu(1)+E^{\mathbb{Q}_{n}}(\int_{0}^{t\wedge\tau_{n}%
}X_{s}(V(s,X_s,\cdot))ds)\\
&
\leq\mu(1)+V_0 E^{\mathbb{Q}_{n}}(\int_{0}^{t\wedge\tau_{n}}X_{s}(1)ds)
\end{align*}
Then by Gronwall's inequality
$E^{\mathbb{Q}_{n}}(X_{t}(1))\leq\mu(1)e^{V_0t}$ and
therefore%
\begin{align*}
&  E^{\mathbb{Q}_{n}}\left(\int_{0}^{t}[\int\int(V(s,X_{s},x)V(s,X_{s},y)+1)Q(X_{s}%
;dx,dy)]ds\right)\\
&  \leq(V_0^{2}+1)\mu(1)e^{V_0t}t+2t=K(t).
\end{align*}
and then by Chebyshev ${\mathbb{Q}_{n}}(\tau_{n}<t)\leq
K(t)/n\rightarrow0$ and $n\rightarrow\infty$.
Then
\bean{} E^{\mathbb{P}_\mu}(Z^V_t)&&\geq E^{\mathbb{P}_\mu}(Z^V_{t\wedge \tau_n}  1_{\tau_n\geq t})
\\&& =E^{\mathbb{P}_\mu}(Z^V_{t\wedge \tau_n})-E^{\mathbb{P}_\mu}(Z^V_{t\wedge \tau_n} 1_{\tau_n< t})\\&&
= 1-\mathbb{Q}_n(\tau_n<t)
\eean
and therefore $E^{\mathbb{P}_\mu}(Z^V_t)=1$ and $Z^V_t$ is a $\mathbb{P}_\mu$-martingale.

Moreover,
\be{} E^{\mathbb{Q}_\mu}\left\langle M^{V}(\phi )\right\rangle
_{t}=E^{\mathbb{Q}_\mu}\left(\int_{0}^{t}[\int\int(\phi(x)\phi(x)Q(X_{s}
;dx,dy)]ds\right) <\infty
\ee
and therefore $M^V_t(\phi)$ is a $\mathbb{Q}_\mu$-martingale.
\end{proof}

\begin{remark}  Ethier and Shiga (2000,2002) \cite{ES-00}, \cite{ES-02} established the Girsanov formula for a Fleming-Viot
process with unbounded selection which arises as the diffusion limit of a model of Tachida (1991) \cite{Ta-91} with type space $\mathbb{R}$, house of cards mutation with mutation source $\nu_0= N(0,\sigma_0^2)$ and fitness function $V(x)=x$. (The model was proposed by Tachida in the context of an ongoing discussion of the roles of genetic drift and weak selection in protein evolution.)

Overbeck, R\"ockner and Schmuland (1995) studied Fleming-Viot processes with interactive selection using Dirichlet forms.

Evans and Perkins (1994) \cite{EP-94} extend the Girsanov formula to the case of  interacting species modeled by two interacting super-Brownian motions with either competition or predation.  In the case of predation in which collisions effect only the prey they establish existence and uniqueness in dimensions one, two and three.

\end{remark}

We will return to a systematic discussion of mutation-selection systems in Chapter 12.

\chapter{Genealogy and History}
\section{Introduction}

In this Chapter we will focus on a neutral Fleming-Viot process with (or without) mutation. If for the moment we ignore mutation then we can focus on the family relations among the members of the population, for example the ancestral relation between a finite random sample from the population at a fixed time.  This was the purpose of the Kingman coalescent that has become a standard tool of population genetics. With mutation it is also of interest to trace the mutational history of an individual and its ancestors. This is the purpose of the historical process.  In the case of the infinitely many sites model the mutational history of an individual is built into the state of the individual and in this context we can explore the genealogy and mutational history in a unified manner.
More generally, the idea is that  giving the individuals labels or coding for certain
genealogical or historical information can be a useful mathematical tool.

We will introduce some important tools in studying this richer structure including the Kingman
coalescent (Kingman (1982) \cite{K-82a}, \cite{K-82b}), the  look-down process (Donnelly-Kurtz \cite{DK-96}), the tree-valued Fleming-Viot process (Greven-Pfaffelhuber-Winter \cite{GPW-08}) and the analogue of the historical process (Dawson-Perkins \cite{DP-91}).

We will return to the question of the genealogical structure of Fleming-Viot processes with both mutation and selection in Chapter 12.

\section{Family Structure of the neutral Fleming-Viot Process}

\subsection{Fleming-Viot with Feller Mutation Semigroup}

We assume that the space of types is a compact set, $E$. $\ $The mutation
process is assumed to be a Feller process with Feller semigroup $\{S_{t}%
:t\geq0\}$ on $C(E)$ and with transition function $p_{t}(x,dy)$. We assume
that the generator, is the closure of $(A,D(A))$ where $A$ is a linear
operator defined on a linear subspace, $D(A),$ of $C(E).$ \ By the theory of
Feller processes, without loss of generality, we can assume that $D(A)$
contains a countable subset that is convergence determining and that there
exists a c\`{a}dl\`{a}g version of the process, $(D([0,\infty),E),(\mathcal{D}%
_{t})_{t\geq0},\{P_{x}:x\in E\}).$

A probability measure \ $\mathbb{P}_{\mu}$ on $C([0,\infty),\mathcal{P}(E))$ is said
to be a solution of the \textit{neutral Fleming-Viot martingale problem}
$\mathbb{MP}_{(A,\gamma Q,0)}$ (wrt $D(A)$) with initial condition $\mu$ and
resampling rate function $\gamma\in C([0,\infty),\mathbb{R}^{+})$ if%
\[
\mathbb{P}_{\mu}(X_{0}=\mu)=1
\]
and for each $\phi\in C_{b}^{+}(E)\cap D(A),$%
\[
M_{t}(\phi):=\left\langle \phi,X_{t}\right\rangle -\left\langle \phi
,X_{0}\right\rangle -\int_{0}^{t}\left\langle A\phi,X_{s}\right\rangle ds
\]
where $M_{t}$ defines a martingale measure $M(ds,dx)$ with covariance%
\begin{align*}
\left\langle M(dx),M(dy)\right\rangle _{t}  &  =\int_{0}^{t}\gamma
(s)Q(X_{s};dx,dy)ds\\
Q(\mu;dx,dy)  &  =\delta_{x}(dy)\mu(dx)-\mu(dx)\mu(dy).
\end{align*}

\subsubsection{Equivalent Martingale Problem}

For each $n\geq1$, define the Feller semigroup $\{S_{t}^{(n)}:t\geq0\}$ on
$C(E^{n})$ by%
\begin{align*}
&  S_{t}^{(n)}f(x_{1},\dots,x_{n})\\
&  :=\int\dots\int f(y_{1},\dots,y_{n})p_{t}(x_{1},dy_{1})\dots p_{t}%
(x_{n},dy_{n})
\end{align*}
and let $A^{(n)}$ denote its generator.

We can also consider integrals with respect to the martingale measure of the
form%
\begin{align*}
M_{t}(f)  &  :=\int\dots\int\gamma(s)f(x_{1},\dots,x_{n})M(ds,dx_{1}%
)M(ds,dx_{2})\dots M(ds,dx_{n}),\;\;\\
f  &  \in C_{sym}(E^{n})
\end{align*}
and using It\^o's lemma verify that%
\[
\left\langle M(f)\right\rangle _{t}=\sum_{1\leq i<j\leq n}\int_{0}^{t}%
\gamma(s)\left(  \left\langle \Phi_{ij}^{(n)}f,X_{s}^{n-1}\right\rangle
-\left\langle f,X_{s}^{n}\right\rangle \right)  ds
\]
where%
\bea{}
&&\Phi_{ij}^{(n)}    :C(E^{n})\rightarrow C(E^{n-1})\\
&&\Phi_{ij}^{(n)}f(x_{1},\dots,x_{n-1})    :=f(x_{1},\dots,x_i,\dots,x_{j-1},x_{i}%
,x_{j},\dots,x_{n-1}).\nonumber
\eea

(Hint: First do this for linear combinations of functions of the form
$f(x_{1},\dots,x_{n})=\prod_{i=1}^{n}\varphi_{i}(x_{i})$ and then take limits
in $\ L^{2}.)$

Now consider the collection, $D(G),$ of subset of $C(M_{1}(E))$ of all linear
combinations of functions of the form%
\be{poly}
F_{f}(\mu)=F(f,\mu):=\left\langle f,\mu^{n}\right\rangle ,\;\;n\in\mathbb{N},\;f\in
C_{sym}(E^{n})\cap D(A^{(n)}).
\ee
A function of the form (\ref{poly}) is called a {\em polynomial of degree }$n$.

(Note that $D(G)$ is an algebra of functions in $C(M_{1}(E))$ that separates
points and is therefore dense by the Stone-Weierstrass theorem.)

Note that since the resampling rate $\gamma(t)$ is not assumed to be constant we have a time-inhomogeneous Markov process.
In this case we can verify that for each $f\in
C_{sym}(E^{n})\cap D(A^{(n)})$,
\be{}
M_{f}(t):=F_{f}(X_{t})-\int_{0}^{t}G_{t}F_{f}(X_{s})ds\quad\text{is a martingale}
\ee
where%
\be{}
G_{s}F_{f}(\mu)=\left\langle A^{(n)}f,\mu^{n}\right\rangle +\sum_{1\leq
i<j\leq n}\gamma(s)(\left\langle \Phi_{ij}^{(n)}f,\mu^{n-1}\right\rangle
-\left\langle f,\mu^{n}\right\rangle ),\quad 0\leq s\leq t.
\ee

We thus obtain a second martingale problem formulation for the Fleming-Viot
process, $\mathbb{MP}_{(D(G),\gamma\Phi)}$. It turns out that this is
equivalent to the $\mathbb{MP}_{(A,\gamma Q,0)}.$

\subsubsection{A Function-valued Dual}

For each $\mu\in M_1(E)$, let
\bea{IHD}&\\
H_{s}F(f,\mu)  &  =\left\langle A^{(n)}f,\mu^{n}\right\rangle +\sum_{1\leq
i<j\leq n}\gamma(t-s)(\left\langle \Phi_{ij}^{(n)}f,\mu^{n-1}\right\rangle
-\left\langle f,\mu^{n}\right\rangle )\quad 0\leq s\leq t,\nonumber
\eea
where we can interpret $H_{s}$ as the generator of a function valued
dual process \be{}\{Y_{s}:0\leq s\leq t\}\ee which has jumps $f\rightarrow\Phi_{ij}^{(n)}f$
at rate $\gamma(s)$ for each pair $0\leq i<j\leq n$, and between jumps evolves
according to the semigroup $f\rightarrow S_{t}^{(m)}f$ $\ $if $f\in C(E^{m}).$

This yields the duality relation%

\[
E_{\mu}[\left\langle F(f,X_{t})\right\rangle ]=E_{f}[F(Y_t,\mu)]
\]
where $E_{f}$ denotes expectation with respect to the law of the
function-valued process, $Y_{t}$, starting at $Y_{0}=f$.

\begin{proof} The proof is analogous to the proof of
Theorem \ref{DR} but modified to take into account the that we are now working with a time inhomogeneous process.
The key step requires that
\be{} G_sF(X(s),Y(t-s))=H_{t-s}(f(X(s),Y(t-s))\ee
which is satisfied by the choice (\ref{IHD}).
\end{proof}

\begin{remark}
This  is a consequence of the fact that the dual process looks backwards in time.
\end{remark}

\subsection{Two Finite Particle Systems and Moment Measures}

Countable exchangeable particle representations of random measures and their genealogical structures have proved to be useful in study the properties
of random measures and measure-valued processes (see \cite{DH-82}).    In the study of Fleming-Viot processes a construction (now known as the {\em look-down process}) of Donnelly and Kurtz (\cite{DK-96})
has become a standard tool in this subject.   This will be described below. For simplicity, in this section we consider the homogeneous case, $\gamma(t)\equiv \gamma$.

\subsubsection{The n-Particle Look-Down Process}

We begin with a graphical construction.  For each $N\in\N$ let $\mathcal{I}_N=\{1,\dots,N\}$ and  consider a collection of independent rate $\gamma$ Poisson point processes:
\be{}
\{(N^{LD}_{j,i}(t))_{t\geq 0}\}_{1\leq i<j\leq N}.
\ee
We consider $N$ levels indexed by $i=1,\dots,N$ and at each jump of the of the process $N_{j,i}$ we draw an vertical arrow from level
$j$ down to level $i$. At  time $t$ a {\em path} from $i$ back to $j$ at time $s$ is a specified  by $n\in \N$,  a sequence $i=i_{n}\geq i_{n-1}>\dots >i_{1}>i_{0}=j$ and times $s\leq u_1\leq u_2\leq \dots \leq u_n\leq t$ where the last jump (back in time) from level $i_{k}$ was to level $i_{k-1}$, and occurred at time $u_k$.

Given $(i,t)$ there is a unique  ancestor at time $0\leq s<t$
\be{AD} A_s(i,t)\in\mathcal{I}_N,\;\;0\leq s\leq t\ee such that there is a path from $(A_s(i,t),s)$ to $(i,t)$.

\begin{figure}[h]
\begin{center}
\includegraphics[scale=0.4]{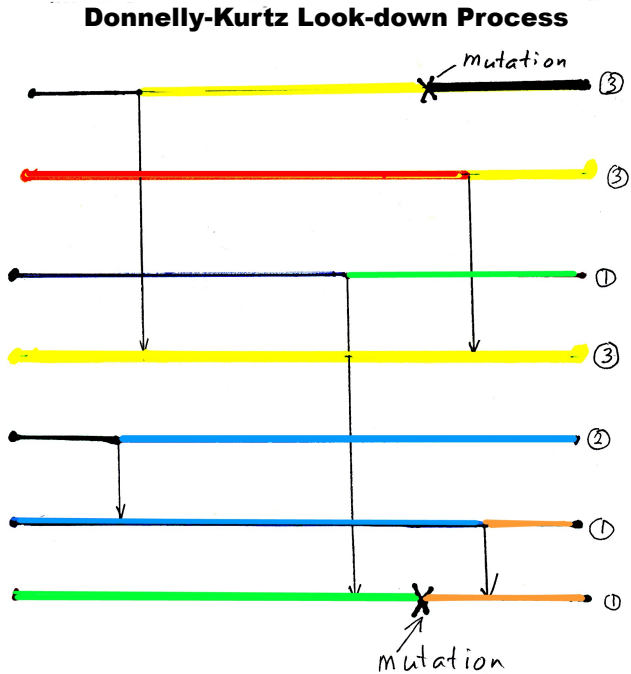}
\caption{Graphical construction of the look-down process}
\end{center}
\end{figure}

We can then define a pseudometric on $\mathcal{I}_N$ (as in \cite{GPW-08})
\be{PM} d_t(i,j):=\left\{\begin{split}
&2(t-\sup\{s\in [0,t]: A_s(i,t)=A_s(j,t)
\})\quad \text{if  } A_0(i,t)=A_0(j,t),\\&
2t+r_0(A_0(i,t),A_0(j,t)),\quad \text{if  } A_0(i,t)\ne A_0(j,t).\end{split}\right.
\ee
This induces as usual a metric space by passing to equivalence classes.

At time $t$ we can decompose $\{1,\dots,N\}$ into equivalence classes where two points $i,j$ belong to the same class if $A_0(i,t)=A_0(j,t)$. Each equivalence class then defines a tree and the set of equivalence classes defines a forrest.
We then obtain a {\em forrest-valued process} that becomes a {\em rooted tree-valued process}  with a constant number,  $N$, of leaves after a finite time.

Now consider a system of $N$ particles $\zeta_1(t),\dots , \zeta_N(t)$ moving in the space $E$ according to a Markov process with generator

\bea{ast0}
&&  \mathcal{C}^{N}f(x_{1},\dots,x_{N})\\
&&  =A^{(N)}f(x_{1},\dots,x_{N})+\sum_{1\leq i<j\leq N}\gamma\{f(\theta
_{ij}(x_{1},\dots,x_{N})-f(x_{1},\dots,x_{N})\}\nonumber
\eea
and where $\theta_{ij}:E^{N}\rightarrow E^{N}$ is defined by $\theta
_{ij}(x_{1},\dots,x_{N}),\;i<j,\;$ is the element of $E^{N}$ obtained from
$(x_{1},\dots,x_{N})$ by replacing the jth component by the ith.

$\mathcal{C}_{t}^{N}$ can be identified as the  generator of
an $N$-particle system (the finite Donnelly-Kurtz look-down process) in which
the dynamics is as follows:

\begin{itemize}
\item at rate $\gamma$  the particle with label $j$ makes a jump
to the location in $E$ of particle $i$ (with $i<j$)

\item between jumps of the previous type the particles perform independent
copies of the mutation process.$\bigskip$
\end{itemize}

\subsubsection{The Moran n-Particle Process}

The second n-particle process in $E$, i.e. $E^{n}$-valued process, $Y^{(n)},$
has generator%
\begin{align*}
&  \tilde{C}^{n}f(x_{1},\dots,x_{n})\\
&  =\frac{1}{2}\sum_{i\neq j}\gamma\{f(\theta_{ij}(x_{1},\dots
,x_{n}))-f(x_{1},,x_{n})\}+A^{(n)}f(x_{1},\dots,x_{n}).
\end{align*}
and define%
\[
\eta_{t}^{(n)}\equiv\frac{1}{n}\sum_{i=1}^{n}\delta_{Y_{i}^{(n)}(t)}.
\]
This is the  continuous time \textit{Moran model} for a
population of size $n$. (We will see below that the sequence $\{\eta_{t}%
^{(n)}:t\geq0\}$ is tight and that every limit point satisfies the
Fleming-Viot martingale problem.)

\begin{lemma}
The empirical measure processes%
\[
\eta_{t}^{(n)}\equiv\frac{1}{n}\sum_{i=1}^{n}\delta_{Y_{i}^{(n)}(t)}%
\overset{\mathcal{L}}{=}Z_{t}^{(n)}:=\frac{1}{n}\sum_{i=1}^{n}\delta
_{\zeta_{i}^{(n)}(t)}%
\]
and the resulting measure-valued process is Markov.
\end{lemma}

\begin{proof}
To verify this note that both satisfy the same martingale problem \ given by
$G|\{F_{f}:f\in C_{sym}(E^{n})\}.$ To verify this it suffices to check
that the generators $C^{n}$ and $\tilde{C}^{n}$ agree on symmetric
functions. But if $f(x_{1},\dots,x_{n})$ is symmetric, then%
\begin{align*}
C^{n}f(x_{1},\dots,x_{n})  &  =\sum_{1\leq i<j\leq n}\gamma\{f(\theta
_{ij}(x_{1},\dots,x_{n})-f(x_{1},\dots,x_{n})\}+A^{(n)}f(x_{1},\dots,x_{n})\\
&  =\frac{1}{2}\sum_{i\neq j}\gamma\{f(\theta_{ij}(x_{1},\dots
,x_{n}))-f(x_{1},,x_{n})\}+A^{(n)}f(x_{1},\dots,x_{n})\\
&  =\tilde{C}^{n}f(x_{1},\dots,x_{n}).
\end{align*}

It then that it suffices to show the solutions to these two martingale
problems both have the same one dimensional marginals. But this follows since
they have the same moment measures%
\[
E_{\mu}\left[  \int\dots\int f(x_{1},\dots,x_{k})\prod_{i=1}^{k}\eta_{t}%
^{(n)}(dx_{i})\right]  =E_{\mu}\left[  \int\dots\int f(x_{1},\dots,x_{k}%
)\prod_{i=1}^{k}X_{t}^{(n)}(dx_{i})\right]
\]
for any $k,n\in\mathbb{N}$ and that these quantities are constant for $n>k$.
But the terms inside the [$\cdot$] is a sum of symmetric functions of
$(x_{1},\dots,x_{n})$ and therefore the expectations agree.
\end{proof}

In other words we have two particle encodings of the same measure-valued process.

\begin{remark}
Donnelly and Kurtz \cite{DK-96} show that one can also construct a coupling of $(\zeta
_{1}(t),\dots,\zeta_{n}(t))$ and $(Y_{1},\dots,Y_{n})$ so that $\zeta^{(n)}$
is a random permutation of $Y^{(n)}$ that is independent of $\eta_{t}^{(n)}$.
\end{remark}

\begin{remark}
The relationship%
\[
GF_{f}(\mu)=GF(f,\mu)=\left\langle \mathcal{C}^{N}f,\mu^{N}\right\rangle
,\forall\,\,f\in D(A^{(N)})\cap C(E^{N})\;\;\forall\;\mu\in M_{1}(E)
\]
implies that for any solution to the Fleming-Viot martingale problem,
$\{X_{t}\}$, and solution $\{x_{1}(t),x_{2}(t),\dots,x_{N}(t)\}$ of the
$\mathcal{C}$-martingale problem and $N\in\mathbb{N},$%

\bea{}&&\\
&& E_{\mu^{N}}\left[  \int\dots\int f(x_{1}(t),\dots,x_{N}(t))\right]\nonumber\\&&=  \int\dots\int E_{x_1,\dots,x_N} f(x_{1}(t),\dots,x_{N}(t))\mu(dx_1)\dots,\mu(dx_N)\nonumber\\&&  =E_{\mu
}\left[  \int\dots\int f(x_{1},\dots,x_{N})\prod_{i=1}^{N}X_{t}(dx_{i}%
)\right]=E_\mu(F(f,X_t)).\nonumber
\eea
(The proof of this uses the fact that both can be represented by the same
function-valued dual $\{Y(t)\}_{t\geq 0}$.)
\end{remark}

\subsubsection{Tree-valued Moran processes}

The tree-valued Moran process is constructed as in Greven-Pfaffelhuber-Winter \cite{GPW-08}.
It is obtained in the as described above for the look-down tree process based on the corresponding
graphical description.


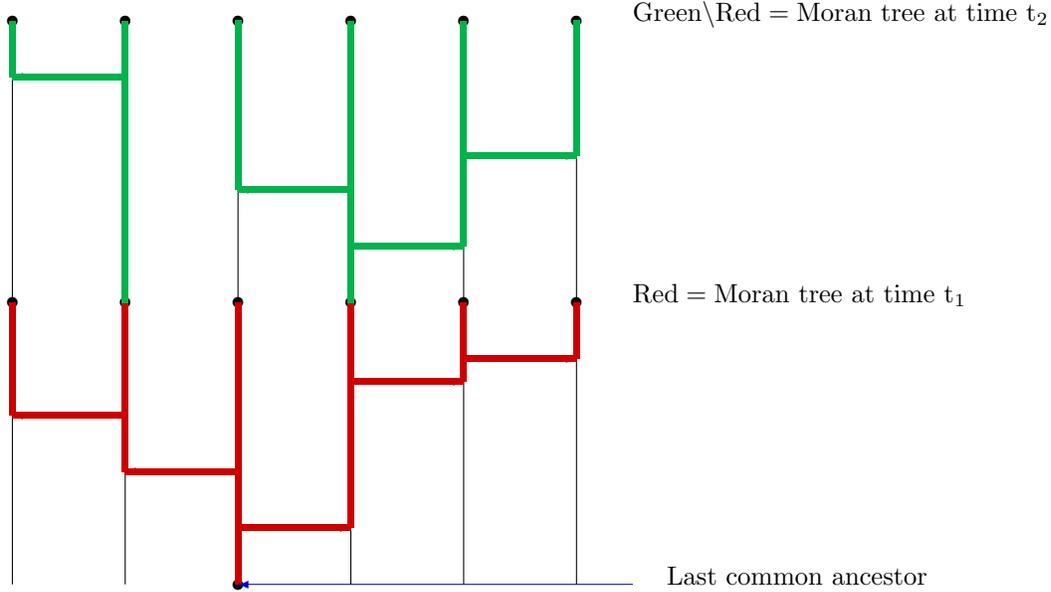
\begin{figure}[h]
\begin{center}
\setlength{\unitlength}{0.15cm}
\begin{picture}(150,60)
\put(0,50){\circle*{1}}
\put(10,50){\circle*{1}}
\put(20,50){\circle*{1}}
\put(30,50){\circle*{1}}
\put(40,50){\circle*{1}}
\put(50,50){\circle*{1}}
\put(0,25){\circle*{1}}
\put(10,25){\circle*{1}}
\put(20,25){\circle*{1}}
\put(30,25){\circle*{1}}
\put(40,25){\circle*{1}}
\put(50,25){\circle*{1}}
\put(20,0){\circle*{1}}
\put(0,0){\line(0,1){50}}
\put(10,0){\line(0,1){50}}
\put(20,0){\line(0,1){50}}
\put(30,0){\line(0,1){50}}
\put(40,0){\line(0,1){50}}
\put(50,0){\line(0,1){50}}
\linethickness{0.8mm}

\put(10,15){\red\vector(-1,0){10}}
\put(20,10){\red\vector(-1,0){10}}
\put(40,20){\red\vector(1,0){10}}
\put(30,18){\red\vector(1,0){10}}
\put(20,5){\red\vector(1,0){10}}
\put(20,5){\red\line(0,-1){5}}
\put(0,25){\red\line(0,-1){10}}
\put(10,25){\red\line(0,-1){15}}
\put(20,25){\red\line(0,-1){25}}
\put(50,25){\red\line(0,-1){5}}
\put(40,25){\red\line(0,-1){7}}
\put(30,25){\red\line(0,-1){20}}
\put(55,50){$\rm{Green\backslash Red }=\text{Moran tree at time }t_2$}
\put(55,25){$\rm{Red }=\text{Moran tree at time }t_1$}
\put(10,45){\green\vector(-1,0){10}}
\put(30,35){\green\vector(-1,0){10}}
\put(40,38){\green\vector(1,0){10}}
\put(30,30){\green\vector(1,0){10}}
\put(0,50){\green\line(0,-1){5}}
\put(10,50){\green\line(0,-1){25}}
\put(20,50){\green\line(0,-1){15}}
\put(30,50){\green\line(0,-1){25}}
\put(40,50){\green\line(0,-1){20}}
\put(50,50){\green\line(0,-1){12}}
\linethickness{0.05mm}
\put(55,0){\blue\vector(-1,0){35}}
\put(58,0){\text{Last common ancestor}}
\end{picture}
\end{center}
\caption{Graphical Representation of the Moran Tree at two times. }
\end{figure}

As before  $\mathcal{I}_N=\{1,\dots,N\}$.
Let  $\{N^{M}_{i,j}:1\leq i,j\leq N\}$ be a realization of a family of rate $\gamma/2$ Poisson point processes.  We say that for
$i,j\in \mathcal{I}_N$  and for $0<s<t<\infty$ there is a path of descent from $(i,s)$ to $(j,t)$ if there exists $n$, $s\leq u_1<u_2\leq\dots <u_n\leq t$ and $i_1,\dots, i_n\in\mathcal{I}_N$ such that for all $k\in \{1,\dots, n-1\}$ (and putting $i_0=i$ and $i_n=j$), $N^{M}_{i_{k-1},i_k}[u_{k-1};u_k]= N^M_{i_{k-1},i_k}\{u_k\}=1$ and $N^M_{m,i_{k-1}}[u_{k-1},u_k)=0$ for all $m\in \mathcal{I}_N$ as well as $N^M_{m,i}[s,u_1)=N^M_{m,j}(u_n,t]=0$.
We define $A_s(i,t)$ and the pseudometric $d_t(i,j)$  as in (\ref{PM}) and (\ref{AD}).

Let $\mathbb{M}$ denote the set of equivalence classes of metric measure spaces (see Appendix I, Section 16.5.1). We call $\mathcal{U}^N=(\mathcal{U}^N_t)_{t\geq 0}$ the {\em tree-valued Moran dynamics} with population size $N$, where for $t\geq 0$
where $\mathcal{U}^N_t\in\mathbb{M}$ is the equivalence class of the metric measure space
\be{}  \mathcal{U}^N_t := \overline{(\mathcal{I},d^N_t,\frac{1}{N}\sum_{i\in\mathcal{I}}\delta_i)}.\ee

\subsection{Extension to Infinite Particle Systems: Some Preliminaries}

In this section we state for convenience some well-known results that are
needed in the next section.

\begin{theorem}
(de Finetti's Theorem)\newline (a) Let $P\in M_{1}(M_{1}(E))$. Then there
exists a sequence $\{Z_{n}\}$ of $E$-valued exchangeable random variables
defined on a probability space $(\Omega,\mathcal{G},P_{dF})$ with
$\Omega=E^{\mathbb{N}},\;$ such that $(Z_{1},\dots,Z_{n})$ has joint
distribution%
\[
P^{(n)}(dx_{1},\dots,dx_{n})=\int_{M_{1}(E)}\mu(dx_{1})\dots\mu(dx_{n}%
)P_{dF}(d\mu),\;n\in\mathbb{N}%
\]
\newline (b) Consider the sequence of $\{Z_{n}\}$ of $E$-valued exchangeable
random variables. Let
\[
X_{n}(\omega)=\frac{1}{n}\sum_{i=1}^{n}\delta_{Z_{i}}.
\]
Then%
\[
X(\omega)=\lim_{n\rightarrow\infty}X_{n}(\omega)\in M_{1}(E)\text{ exists for
}P_{dF\text{ }}\text{a.e. }\omega
\]
where the limit is taken in the weak topology on $M_{1}(E)$ and $X$ has
probability law $P$.\newline (c) Let $\mathcal{G}_{n}$ =$\sigma(X_{n}%
,Z_{n+1},Z_{n+2},\dots).$Then $\mathcal{G}_{\infty}:=\cap\mathcal{G}_{n}$\ is
the $\sigma$-algebra of exchangeable events. Then $X$ is $\mathcal{G}_{\infty
}$-measurable and conditioned on $\mathcal{G}_{\infty},\;\{Z_{n}\}$ is a
sequence of i.i.d. random variables with marginal distribution $X$.
\end{theorem}
\begin{proof}
See \cite{D-93}, Theorem 11.2.1.\emph{\emph{}}
\end{proof}

\begin{corollary}
\label{defincor}Given $f\in C(E),$ and $\varepsilon>0,$ there exists $C>0$ and
$\eta>0$ (both depending only on $\varepsilon$ and $\Vert f\Vert$) such that%
\[
P\left(  \left|  \int f(x)X_{n}(dx)-\int f(x)X(dx)\right|  \geq\varepsilon
\right)  \leq Ce^{-n\eta}.
\]
\end{corollary}

\begin{proof}
This follows immediately from the de Finetti disintegration and Azuma's inequality (see Appendix Lemma \ref{cramer}).
\end{proof}

\subsection{The Countable Donnelly-Kurtz Look-Down Process}

In this section we describe the look-down process of Donnelly and Kurtz (\cite{DK-96}). This
process provides a representation of the Fleming-Viot process in terms of an
exchangeable infinite particle system.

Note that the processes with generators $\{\mathcal{C}^{n}:n\in\mathbb{N\}}$
are consistent. Then taking the projective limit we obtain an infinite
particle system described as follows. For any $n\in\mathbb{N}$ and $f\in
D(A^{(n)})$, let

\bigskip%
\bea{ast02}
&  \mathcal{C}f(x_{1},\dots,x_{n})\\
&  =A^{(n)}f(x_{1},\dots,x_{n})+\sum_{1\leq i<j\leq n}\gamma\{f(\theta
_{ij}(x_{1},\dots,x_{n})-f(x_{1},\dots,x_{n})\}\nonumber
\eea
and where $\theta_{ij}:E^{n}\rightarrow E^{n}$ is defined as above.

$\mathcal{C}$ can be identified as the (time-dependent) generator of an
$\infty$-particle system, $(\zeta_{1}(t),\zeta_{2}(t),\dots)_{t\geq0}$, (the
Donnelly-Kurtz look-down process) in which the dynamics is as follows:

\begin{itemize}
\item at rate $\gamma$ at time $t$ the particle with label $j$ makes a jump to
the location of particle $i$ (with $i<j$)

\item between jumps of the previous type the particles perform independent
copies of the mutation process with generator $A$.$\bigskip$
\end{itemize}

\begin{remark}
The relationship%
\[
GF_{f}(\mu)=\left\langle \mathcal{C}^{n}f,\mu^{n}\right\rangle ,\forall
\,\,f\in D(A^{(n)})\cap C(E^{n})\;\;\forall\;\;\mu\in M_{1}(E)
\]
implies that for any solution to the Fleming-Viot martingale problem,
$\{X_{t}\},\,$and the solution $\{\zeta_{1}(t),\zeta_{2}(t),\dots\}$ of the
$\mathcal{C}$-martingale problem and $n\in\mathbb{N},$%
\[
E_{\mu^{n}}\left[  \int\dots\int g(\zeta_{1}(t),\dots,\zeta_{n}(t))\right]
=E_{\mu}\left[  \int\dots\int g(x_{1},\dots,x_{n})\prod_{i=1}^{n}X_{t}%
(dx_{i})\right]  .
\]
(The proof of this uses the fact that both can be represented by the same
function-valued dual.)
\end{remark}

\begin{lemma}
If $(\zeta_{1}(0),\zeta_{2}(0),\dots)$ is an exchangeable \ sequence, then for
any fixed $t$, $(\zeta_{1}(t),\zeta_{2}(t),\dots)$ is an exchangeable sequence,
\end{lemma}

\begin{proof}
From the above we have%
\[
\int_{M_{1}(E)}E_{\mu}[\left\langle f,X_{t}^{n}\right\rangle ]\nu(d\mu
)=\int_{M_{1}(E)}E_{\mu^{\infty}}[f(\zeta_{1}(t),\dots,\zeta_{n}(t)]\nu
(d\mu).
\]
The \ left side is the expectation for a Fleming-Viot process with initial
distribution $\nu$ and the right side is the expectation for the particle
system under the assumption that $(\zeta_{1}(0),\zeta_{2}(0),\dots)$ is an
exchangeable sequence with%
\[
P(\zeta_{1}(0)\in B_{1},\dots,\zeta_{n}(0)\in B_{n})=\int_{M_{1}(E)}%
\prod_{i=1}^{n}\mu(B_{i})\nu(d\mu)
\]
Then%
\[
P(\zeta_{1}(t)\in B_{1},\dots,\zeta_{n}(t)\in B_{n})=\int\int\prod_{i=1}%
^{n}X_{t}(B_{i})P_{\mu}(dX)\nu(d\mu)
\]
for all $t\geq0$ where $P_{\mu}$ is the law of the Fleming-Viot process
starting at $\mu$ and $\nu$ is the law of $X_{0}.$ Hence if $(\zeta
_{1}(0),\zeta_{2}(0),\dots)$ is an exchangeable \ sequence, then for fixed
$t$, $(\zeta_{1}(t),\zeta_{2}(t),\dots)$ is an exchangeable sequence.
\end{proof}

\bigskip Since $(\zeta_{1}(t),\zeta_{2}(t),\dots)$ is an exchangeable
sequence, the corresponding de Finetti measure
\[
Z_{t}:=\lim_{n\rightarrow\infty}\frac{1}{n}\sum_{i=1}^{n}\delta_{\zeta_{i}(t)}%
\]
exists a.s. in the weak topology and has the same distribution as $X_{t}$. Let
$\mathcal{G}_{t}^{n}=\sigma\{Z^{n}(s),\zeta_{n+1}(s),\zeta_{n+2}(s),\dots
\}$,\ $\mathcal{G}_{t}=\cap_{n}\mathcal{G}_{t}^{n}.$

It also follows from de Finetti's theorem that%
\be{astast}
E[f(\zeta_{1}(t),\dots,\zeta_{k}(t))|\mathcal{G}_{t}]=\left\langle
f,Z^{(k)}(t)\right\rangle .
\ee
\bigskip

In fact, we will show that the process $\{Z_{t}:t\geq0\}$ is a version of the
Fleming-Viot process. Given a convergence determining class of functions
$\{f_{n}\}$ (with $\Vert f_{n}\Vert\leq1$ for each $n)$ on E consider the
metric $\rho$defined by
\[
\rho(\mu,\nu)=\sum_{n}\frac{1}{2^{n}}|\left\langle f_{n},\mu\right\rangle
-\left\langle f_{n},\nu\right\rangle | .
\]

\beT{DK-96T}(Donnelly-Kurtz \cite{DK-96})
(a) Let $\zeta=(\zeta_{1},\zeta_{2},\dots)$ be a Markov process in $E^{\infty
}$ with generator $C$ and suppose that $(\zeta_{1}(0),\zeta_{2}(0),\dots)$ is
exchangeable. Then for each $t>0$, $(\zeta_{1}(t),\zeta_{2}(t),\dots)$ is
exchangeable, and the process given by the de Finetti measure
\begin{align*}
Z_{t}  &  =\lim_{n\rightarrow\infty}Z_{t}^{(n)}\\
Z_{t}^{(n)}  &  :=\frac{1}{n}\sum_{i=1}^{n}\delta_{\zeta_{i}(t)}%
\end{align*}
is a continuous $M_{1}(E)$-valued process that is a solution to the
Fleming-Viot martingale problem with mutation operator $A$.\newline (b) With
probability one, $Z_{t}^{(n)}$ converges uniformly (in t) in the weak topology
on $M_{1}(E)$.
\end{theorem}

\begin{proof}
(a) We first note that the sequence $\{Z_{t}^{(n)}\}$ is tight in
$D([0,\infty),M_{1}(E))$. Since there is a countable convergence determining
class of functions in $D(G)$ it suffices to prove that for $\{F_{f}%
(Z_{t}^{(n)})\}$ is tight in $D([0,\infty),\mathbb{R}),$ with $f\in
D(A^{(n)})$ for some $n\in\mathbb{N}$. But since
\[
M_{f}(t):=F_{f}(Z_{t}^{(m)})-\int_{0}^{t}GF_{f}(Z_{s}^{(m)})ds
\]
is a bounded martingale for $m>n$, we can verify this by applying
Lemma (Appendix I, \ref{L.tightness}). Moreover since the jump sizes go uniformly to zero, this
implies that $\{f(Z_{t}^{(\infty)})\}$ is continuous, a.s. (e.g. Theorem 10.2,
Chapt. 3, Ethier and Kurtz \cite{EK-86}).

Moreover from the above we know that the moment equations are the moment
equations of the Fleming-Viot process. (This can be extended o joint moments
at a finite set of times $t_{1}<t_{2}<\dots<t_{k}.)$ This gives the existence
of a solution to the Fleming-Viot martingale problem. By the uniqueness proved
above this means that that $Z$ is a version of Fleming-Viot.

(b) \ By Corollary (\ref{defincor}) we have%
\[
P\left\{  \left|  \int f(x)Z_{n}(t,dx)-\int f(x)Z(t,dx)\right|  \geq
\varepsilon\right\}  \leq Ce^{-n\eta}%
\]
where $C$ and $\eta$ depend only on $\varepsilon$ and $\Vert f\Vert$.We will
show that in fact the processes $Z_{t}^{(n)}$ converges uniformly in $t$ to a
solution of the Fleming-Viot martingale problem. \ Since there is a countable
convergence determining class of functions in $D(A)$ it suffices to prove this
for $f(Z_{t}^{(n)}),\;f\in D(A).$

Let $R_{i}(t,h)=1$ if $\zeta_{i}$ ``looks down'' during the time interval
$(t,t+h]$ and 0 otherwise. Note that $P(R_{i}(t,h)=0)=e^{-(i-1)h}.$ For
$\varepsilon>0$%
\bea{ast}
&  P\left(  \sup_{t\leq s<t+h}\left|  \int f(x)Z_{{}}^{(n)}(s,dx)-\int
f(x)Z^{(n)}(t,dx)\right|  \geq\varepsilon\right) \\
&  =P\left(  \sup_{t\leq s<t+h}\left|  \frac{1}{n}\sum_{i=1}^{n}f(\zeta
_{i}(s))-\frac{1}{n}\sum_{i=1}^{n}f(\zeta_{i}(t))\right|  \geq\varepsilon
\right) \nonumber\\
&  \leq P\left(  \sup_{t\leq s<t+h}\left|  \frac{1}{n}\sum_{i=1}^{n}\left[
(f(\zeta_{i}(s))-f(\zeta_{i}(t))-\int_{t}^{s}Af(\zeta_{i}(u))du)\right]
\right.  \right. \nonumber\\
&
\;\;\;\;\;\;\;\;\;\;\;\;\;\;\;\;\;\;\;\;\;\;\;\;\;\;\;\;\;\;\;\;\;\;\;\;\;\;\;\;\;\;\;\;\;\;\left.
\left.  1_{\{R_{i}=0\}})\bigskip\right|  >\frac{\varepsilon}{4}\right)
\nonumber\\
&  +P\left(  (2\Vert f\Vert+h\Vert Af\Vert)\left|  \frac{1}{n}\sum_{i=1}%
^{n}(R_{i}(t,h)-1+e^{-(i-1)h}\right|  >\frac{\varepsilon}{4}\right)
\nonumber\\
&  +P\left(  \frac{1}{n}\sum_{i=1}^{n}\int_{t}^{t+h}|Af(X_{i}%
(u))|du>\frac{\varepsilon}{4}\right) \nonumber\\
&  +P\left(  (2\Vert f\Vert+h\Vert Af\Vert\frac{1}{n}\sum_{i=1}^{n}%
(1-e^{-(i-1)h})\geq\frac{\varepsilon}{4}\right) \nonumber
\eea

(The first and third terms come from $1_{\{R_{i}=0\}}.$The second and fourth
terms comes from observing that on $1_{\{R_{i}=1\}}$ the increment is bounded
by ($2\Vert f\Vert+h\Vert Af\Vert).$

The independence of the $R_{i}$ from the evolution of the $\zeta_{i}$ between
look-downs implies that the process in the first term on the right is a
martingale. By Doob's inequality this is less that or equal to
\[
\inf_{\lambda>0}\frac{1}{\Phi(\lambda)}E\left\{  \Phi\left(  \lambda
\times\frac{1}{n}\sum_{i=1}^{n}\left[  (f(\zeta_{i}(t+h))-f(\zeta_{i}%
(t))-\int_{t}^{t+h}Af(\zeta_{i}(u))du)\right]  (1-R_{i}(t,h))\right)
\right\}
\]
for any convex function $\Phi.$ Then by part (b) of the large deviation Lemma \ref{cramer}
(for a sum of bounded independent zero mean r.v's) that there exits $C$ and
$\eta>0$ (depending only on $\varepsilon$ and $\ (\Vert f\Vert+h\Vert
Af\Vert)$) such this term is bounded by $Ce^{-n\eta}.$

The second term is bounded by a similar expression. The third and fourth terms
are zero if $h\Vert Af\Vert<\frac{\varepsilon}{4}$ and $(2\Vert f\Vert+h\Vert
Af\Vert)(1-e^{-h(n-1)})<\frac{\varepsilon}{4}.$ Therefore $C$ and $\eta$ may
be selected depending only on $\varepsilon,\Vert f\Vert$ and $\Vert Af\Vert$
such that for $h$ sufficiently small,
\[
P\left(  \sup_{t\leq s<t+h}\left|  \int f(x)Z_{{}}^{(n)}(s,dx)-\int
f(x)Z^{(n)}(t,dx)\right|  \geq\varepsilon\right)  \leq Ce^{-\eta n}.
\]
Let $h_{n}\rightarrow0$ slowly enough so that $\sum e^{-\eta n}/h_{n}<\infty$
for every $\eta>0$ and fast enough so that $nh_{n}\rightarrow0$ (e.g.
$h_{n}=n^{-2})$. For $T>0,$ let $H_{T,n}=\{kh_{n}:k\leq T/h_{n}\}$. Let $f\in
D(A)$, recall that $\int f(x)Z(\cdot,dx)$ is continuous and define%
\[
D_{n}=\left\{  \sup_{t\leq T}\sup_{t\leq s\leq t+h_{n}}\left|  \int
f(x)Z(s,dx)-\int f(x)Z(t,dx)\right|  <\varepsilon\right\}  .
\]
Note that $D_{n}\subset D_{n+1},$ and by continuity of $\int f(x)Z(\cdot
,dx),\;P(D_{n})\rightarrow1.$ For $n\;$sufficiently large ($h_{n}$
sufficiently small),%
\begin{align*}
&  P\left(  \left\{  \sup_{t\leq T}\left|  f(x)Z_{n}(t,dx)-\int
f(x)Z(t,dx)\right|  \geq3\varepsilon\right\}  \cap D_{n}\right) \\
&  \leq\sum_{t\in H_{T,n}}P\left(  \left\{  \left|  f(x)Z_{n}(t,dx)-\int
f(x)Z(t,dx)\right|  \geq\varepsilon\right\}  \right) \\
&  +\sum_{t\in H_{T,n}}P\left(  \sup_{t\leq s\leq t+h_{n}}\left\{  \int
f(x)Z_{n}(s,dx)-\int f(x)Z_{n}(t,dx)\geq\varepsilon\right\}  \right) \\
&  \leq\frac{2CT}{h_{n}}e^{-\eta n}.
\end{align*}
where we have used the Corollary to de Finetti's theorem for the first term.
Summing over $n$, the right side converges, and Borel-Cantelli and the
properties of $D_{n}$ ensure that%
\[
\lim_{n\rightarrow\infty}\sup_{t\leq T}\left|  \int f(x)Z_{n}(t,dx)-\int
f(x)Z(t,dx)\right|  =0
\]
and the proof is complete.
\end{proof}

\begin{corollary}
Let $\tau$ be a finite $\{\mathcal{G}_{t}\}$ stopping time. Then $\{\zeta
_{1}(\tau),\zeta_{2}(\tau),\dots\}$ is exchangeable and
\[
Z(\tau)=\lim_{n\rightarrow\infty}\frac{1}{n}\sum\delta_{\zeta_{i}(\tau)}.
\]
\end{corollary}

\begin{proof}
If $\tau$ is discrete, then (\ref{astast}) implies the exchangeability and $Z(\tau
)=\lim_{n\rightarrow\infty}\frac{1}{n}\sum\delta_{\zeta_{i}(\tau)}$.  For the
general case, let $\tau_{n}$ be a decreasing sequence of $\{\mathcal{G}_{t}\}$
stopping times converging to $\tau$. By the right continuity of the
$\zeta,\,(\zeta_{1}(\tau_{n}),\zeta_{2}(\tau_{n}),\dots)\rightarrow(\zeta
_{1}(\tau),\zeta_{2}(\tau),\dots)$ exchangeability follows. This together with
the uniform convergence gives $Z(\tau)=\lim_{n\rightarrow\infty}\frac{1}%
{n}\sum\delta_{\zeta_{i}(\tau)}.$
\end{proof}

\subsubsection{ Explicit Construction of the Look-Down Process}

\begin{theorem}
Let $\{S_{t}\}$ be a Feller semigroup on $C(E)$ with $E$ compact. Then
\newline (a) for each $x\in E$ there exists a probability measure $P_{x}$ on
$\mathcal{B}(D_{E})\;$satisfying%
\[
P_{x}(\omega(0)=x)=1,\text{ }%
\]
and for $s\leq t$%
\[
P_{x}(f(\omega(t))|\sigma(\omega(u):u\leq s)=(S_{t-s}f)(\omega(s)),\;P_{x}%
\text{-a.s. }\forall\;f\in C(E)
\]
\newline (b) There exists a standard probability space $(\Omega^{A}%
,\mathcal{F}^{A},Q_{A})$ and a measurable mapping $\zeta:(E\times\Omega
^{A},\mathcal{E}\otimes\mathcal{F}^{A})\rightarrow(\mathcal{D}_{E}%
,\mathcal{B}(D_{E}))$ such that for each $x\in E$
\[
Q_{A}(\{\omega:\zeta(x,\omega)\in B\})=P_{x}(B)\text{\ }\forall\;B\in
\mathcal{B}(D_{E})
\]
Furthermore, $\zeta(.,\omega)$ is continuous at $x$ for $Q_{A}$-a.e. $\omega$,
for each $x\in E$.
\end{theorem}

\begin{proof}
(a) It is well -known that a Feller process has a c\`{a}dl\`{a}g version.
(This can be verified by obtaining it as the weak limit of the jump processes
associated with the Yosida approximation, $A_{n}=A(I-\frac{1}{n}A)^{-1})$ of
$A$.

(b) The mapping $x\rightarrow P_{x}$ from $E$ to $M_{1}(D_{E})$ is continuous
if the latter is given the weak topology. To see this consider first the
finite dimensional distributions, $0<t_{1}<\dots<t_{n}$%
\[
E_{x}(f_{n}(\omega(t_{n}))\dots f_{1}((t_{1})))=(S_{t_{1}}f_{1}(\dots
S_{t_{n-1}-t_{n-2}}(f_{n-1}S_{t_{n}-t_{n-1}}f_{n})))(x)
\]
By the Feller property, $S_{t}f(x)$ is continuous in $x$ for all choices of
$f_{1},\dots,f_{n}$. This implies that the finite dimensional distributions
are continuous in $x$. It now suffices to show that the measures $\{P_{x}:x\in
E\}$ are relatively compact in $D([0,\infty),E)$. \ To get tightness of paths
in $D([0,\infty),E)$ it suffices to show tightness of $\{f(x(t)):t\geq0\}$ in
$D([0,\infty),\mathbb{R)}$ for each $f\in D(A).$ Since $f(x(t))-\int_{0}%
^{t}Af(x(s))ds$ is a bounded martingale for each $f\in D(A),$ we can verify
the latter condition by using the tightness lemma. Hence we get that if
$x_{n}\rightarrow x$, then $P_{x_{n}}\Longrightarrow P_{x}$. Since the mapping
$x\rightarrow P_{x}$ is continuous, existence of a representation $(\Omega
^{A},\mathcal{F},Q_{A},\{\xi(x)\}_{x\in E}),\;\xi:\Omega^{A}\times
E\rightarrow D_{E}$ (for each $x\in E,\;$ $\xi(\cdot,x)$ is measurable on
$\Omega^{A})$ follows from the extension of Skorohod's almost sure
representation theorem due to Blackwell and Dubins (1983) \cite{BD-83}. In their
representation for each $x\in E$, $\xi(\cdot,\cdot)$ is almost surely
continuous at $x$. It remain to show that there exists a jointly measurable
version. \ We will construct a jointly measurable function $\zeta(\cdot
,\cdot)$ such that at each for each $x\in E$,
\[
\zeta(\omega,x)=\xi(\omega,x)\text{ a.e. }\omega
\]
that is, $\zeta$ is a version of $\xi$. To do this let $\{x_{m}:m=1,2,\dots\}$
be an enumeration of a countable dense set in $E$. Let $\rho$ be a complete
separable metric on $D([0,\infty),E)$. Consider the finite measurable
partitions of $E^{(n)}=\cup_{m}E_{m}^{(n)}$ where
\[
x\in E_{m}^{(n)}\text{ if }x_{m}\text{ is the closest among }\{x_{m}\text{,
}m\leq n\}\text{ to }x
\]
and in the case of ties $x$ is assigned to the smallest such $x_{m}$. We then
define the jointly measurable functions%
\[
\tilde{\zeta}_{n}(\omega,x)=\xi(\omega,x_{m})\text{ }\in D([0,\infty),E)\text{
if }x\in E_{m}^{(n)}.
\]
In particular \ $\zeta_{n}(\omega,x_{k})=\xi(\omega,x_{k})$ for all
sufficiently large $n$. Now define the jointly measurable function
\[
\eta(\omega,x)=\lim_{n\rightarrow\infty}\max_{n^{\prime},n^{\prime\prime}\geq
n}\{\rho(\tilde{\zeta}_{n^{\prime}}(\omega,x),\tilde{\zeta}_{n^{^{\prime
\prime}}}(\omega,x))\}
\]
and the jointly measurable function%
\[
\zeta(\omega,x):=1_{\eta=0}(\omega,x)\lim_{n\rightarrow\infty}\tilde{\zeta
}_{n}(\omega,x)+1_{\eta>0}(\omega,x)\zeta_{x}^{0}%
\]
where $\zeta_{x}^{0}$ is the constant function $\zeta_{x}^{0}\equiv x.$ Note
that for each $x\in E$%
\[
\zeta(\omega,x)=\xi(\omega,x)\text{ }Q_{A}\text{-a.e. }\omega
\]
since $\xi(\omega,x_{m})\rightarrow\xi(\omega,x)$ if $x_{m}\rightarrow x$ for
a.e. $\omega$ by the defining property of $\xi$. This means that $\zeta$ is a
version of $\xi$.
\end{proof}

\begin{example}
For Brownian motion we can simply take $\zeta(\omega,x)=x+W(\cdot)$ where
$W(\cdot)$ is a standard Brownian motion starting at $0$.
\end{example}

\begin{theorem}
Let $(\Omega^{\mathcal{P}},\mathcal{F}^{\mathcal{P}},Q_{\mathcal{P}})$ denote
a probability space on which there is defined a rate 1 Poisson process, $N$
and $(\Omega^{A},\mathcal{F}^{A},Q_{A})$ the probability space on which we
have defined the mutation process, $\zeta,$ as above. Given $\mathbf{\zeta
(0)=}\{\zeta_{1}(0),\zeta_{2}(0),\dots\}$ there exists a measurable process,
$\{\zeta_{i}(t):i\in\mathbb{N\}}$ on the probability space%
\[
\Omega^{LD}:=((\Omega^{A},\mathcal{F}^{A},Q_{A})^{\mathbb{N}}\times(\Omega
^{A},\mathcal{F}^{A},Q_{A})^{\mathbb{N}^{3}}\times(\Omega^{\mathcal{P}%
},\mathcal{F}^{\mathcal{P}},Q_{\mathcal{P}})^{\mathbb{N}^{2}})
\]
with law given by that of the look-down process started at $\mathbf{\zeta}(0).$
\end{theorem}

\begin{proof}
$\omega\in\Omega^{LD}$ has the form $\omega=((\omega^{1})_{i\in\mathbb{N}%
},(\omega^{2})_{ijk\in\mathbb{N}^{3}},(\omega^{3})_{ji\in\mathbb{N}^{2}})$.
\ For $i\in\mathbb{N}$, let $U_{i0}=\omega_{i}^{1},$ for $1\leq j<i<\infty$
and $k\geq1$ let $U_{jik}=\omega_{jik}^{2}$ and for $1\leq j<i<\infty,$ let
$N_{ji}=\omega_{ji}^{3}$. . Thus the $\{U_{jik},U_{i0}\}$ are independent
copies of the mutation process and $\{N_{ij}\}$ are independent Poisson
processes. \ Put $N_{i}=\sum_{j:j<i}N_{ji}$. The dynamics of the system
$\{\zeta_{i}(\cdot)\}_{i\in\mathbb{N}}$ is as follows.

\begin{itemize}
\item Until the first jump in $N_{i}$, $\zeta_{i}$ evolves according to
$U_{i0}(\zeta_{i}(0),.).$

\item If the kth jump of $N_{ji}$, occurs at time $\tau_{ji}^{k},$ then
$\zeta_{i}$ assumes the value of $\zeta_{j}$ at time $\tau_{ji}^{k}$ and then
evolves according to $U_{jik}(\zeta_{i}(\tau_{ji}^{k}),\cdot)$ until the next
jump of $N_{i}$. \
\end{itemize}

It is then easy to verify that the system $\{\zeta_{i}(\cdot)\}_{i\ni
\mathbb{N}}$ is a version of the look-down process.
\end{proof}

\begin{remark}
The process satisfies the following system of stochastic integral equations.
\ For any $f\in C(E),$%
\begin{align*}
f(\zeta_{i}(t))  &  =f(U_{i,0}(\zeta_{i}(0),t))1(\sum_{j<i}N_{ji}(t)=0)\\
&  +\sum_{j=0}^{i-1}\int_{0}^{t}f(U_{j,i,N_{ji}(s)}((\zeta_{j}%
(s),t-s)))-f(\zeta_{i}(s-)))dN_{ji}(s).
\end{align*}
\end{remark}

\subsection{Genealogy  and the Kingman Coalescent}

In this section we describe the embedding of the genealogical tree in \ the
countable particle system.

If the mutation process is stationary, then we can consider the look-down
process on the time interval $(-\infty,\infty)$ and assume that $\{\zeta
_{1}(t):-\infty<t<\infty\}$ is stationary. In this stationary case we can
trace the ancestry of a particle by following the process backward in time.
For $s<t$ we define $a_{j}(s,t)$ to be the level of the ancestor at time $s$
of the jth level particle at time $t$. To be precise, for $s<t$ let
$N_{j}(s,t]=\sum_{i:i<j}N_{ij}(s,t].$ Define $\gamma_{i}(t)=\sup
\{u<t:N_{j}(u,t]>0\}$ (last time before $t$ that the particle $j$ passed on
its type to a particle at a lower level) . Let $\alpha_{j}(\gamma_{j}(t))$ be
the index $i$ such that $\gamma_{j}(t)\in N_{ij},$ that is, the index of .
Define $a_{j}(s,t)=j$ for $\gamma_{j}(t)\leq s<t$ and $a_{j}(s,t)=\alpha
_{j}(\gamma_{j}(t))$ for $\gamma_{\alpha_{j}(\gamma_{j}(t))}(\gamma
_{j}(t))\leq s<\gamma_{j}(t)$ and extend the definition of $a_{j}(s,t)$ to all
$s<t$ in the obvious way.

For $0<s<t<\infty,$ let
$\Gamma_{n}(s,t):=\{a_{j}(s,t):j=1,\dots,n\}$, that is,
$\Gamma_{n}(s,t)$ \ is the set of indices of particles at time $s$
that have descendants among the 1st $n$ particles at time $t$.
Then letting $N_{t,n}(s):=|\Gamma_{n}(s,t)|$ to denote the
cardinality of $\Gamma_{n}(s,t)$ it follows that
$D_{n}(u)=|\Gamma_{n}((t-u)-,t)|$ is a pure death process with
transition intensity $\left(
\begin{array}
[c]{c}%
k\\
2
\end{array}
\right)  $ from state $k$ and $D_{n}(0)=n.\,$Let $\tau_{n,k}=\inf
\{u\geq0:D_{n}=k\}$. Then%
\[
E[\tau_{n,k}]=\sum_{m=k+1}^{n}\frac{1}{\left(
\begin{array}
[c]{c}%
m\\
2
\end{array}
\right)  }=\frac{2}{k}-\frac{2}{n}%
\]
converges as $n\rightarrow\infty$. Then $D(u)=\lim_{n\rightarrow\infty}%
D_{n}(u)<\infty$ for all $u>0$.

Let $N_{t,n}(s)$ denote the number of distinct individuals
(ancestors) at time $s$ that have descendants among
$(\zeta_{1}(t),\dots\zeta_{n}(t))\;$and let $\Gamma_{n}(s,t)$
denote the collection of indices of these $N_{t,n}(s)$ particles.
Since $N_{t,n}(s)$ is monotone increasing in $s$, we can associate
a binary branching Feller process in a natural way. \ For $u>0$
let $R^{n}(u)$ denote the equivalence relation of $\{1,\dots,n\}$
where $i$ and $j$ are in the same equivalence class iff they have
the same ancestors at time $t-u$,
that is, $i\sim j$ if $a_{i}((t-u)-,t)=a_{j}((t-u)-,t).$ Let $\mathfrak{C}%
^{n}$ denote the set of equivalence classes on $\{1,\dots,n\}$. Let
$D_{t,n}(u):=N_{t,n}(t-u)-),\;0\leq u\leq t,$ the number of equivalence
classes in $R^{n}(u)$ (we take right continuous versions of all processes.)

\begin{theorem} (Kingman)
(a) The $\mathbb{N}$-valued process $\{D_{t,n}(u):n\geq0\}$ is a pure death
process with death rates $d=\left(
\begin{array}
[c]{c}%
k\\
2
\end{array}
\right)  $ and $N_{t,n}(t-0)=n$.\newline (b) Let $N_{t}(s):=\lim
_{n\rightarrow\infty}N_{t,n}(s)$ denote the number of distinct ancestors at
time $s$ of the infinite set of particles $\{\zeta_{1}(t),\zeta_{2}%
(t),\dots\}$. Then for $s<t$, $N_{t}(s)<\infty,$ a.s.\newline (c) Let
$D_{t}(u):=N_{t}((t-u)-).$ Then $\{D_{t}(u):u>0\}$ is a Markov pure death
process started from an entrance boundary at $\infty$ with death rates
$d_{k}=\left(
\begin{array}
[c]{c}%
k\\
2
\end{array}
\right)  .$\newline
\end{theorem}

\begin{proof}
(a) The times between jumps in the $n$-particle look-down process are i.i.d.
exponential r.v.'s with mean $\frac{2}{n(n-1)}.$ Therefore the time since the
last look-down is exponential with mean $\frac{2}{n(n-1)}.$ To obtain the
second-to-last look-down time we then consider the resulting $(n-1)$-particle
system and the distribution of its last look-down is exponential with mean
$\frac{2}{(n-1)(n-2)}.$ Continuing in this way we get that the time between
the $(k-1)$st last look-down and kth last look-down is an exponential r.v.
with mean $\frac{2}{(n-k+1)(n-k)}$. Since the times between these look-downs
are also independent we conclude that $\{D_{t,n}(s):s\geq0\}$ is a pure death
process with death rates $d_{k}=\frac{k(k-1)}{2}$.\newline (b) Let $\tau
_{n,k}:=\inf\{s:D_{t,n}(s)=k\}$. Then $\tau_{n,k}=2E_{1}/(n(n-1))+\tau
_{n-1,k}$ where $E_{1}$ and $\tau_{n-1,k}$ are independent r.v. and $E_{1}$ is
exponential with mean $1$. \ From this we obtain the representation
$\tau_{n,k}=2E_{1}/(n(n-1))+\dots+2E_{n-k}/(k+1)k)$ where $\{E_{m}\}$ are iid
Exp(1) r.v.'s. Since $\sum_{j=k}^{\infty}\frac{1}{j(j+1)}<\infty$, we conclude
that $\lim_{\rightarrow\infty}\tau_{n,k}<\infty$ a.s. Consequently,
$N_{t}(s)<\infty$ a.s. if $s<t$.\newline (c) This follows from (a), the
consistency of the processes $\{N_{t,n}(\cdot):n\in\mathbb{N\}}$ and the
construction of $N_{t}$ as the projective limit of the $\{N_{t,n}(\cdot
):n\in\mathbb{N\}}$.
\end{proof}

Note that the processes $\{R^{n}(u)\}_{n\in\mathbb{N}}$ are
consistent and we can take the projective limit, $\{R(u)\}$. This
process can be described as follows. Let $\Gamma(s,t)$ be the
collection of indices of particles at time $s<t$ that have a
descendent at time $t$ in the infinite look-down processes. By the
last Theorem (b) $\Gamma(s,t)$ is a.s. finite and is therefore
associated with an equivalence relation on $\mathbb{N}$ having a
finite number of equivalence classes which we denote by $R(t-s)$.
In other words, $R(u)$ is the equivalence relation on $\mathbb{N}$
in which $i$ and $j$ belong to the same equivalence class iff they
have the same ancestors at time $t-u$. Let
$\mathfrak{C}\subset2^{\mathbb{N\times N}}$ denote the set of
equivalence relations on $\mathbb{N}$ with the subspace topology
when $2^{\mathbb{N\times N}}$ is given the product topology. Then
$\mathfrak{C}$ is a compact metrizable space. A probability
measure on $\mathfrak{C}$ is called exchangeable if it is
invariant under permutations on $\mathbb{N}$.

From the limiting argument \ above we conclude that $R(s)$ is a $\mathfrak{C}%
$-valued continuous time Markov chain called \textit{Kingman's}
\textit{coalescent} which is characterized by the property that its
restriction to $\{1,\dots,n\}$ is the coalescent described above.

We will now consider the genealogical development in ``forward time''. Note
that $t-\tau_{1}=\sup\{s:$ all particles at time $t$ have a common ancestor at
time $s\}.$ Define $\bar{D}(s)=D((\tau_{1}-s)-),\;\bar{R}(s)=R((\tau
_{1}-s)-).$ Then $\bar{D}$ is a pure birth process with $\bar{D}(0)=2$ and
birth rates $\frac{k(k-1)}{2}$ and a.s. finite explosion time $\hat{\tau
}_{\infty}:=\lim_{k\rightarrow\infty}\hat{\tau}_{k}$ where $\hat{\tau}%
_{k}:=\inf\{s:\bar{D}(s)=k\}$. We denote by $\bar{D}^{t},\bar{R}^{t}$ the
corresponding processes conditioned on $\{\hat{\tau}_{\infty}=t\}$.

\begin{figure}[h]
\begin{center}
\includegraphics[scale=0.40]{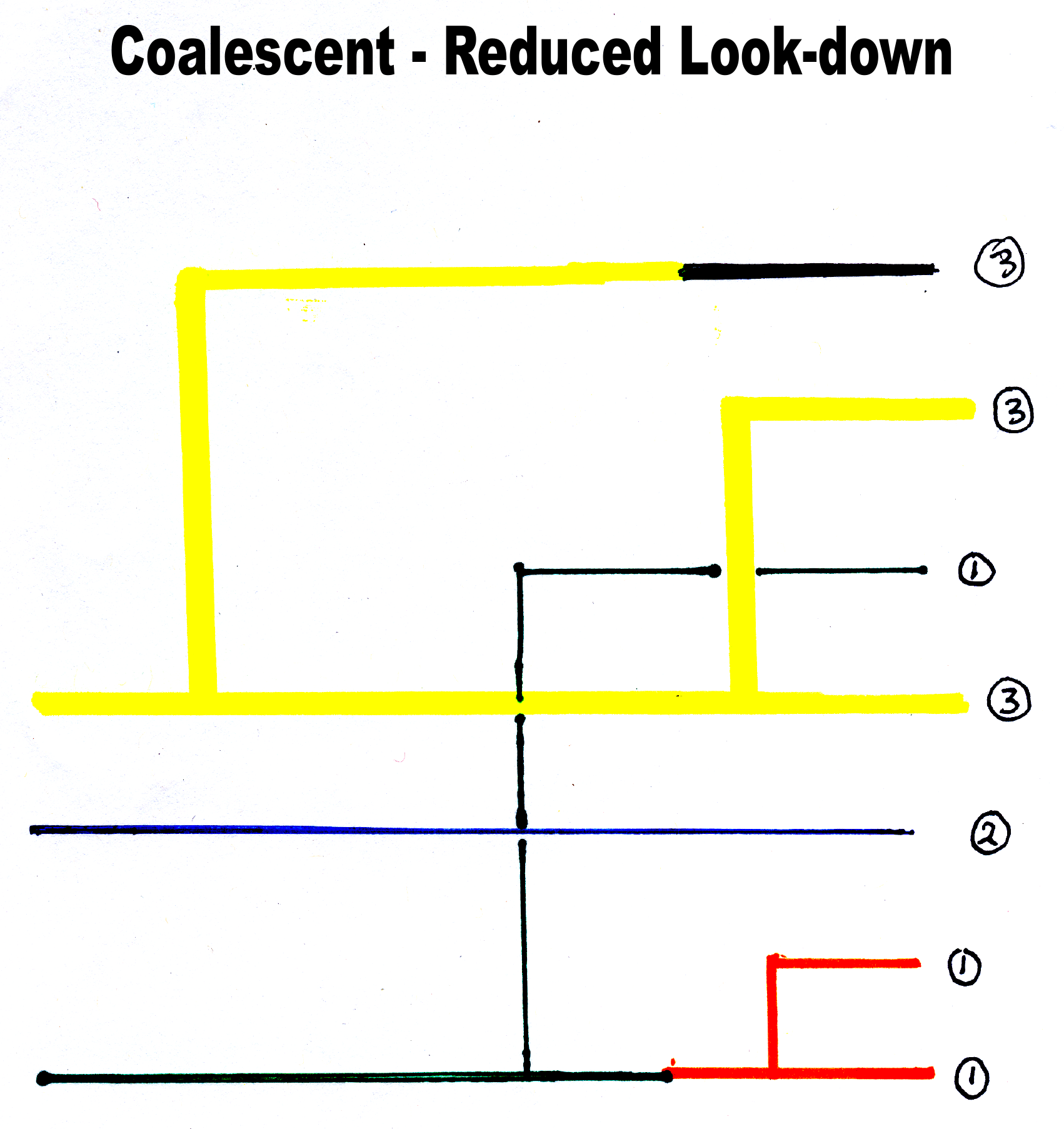}
\caption{Coalescent}
\end{center}
\end{figure}

\subsection{Poly\'{a} Urn Scheme}

\begin{theorem}
(Poly\'{a} Urn Representation) (a) For each $\ j,$ the limit%
\[
\tilde{M}_{j}(s)=\lim_{u\uparrow t}\frac{N_{t}(s,u,j)}{\sum_{k}N_{t}%
(s,u,k)},\;j=1,\dots,N_{t}(s)
\]
exists a.s. where $N_{t}(s,u,j)$ is the number of atoms at time $u$ with
common ancestor $j$ at time $s$.\newline (b) At each time $s,$ the random
vector $(\tilde{M}_{1}(s),\dots,\tilde{M}_{N_{t}(s)})$ is uniformly
distributed over the simplex $\Delta_{D(s)-1},$ that is, it is distributed as
the Dirichlet $D(1,\dots,1)$.
\end{theorem}

\begin{proof}
(a) For fixed $s>0$ we consider an urn model involving $N_{t}(s)$ types. At
each jump time, $u$, $N_{t}(.),$ $s<u<t,$ a particle of type $j$ is added with
probability%
\[
\tilde{M}_{j}(s)=\lim_{u\uparrow t}\frac{N_{t}(s,u,j)}{\sum_{k}N_{t}%
(s,u,k)},\;j=1,\dots,N_{t}(s)
\]
that is, it is an $N_{t}(s)$-type Poly\'{a} urn. \ It follows from the results
of Blackwell and Kendall (1964) \cite{BK-64} on the Martin boundary of the Poly\'{a} urn
process (and the fact that we have conditioned on $\{\tau=t\}$) that the
vector $(\tilde{M}_{j}(s),\dots,\tilde{M}_{N_{t}(s)})$ is uniformly
distributed over the simplex $\Delta_{N_{t}(s)-1}.$ Moreover \ the same result
implies that at the time of a split, the mass $\tilde{M}_{j}(u)$ is divided
into two equivalence classes of masses $U\tilde{M}_{j}(u)$ and $(1-U)\tilde
{M}_{j}(u)$ where $U$ is an independent $U[0,1]$ r.v. Recall that if an
$m$-type Polya urn is started with $n_{j}$ initial particles of type $j$, then
the joint distribution of the limiting proportions $(x_{1},\dots,x_{m})$ in
$\Delta_{m-1}$ \ is distributed via the Dirichlet $D(n_{1},\dots,n_{m})$.
Using this we will show that if a split occurs at time $s$ then%
\[
P(jth\text{ class splits%
$\vert$%
}\tilde{M}_{1}(\tau_{k}-),\dots\tilde{M}_{k-1}(\tau_{k}-))=\tilde{M}_{j}%
(\tau_{k}-),\;j=1,\dots,k-1.
\]
Taking advantage of the fact that merging types in a Polya urn yields a new
P\'{o}lya urn, it suffices to prove this for a two type urn in which case we
obtain%
\[
f_{m,n}(x,1-x)=\frac{\Gamma(m+n)}{\Gamma(m)\Gamma(n)}x^{m-1}(1-x)^{n-1}%
,\;x\in\lbrack0,1]
\]
Then
\begin{align*}
P(1st\text{ class splits%
$\vert$%
}\tilde{M}_{1}  &  =x,\tilde{M}_{2}=1-x)\\
&  =\frac{\frac{m}{n+m}f_{m+1,n}(x,1-x)}{\frac{m}{n+m}f_{m+1,n}%
(x,1-x)+\frac{n}{n+m}f_{m,n+1}(x,1-x)}=x.
\end{align*}
Finally it is clear from the above construction that the $\{\hat{M}_{j}%
(\tau_{k}-):j=1,\dots,k-1;k\in\mathbb{N\}}$ \ is independent of $\{\tilde
{D}^{t}(u):u>0\}$. This completes the proof that the P\'{o}lya urn scheme and
coalescent yield the same probabilistic mechanism, that is,%
\[
\mathit{Law}(\tilde{M}_{1}(s),\dots,\tilde{M}_{N_{t}(s)}:s<t)=\mathit{Law}%
(M_{1}(s),\dots,M_{N_{t}(s)}:s<t)\text{.}%
\]
(b) follows from the result of Blackwell and Kendall \cite{BK-64}.
\end{proof}

\begin{corollary}
(a) At the time \ of a split, one equivalence class of mass $M_{C}$ is split
into two equivalence classes of masses, $M_{C^{\prime}}=UM_{C}$ and
$M_{C^{\prime\prime}}=(1-U)M_{C}$ where $U$ is an independent uniform $(0,1)$
random variable. \newline (b) For each $C$, the probability that $C$ splits is
given by $M_{C}$.
\end{corollary}

\begin{proof}
Kingman (1982) \cite{K-82a}.
\end{proof}

\section{Dynamics of   population structure and history}
\subsection{The tree-valued Fleming-Viot process}

Recently, Greven-Pfaffelhuber-Winter (2008)  \cite{GPW-08} have identified the analogue of the Fleming-Viot limit of Moran processes in the enriched framework of Moran tree processes (without mutation).  This required some new concepts and techniques, in particular the notion of metric measure space and the Gromov-weak topology (see Appendix I, \ref{s.AGH}).

Let $\mathbb{M}$ denote the set of equivalence classes of metric measure space with the Gromov-weak topology. Let $\mathbb{M}_c$ denote the subset of compact metric measure spaces.

The resulting {\em Greven-Pfaffelhuber-Winter tree-valued Fleming-Viot process} has as state space space the set $\mathbb{U}$ of ultrametric measure spaces furnished with the Gromov-weak topology. The process is characterized by a martingale problem analogous to the Fleming-Viot martingale problem and again the proof of uniqueness is based on duality. We briefly outline some of the main ingredients.

We first define a class of test functions needed to define the generator.

A function $\Phi=\Phi^{n,\phi}:\mathbb{M}\to\mathbb{R}$ is called a polynomial of degree $n$ if $n$ is the minimal number such that there exists a function $\phi: \mathbb{R}_+^{\scriptsize{\left(\begin{array}{c} \N\\2\end{array}\right)}}\la \mathbb{R}$ but  and $\phi $ depends only on $(r_{i,j})_{1\leq i<j\leq n}$ such that
if $\chi= \overline{(X,r,\mu)}$,
\bea{}&& \Phi(\chi)=\langle \nu^{\chi},\phi\rangle :=\int_{\mathbb{R}_+^{\scriptsize{\left(\begin{array}{c} \N\\2\end{array}\right)}}}\nu^\chi(d\underline{\underline{r}})\phi(\underline{\underline{r}})\quad \text{where } \underline{\underline{r}}:=(r_{i,j})_{1\leq i<j}\\&&
=\int_{X^n} \phi(\{r(x_i,x_j)\}_{1\leq i<j\leq n})\mu(dx_1)\dots\mu(dx_n).\nonumber\eea

Let $\Pi^1:=\{\Phi^{n,\phi}:n\in\N,\phi\in C_b(\mathbb{R}_+^{\scriptsize{\left(\begin{array}{c} n\\2\end{array}\right)}})\}$

For $\Phi\in\Pi^1$,  The define
\be{}\Omega ^{\uparrow}\Phi = \Omega ^{\uparrow,\rm{grow}}\Phi+ \Omega ^{\uparrow,\rm {res}}\Phi\ee

\be{} \Omega ^{\uparrow,\rm{grow}}\Phi(\upsilon) :=\langle \nu^\upsilon,\rm{div}(\phi)\rangle,\ee
\be{} \rm{div}(\phi):=2 \sum_{1\leq i<j\leq n}\frak{\partial \phi}{\partial r_{i,j}}.\ee

\be{} \Omega ^{\uparrow,\rm{res}}\Phi (\upsilon):= \frac{\gamma}{2} \sum_{1\leq k,\ell\leq n} \left(\langle\nu^\upsilon,\phi\circ \theta_{k,\ell}\rangle -\langle \nu^\upsilon\rangle\right),\ee
\be{} \theta_{k,\ell}((r_{i,j})_{1\leq i<j}\rangle:= \left\{\begin{split} & r_{i,j} \qquad& i,j\ne 1\\
&r_{i\wedge k,i\vee k} & j=\ell\\ & r_{j\wedge k,j\vee k}& i=\ell.\end{split}\right.
\ee

\beT{} (\cite{GPW-08}, Theorem 1) Let $\mathbf{P}_0\in\mathcal{P}(\mathbb{U})$. The $(\mathbf{P}_0,\Omega^{\uparrow},\Pi^1)$-martingale problem has a unique solution $\mathcal{U}=(\mathcal{U}_t)_{t\geq 0}$. Moreover a.s.

(i) $\mathcal{U}$ has continuous sample paths in $\mathbb{U}$

(ii) For all $t>0$ $\mathcal{U}_t\in \mathbb{U}_c =\mathbb{U}\cap \mathbb{M}_c $

(iii) For all $t>0$, $\nu^{\mathcal{U}_t}((0,\infty)^{\scriptsize{\left(\begin{array}{c} \N\\2\end{array}\right)}})=1$ and
\be{} \{ t\in[0,\infty):\nu^{\mathcal{U}_t}((0,\infty)^{\scriptsize{\left(\begin{array}{c} \N\\2\end{array}\right)}})<1 \}\ee
has Lebesgue measure 0.
\end{theorem}

Finally, they establish the weak convergence of the tree-valued Moran processes to the tree-valued Fleming-Viot process.
\beT{} (Greven-Pfaffelhuber-Winter \cite{GPW-08}, Theorem 2) (Convergence of the tree valued  Moran to Fleming-Viot dynamics).  For $N\in\N$, let $\mathcal{U}^N$ be the tree-valued Moran dynamics with population size $N$, and let $\mathcal{U}=(\mathcal{U}_t)_{t\geq0}$ be the tree-valued Fleming-Viot dynamics.  If  $\mathcal{U}^N_0\Rightarrow \mathcal{U}_0$ as $N\to\infty$, weakly in the Gromov-weak toplogy (see Appendix ?), then
\be{} \mathcal{U}^N \Nto \mathcal{U}\ee weakly in the Skorohod topology on $\mathcal{D}_{\mathbb{U}}([0,\infty))$.
\end{theorem}

Uniqueness is again proved by duality where the dual is a tree-valued Kingman coalescent.
\subsection{The Historical Process}

We have incorporated mutation above in the explicit construction of the
look-down process in the case of a Feller mutation process.\ We can also allow
for more general mutation processes including time inhomogeneous processes.

An important special case is the \textit{historical process} (Dawson-Perkins (1991) \cite{DP-91}). In this case the
mutation process is a path-valued process $(t,Y^{t}):=(t,Y(\cdot\wedge
t))\in\mathbb{R}_{+}\times$ $D([0,\infty),E).$ The state space is
$\hat{E}=\{(t,Y^{t}):t\geq0,y\in D([0,\infty),E)\}$ with the subspace topology
from $\mathbb{R}_{+}\times D([0,\infty),E).$ If $y,w\in D([0,\infty),E)$ and
s$\geq0,$ let%
\[
(y/s/w)(t)=\left\{
\begin{array}
[c]{c}%
y(t)\;\;\;\;\;t<s\\
w(t-s)\;\;\;t\geq s
\end{array}
\right.
\]
Note that $\hat{E}$ is Polish since it is a closed subset of the Polish space
$\mathbb{R}_{+}\times$ $D([0,\infty),E).$

\begin{definition}
Let $W_{t}:D([0,\infty),\hat{E})\rightarrow\hat{E}$ denote the coordinate maps
and for $(s,y)\in\hat{E}$, define $\hat{P}_{s,y}$ on $D([0,\infty),\hat{E})$
with its Borel $\sigma$-algebra $\mathcal{\hat{D}}$, by%
\[
\hat{P}_{s,y}(W_{.}\in A)=P_{y(s)}((s+\cdot,y/s/Y^{.})\in A),
\]
i.e. under $\hat{P}_{s,\xi}$ we run $\xi$ up to some $s$ and then tag a copy
of $Y$ starting at $y(s).$
\end{definition}

\begin{lemma}
$(W,(\hat{P}_{s,y})_{(s,y)\in\hat{E}})$ is a strong Markov process with
semigroup%
\[
\hat{P}_{t}:C_{b}(\hat{E})\rightarrow C_{b}(\hat{E}).
\]
\end{lemma}

\begin{proof}
See Perkins.
\end{proof}

\begin{remark}
This Lemma remains true if we simply assume that the mutation process
satisfies: $x\rightarrow P_{x}$ from $E$ to $M_{1}(D_{E})$ is continuous. This
was proved above for a Feller process.
\end{remark}

If we assume condition%
\be{}
(s,y)\rightarrow\hat{P}_{s,y}\text{ from }\hat{E}\text{ to }M_{1}%
(D([0,\infty),\hat{E})\text{ is continuous}%
\ee
then we may carry out the the explicit construction of the look-down process
exactly as above. The resulting process is the \textit{historical look-down
process}.

\begin{remark}
The historical process can be defined even if the mutation process is
non-Markovian. For example we can consider the case in which the mutation
process is a fractional Brownian motion. \ In this case we modify the above
definition as follows:%
\[
\hat{P}_{s,y}(W_{.}\in A)=P((s+\cdot,y/s/Y^{.})\in A|Y(u)=y(u),u\leq s),
\]
\end{remark}

\section{Some Applications of the Coalescent and Look-Down Process}

\subsection{Ergodicity for the Fleming-Viot Process}

As a first application of the look-down process we establish ergodicity for
the Fleming-Viot process under the assumption that the mutation process is ergodic.

\begin{theorem}
Assume that $\theta >0$ and the mutation process has a unique stationary distribution
$\pi(dx)$ on $E$. Then the Fleming-Viot process, $X_{t}$ has a unique
stationary distribution on $M_{1}(E).$
\end{theorem}

\begin{proof}
Note that $\zeta_{1}(t)=U_{10}(\zeta_{1}(0),t)$ for all $t\geq0.$ Therefore if
the mutation process has a stationary distribution, then we can assume that
$\zeta_{1}(t)$ is defined for $-\infty<t<\infty$ and otherwise define the
process as above on $(-\infty,\infty)$. It is then easy to verify iteratively
that $\{\zeta_{1}(t),\dots,\zeta_{n}(t)\}$, $n=1,2,\dots$ are stationary. We
will show below that w.p.1 all particles have a common ancestor at a finite
time in the past. It is then easy to check that the resulting $E^{\infty}%
$-valued process will be stationary, as will the Fleming-Viot process.
\end{proof}

\begin{remark}
For the infinitely many alleles model it can be shown (cf. Ethier)  that%
\[
\Vert P_{\mu,t}(\cdot)-\Pi_{\theta}\Vert\leq1-P(D(t)=0)
\]
where $P_{\mu,t}$ is the law at time $t$ of the infinitely many alleles model
with parameter $\theta$ at time $t$ and $\Vert\cdot\Vert$ is the total
variation norm, and $D(t)$ is the pure death process starting at $\infty$ with
death rates%
\[
d_{k}=\frac{1}{2}k(k-1+\theta),\text{\ \ }k\geq0.
\]
\end{remark}

\subsection{Atomic Structure of the Infinitely Many Alleles Model}

Recall that in this case the mutation operator is bounded. \ Note if there is
no mutation, then at time $t$ the infinite collection of particle have only
finitely many ancestors at time zero and therefore only finitely many types
(i.e. a finite set of atoms). \ In the case of a bounded mutation rate, at
most countably many mutations occur (finitely many along each $\zeta
_{i}(s),0\leq s\leq t$). \ Therefore there are at most countably many types
(or countably many atoms). To show that there are actually countably many
types note that the mutations occur according to independent Poisson processes
for the countably many particles in the representation. \ Now consider the
sequence of time of jumps in the Kingman coalescent process $D_{n}.$ Note
that
\[
\int_{0}^{t}D(u)du=\infty\text{\ w.p.1}%
\]
Therefore the Poisson processes running along the genealogical tree ending at
time $t$ has infinitely many jumps. Hence there are infinitely many types
present at a fixed time $t$ w.p.1.\bigskip

\subsubsection{Schmuland's Theorem}

For an infinite-alleles model, with probability 1, there will be times at
which the Fleming-Viot measure consists of a single atom iff $\theta<1.$ If
$\theta\geq1$, then there will always be an infinite number of atoms.

\begin{proof} This was first proved by Schmuland \cite{Sch-91} using Dirichlet forms. (Here we give a proof due to Donnelly and Kurtz.)

Let $S_{1}(t)$ denote the size of the largest atom in $Z(t)$ and define
$\tau_{1}=\inf\{t:S_{1}(t)=1\}.$ Define recursively, $\alpha_{1}=\inf
\{t:S_{1}(t)\geq\frac{3}{4}\},\,\ \beta_{k}=\inf\{t>\alpha_{k}:S_{1}%
(t)\leq\frac{1}{2}\}$ and $\alpha_{k+1}=\inf\{t>\beta_{k}:S_{1}(t)\geq
\frac{3}{4}\}$. Fix a time interval $[\alpha_{k},\beta_{k})$ and define
$\tilde{S}(t)=S_{1}(\alpha_{k}+t).$ By the strong Markov property, we can let
$g$ be the indicator of the location of the largest atom at time $\alpha_{k}$.
\ Noting that this location does not change during the time interval
$[\alpha_{k},\beta_{k})$, therefore (by the martingale problem) $\tilde{S}$
satisfies%
\[
\tilde{S}(t)=S_{1}(\alpha_{k})+\int_{0}^{t}\sqrt{\tilde{S}(s)-\tilde{S}%
^{2}(s)}dW_{k}(s)-\int_{0}^{t}\frac{\theta}{2}\tilde{S}(s)ds,\;t<\beta
_{k}-\alpha_{k}%
\]
for some standard Brownian motion $W_{k}.$ This corresponds to a Wright-Fisher
diffusion with generator%
\[
G_{\theta}f(x)=\frac{1}{2}x(1-x)f^{\prime\prime}(x)-\frac{\theta}{2}%
xf^{\prime}(x)
\]
This process has ``speed measure'' given by
\begin{align*}
m(x)  & =\frac{2}{x(1-x)}e^{\int^{x}\frac{\theta y}{y(1-y)}dy}\\
& =\frac{2}{x(1-x)}\frac{1}{(1-x)^{\theta}}.
\end{align*}

Recall that a boundary point $1$ of a diffusion in natural scale is accessible
(cf. Breiman) \ iff%

\[
\int^{1}(1-x)m(x)dx <\infty.
\]
Noting that for the Wright-Fisher diffusion  the scale function, $s(x)\sim
const\cdot x$ near $1$, this  becomes%
\[
\int_{\frac{1}{2}}^{1}\frac{2}{x}\frac{1}{(1-x)^{\theta}}dx<\infty,
\]
that is, $\theta<1$ (note that the scale function, $s(x)\sim const\cdot x$
near $1$.  \ That is
\[
P(\tilde{S}(t)=1\text{ for some }s\in\lbrack\alpha_{k},\beta_{k}))>0
\]
iff $\theta<1.$  $S_{1}$ can reach $1$ iff $\theta<1$. If $\theta<1,$ a
renewal argument shows that $\tau_{1}<\infty,\,a.s.$. If $\theta\geq
1,\;\tau_{1}>\beta_{k}$ for every $k$. Since $(\beta_{k}-\alpha_{k})$ is iid
by the strong Markov property for $S_{1},$ we must have $\lim_{k\rightarrow
\infty}\beta_{k}=\infty$ and hence $\tau_{1}=\infty.$

The proof that there are infinitely many atoms if $\theta\geq1$ proceeds by
considering the sum of the sizes of the two largest atoms and showing that there
will always be at least three, etc.
\end{proof}

\begin{remark}
\begin{remark}
This was originally proved by B. Schmuland in the Fleming-Viot case using
Dirichlet forms and calculating capacities.
\end{remark}
\end{remark}

\subsection{The Infinitely Many Sites Model}

The infinitely-many-sites model has $E=[0,1]^{\mathbb{Z}_{+}}$ and mutation
kernel%
\[
P(\mathbf{x},.)=\int_{0}^{1}\delta_{(\xi,\mathbf{x)}}(\cdot)d\xi.
\]

Here we interpret $[0,1]\,$as the sites in a DNA string and $\mathbf{x}%
=(x_{1},x_{2},\dots,0,0,0,\dots)$ where $x_{1}$,$x_{2},\dots$ denotes the
sites at which mutations have occurred with $x_{1}$ denotes the latest
mutation, $x_{2}$ the site of the second most recent mutation, etc. Note that
we can identify this with the historical infinitely many types model by
reinterpreting the jump to $x\in\lbrack0,1]$ as the site of the most recent
mutation (rather than as a label for a new type). This model is used in the
analysis of large sets of DNA sequence data. \ In particular, given a finite
population the analysis begins by by identifying the sites at which at least
two members of the population differ this indicating a mutation has occurred
in one of their ancestors (after the most recent common ancestor at which all
members of the population had identical DNA sequences).

\beT{} (Ethier-Griffiths \cite{EG-87}, Theorem 2.3)  The infinitely many sites process $X(t)$ has a unique stationary distribution $P_{st}\in\mathcal{P}(\mathcal{P}(E))$ and is ergodic $X(t)\Rightarrow X(\infty)$.
\end{theorem}

We can now consider  a random sample of size $n$, that is, a point in $E^n$ chosen according to $(X(\infty))^{\otimes n}$ and the corresponding moment measures. A site $z\in [0,1]$ is said to be {\em segregating} with respect to the sample if it appears in at least one but not all of the $n$ sequences.

Given an ordered k-tuple $(\mathbf{x}_{1},\dots,\mathbf{x}_{k})\in E^{k}$ it forms
a {\em tree } if \begin{description}
                   \item[a] the coordinates of\textit{ }$\mathbf{x}_{1}$
are distinct
                   \item[b]   if $i,i^{\prime}\in\{1,\dots,k\},$ and $j,j^{\prime}%
\in\mathbb{Z}_{+}$, and $x_{ij}=x_{i^{\prime},j^{\prime}}$, then $x_{i,j+\ell
}=x_{i^{\prime}j^{\prime}+\ell}$ for all $\ell\geq1$,
                   \item[c] there exists
$j_{1},\dots,j_{d}\in \mathbb{Z}_{+}$ such that $x_{1,j_{1}}=\dots x_{kj_{k}}$, that is, they have a common ancestor
                 \end{description}

\begin{figure}[h]
\begin{center}
\includegraphics[scale=0.75]{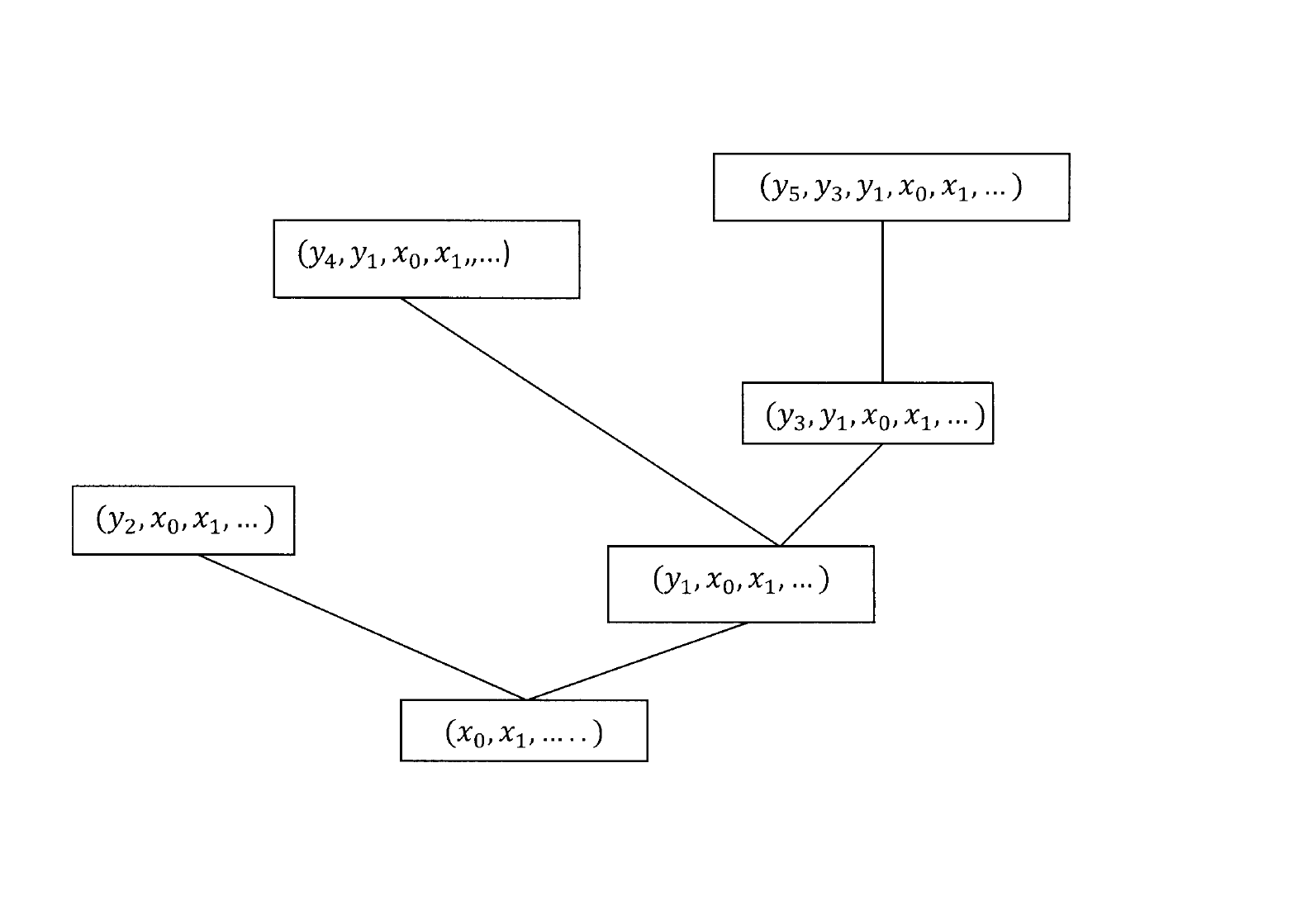}
\caption{Infinitely Many Sites - Segregating sites}
\end{center}
\end{figure}

Let \be{}\mathcal{P}^0_a(E):= \{\mu\in\mathcal{P}_a(E):\mu^n(\mathcal{T}_n)=1\;\forall\;n\in\N\}.\ee
Then
\be{}  P(X(t)\in \mathcal{P}^0_a(E)\;\forall\; t>0)=1,\quad\text{and  } P_{st}(\mathcal{P}^0_a)=1.\ee

We can classify the tree structures into equivalence classes where two trees
are equivalent if they are equal after a relabeling of $[0,1]$.

For $i,j\in\N$ let $T_{ij}$ be the equivalence class of trees of the form
\be{} ((\mathbf{x}_0,\dots, \mathbf{x}_{i-1},\mathbf{z}_0,\mathbf{z}_1,\dots),\;((\mathbf{y}_0,\dots, \mathbf{y}_{j-1},\mathbf{z}_0,\mathbf{z}_1,\dots))\ee
where the $\mathbf{x}_{0},\dots,\mathbf{x}_{i-1},\mathbf{y}_{0},\dots,\mathbf{y}_{j-1}, \mathbf{z}_0,\dots$ are distinct.

Watterson (1975) obtained the distribution of the number of segregating sites as follows.
Let
\be{} p_{i,j}=\int \mu^2 (T_{i,j})P_{st}(d\mu).\ee
Then it can be shown that
\be{} p_{i,j}=\left(\begin{array}{ll}i+j\\\;\; i \end{array}\right)\left(\frac{\theta}{2(1+\theta)}\right)^{i+j}\frac{1}{1+\theta}.\ee
Let $S_k$ denote the number of segregating sites in a random sample of size $k$. Then the distribution of $S_2$ is geometric,
\be{} P(S_2=n)=\sum_{i=0}^{n} \left(\begin{array}{ll}n\\i \end{array}\right)\left(\frac{\theta}{2(1+\theta)}\right)^{n}\frac{1}{1+\theta}
=\left(\frac{\theta}{1+\theta}\right)^n\frac{1}{1+\theta}.\ee

Tavar\'e (1984) \cite{T-84} proved that
\be{} P(S_n=s)=\frac{n-1}{\theta}\sum_{j=1}^{n-1}(-1)^{j-1}\left(\begin{array}{ll}n-2\\j-1 \end{array}\right)\left(\frac{\theta}{j+\theta}\right)^{s+1}.\ee

Given $k$
sequences $\mathbf{x}_{1},\dots,\mathbf{x}_{k}$ in $E=[0,1]^{\mathbb{Z}_{+}}$
and a vector of multiplicities $\mathbf{n=(n}_{1},\dots,\mathbf{n}_{k}%
)$, consider the equivalence class of trees containing \be{}(\mathbf{x}_{1}%
,\dots\mathbf{x}_{1},\mathbf{x}_{2},\dots,\mathbf{x}_{2},\dots,\mathbf{x}%
_{k},\dots,\mathbf{x}_{k}).\ee Then an analogue of the Ewens sampling formula
is to determine \ the distribution of a finite sample taken from the
stationary distribution of the infinitely many sites model. In particular
given an ordered random sample of size $n,$ what is the probability that the
sample has tree structure $T$ with multiplicities $(n_{1},\dots,n_{k}).$
\ These questions have been studied by Griffiths (1982) \cite{G-82}),  Wakeley (1998) \cite{W-98}, Ethier and Griffiths \cite{EG-87}, and Griffiiths
and Tavar\'{e} \cite{GT-95}.

\subsection{Wandering Distributions}

In 1973 Ohta and Kimura introduced the \textit{stepwise mutation model} to
describe electrophoretically detectable alleles in a population. \ In this
case the different alleles are represented by points on $\mathbb{Z}^{1}$ and
the mutation semigroup is that of simple random walk on $\mathbb{Z}^{1}$. \ In
simulations they discovered the tendency for the entire population to be
somewhat spread out but essentially to wander around as a loose clump.   Moran investigated this in \cite{Mo-75}.\ The
explanation for this was given by Kingman using the coalescent \cite{K-82a}, (see \cite{K-00} for an interesting history of this and its role in the origins of the coalescent).

\begin{figure}[h]
\begin{center}
\includegraphics[scale=0.4]{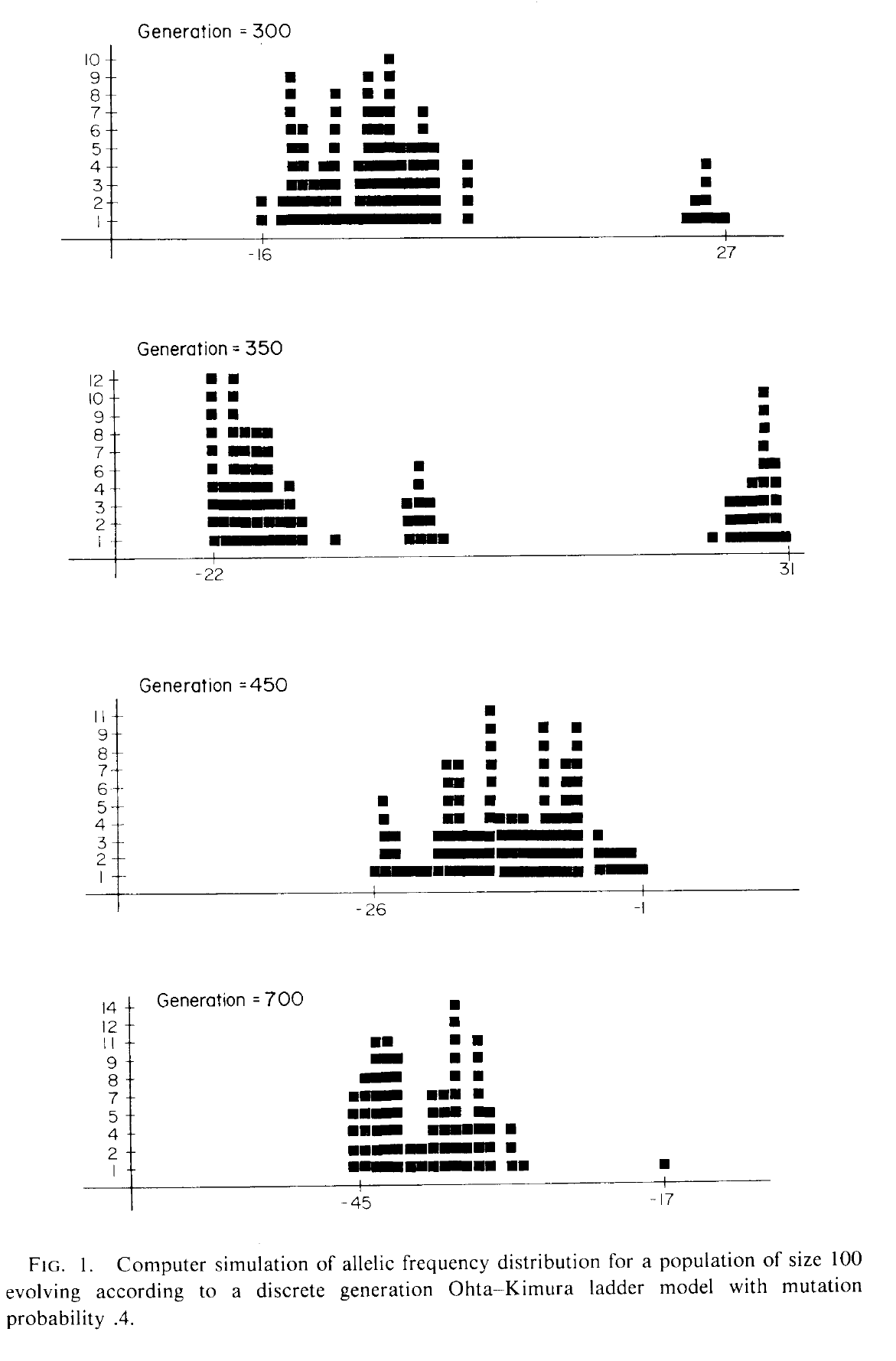}
\caption{Wandering distribution}
\end{center}
\end{figure}

 The analogous
phenomenon with continuous types was given by Dawson
and Hochberg (1982) \cite{DH-82} using the infinite particle representation. \ In fact consider
the particle $\zeta_{1}$ in the lookdown process. \ Then it follows a simple
random walk on $\mathbb{Z}^{1}$. \ Morover the entire population branches off
from $\zeta_{1}$ at a finite time in the past. \ It can be shown that the
relative to the particle $\zeta_{1}$, the evolving population cloud approaches
an equilibrium thus describing what was called the ``wandering distribution''.

\subsection{Support Properties.}

Consider the Fleming-Viot process in which the mutation semigroup is the
standard Brownian motion semigroup in $\mathbb{R}^{d}$. This model, in the
case of $\mathbb{R}^{1},$ arises as the limit of the stepwise mutation model
on $\varepsilon\mathbb{Z}^{1}$ when $\varepsilon\rightarrow0$. \ In the
general case it serves as a model of a population described by $d$ continuous characteristics.

\beT{Compact}
Consider the Fleming-Viot process with Brownian motion in $\mathbb{R}^{d}$ as
the mutation process. \ Then at a fixed time (a) the closed support of $X_{t}$
is compact with probability one\newline (b) if $d>2,$ then the measure is
supported on a set of Lebesgue measure zero.
\end{theorem}

\begin{proof}
We begin by determining the distribution of the times $\tau_{\infty,n}$ in the
look down process. Let $\tau_{\infty,n}:=\inf\{s:|\Gamma(s,t)|\geq n\}$ and
$B_{n}:=n(t-\tau_{\infty,n}).\,\;$Then%
\[
B_{n}:=n\sum_{k=n+1}^{\infty}\frac{2E_{k}}{k(k+1)}.
\]
where $\{E_{k}\}$ are i.i.d. exponential one r.v.'s. \ Using standard
exponential estimates (see Dawson and Vinogradov (1992)) it can be verified
that \ for $\ 0<\nu<\frac{1}{2},$%
\[
P(|B_{n}-2|>n^{-\nu})\leq2e^{-3n^{1-2\nu}/4}+c/n^{1+\nu}%
\]
Therefore, by Borel-Cantelli there exists a random $n(\nu)$ such that%
\begin{align*}
|B_{n}-2| &  \leq n^{-\nu}\;\;\forall\;\;n\geq n(\nu),\;\;P\text{-a.s.\ and}\\
B_{n} &  \leq3\;\;\forall\;\;n\geq n(\nu),\;\;P\text{-a.s.}%
\end{align*}
and%
\[
(t-\tau_{\infty,n})\leq\frac{3}{n}\text{ if }n\geq n(\nu),\;P-a.s.
\]
We then observe that for the Brownian mutation process, in the time interval
$(\tau_{\infty,n},t]$ a particle has a displacement that is normally
distributed with variance $B_{n}/n$ (and different particles have independent
displacements)$.$ Let $A_{n+1}$ be the event that some particle with index in
$\Gamma(\tau_{\infty,2^{n+1}},t)$ is a distance greater than $\varepsilon_{n}$
from its ancestor at time $\tau_{2^{n}}.$ Then $P(A_{n+1})\leq c2^{n+1}%
e^{-\frac{2^{n}.\varepsilon_{n}^{2}}{3}}$ for some constant $c$. Taking
$\varepsilon_{n}=2^{-n/(2+\eta)}$ with $0<\eta<\frac{1}{2},$ we get
\begin{align*}
&  \sum_{n=1}^{\infty}P(A_{n+1})\\
&  \leq\sum_{n=1}^{\infty}c2^{n+1}e^{-\frac{2^{\frac{n\eta}{2+\eta}}}{3}}\\
&  <\infty.
\end{align*}
Then by the Borel-Cantelli only finitely many of the $A_{n}$ occur.

Let $\delta_{n}=\sum_{k\geq n}\varepsilon_{k}=\frac{2^{-n/(2+\eta)}%
}{1-2^{-1/(2+\eta)}}.$

(a) Then for all $n$ sufficiently large,%
\[
\text{supp}(X_{t})\subset\cup_{j\in\Gamma(\tau_{\infty,2^n},t)}B(\delta
_{n},\zeta_{j}(\tau_{\infty,2^n})).
\]
where $B(\delta,\zeta)$ denotes a ball of radius $\delta$ centered at $\zeta$.
This implies that supp$(X_{t})$ is contained in bounded subset of
$\mathbb{R}^{d}$ and is therefore compact.

(b) This implies that supp$(X_{t})$ is contained in a set of Lebesgue measure
less that or equal to%
\[
c2^{n+1}\left(  \frac{2^{-n/(2+\eta)}}{1-2^{-1/(2+\eta)}}\right)
^{d}\rightarrow0\text{ as }n\rightarrow\infty
\]
if $d>2.$
\end{proof}

\section{Generalizations of the Kingman coalescent}

\subsection{Coalescent with time varying population size}

Since constant population size is unrealistic it is important to determine to what extent the results obtained from coalescent
 theory are robust, that is, what happens to the ancestral structure if the population size is randomly-varying in time.
  Under the assumption that the individuals in the population are exchangeable and the rescaled backward population size process
   converges to a continuous time Markov chain Kaj and Krone \cite{KK-03}  show that the ancestral process is a stochastic time change
   of the Kingman coalescent.

\subsection{Generalized coalescent}

In the Kingman coalescent the jumps correspond to the coalescence of exactly two clusters.
In (1999) Pitman \cite{P-99}  and Sagitov \cite{S-99} introduced coalescents with multiple collisions  called $\Lambda$-coalescent
in which many clusters can merge simultaneously  into a single cluster. The relation between these Lambda-coalescents and a class of discontinuous Fleming-Viot processes was established in the 2005 paper of Birkner, Blath, Capaldo, Etheridge, M\"ohle, Schweinsberg and Wakolbinger
\cite{BBCEM-05}.  Here we just briefly describe these objects.

A generalized coalescent process is a Markov process $\{\Pi(t)\}_{0\leq t\leq T}$ with state space given by the space of partitions of $\mathbb{N}$ and such that the law is invariant under permutations of $\mathbb{N}$ - see \cite{Ber}.

Consider a Markov  process with state space $\mathcal{P}([0,1])$ and with generator
\be{} GF(\mu)= \int_{(0,1]} y^{-2}\Lambda(dy) \int\mu(da)(F((1-y)\mu+y\delta_a)-F(\mu)\ee
where $\Lambda$ is a finite measure on $[0,1]$ with $\Lambda(\{0\})=0$.

Consider functions of the form $F(\mu)=\int \dots\int f(a_1,\dots,a_p)\mu(da_1)\dots \mu(da_p)$.  Then the generator has the form
\[{G}F(\mu)=\sum_{J\subset\{1,\dots,p\},|J|\geq2}\beta^\Lambda_{p,J}\int\mu(da_1)\dots\mu(da_p)
(f(a_1^J,\dots,a^J_p)-f(a_1,\dots,a_p))\]
\[\beta^\Lambda_{p,J}=\int_{[0,1]}y^{j-2}(1-y)^{p-j}\Lambda(dy)\]
where $a_1^J,\dots,a^J_p$ denotes the coalescence of the $a_j\in J$. This process is called a {\em generalized Fleming-Viot process}.

To set up the connection between the Fleming-Viot process and the coalescent process, fix a time $T$ and pick a sequence of individuals labelled $1,2,3,\dots$ independently and uniformly on $[0,1]$. Then for $t\leq T$ let  $\Pi(t)$ denote collecting together individuals having the same ancestor at time $T-t$.  This results in The $\Lambda$-coalescent process.  Then by Kingman's theory of exchangeable partitions, for every $t\geq 0$ each block of $\Pi(t)$ has an asymptotic frequency and the ranked sequence of these frequencies yields a Markov process called the mass-coalescent.  The has the same distribution as the ranked sequence of jump sizes of the Fleming-Viot process.

 Pitman \cite{P-99} proved that the $\Lambda$-coalescent has proper ordered frequencies, that is, \be{pof}\sum_i f(\pi_i)=1 \text{  if and only if  }\int_0^1 \Lambda(dx) x^{-1}=\infty.\ee

Consider the case where  $\nu([\ve,1])$ is regularly varying with index $-\gamma$ as $\ve\to 0$.
If $1<\gamma<2$, then the coalescent {\em comes down from infinity}, that is, $\Pi(t)$ with $t>0$ has finitely many blocks.

\begin{example}
The classical Fleming-Viot corresponds to $\Lambda=\delta_0$.
\end{example}
\begin{example}
Another special case, the Bolthausen-Sznitman coalescent
is given by $\Lambda=U([0,1])$. This has interesting connections to the random energy model of Derrida which arises in statistical physics and Neveu's CSBP.
 In contrast to the Kingman coalescent, the Bolthausen-Sznitman coalescent  {\em does not come down from infinity} but has infinitely many clusters at all times.
\end{example}
\begin{example}

If $\Lambda= \text{ Beta}(2-\alpha,\alpha)$ for $0<\alpha <2$,  the Fleming-Viot
process $Y_t$ is associated to a time-changed {\em $\alpha$-stable
continuous state measure-valued branching process} $X_t$ as follows:
\be{} Y_t=  \frac{X_{\tau(t)}}{X_{\tau(t)}([0,1])}.\ee

\end{example}

In the case (\ref{pof}), Bertoin and LeGall \cite{BL-03}, \cite{BL-06} also obtain the generalized Fleming-Viot process as the solution for $x\in [0,1]$ of the stochastic equation
\be{}  Y_t(x)=x+\int_{[0,1]\times(0,1)\times(0,1]}N(ds,du,dr)r(1_{\{u\leq Y_{s-}(x)\}}-Y_{s-}(x))\ee
where $N$ is an $\mathcal{F}_t$ Poisson point process with intensity $dt\otimes du\otimes \nu(dr)$ with $\nu(dx)=\frac{\Lambda(dx)}{x^2}$. This has been generalized by Dawson and Li \cite{DLi} to include the case in which $\Lambda(\{0\})>0$.

\chapter[Spatially Structured Models]{Spatially Structured and Measure-valued Models}

\section{Introduction}

In this chapter we consider spatially structured population systems in the context of measure-valued processes on a discrete (countable) space $S_1$ or more generally on a Polish space $S$.  In addition an objective is consider the scaling limit of particle systems so that we also consider the situation of parametric families $\{S_\ve\}_{\ve\in(0,1]}\text{  or  }\{S_N\}_{N\in\N}\subset S$.
We denote by $\mathcal{N}_f(S_1)$, $\mathcal{M}_f(S)$ the space of finite counting measures, respectively, finite  Borel measures on $S$.

\subsection{Spatial dynamics}

The spatial systems we consider involve two basic types of dynamical mechanism:
\begin{itemize}
\item Migration between sites usually described by a random walk in the discrete case or Markov process in the general case,
\item Local interactions at a site such as reproduction, competition, etc.

\end{itemize}

In the  discrete case we  consider a finite or
countable set of sites, $S_1$,  together with an irreducible
continuous time Markov chain on $S_1$. The Markov chain describes
the migration of individuals between sites.

\subsection{Random Walks}

We often  consider the special case in which the migration Markov chain is a random walk on a countable (additive) abelian group
 $S_1$  with transition kernel $p(x,y)=p(x-y)$
where $p(\cdot)$ is a  symmetric probability kernel on $S_1$ with $p(0)=0$.
The corresponding continuous time
random walk with transition rate $\gamma$ is usually denoted by $q^\gamma_{x,y}=\gamma p(x-y)$.

\begin{example} $S_1=\mathbb{Z}^d\subset \mathbb{R}^d =S$.

Let   $p(\cdot)$ is a finite range kernel on $\mathbb{Z}^d$ which satisfies
\be{rw0}   \sum_{x\in\mathbb{Z}^d}x^ip(x)=0, \text{ for } i=1,\dots,d.\ee
 \be{rw1}  \sum_{x\in\mathbb{Z}^d}x^ix^jp(x)=\delta_{i,j}\sigma^2,\quad i,j=1,\dots, d, \ee
 where $\delta_{i,j}$ is Kronecker's delta.

 \beP{}
 (a) The random walk is transient if and only if $d\geq 3$.

(b)  Consider the rescaled lattice $S_{\ve}= \sqrt{\ve}\, \mathbb{Z}^d$
 and random walk with kernel $p_\ve(x)=p(x/\sqrt{\ve})$ and jump rate $\frac{\gamma}{\ve}$.  Then the {\em scaling limit} of the random walk is  $\{ \sqrt{\gamma}\,\sigma B_t\}_{t\geq 0}$ where $B_t$ is a standard Brownian motion on $\mathbb{R}^d$.
\end{proposition}

\begin{remark} We can also consider infinite range random walk kernels for which the scaling limit is a $\alpha$-symmetric stable process, $0<\alpha<2$.

\end{remark}

\end{example}
\strut
\begin{example} Island Model.  For the finite group $S_N=\{0,\dots,N-1\}$ (with addition modulo $N$) we consider the random walk
\be{}
p(x)=\frac{1}{N-1}\quad \text{ if } x\ne0,\quad p(0)=0
\ee

\end{example}
\strut

\begin{example}For $N\in\N$ we denote by
 $S_N$  {  the  Hierarchical Lattice
 }  $\Omega^0_N$ and the corresponding rescaled lattices $\Omega_N^j$ defined by:
{
\begin{align*}
& \Omega_{N}^{j}=
\{(\xi_{\ell})_{\ell\in\mathbb{Z},\ell\geq-j}:\xi_{\ell}\in\{0,1,\dots
,N-1\},\\&\qquad \qquad \qquad\qquad\qquad \exists\ell_{0},
\xi_{\ell}=0\ \forall\ \ell\geq\ell_{0}.\}
\end{align*}%
The group operation is componentwise addition modulo $N$.

The hierarchical distance between two points $\xi$ and $\eta$ is
\[
|\xi-\eta|:=\min\{k\in\mathbb{Z}:\mathbb{\xi}_{m}=\eta_{m}%
\ \ \forall\ \ m\geq k\}.
\]%
We also introduce the metric
\begin{align*}
d_{j}(\xi,\eta)  &  =N^{|\xi-\eta|+1}\text{ if }-j<|\xi-\eta
|.
\end{align*}
Then  $(\Omega_{N}^{\infty},d_\infty)\ $is a totally
disconnected, locally compact abelian group. (See Evans \cite{E-89}).
}

\strut

Let $c,d,\alpha \in (0,\infty)$.   The rescaled  $(c,d,\alpha)$-random walk on $\Omega_{N}^{j}$ has jump rates  %
\be{hrw}
q^{(j)}_{\xi,\eta}=\sum_{k=-j+1}^{\infty}\frac{c^{k-1}N^{(k-1)(1-\frac{\alpha}{d})}}{N^{k-1}}\frac{1}{N^{j-|k|}}1_{B_{k}%
^{j}(\xi)}(\eta),%
\ee
that is, a  jump to a point in ball $B^j_k=\{\zeta\in \Omega^j_N: |\zeta-\xi|\leq k\}$   is taken at rate $\frac{c^{k-1}}{N^{k-1}}%
,\ k\geq-j+1 $ and the point to which it jumps is  chosen at random in the ball
$B_{k}^{j}(\xi)$ .
\end{example}

\beP{} (a) The $(c,d,2)$ random walk is transient if and only if $d>2$ or $d=2,\; c>1$. (Sawyer-Felsenstein \cite{SF-83})

(b)
The {scaling limit of the rescaled $(c,d)$-random walks on $\Omega_N^j$   as $ j\to\infty$ } { is the
Evans-L\'{e}vy process on $\Omega_{N}^{\infty}$ with parameters
$(c,d)$.} (See Evans \cite{E-89}.)

(c) (Evans-Fleischmann \cite{EF-96}, Prop. 13) The L\'evy process has a jointly continuous local time if $c<1$ (this corresponds to dimension $d=2-$.
\end{proposition}
\strut

The $(c,d,\alpha)$-random walks also mimic $\alpha$-stable processes in terms of the potential operators (Dawson-Gorostiza-Waklobinger) \cite{DGW-01}, \cite{DGW-04}. The potential operator of this hierarchical random walk has a kernel of the form
$$G_{N}(x)={\rm const.} N^{-|x|(1-1/(\gamma+1))},$$
where $\gamma=\frac{d}{\alpha}-1$   (hence $d>\alpha)$, this can be written as
\be{caricature}G_{N,\gamma}(x)={\rm const.} \rho(x)^{-(d-\alpha )}\ee
where
$$\rho(x)=N^{|x|/d}.$$%
$\rho(x)$ is the ``Euclidean radial distance'' of $x$ from $0$, so that
 the volume of a ball of radius $\rho$ grows like $\rho^d$. Therefore the potential operator of the $(c,d, \alpha)$-random walk and that  for the $\alpha$-stable process in $\mathbb{R}^d$  have the same asymptotic decay.

 In the next three sections we describe  some of the basic spatial systems including interacting particle systems, interacting diffusions and measure-valued processes.

\section{Branching random walk and branching Brownian motions}

\bigskip

Branching random walks and branching diffusions have a long history.  A general  theory of branching Markov processes was developed
in a series of three papers by Ikeda, Nagasawa and Watanabe in 1968, 1969 \cite{INW-68}.  The application of  branching random fields to genetics was introduced by Sawyer (1975) \cite{S-75}.

We consider a
branching random walk (BRW).
   The dynamics are given by:

\begin{itemize}
\item Birth and death at rate ${\gamma}$:
\begin{eqnarray*}
&\delta _x \rightarrow  (k\text{  particles}) \;\delta _x+ {\cdots}+ \delta _x \text{  w.p. }p_k, \quad   \delta _x
\rightarrow \emptyset\quad \text{w.p. }p_0,\\&
\mathcal{G}(z)=\sum_{k=0}^\infty z^kp_k\quad \text{offspring distribution generating function}.
\end{eqnarray*}
\item Spatial random walk in $S_1$ with kernel  $p(\cdot)$
\begin{eqnarray*}
\delta _x \rightarrow &\delta_{y } \text{  with  rate  } p(y-x)
\end{eqnarray*}
\end{itemize}
The BRW is critical, subcritical, supercritical depending on $m=\sum_k kp_k$ $=1, <1,>1$, respectively.

We can write the generator of the branching rate walk as follows:
$\mathcal{D}=\{F:F(\mu)=f(\mu(\phi))=f(\langle \phi,\mu\rangle),\; \phi\in\mathcal{B}_b(S_1),\;f\in C(\mathbb{R})$
and for $F\in \mathcal{D}$,
\bean{} GF(\mu)=&&\sum_x\mu(x) \sum_y p(y)[F(\mu+\delta_{x+y}-\delta_x)-F(\mu)]\\&&+\kappa\sum_x\mu(x)\sum_{k=0}^\infty p_k[F(\mu+(k-1)\delta_x)-F(\mu)]\\&& =
\sum_x\mu(x)\sum_y p(y)[f(\mu(\phi)+\phi(x+y)-\phi(x))-f(\mu(\phi))]\\&&
\kappa\sum_x\mu(x) \sum_{k=0}^\infty p_k[f(\mu(\phi)+(k-1)\phi(x))-f(\mu(\phi))]
\eean

Let $\{S_t:t\geq 0\}$ denote the semigroup acting on $\mathcal{B}_b(S_1)$ associated to the random walk.
Now define the Laplace functional
\be{} u(t,x)=P_{\delta_x}(e^{-X_t(\phi)}).\ee

Then conditioning at the first birth-death event we obtain
\be{BRWLF} u(t,x)=(S_t e^{-\phi})(x)e^{-\kappa t}+\kappa\int_0^te^{-\kappa s}(S_s \mathcal{G}(u(t-s,\cdot)))(x)ds.\ee
Note that this is also valid if we replace the random walk by a L\'evy process on a locally compact abelian group (for example Brownian motion on $\mathbb{R}^d$  with semigroup $\{S_t:t\geq t\}$).

\beP{P.BRWLF} The martingale problem for $G$ is well posed and the Laplace functional of the solution is the unique solution of equation (\ref{BRWLF}).
\end{proposition}

The system of branching Brownian motions (BBM) is defined in the same way with $S=\mathbb{R}^d$ with offspring produced at the location of the parent and between branching the particles perform independent Brownian motions. (For non-local branching see Z. Li \cite{Li-06}.)

\begin{remark}  We sometimes  combine the reproduction and spatial jump by replacing the reproduction and migration by a single mechanism in which an  offspring produced by a birth immediately moves to a new location obtained by taking a jump with kernel $p_\ve$, that is, $\delta_x \to \delta_x+\delta_y$.
\end{remark}

We also consider the   $ \mathcal{N}(S)$-measure-valued process $\{X_t\}$ in which
each particle has mass $\eta$, that is,
\[
X_{t}(A) = \eta\sum_{i=1}^{N(t)}\delta_{x_i(t)},\quad A\subset S
\]
where $x_i(t)$ denotes the location of the ith particles at time $t$.

\subsubsection{Supercritical BRW and BBM}

There is an important relation between supercritical branching Brownian motions and the  Fisher-KPP equation.  This relation was developed by
McKean \cite{M-75} and  Bramson \cite{B-83}.

A basic question concerns the geometrical properties of the supercritical branching random walk.  Biggins \cite{B-78} has proved that
the set $\mathcal{I}^{(n)}$ of positions occupied by nth generation individuals rescaled by a factor $\frac{1}{n}$ has asymptotic shape $\mathcal{I}$ where $\mathcal{I}$ is a convex set.

\section{Interacting particle systems}

\subsection{Coalescing random walks}

Let $S_1$ be a countable abelian group and $p(x)$ a symmetric
finite range kernel on $S_1$, $p(0)=0$.  The state space is $\{0,1\}^{S_1}$.

Consider a collection of particles on $S_1$ undergoing  random walks
which are independent  up to a collision time.  If two particles collide (occupy the same site), then the two particles instantaneously coalesce and are replaced by a single particle.

Now consider the associated a set-valued process $A(t)$ corresponding to the set
of points occupied by  a coalescing random walk. The state space is the collection of finite subsets of $S_1$, $\{\eta\in\{0,1\}^S,\;\eta(x)=0,\;a.a.\; x\}$. Then the transitions are
\be{} A\to A\cup\{y\}\backslash \{x\} \text{ at rate  }q_{x,y}=p(x-y). \ee

Now consider the family of functions on
 \bea{tfvm}&& F(A,\eta)=\prod_{x\in A}\eta(x),\\&& \;A =\text{ finite subset of  }
 S_1,\quad \eta\in \{0,1\}^{S_1}\nonumber\eea

We denote the generator by $H$.
We observe that the generator $H$ applied to functions of the form  $F$ (defined in (\ref{tfvm})) satisfies

\bea{}&& HF(\{x_1,\dots,x_n\},\eta) \\&&=\sum_{i=1}^n\sum_y q_{x_i,y}
(F(\{x_1,\dots,x_n\}\backslash x_i\cup
y,\eta)-F(\{x_1,\dots,x_n\},\eta))\nonumber\eea

\subsubsection{Coalescing random walks with delay}
A variation of the system of coalescing walks is the coalescing random walk with delay.  In this case when a pair of particles occupy the same site, then they coalesce at a random time with exponential distribution, that is, the particles coalesce with finite rate $\kappa$.  These exponential random variable at different sites are independent.

An important quantity is the probability that two particles starting at the same site eventually coalesce.  To determine this let $q_0$ be the rate at which a jump from $0$ occurs,  $\tau_0=\inf\{t>0:Z_t\ne 0\}$ and  $\tau_1:=\inf\{t>\tau_0:Z_t=0\}$.  Then $E[\tau_0]=\frac{1}{q_0}$.  Let $q_e=1-P(\tau_1<\infty)$ (escape probability).

Starting two particles at $0$ we compute the probability that they do not eventually coalesce, $\gamma_e$,
as
\be{} \gamma_e= \frac{2q_0}{2q_0+\kappa}\cdot [q_e+(1-q_e)\gamma_e]\ee
so that
\be{ncp}  \gamma_e= \frac{2q_0q_e}{\kappa+2q_0q_e}.\ee

\subsection{Voter Model}

The voter model was independently introduced by
Clifford, Sudbury (1973) \cite{CS-73}  and  Holley and Liggett (1975) \cite{HL-75}.

Again, let $S_1$ be a countable abelian group and $p(x)$ a symmetric
finite range kernel on $S_1$.  The state space is $\{0,1\}^{S_1}$. The voter model $\xi_t$ is defined by the transitions
\begin{align*}
&\eta \to \eta_{x,y},\;x\ne y,  \qquad \text{ at rate } q_{x,y}\text{  where } \\&
\eta_{x,y}(z)=\eta(z)\; \text{  if  }z\ne x\\& \eta_{x,y}(x)=\eta(y).
\end{align*}

Acting on functions of finite support, the generator can be written

\be{}  Gf(\eta)=\sum_x c(x,\eta)[f(\eta_x)-f(\eta)]\ee where

\bea{}
\quad\;\;\eta_x(y)=\left\{%
\begin{array}{ll}
    \eta(y), & \hbox{ if } x\ne y \\
    1-\eta(y), & \hbox{ if } x=y\\
\end{array}%
\right.\eea and \be{}  c(x,\eta)=\left\{%
\begin{array}{ll}
    \sum_y q_{x,y}\eta(y), & \hbox{ if }\eta(x)=0 \\
    \sum_y q_{x,y}[1-\eta(y)] & \hbox{ if }\eta(x)=1 \\
\end{array}%
\right.\ee

Let  $f_x(\eta):= 1(\eta(x)=1)$.  Now consider the functions
 \bea{tfvm2}&& F(A,\eta)=\prod_{x\in A}\eta(x),\\&& \;A \text{ is a finite subset of  }
 S,\quad \eta\in \{0,1\}^S\nonumber\eea
Noting that $(f(\eta_x)-f(\eta))= 1-2f_x(\eta)$ and a change
occurs at site $x$ with rate $\sum_y
q_{y,x}[f_x(\eta)(1-f_y(\eta))+f_y(\eta)(1-f_x(\eta))]=\sum_yp(y,x)[f_x(\eta)+f_y(\eta)-2f_x(\eta)f_y(\eta)]$
Then \bea{} && GF(\{x_1,\dots,x_n\},\eta)=
G\left(\prod_{i=1}^nf_{x_i}(\eta)\right)\\&&=\sum_{i=1}^n\left(\sum_{y\ne
x_1,\dots,x_n} q_{y,x_i}(f_y-f_{x_i})\prod_{{j\ne i}}f_{x_j}
+\sum_{k=1}^n p(x_k,x_i)(f_{x_k}-f_{x_i})\prod_{j\ne
i}f_{x_j}\right)\nonumber\\&& =\sum_{i=1}^n\sum_y q_{y,x_i}
(F(\{x_1,\dots,x_n\}\backslash x_i\cup
y,\eta)-F(\{x_1,\dots,x_n\},\eta))\nonumber\eea

Now consider a set valued process $A(t)$ corresponding to the set
of points of a coalescing random walk. Then the transitions are
\be{} A\to A\cup\{y\}\backslash \{x\} \text{ at rate  }q_{x,y}. \ee
We denote the generator by $H$.
We observe that the generator $H$ applied to functions of the form  $F(\cdot,\eta),\;\eta\in \{0,1\}^S$ (defined in (\ref{tfvm})) satisfies

\be{} HF(\{x_1,\dots,x_n\},\eta)= GF(\{x_1,\dots,x_n\},\eta).\ee

We then apply the dual representation Theorem \ref{DR}.

\begin{remark}
In the case $S=\mathbb{Z}$, a direct relation between the voter
model and its dual can be demonstrated by a graphical
construction.


\begin{figure}[h]
\begin{center}
\setlength{\unitlength}{0.15cm}
\begin{picture}(150,60)
\put(55,0){$\rm{Time }\;\;t=  0$}
\put(55,50){$\rm{Time }\;\;t=  T$}
\put(0,0){\line(0,1){50}}
\put(10,0){\line(0,1){50}}
\put(20,0){\line(0,1){50}}
\put(30,0){\line(0,1){50}}
\put(40,0){\line(0,1){50}}
\put(50,0){\line(0,1){50}}
\linethickness{0.5mm}

\put(10,15){\red\vector(-1,0){10}}
\put(20,10){\red\vector(-1,0){10}}
\put(40,20){\red\vector(1,0){10}}
\put(30,18){\red\vector(1,0){10}}
\put(20,5){\red\vector(1,0){30}}
\put(10,45){\red\vector(-1,0){10}}
\put(20,40){\red\vector(-1,0){10}}
\put(30,35){\red\vector(-1,0){20}}
\put(40,38){\red\vector(1,0){10}}
\put(30,30){\red\vector(1,0){10}}
\end{picture}
\end{center}
\caption{Graphical Representation of the Voter Model. }
\end{figure}
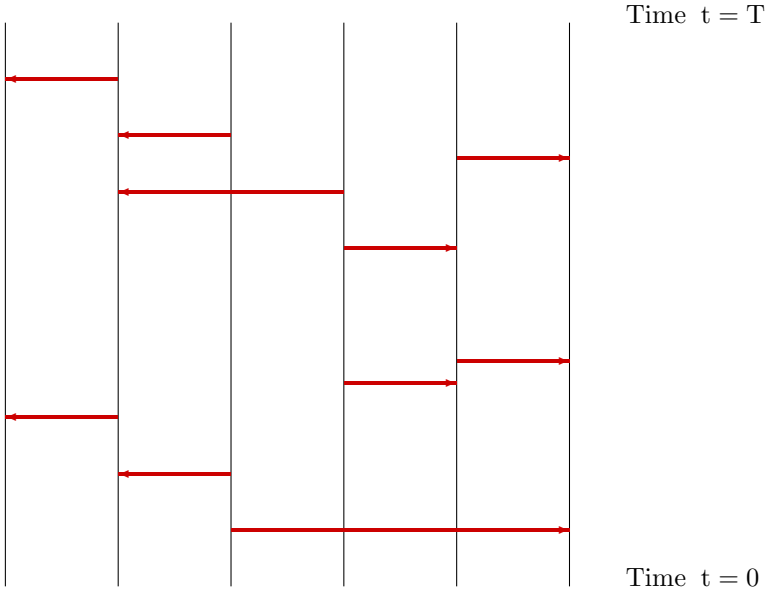

\end{remark}

As in the case of birth and death processes we can also obtain a representation of the voter model by a stochastic integral equation (cf.
Mueller-Tribe (1995) \cite{MT-95}, Kurtz-Protter (1996) \cite{KP-96}) as follows.

Let $\{\Lambda_t(x,y):x,y\in S\}$ be a family of independent
Poisson processes  with rates $p(y-x)$. Consider the system of
stochastic integral equations

\be{vmsie}
\xi_t(x)=\xi_0(x)+\sum_y\int_0^t[\xi_{s-}(y)-\xi_{s-}(x)]d\Lambda_s(x,y)\ee
with $x\in S,\;t\geq 0$, \be{}\xi_0(x)=0,\;\text{a.a.}x.\ee

\beP{} The system (\ref{vmsie}) has a unique solution and has
pregenerator $G$.

\end{proposition}
\begin{proof}
This is a special case of a result in Kurtz and Protter \cite{KP-96}, Chap. 9.
\end{proof}

A key property of the voter model is that the system started with a single non-zero site always dies out.  This follows since $\sum_{x\in S}\xi_t(x)$ is a martingale.

\subsection{Biased voter model}

The biased voter model is a modification of the voter model that arose from  the Williams-Bjerknes (1972) tumour
growth model.  Basic results on this model were obtained by
 Bramson and Griffeath (1980,1981) \cite{BG-80a}, \cite{BG-81} and Lanchier-Neuhauser (2007) \cite{LN-07}.

It is a spin system with generator

\be{}  Gf(\eta)=\sum_x c(x,\eta)[f(\eta_x)-f(\eta)]\ee where
\bea{}\qquad\quad
\eta_x(y)=\left\{
\begin{array}{ll}
    \eta(y), & \hbox{ if } x\ne y \\
    1-\eta(y), & \hbox{ if } x=y.\\
\end{array}\right.%
\eea
The biased voter model spin rates are given by:
\be{}  c(x,\eta)=\left\{%
\begin{array}{ll}
    \beta\sum_y q_{x,y}\eta(y), & \hbox{ if }\eta(x)=0 \\
    \sum_y q_{x,y}[1-\eta(y)] & \hbox{ if }\eta(x)=1 \\
\end{array}
\right.\ee where $ \beta \geq 0$.  If $\beta >1$, the dual process
is a coalescing branching random walk.  If $\beta>1$ the opinion {1} is favoured and growing clusters of $1$'s can form.
In fact Bramson and Griffeath  \cite{BG-81} prove conditioned on non-extinction that there is a growing region whose radius
 grows linearly and that the occupied regions has an asymptotic shape.

\subsection{Oriented percolation}

Consider the lattice $\mathbb{Z}_+\times \mathbb{Z}^d$. Points
$(n_1,x_1)$ and $(n_2,x_2)$ are neighbours iff $n_2=n_1+1$ and
$x_1,x_2$ are neighbours in $\mathbb{Z}^d$. We consider the bonds
between such neighbours and designate them open with probability
$p$ and closed with probability $1-p$. Percolation occurs by the
flow through open bonds.  We say that $(n,x)$ can be reached from
$(m,y)$ is there is a sequence of open bonds joining them.  We
consider the cluster $C_{(0,0)}$ consisting of points that can be
reached from  $(0,0)$.  Also let \be{}
\xi^0_n=\{x:(0,0)\to(n,x)\}.\ee We define \be{}
p_c=\inf\{p:P(\xi^0_n\ne\empty\;\forall\;n)>0\}.\ee

In the case $d=1$, Liggett \cite{L-95} proved that $p_c\leq \frac{2}{3}$ and it is known that
$p_c\geq .6446 $ (cf. Durrett (1985) \cite{D-85}).

See Durrett (1984) \cite{D-84} and Durrett-Tanaka (1989) \cite{DT-89} for the basic properties.

\subsection{Contact process}

The contact process was introduced by Mollison (1977) \cite{M-77}, and basic results were obtained by  Griffeath (1981) \cite{G-81}, Durrett-Griffeath (1982) \cite{DG-82}, and  Durrett-Schonmann (1987)  \cite{DS-87}.

The contact process on $S_1=\mathbb{Z}^d$ has transition rates:

\be{}  c(x,\eta)=\left\{%
\begin{array}{ll}
    \lambda  \sum_{|y-x|=1}\eta(y), & \hbox{ if }\eta(x)=0 \\
    1 & \hbox{ if }\eta(x)=1 \\
\end{array}
\right.\ee
where $\lambda >0$.

There is a critical value $\lambda_c$ of $\lambda$ such that there is a positive probability that the contact process on $\mathbb{Z}^d$ does not die out for $\lambda >\lambda_c$  and below which dies out with probability 1. See Durrett (1988) \cite{D-88} for an introduction to contact process on $\mathbb{Z}^d$. Recent results on the  contact process on hierarchical group are given in Athreya and Swart \cite{AS-08}.

\subsection{Interacting birth and death processes}
We next consider a class of processes which have been used to model chemical reaction diffusion systems.
The state space is $({\mathbb{Z}_+})^{S_1}$ and the generator has the form
\bea{} &\\Gf(x)&=\sum_{\xi\in S_1}\left\{ \lambda_1(x_\xi)[f(x+e_\xi)-f(x)]+\lambda_2(x_\xi)[f(x-e_\xi)-f(x)]\right\}\nonumber\\
&+\sum_{\xi\ne \xi'} x_\xi q_{\xi,\xi'}[f(x-e_\xi+e_{\xi'})-f(x)]\nonumber\eea
where $\lambda_1(\cdot),\;\lambda_2(\cdot) \geq 0$.
If $\lambda_1(k)= \beta_0 +\beta_1k,\qquad \lambda_2(k)= \delta_1 k +\delta_2 k^2$,
the construction and uniqueness for such systems has been established by M.-F. Chen (see for example \cite{C-04}). In the hydrodynamic limit they given rise to reaction-diffusion equations.

 If $\beta_0,\beta_1,\delta_1,\delta_2>0$, this is known as Schl\"ogl's first model. If $\beta_0=0, \beta_1>0, -\delta_1=\delta_2>0$, this corresponds to branching coalescing model (BC-model).

\section{Interacting diffusions}\label{intdif}

We now consider processes $X_t$ with local state space $[0,\infty)^M,\;M\in\N$ and configuration space $([0,\infty)^M)^{S_1}$ where $S_1$ is a countable abelian group.

\begin{itemize}
\item  $g:[0,\infty)\to [0,\infty)\qquad$ is locally Lipschitz continuous
\item $g^{-1}((0,\infty))=(0,b)\qquad $for some $b\in (0,\infty]$,
\item $g(z)\leq C(1+z^2)\qquad$ for some $C<\infty.$

\end{itemize}

Consider the system of stochastic differential equations
\begin{align*}
{d{ X_t^{(i)}(x)}}  &  = [\sum q_{x,y}({
X}_t^{(i)}(y)-{ X}_t^{(i)}(x))]dt  + \sqrt{{ g(
X^{(i)}(x))}} \;dW_t^{(i)}(x)\\&\\&i=1,\dots,M,
\quad x\in S_1
\end{align*}
where $\{(W_t^{(i)}(x))_{i=1,\dots,M}\}_{x\in S_1} $ are independent M-dimensional Wiener
processes and
 $\{q_{x,y}\}_{x,y\in S_1}$  are the transition rates for a symmetric random walk
 on $S_1$.

Assume that there exists $\{\gamma_i:i\in S_1\}$, a positive summable reference measure satisfying
\be{}  \sum_i \gamma_i q_{ij} \leq \Gamma \gamma_j,\quad j\in S_1,\quad\text{for some constant }\Gamma .\ee
The the Spitzer-Liggett space is defined as  $E:= \{z\in [0,b]^{S_1}: \| z\| <\infty\},\quad \| z\| =\sum_i \gamma_i|z_i|$ with the topology of componentwise convergence.

Conditions  for existence and uniqueness of solutions to these systems were obtained by Shiga and Shimizu (1980) \cite{SS-80}.

\subsection{Interacting Feller diffusions.}

The case $M=1$, $g(x)=\gamma x,\; \gamma >0,$ describes a system of interacting Feller CSBP processes.

\subsection{The Wright-Fisher stepping stone model.}

The special case $M=2$, $X^1,X^2\geq 0,\; X^1+X^2\equiv 1$, setting $X_t=X^1_t$, satisfying

\begin{align*}
dX_t(x)  &  =\sum_{y\in S_1}p_{y-x}(X_t(y)
-X_t(x)dt\\&+sX_t(x)(1-X_t(x))dt+\sqrt{2\gamma X_t(x)(1-X_t(x))}dW_t(x)\\
x_0(x)  &  = \theta\in [0,1] \quad\forall x\in S_1
\end{align*}
with $s=0$ describes the  stepping stone diffusion approximation  to the neutral stepping stone
model introduced by Mal\'{e}cot (1948) \cite{M-48}, and studied by Kimura (1953) \cite{K-53},  Kimura-Weiss (1964) \cite{KW-64}, Nagylaki (1974) \cite{N-74}, and  Sawyer (1976) \cite{S-76a}.

\begin{remark}
The voter model corresponds to the limiting case $s=0,\;\gamma\to
\infty$.
The additional term with coefficient $s$ represents the non-neutral case in which type 1 has fitness $s\ne 0$.
\end{remark}

Refer to Section \ref{s.NSSM} for the development of this model.

\subsection{Mean-field limit of exchangeable interacting diffusions}

Consider the system of exchangeable diffusions on $S=\{0,\dots,N-1\}$

\begin{eqnarray*}
&&dx_{\xi}(t)  = c(x_{\xi[1]}%
(t)-x_{\xi}(t))dt+\sqrt{2g(x_{\xi}(t)))}dw_{\xi}(t)\\&&x_\xi(0)=\theta_0\quad\forall
\xi\\&&{x_{\xi_{[1]}}(t):=\frac{1}{N}\sum_{\xi=1}^Nx_{\xi}(Nt) \qquad\qquad {\textbf{mean-field process}}}
\end{eqnarray*}

 Renormalized system:
 \be{} Z_1^N(t)=  x_{\xi_{[1]}}(Nt),\ee
\[
dZ^N_{1}(t)=-cZ^N_{1}(t))dt+\sqrt{\frac{2}{N}\sum_{\xi =1}^Ng(x_\xi(Nt))}\;dw(t).\]

Stationary measures:
\[
\Gamma^{g}_\theta(A):=\frac{1}{Z(g)}\int_A\frac{1}{g(x)}\exp\left[\int_\theta^x\frac{\theta-y}{g(y)}dy\right]dx\]

\[{\mathcal{F}(g)(\theta):=\int g(x)\Gamma^{g}_\theta(dx)}\]

\begin{theorem} As $N\to\infty $ \be{}Z^N_1(\cdot)\Rightarrow Z_1(t)\ee where $Z_1(t)$ satisfies
\[
dZ_{1}(t)=-cZ_{1}(t))dt+\sqrt{2g_1(Z_1(t))}\;dw(t).\]
where
\[{g_{1}=\mathcal{F}(g)}.\]
\end{theorem}
\begin{proof}
See Dawson-Greven (1993) \cite{DG-93}.
\end{proof}

\section{Measure-valued branching processes}

\subsection{Super-Brownian motion}
\subsubsection{Introduction}

Super-Brownian motion (SBM) is a measure-valued branching process which
generalizes the Jirina process.  It was constructed by S. Watanabe (1968)
\cite{W-68} as a continuous state branching process and Dawson (1975)
\cite{D-75} in the context of SPDE.  The lecture notes by Dawson (1993) \cite{D-93} and Etheridge (2000) \cite{E-00} provide introductions to measure-valued processes. The books of Dynkin \cite{Dy-94}, \cite{Dy-02}, Le Gall \cite{LG-99},  Perkins \cite{P-02} and Li \cite{Li-09} provide comprehensive developments of various aspects of measure-valued branching processes.
In this section we begin with  a brief introduction and then survey some aspects of superprocesses which are important for the study of stochastic population models.  Section 9.5 gives a brief survey of the small scale properties of SBM and Chapter 10 deals with the large space-time scale properties.

Of special note is the discovery in recent years that super-Brownian motion arises as the scaling limit of a number of models from particle systems and statistical physics. An introduction to this class of {\em SBM invariance principles} is presented in Section 9.6 with emphasis on their application to the voter model and interacting Wright-Fisher diffusions. A discussion of the invariance properties of Feller CSB in the context of a renormalization group analysis is given in Chapter 11.

\subsubsection{The SBM Martingale Problem}

Let $(D(A),A)$ be the generator of a Feller process on a locally compact metric space $(E,d)$ and  $\gamma\geq 0$.
The  probability laws $\{P_{\mu}:\mu\in M_{f}(E)\}$ on
$C([0,\infty
),M_{f}(E))$ of the  superprocess associated to $(A,a,\gamma)$  can be
characterized as the unique solution of the following martingale
problem:

\[
M_{t}(\varphi):=\left\langle \varphi,X_{t}\right\rangle -\int_{0}%
^{t}\left\langle  A\varphi,X_{s}\right\rangle ds
\]
is a $P_{\mu}$-martingale with increasing process%
\[
\left\langle M(\varphi)\right\rangle _{t}=\int_{0}^{t}\gamma \left\langle  \,
\varphi ,X_{s}\right\rangle ds
\]
for each $\varphi\in D(A)$.

Equivalently, it solves the martingale problem
\begin{eqnarray*}
GF &=&\int A \frac{\delta F}{\delta \mu (x)}\mu (dx) \\
&&+\frac{\gamma}{2}\iint \frac{\delta ^{2}F}{\delta \mu (x)\delta \mu
(y)}\delta _{x}(dy)\mu (dx)
\end{eqnarray*}
\[
{D}({G}):=\{F(\mu )=e^{-\mu (\varphi )},\,\,\varphi \in {\cal B}_{+}(%
{\Bbb R}^{d})\}
\]

The special case $E=[0,1]$, $Af(x)=[\int f(y)\nu_0(dy)-f(x)]dy$ is the Jirina process. The special case $E=\mathbb{R}^d$  $A=\frac{1}{2}\Delta$ on $D(A)=
C_b^{2}(\mathbb{R}^{d})$ is called super-Brownian motion.

\subsection{Super-Brownian Motion as the Limit of Branching Brownian Motion}

Given a system of branching Brownian motions  on $S={\mathbb R}^d$ and  $\ve>0$ we consider the measure-valued process, $X^\ve$, with particle mass  $
m_\varepsilon =\varepsilon $ and branching rate $\gamma _\varepsilon =\frac \gamma
\varepsilon $, that is,
\be{} X^\ve(t)=  m_\ve\sum_{j=1}^{ N(t)}\delta_{x_j(t)}\ee
where $N(t)$ denotes the number of particles alive at time $t$ and $x_1(t),\dots,x_{N(t)}$ denote the locations of the particles at time $t$.
Given an initial set of particles, let  $\mu_\ve=m_\ve\sum_{j=1}^{N(0)}\delta_{x_j(0)}$, let $P^\ve_{\mu_\ve}$ denote the probability law of $X^\ve$ on $D_{M_F(\mathbb{R}^d)}([0,\infty))$. Let $\{\mathcal{F}_t\}_{t\geq 0}$ be the canonical filtration on $D([0,\infty ),M_F(\mathbb{R}^d))$.

\begin{notation}
$\mu (\phi )= \langle \phi,\mu\rangle =\int \phi d\mu $.
\end{notation}
Let $C(M_F(\mathbb{R}^d))\supset D(G_\ve):=\{F(\mu )=f(\langle \phi,\mu \rangle )):f\in C^2_b(\mathbb{R}),\,\phi \in C^2_b({\mathbb R}^d)\}$.
Then $D(G_\ve)$ is measure-determining on $M_F(\mathbb{R}^d) $ (\cite{D-93}, Lemma 3.2.5.).

Then using It\^o's Lemma, it follows that $P^\ve_{\mu_\ve}\in \mathcal{P}(D([0,\infty ),M_F(\mathbb{R}^d)))$ satisfies the $G^\ve$-martingale problem
where for
 $F\in D(G^\varepsilon ),$%
\begin{eqnarray*}
G^\varepsilon F(\mu ) &=&f^{\prime }(\mu (\phi ))\mu (\frac 12\Delta \phi
)+\frac \varepsilon 2f^{\prime \prime }(\mu (\phi ))\mu (\nabla \phi \cdot
\nabla \phi ) \\
&&+\frac \gamma {2\varepsilon ^2}\int [f(\mu (\phi )+\varepsilon \phi
(x))+f(\mu (\phi )-\varepsilon \phi (x))-2f(\mu (\phi ))]\mu (dx).
\end{eqnarray*}

We can also obtain the Laplace functional of $X^\ve(t)$ using  Proposition \ref{P.BRWLF} with $\{S_t:t\geq 0\}$  the Brownian motion semigroup on $C(\mathbb{R}^d)$ and $\mathcal{G}(z)=\frac{1}{2}+\frac{1}{2}z^z$.

\begin{theorem}  Assume that $X^\varepsilon (0)=\mu_\varepsilon
\Rightarrow \mu$ as $\ve\to 0$.

Then\\
(a) $P_{\mu_\varepsilon }^\varepsilon \stackrel{\varepsilon
\rightarrow 0}{\Longrightarrow }{\mathbb P}_\mu\in \mathcal{P}(C_{M_F(\mathbb{R}^d)}([0,\infty))$
and ${\mathbb P}_\mu$ is the unique solution to the martingale problem: for all $%
\phi \in C^2_b(\mathbb{R}^d),$%
\be{SBMMP}
 M_t(\phi ):=X_t(\phi )-\mu(\phi )-\int_0^tX_s(\frac{1}{2}\Delta\phi )ds
\ee
is an $({\mathcal F}_t^X)-$martingale starting at zero with increasing process
\[
\langle M(\phi )\rangle_t=\gamma \int_0^tX_s(\phi ^2)ds.
\]

(b) The Laplace functional of $X_t$ is given by
\be{LF0}
{\mathbb P}_\mu \left(e^ {\left(-\int\phi(x)\mathbb{X}_t(dx)\right)}\right)=e^{ \left(-\int v_t(x)\mu(dx)\right)}.
\ee
where  $v(t,x) $ is the unique solution of
\be{LLE0}
\frac{\partial v(t,x)}{\partial t}= \frac{1}{2}\Delta v(t,x)-\frac \gamma 2v^2(t,x),\quad v_0=\phi
\in C^2_{+,b}(\mathbb{R}^d).
\ee

(c)  The total mass process $\{X_t(\mathbb{R}^d\}_{t\geq 0}$ is a Feller CSBP.
\end{theorem}

\begin{proof}

Step 1. Tightness of probability laws of $X^\ve$ on $D_{M_F(\mathbb{R}^d)}([0,\infty$  and a.s. continuity of limit points. In order to prove
tightness it suffices to prove that {for } $\delta >0$ there exists a compact subset $K\subset\mathbb{R}^d$ and $0<L<\infty$
 such  that
\bea{TC1}
&  \\
& P_{\mu_\varepsilon }^\varepsilon (\sup_{0\leq t\leq T}X_t(K^c)>\delta
)<\delta ,\,\,\,P_{\mu_\varepsilon }^\varepsilon (\sup_{0\leq t\leq
T}X_t(1)>L)<\delta\nonumber\eea
and
\bea{TC2}
& P_{\mu_\varepsilon }^\varepsilon \circ (X_t(\phi ))^{-1}\text{ is tight
in  }%
D_{\mathbb{R}}([0,\infty )) \text{ for } \phi \in C^2_c(\mathbb{R}^d).
\eea
This can be checked by standard moment and martingale inequality arguments. For example for (\ref{TC1}) it suffices to show that
\be{}\sup_{0< \ve\leq 1}\sup_{\delta>0} E(\sup_{0\leq t\leq T} \langle e^{-\delta \|x\|}(1+\|x\|^2),\X^\ve(t)\rangle )<\infty,\ee
and (\ref{TC2}) can be verified using the Joffe-M\'etivier criterion (see Appendix, (\ref{JMC})).
The  a.s. continuity of any limit point then follows from Theorem \ref{DTC} since the maximum jump size in $X^\ve$ is $\ve$.

Moreover, if ${\mathbb P}_\mu$ is a limit point, it is also easy to check (cf. Lemma \ref{WCM}) that for $\phi$
in $C^2_b(\mathbb{R}^d)$, $M_t(\phi )$ is a ${\mathbb P}_\mu$-martingale and ($F_1(\mu )=\mu
(\phi ),\,\,F_2(\mu )=\mu (\phi )^2$)
\begin{eqnarray*}
\langle M(\phi )\rangle _t &=&\lim_{\varepsilon \rightarrow 0}\int_0^t(G_\varepsilon
F_2(X_s)-2F_1(X_s)G_\varepsilon F_1(X_s))ds \\
&=&\gamma \int_0^tX_s(\phi ^2)ds.
\end{eqnarray*}
\newline
As pointed out above,  (\ref{SBMMP})  and  Ito's formula yields  an equivalent
formulation of the martingale problem, namely: for $f\in C_b^2({\mathbb R})$, $\phi
\in C^2_b(\mathbb{R}^d)$,  and $F(\mu )=f(\mu (\phi )),$%
\be{MP2}
F(X_t)-\int_0^tGF(X_s)ds\,\,is\,\,\,a\,\,{\mathbb \,P}_\mu-martingale
\ee
where
\[
GF(\mu )=f^{\prime }(\mu (\phi ))\mu (\frac{1}{2}\Delta\phi )+\frac \gamma 2\,f^{\prime
\prime }(\mu (\phi ))\mu (\phi ^2).
\]

Step 2. (Uniqueness)
In order to prove (b) we first verify that
${\mathbb P}_\mu$ also solves the following time dependent martingale problem.
 Let  $\psi :[0,\infty)\times
E\rightarrow [0,\infty)$ such that   $ \psi$, $\frac \partial {\partial s}\psi$   and $\Delta \psi$  are bounded and
strongly continuous in $C_b(\mathbb{R}^d)$.  Assume that
\be{smooth}
\left\| \frac{\psi (s+h,\cdot )-\psi (s,\cdot )}h-\frac \partial {\partial
s}\psi (s,\cdot )\right\| _\infty {\rightarrow }0\quad\text{as } h\to 0.
\ee
Then
\be{IMP}
\exp (-X_t(\psi _t))+\int_0^t\exp (-X_s(\psi _s))X_s((A+\frac \partial
{\partial s})\psi _s)ds-\frac \gamma 2\int_0^t\exp (-X_s(\psi _s))X_s(\psi
_s^2)ds
\ee
is a ${\mathbb P}_\mu$-martingale.  Let ${\mathbb P}_{\mu}^{{\mathcal F}_t}$ denote the conditional
expectation with respect to ${\mathcal F}_t$ under ${\mathbb P}_{\mu}.$

To prove (\ref{IMP}) first note that applying (\ref{MP2}) to $\exp (-\mu (\phi ))$ with $\phi \in
C^2_b(\mathbb{R}^d)$, we obtain
\be{EMP}
{\mathcal E}_t(\phi )=\exp (-X_t(\phi ))+\int_0^t\exp (-X_s(\phi ))X_s(A\phi
)ds-\frac \gamma 2\int_0^t\exp (-X_s(\phi ))X_s(\phi ^2)ds
\ee
is a ${\mathbb P}_\mu$-martingale.

Next take \bea{} && u(s,X_t)=\exp (-X_t(\psi _s)),\quad v(s,X_t)=\exp (-(X_t(\psi
_s))X_t(\frac \partial {\partial s}\psi _s),\text{   and}\\&&  w(s,X_t)=\exp (-X_t(\phi
))(X_t(A\psi _s))\nonumber\eea
so that
for $t_2>t_1$
\be{} u(t_2,X_{t_2})-u(t_1,X_{t_2}) = -\int_{t_1}^{t_2}v(s,X_{t_2})ds. \ee
Then using (\ref{EMP}) we have
\be{1.1}\begin{split}
&&{\mathbb P}_{\mu}^{{\mathcal F}_{t_1}}[u(t_1,X_{t_2})-u(t_1,X_{t_1})]
 = -{\mathbb P}_{\mu}^{%
{\mathcal F}_{t_1}}[\int_{t_1}^{t_2}w(t_1,X_s)ds]    \\
&  &+\frac \gamma 2{\mathbb P}_{\mu}^{{\mathcal F}_{t_1}}[\int_{t_1}^{t_2}u(t_1,X_s)
X_s(%
\psi _s^2)ds].\end{split}
\ee
Let $\Lambda ^n$ be a partition of $[t_1,t_2]$ with mesh($\Lambda
^n)\rightarrow 0$ and
\begin{eqnarray*}
\psi ^n(s,x) &:&=\sum_{i=1}^n\psi (t_i^n,x)1_{[t_i^n,t_{i+1}^n)}(s) \\
X^n(s) &:&=\sum_{i=1}^nX_{t_{i+1}^n}1_{[t_i^n,t_{i+1}^n)}(s)
\end{eqnarray*}
Let $u^n(t,X_t):=\exp (-X_t(\psi _t^n))$.

Then by (\ref{1.1})
\begin{eqnarray*}
{\mathbb P}_{\mu}^{{\mathcal F}_{t_1}}[u^n(t_2,X_{t_2})-u^n(t_1,X_{t_1})] &=&-{\mathbb P}%
_m^{{\mathcal F}_{t_1}}[\int_{t_1}^{t_2}\exp (-X_s^n(\psi _s))X_s^n(\frac
\partial {\partial s}\psi _s)ds] \\
&&-{\mathbb P}_{\mu}^{{\mathcal F}_{t_1}}[\int_{t_1}^{t_2}\exp (-X_s(\psi
_s^n))X_s(A\psi _s^n)ds] \\
&&+\frac \gamma 2{\mathbb P}_{\mu}^{{\mathcal F}_{t_1}}[\int_{t_1}^{t_2}\exp (-X_s(\psi
_s^n))X_s((\psi _s^n)^2)ds].
\end{eqnarray*}
Standard arguments show that this converges to
\begin{eqnarray*}
{\mathbb P}_{\mu}^{{\mathcal F}_{t_1}}[u(t_2,X_{t_2})-u(t_1,X_{t_1})] &=&-{\mathbb P}_m^{%
{\mathcal F}_{t_1}}[\int_{t_1}^{t_2}\exp (-X_s(\psi _s))X_s(\frac \partial
{\partial s}\psi _s)ds] \\
&&-{\mathbb P}_{\mu}^{{\mathcal F}_{t_1}}[\int_{t_1}^{t_2}\exp (-X_s(\psi _s))X_s(A\psi
_s)ds] \\
&&+\frac \gamma 2{\mathbb P}_{\mu}^{{\mathcal F}_{t_1}}[\int_{t_1}^{t_2}\exp (-X_s(\psi
_s))X_s((\psi _s)^2)ds]
\end{eqnarray*}
which completes the proof of (\ref{IMP}).

Now let $v_t=V_t\phi $ be the unique solution (see \cite{P-83}) of
\be{LLE}
\frac{\partial v_t}{\partial t}=Av_t-\frac \gamma 2v_t^2,\quad v_0=\phi
\in C^2_{+,b}(\mathbb{R}^d).
\ee

Then $v_t$ satisfies (\ref{smooth}).  Applying (\ref{IMP}) we deduce that $\{\exp
^{(-X_s(v_{t-s}))}\}_{0\leq s\leq t}$ is a martingale. Equating mean values at $%
s=0$ and $s=t$ we get the fundamental equality

\be{LF}
{\mathbb P}_\mu{\mathbb (}\exp (-X_t(\phi ))=\exp (-\mu(V_t\phi )).
\ee
The extension from $\phi \geq 0$ in $D(A)$ to $\phi \geq 0$ in $b{\mathcal E}$
follows easily by considering the ``weak form''
\[
V_t\phi =P_t\phi -\frac \gamma 2\int_0^tP_{t-s}(V_s\phi )^2ds
\]
and then taking bounded pointwise limits.

(\ref{LF})  proves the uniqueness of the law ${\mathbb P}_\mu(X_t\in \cdot )$ for any
solution of (MP) and hence the uniqueness of ${\mathbb P}_\mu$ (see  \cite{SV-79} or \cite{EK-86}, Chapt. 4, Theorem 4.2).

(c) follows by taking $\phi\equiv 1$ and comparing the Laplace transforms of the transition measures.

\end{proof}

\begin{corollary}
(Infinite divisibility)
$ {\mathbb P}_\mu$ is an infinitely divisible probability measure with canonical representation
\bea{}&&
-\log ({\mathbb P}_\mu(\exp^{ (-X_t(\phi ))})\\&& = -\int \log ({\mathbb P}_{\delta _x}(\exp
^{(-X_t(\phi ))} )\mu(dx) = \int_{M_F(\mathbb{R}^d)\backslash \{0\}} (1-e^{\nu(\phi)})R_t(x,d\nu)\nonumber
\eea
where the canonical measure $R_t(x,d\nu)\in M_F(M_F(\mathbb{R}^d)\backslash \{0\})$ satisfies\\
$R_t(x, M_F(\mathbb{R}^d)\backslash \{0\})= \frac{2}{\gamma t}$ and the normalized measure is an exponential probability law with mean $\frac{\gamma t}{2}$.

\end{corollary}

The infinite divisibility of SBM allow us to use the the theory of infinitely divisible random measure (see e.g. Dawson (1992) (1993), \cite{D-92}, \cite{D-93})
to obtain detailed properties of the process.

\subsubsection{Weighted occupation time}

If $\{X_t:t\geq 0\}$ is a super-Brownian motion, then we defined the associated {\em weighted occupation time} $Y_y:t\geq 0$ as
\be{}  Y_t(A)=\int_0^t X_s(A)ds,\quad A\in\mathcal{B}(\mathbb{R}^d).\ee

\beT{} (Iscoe (1986) \cite{I-86a})\\
Let $\mu\in M_F(\mathbb{R}^d)$ and $\phi,\psi\in C^2_{c,+}(\mathbb{R}^d)$.  Then the joint Laplace functional of $X_t$ and $Y_t$ is given by
\be{} E_\mu\left[ e^{-\langle \psi,X_t\rangle -\langle \phi,Y_t\rangle}\right] =e^{-\int (U^\phi_t\psi)(x)\mu(dx)}\ee
where $(U^\phi_t\psi((x)))$ is the unique solution of
\bea{uphie} &&\frac{\partial u(t,x)}{\partial t}=\frac{1}{2} \Delta u(t,x) - \frac{\gamma}{2}u^2(t,x)+\phi(x),\\&&
u(0,x)=\psi(x). \nonumber \eea

\end{theorem}
\begin{proof}
Now consider two semigroups on $C_{0,+}\mathbb{R}^d)$:\begin{itemize}
\item $V_t\psi $ given by the solution of
\be{} v(t)=S_t\psi-\int_0^t S_{t-s}(v(s))^2ds,\quad \ee
\item $W_t$ given by
\be{} W_t\psi =\psi +t\phi,\qquad \dot{W}_t\psi =\phi.\ee
\end{itemize}
Using the iterated conditioning and the Markov property at the time points $\{\frac{N-1}{N}t,\frac{N-2}{N}t,\frac{N-3}{N}t,\dots,\frac{1}{N}t\}$ we obtain
\bea{} &&E_\mu\left(\exp\left[-\langle\psi,X_y \rangle -\int_0^t\langle \phi,X_s\rangle ds\right]\right)\\&&
= \lim_{N\to\infty}\exp\left[ -\langle(V_{\frac{t}{N}}W_{\frac{t}{N}})^N\psi,\mu\rangle \right]\nonumber
\eea
Then by the Trotter-Lie product formula (cf. Chorin et al \cite{CHMM-78})
\be{uphi} U_t=\lim_{N\to\infty}\left(V_{\frac{t}{N}}W_{\frac{t}{N}}\right)^N\qquad\text{on  } C_{0,+}(\mathbb{R}^d).\ee
and therefore
\be{}E_\mu\left(\exp\left[-\langle\psi,X_t \rangle -\int_0^t\langle \phi,X_s\rangle ds\right]\right)
=\exp\left(-\langle U_t\psi_t,\mu\rangle\right)\ee
noting that the interchange of limit and integral is justified by dominated convergence since
\be{} (V_{\frac{t}{N}}W_{\frac{t}{N}})^N\psi\leq (S_{\frac{t}{N}}W_{\frac{t}{N}})^N\psi \leq S_t\psi +(1+t)\|\phi\|.\ee
Finally, we can verify that the semigroup $U^\phi_t$ defined by (\ref{uphi}) satisfies (\ref{uphie}).

\end{proof}

As an application of this Iscoe established the {\em compact support property } of super-Brownian motion:

\beT{} Let  $\{ X_t:t\geq 0\}$ with be super-Brownian motion with initial measure $\delta_0$. Then

\be{} P_{\delta_0}(\sup_{0\leq t<\infty}  X_t(\mathbb{R}^d\;\backslash\; \overline{B(0,R)}\, )\,>0)=1-e^{-\frac{u(0)}{R^2}},\ee
where $u(\cdot)$ is the solution of
\bea{} &&\Delta u(x)=u^2(x),\; \; x\in B(0,1),\\&&
u(x)\to \infty \text{  as  } x\to \partial B(0,1).\nonumber\eea
\end{theorem}
\begin{proof}
Iscoe (1988) \cite{I-88}.
\end{proof}

\subsubsection{Convergence of renormalized BRW and interacting Feller CSBP}

A number of different variations of rescaled branching systems on $\mathbb{R}^d$ can be proved to converge to SBM. For example,
the following two results can be proved following the same basic steps.
\beT{}
Let  $\ve=\frac{1}{N}$. Consider a sequence of branching random walks $X^{\ve}_t$ on $\ve\mathbb{Z}^d$ with random walk kernel $p^\ve(\cdot)$ which satisfies (\ref{rw1}), particle mass $m_\ve=\ve$, branching rate $\gamma^\ve =\frac{\gamma}{\ve}$ and assume that $X^\ve_0\Rightarrow X_0$ in $\mathcal{M}_F(\mathbb{R}^d)$.  Then $\{X^\ve_t\}_{t\geq 0}\Rightarrow \{X_t\}_{t\geq 0}$
where $X_t$ is  a super-Brownian motion on $\mathbb{R}^d$ with $A=\frac{\sigma^2}{2}\Delta$ and branching rate $\gamma$.

\end{theorem}

\begin{remark}
The analogue of these results for more general branching mechanisms with possible infinite second moments are established in  Dawson (1993), \cite{D-93} Theorem 4.6.2.

\end{remark}
\beT{} Consider a sequence of interacting Feller CSBP in which the rescaled random walks converge to Brownian motion.  Then the interacting Feller CSBP converge to SBM.

\end{theorem}

\subsection{The Poisson Cluster Representation}

Let $X$ be an infinitely divisible random measure on a Polish space $E$ with finite expected total mass
$E[X(E)]<\infty$.  Then (cf.  \cite{D-93}, Theorem 3.3.1) there exists a measure $X_D\in
M_F(E)$ and a measure $R\in M(M_F(E)\backslash \{0\})$ satisfying
\[
\int (1-e^{-\mu (E)})R(d\mu )<\infty
\]
and such that
\[
-\log (P(e^{-X(\phi )}))=X_D(\phi )+\int (1-e^{-\nu (\phi )})R(d\nu ) .
\]

$X_D$ is called the deterministic component and $R$ is called the {\it %
canonical measure}$.$ For example, for a Poisson random measure with
intensity $\Lambda ,$ $R(d\nu )= \int \delta _{\delta _x}(d\nu )\Lambda %
(dx).$
If
we replace each Poisson point, x, with a random measure ({\it cluster}) with
probability law $R(x,d\nu )$ then we obtain
\[
R(d\nu )= \int \Lambda (dx)R(x,d\nu ).
\]
In the case of super Brownian motion with $X_0=\mu \in M_F(\mathbb{R}^d)$, $X_t$ for $t>0$ is infinitely
divisible with $X_D=0$ and canonical measure $\,R_t(d\nu )=\int
\mu(dx)R_t(x,d\nu )$ where
\[
V_t\phi (x)=\int (1-e^{-\nu (\phi )})R_t(x,d\nu ),\,\,\,\phi \geq 0
\]
(see e.g. Dawson and Perkins (1991) \cite{DP-91}). Since $(V_t  \theta )(x)
\, = \, \frac{2 \theta }{(2 \, + \, \theta \, \gamma  \,t)}, \\$
$R_t(x,\, M_F({\mathbb R}^d)\, \backslash \,
\{0\}) \, = \lim_{\theta \rightarrow \infty }(V_t \, \theta )(x)
=\frac 2{\gamma t}$
and
\[
\int e^{-\theta \nu (1)}R_t(x,d\nu )=\frac{(\frac 2{\gamma t})^2}{(\frac
2{\gamma t})+\theta }.
\]

Hence $R_t(x,\nu (1)\in \cdot )$ is $\frac 2{\gamma t}$ times an exponential
law with mean $\frac{\gamma t}2.$ Then using the above and  (\ref{EXT}) one can check that
\[
\int e^{-\theta \nu (\phi )}R_t(x,d\nu )=\lim_{\varepsilon \downarrow
0}\varepsilon ^{-1}{\mathbb P}_{\varepsilon \delta _x}(e^{-\theta X_t(\phi
)}1(X_t>0)).
\]
Hence $R_t(x,\cdot )$ can be interpreted as the unnormalized distribution of
$X_t$ starting with infinitesimal mass at $x$ (when it is nonzero).

\subsection{The Palm measure}
We now consider the locally size-biased law which is given by the Campbell measure
\be{} \bar R_t(x,A\times B)=\nu(A)1_B(\nu)R_t(x,d\nu),\qquad A\in\mathcal{B}(\mathbb{R}^d),\; B\in\mathcal{B}(M_F(\mathbb{R}^d))
\ee
The {\em Palm measure at $y\in\mathbb{R}^d$} is defined by
\be{}
(R_t(x,d\mu))_y=
\frac{\bar R_t(x,dy\times d\mu) }{I(dy)},\qquad I(A)=\int \mu(A) R_t(x,d\mu)
\ee
\begin{remark}
In the case of a single point $X_t=x_t\delta_e$ (instead of $\mathbb{R}^d$), \be{}
(R_t(x,B))_y=
\int_B \frac{\mu(e) R_t(x, d\mu) }{I(e)}
\ee
we obtain the size-biased distribution of the mass at the point.
\end{remark}

The Laplace functional of the Palm measure is given by  (see \cite{DP-91})
\be{}
\int e^{-\mu(\phi)}(R_t(x,d\mu))_y =P_y\left( e^{-\gamma\int_0^t(V_r\phi)(w_r)dr}|w_t=x\right)
\ee where $P_y$ is the law of Brownian motion started at $y$  and $V_t\phi$ satisfies
\be{} V_t\phi(x)=S_t\phi(x)-\frac{\gamma}{2}\int_0^tS_u(V_{y-u}\phi)^2du.\ee

\subsection{Excursions}

$w$ $\in C([0,\infty ),M_F(E))$ is an $M_F(E)-$valued with lifetime  $%
\tau $ starting at x if

\begin{eqnarray*}
&&w_0=0,\,\,\tau (w)>0,\,\,w_t\ne 0,\,\,0<t<\tau
(w),\,\,w_t=0,\,\,t>\tau (w) \\
&&w_t(1)^{-1}w_t\Rightarrow \delta _x\text{ as}\ t\rightarrow 0+.
\end{eqnarray*}

Let $C_x^0=C_x^0([0,\infty ),M_F(E)\backslash \{0\})$ denote the set of all
excursion paths.

\begin{theorem} (Super-excursions and canonical measure)
(El Karoui and Roelly \cite{EKR-91}, Li and Shiga \cite{LS-95})
There is a unique $\sigma -$finite kernel $%
R(x,\cdot )$ from $E$ to $C_x^0$ such that\newline
(i) for each $x\in E$ $ R(x,\cdot )$ is supported by $C_x^0$ \newline
(ii) $\int w_t(1)R(x,dw)\rightarrow 1$ as $t\rightarrow 0+,$ and \newline
(iii) The measure $R(x,\cdot )$ is Markov with the same transition laws as
super-Brownian motion and has t-marginal distribution equal to $R_t(x,\cdot
).$\newline
(iv)
\[
E_m(e^{-X_t(\phi )})=e^{ \left(-\int \int 1_{\tau (w)\geq t}(1-e^{-w_t(\phi
)})R(x,dw)m(dx)\right)}.
\]

\end{theorem}

\begin{remark}
$R(x,\{\tau >t\})=\frac 2{\gamma t}$ and we can verify that
\[
R(x,B\cap \{\tau >t\})=\lim_{\varepsilon \rightarrow 0+}\frac{%
P_{\varepsilon \delta _x}(B\cap \{X_t>0\})}\varepsilon
\]
\end{remark}

The {\em integrated super-excursion} (ISE) is defined to be the random measure on $\mathbb{R}^d$ given by
\be{} \mathcal{I}(B):= \int_0^\infty w_s(B)ds\quad \text{
conditioned on } \mathcal{I}(\mathbb{R}^d)=1.\ee It has been proved that ISE also describes the scaling
limit of random trees in high dimensions (see Derbez and Slade (1997) \cite{DS-97} and is  conjectured to arise in the scaling limit of critical
percolation clusters in high dimensions (see Slade (2002), (2006)  \cite{S-02}, \cite{S-06} for a complete  exposition of these results.)

\subsection{The Historical Branching Process}

An enrichment of SBM called the  the Historical
Brownian motion (HBM) was introduced by Dawson and Perkins (1991) \cite{DP-91} and studied by Dynkin (1991) \cite{Dy-91b},  Le Gall (1991) \cite{LG-91}
and Donnelly and Kurtz (1996) \cite{DK-96}.

Given the Brownian motion $\xi $, the {\it path-valued process} $\bar{\xi}%
(t)=\xi (\cdot \wedge t)$ is a time-inhomogeneous continuous strong Markov
process taking values in the Polish space $C=C({\mathbb R}_{+},E).$ We let $%
P_{r,y}\,\,($here $y(\cdot \wedge r)=y)$ denote the law of $\bar{\xi}$
starting at $y$ at time $r$. Historical Brownian motion is obtained by
replacing $\xi $ with the path-valued Brownian motion $\bar{\xi}$ in the
general superprocess construction outlined above.

We can consider the corresponding empirical process for the BBM
\[
H_t^\varepsilon =\sum_{{\bf \alpha }\sim t}m_\varepsilon \delta _{\bar{%
\xi}_\alpha (t)}.
\]
By a limiting procedure analogous to the above (but involving additional technical considerations), we get the $\xi $%
-historical process $\{H_t\}$ and law ${\mathbb Q}_{\tau ,m}$ on $\Omega _H$.
Let $y^t:=y(\cdot \wedge t)$, $C^t:=\{y:y=y^t\}$, $M_F(C)^t=\{m\in
M_F(C):y^t=y\,\,\,\,m-a.e.\,\,\,\,y\},$ . ${\mathcal C}$ is the Borel $\sigma -$%
field of $C$ and ${\mathcal C}_t$ is the sub-$\sigma -$field generated by the
paths up to time $t.$
\begin{eqnarray*}
\Omega _H &=&\{H_{\cdot }\in C({\mathbb R}_{+},M_F(C)):H_t\in
M_F(C)^t\,\,\,\,\forall t\geq 0\}, \\
{\mathcal H}[s,t] &=&\sigma (H_u:s\leq u\leq t),\,\,{\mathcal H}[s,\infty )=\sigma
(H_u:u\geq s).
\end{eqnarray*}

Foe $E=\mathbb{R}^d$ let $C_{0,\infty \text{ }}^d$ denote the continuous functions from $%
[0,\infty )$ to ${\mathbb R}^d$ and $D(\bar{A})$ denote the subset of $f\in
C_{0,\infty }^d$ of the form
\[
f(y)=g(y(t_1),\dots ,y(t_n)),\,\,\,0\leq t_1\leq t_2\dots \leq t_n,
\]
where $g$ is infinitely differentiable and constant outside a compact set.

For $f\in D(\bar{A})$ and $t>0$, let
\[
\bar{A}_tf(y):=\frac 12\sum_{i=1}^d\sum_{k=0}^{n-1}\sum_{\ell
=0}^{n-1}1(t<t_{k+1}\wedge t_{\ell +1})g_{x_{kd+i\,\,},x_{\ell
d+i}}(y^t(t_1),\dots ,y^t(t_n)).
\]
where $g_{x_{kd+i\,\,},x_{\ell
d+i}}$ denote the second partial derivatives of $g$ with respect to $x_{kd+i\,\,},x_{\ell
d+i}$,  $k,\ell =1,\dots, n-1$, $i=1,\dots,d$.
\beT{HMP1}
(Perkins (1995) \cite{P-95}) If $\tau \geq 0$ and $m\in M_F(C)^\tau $, ${\mathbb {Q}}_{r ,m}$
is the unique law on $(\Omega _H,{\mathcal H}[\tau ,\infty ))$ such that

\be{HMP}\begin{split}
&\text{for each }f\in D(\bar{A}),\,\,\,M_t(f)=H_t(f)-m(f)-\int_r ^tH_s(
\bar{A}_sf)ds,\text{ } t\geq r
\\
 &\text{ is an }({\mathcal H}[r ,\infty ))\text{-martingale, starting at zero
when }t=r, \\
& \text{and with quadratic variation}\\
&\langle M(f)\rangle _t=\int_r^t\gamma H_s(f^2)ds.\end{split}
\ee
\end{theorem}

The Laplace functional characterizing $H$ has the form
\begin{eqnarray*}
{\mathbb Q}_{r,m}(\exp (-H_t(\phi ))) &=&\exp (-m(V_{r,t}\phi )) \\
V_{r,t}\phi (y) &=&S_{r,t}\phi (y)-\frac \gamma 2\int_r^tS_{r,s}(V_{s,t}\phi
)^2ds,\,\,\,y\in C^t \\
S_{r,t}\phi (y) &=&P_{r,y(r)}\phi (y/r/\xi (.\wedge t))
\end{eqnarray*}
where if $y,\xi \in C_{0,\infty }^d,$ $r\in [0,\infty ),$then $(y/r/\xi )\in
C_{0,\infty }^d$ is defined by
\[
(y/r/\xi )(u):=\left\{
\begin{array}{l}
y(u)\text{ if }u<r \\
\xi (u-r)\text{ if }u\geq r
\end{array}
.\right.
\]

\begin{corollary}
Let $\pi _t:C\rightarrow E$, $\pi _t(y):=y(t).$ Then .
\[
\pi _t(H_t)\text{ is a version of the SBM }X_t
\]
\end{corollary}

\begin{remark} \label{Remark2.1}
If $A_1,\dots ,A_n$ are disjoint sets in ${\mathcal C}_r $ and $%
f_i(t,y)=1_{A_i}(y),$ then clearly $f_i\in D(A)$ and $A_{\tau ,m}f_i=0$ .
(HMP) shows that under ${\mathbb {Q}}_{r ,m},\,\,$ $X_t^i=H_{r +t}(A_i),$
$i=1,\dots ,n$, are independent FB processes starting at $m(A_i),$ $%
i=1,\dots ,n$. Therefore
\be{FB1}
{\mathbb {Q}}_{r ,m}(\exp (-\sum_i^n\lambda _iH_{\tau +t}(A_i)))=\exp
(-\sum_{i=1}^n2\lambda _im(A_i)/(2+\lambda _i\gamma t)),\,\,\,\lambda _i\geq
0,
\ee

and
\be{FB2}
{\mathbb {Q}}_{r ,m}(H_{\tau +s}(A_i)=0\,\,\forall s\geq t)=\exp
(-2m(A_i)/(\gamma t)).
\ee
\end{remark}

\begin{proposition}
(Historical Cluster Representation) Let $r \leq s<t$. Under ${\mathbb {Q}}
_{r ,m}$, the conditional distribution of $H_t(\{y:y^s\in \cdot \})$
given ${\mathcal H}[r ,s]$ is the law of $\sum_{i=1}^M\delta _{y_i}m_i$,
where $\{y_1,\dots ,y_M\}$ are the points of a Poisson point process with
intensity $H_s(\cdot )2\gamma ^{-1}(t-s)^{-1}$ and $\{m_1,\dots ,m_M\}$ are
independent exponential masses with mean $\gamma (t-s)/2$.
\end{proposition}

{\bf Proof}. As above we may take $s=r $ and argue unconditionally. An easy
calculation shows the Laplace functional of the above Poisson cluster random
measure, $\Xi $, is
\be{HCR}
{\mathbb {P}(}\exp \bigl( -\Xi (\phi )\bigr) )=\exp \left( -\int 2\phi (y){%
\bigl( 2+\phi (y)\gamma (t-r )\bigr)}^{-1}\,m(dy)\right) .
\ee
If $\phi (y)=\sum\limits_1^n\lambda _i1_{B_i}(y)$ for $B_1,\ldots ,B_n$
disjoint Borel sets in $C$, (\ref{FB1}) implies \newline
${\mathbb Q}_{r ,m}\biggl( \exp (-\int \phi \bigl( \overline{y}(\tau )%
\bigr)
H_t(d\bar{y}))\biggr) $ is also given by the right side of (\ref{HCR}). By taking
limits of simple functions we see that the Laplace functionals of the random
measures in question are equal and hence so are their laws.
$\Box$

\begin{remark} Perkins developed the {\em historical stochastic calculus } which is an analogue of It\^o stochastic
calculus for the historical superprocess (see Perkins (1992),\cite{P-92}, \cite{P-95}) and in Evans-Perkins (1998)
 \cite{EP-98} used it to obtain a model of competing superprocesses with killing according to a collision local time.
 See Perkins (2002) \cite{P-02} for a comprehensive exposition of these developments.

\end{remark}

\subsection{A skew product construction}
\label{ss.skew}
Etheridge and March \cite{EM-91} proved that between normalized SBM conditioned to have constant mass is a Fleming-Viot process. Perkins then
proved that (unconditioned) normalized SBM could be viewed as a Fleming-Viot process with time-dependent  resampling rate inversely  to the total mass
process for the SBM (Perkins Disintegration Theorem \cite{P-91}).   We leave the precise statement and proof of these results to Chapter 12 but using these ideas
now give an informal construction of SBM from the Feller CSB and the Fleming-Viot process.

Consider the following $\mathbb{R}^{+}\times M_{1}(E)$-valued martingale
problem for
$(x(t),Y(t))$%
\bea{}
&  x(t)\text{ is a martingale with increasing process}\\
&\left\langle x\right\rangle _{t}   =\int_{0}^{t}x(s)ds\nonumber\eea
Let
\be{}
\tau   :=\inf\{t:x(t)=0\}.\ee
In other words $x(t)$ is a critical Feller CSBP process.

Perkins \cite{P-91} introduced a Fleming-Viot process $Y_t$ with time varying resampling rate $1/x(t)$ that satisfies the following
martingale problem:
\bean{}
M_{t}(\varphi)  &  :=\{\left\langle \varphi,Y_{t}\right\rangle -\int_{0}%
^{t}\left\langle A\varphi,Y_{s}\right\rangle ds\}_{0\leq
t<\tau}\text{ is a
martingale }\\&
\left\langle M(\varphi)\right\rangle _{t}    =\int_{0}^{t}\frac{1}%
{x(s)}[\left\langle \varphi^{2},Y_{s}\right\rangle -\left\langle
\varphi ,Y_{s}\right\rangle ^{2}]ds
\eean{}
Now consider the process%
\[
X(t):=x(t)Y_{t}%
\]
Let us verify that $X$ is a solution to the super-A martingale
problem.
Applying Ito's Lemma to the semimartingale $\left\langle \varphi
,X(t)\right\rangle $ we obtain%
\begin{align*}
&  d\left\langle \varphi,X(t)\right\rangle  \\
&  =x_{t}\left\langle A\varphi,Y(t)\right\rangle
+x_{t}dM_{t}(\varphi
)+\left\langle \varphi,Y_{t}\right\rangle dx_{t}\\
&  =\left\langle A\varphi,X(t)\right\rangle +\tilde{M}_{t}(\varphi)\\
\tilde{M}_{t}(\varphi)  &  :=x_{t}dM_{t}(\varphi)+\left\langle
\varphi
,Y_{t}\right\rangle dx_{t}\\
\left\langle \tilde{M}(\varphi)\right\rangle _{t}  &  =\int_{0}^{t}%
\{\frac{x_{s}^{2}}{x_{s}}[\left\langle
\varphi^{2},Y_{s}\right\rangle -\left\langle
\varphi,Y_{s}\right\rangle ^{2}]+\left\langle \varphi
,Y_{s}\right\rangle ^{2}x_{s}\}ds\\
&  =\int_{0}^{t}\{\left\langle \varphi^{2},X_{s}\right\rangle -x_{s}%
\left\langle \varphi,Y_{s}\right\rangle ^{2}+\left\langle \varphi
,Y_{s}\right\rangle ^{2}x_{s}\}ds\\
&  =\int_{0}^{t}\left\langle \varphi^{2},X_{s}\right\rangle ds.
\end{align*}


\subsection{The Donnelly-Kurtz countable particle representation}

Donnelly-Kurtz (1999) \cite{DK-99b} developed a countable particle representation that includes both the Fleming-Viot process and superprocess in the same spirit as the lookdown process construction of the Fleming-Viot process but we do not consider this here.

\subsection{Brownian Snake Representation of SBM
$X_{t}$}

In (1991) Le Gall \cite{LG-91}  developed the {\em Brownian snake representation} of super-Brownian motion.  This is based on a deep connection between
Brownian excursions and Feller CSB that was suggested by the discovery by  Neveu and Pitman (1989) \cite{NP-89} of a natural branching process structure
within a Brownian excursion. We have already met this in the discussion of the continuum random tree.

The Brownian snake is based on reflecting  Brownian motion $\{\zeta_t\}_{t\geq 0} \in \mathbb{R}_+$ which serves as the {\em lifetime process}.
 Let $C_s([0,\infty),\mathbb{R}^d)$  denote the set of continuous paths in $\mathbb{R}^d$ stopped at time $s$.
Then given $\{\zeta_t\}_{t\geq 0}$ the Brownian snake  $ W_{s} \in C_s([0,\infty),\mathbb{R}^d)$ and is characterized as follows: for
$s_{1}<s_{2}$%
\begin{align*}
&  { W_{s_{1}}(u)=W_{s_{2}}(u)},\,\ \ u\leq{\min_{s\in\lbrack s_{1},s_{2}]},%
\zeta_{s}}\\
&{ \text{the continuations are given by independent Brownian paths  }}{ u\geq \min_{s\in\lbrack
s_{1},s_{2}]}\zeta
_{s}}
\end{align*}
The excursion measure of the Brownian snake starting at $x$ is defined by
\be{}  N_x(df d\omega)=n(df) \Theta^f_x(d\omega)\ee
where $n(df)$ is the It\^o excursion measure (on the set of positive excursions)
 and $\Theta^f_x$ is the law of the  Brownian snake started at $x$ with lifetime process $f$.


Let
\be{}\mathcal{W}=\cup_{t\geq 0} C([0,t],\mathbb{R}^d),\quad \zeta_w=t\text{  if  } w\in C([0,t],\mathbb{R}^d).\ee

\[
\int_{\mathbb{R}^d} \varphi(x) X_{t}(dx) =\int_{\mathcal{W}\times
R^d}\left[
\int_{0}^{\sigma(W)}\varphi(W_{s}(\zeta_{s}))L_{ds}^{t}\right]
\Pi(dw\ dx)
\]
$ (L^a_s)_{s\geq 0} $ local time of $\zeta$ at $a$  and
$\sigma(W)$ is the lifetime of the excursion and
 where $\Pi$ is a Poisson field with intensity%
\[
\pi(dw~dx)=\int_{\mathbb{R}^d}N_{y}(dw)\otimes\delta_{y}(dx)X_0(dy)
\]
and $N_{y}$ is the  It\^o excursion measure for the Brownian
snake from the trivial path  at $y$. The corresponding measures
on paths $\{W_u\}_{0\leq u\leq \zeta_s}$ gives the {\em historical
process}.

\subsection{Catalytic SBM}

The branching rate $
\gamma$ of SBM is assumed to be constant. A superprocess in which the branching rate  $\gamma$ is
 replaced by  a non-negative function or measure is called catalytic SBM. See
Dawson-Fleischmann \cite{DF-00}, \cite{DF-02} for an exposition.



\section{Local properties of critical branching systems}

\subsection{Support Properties of SBM}

We first derive some properties using the skew product representation and the lookdown process for the Fleming-Viot process.

\begin{theorem}
(Iscoe's Compact Support Property) Consider a super Brownian
motion in $\mathbb{R}^{d}$ with compactly supported initial
measure. \ Then at any fixed time, $t\geq0,$ the closed support of
$X_{t}$ is compact with probability one.
\end{theorem}

\begin{proof}
Using the above theorem, this follows immediately from the compact
support property for the Fleming-Viot process, Theorem \ref{Compact}.
\end{proof}

\begin{theorem}
(R. Tribe) Consider super-Brownian motion $\{X_t\}_{t\geq 0}$. \ Let $\tau:=\inf\{t:X_{t}%
(1)=0\}.\,\ $Then%
\[
\lim_{t\rightarrow\tau-}Z(t)=\delta_{\zeta_{1}(\tau)}\text{ a.s.}%
\]
\end{theorem}

\begin{proof}
Note that $\tau$ is predictable and
$(\zeta_{1}(\tau-),\zeta_{2}(\tau -),\dots)$ is exchangeable. It
is known that for the Feller branching
diffusion $\{x(t):t\geq0\}$ that%
\[
\int_{0}^{\tau-}\frac{1}{x(s)}ds=\infty\text{.}%
\]

From this it follows that each $\zeta_{i}$ has lookdowns to
$\zeta_{0}$ at times arbitrarily close to $\tau-$. The result then
follows from the H\"{o}lder continuity of the Brownian paths.
\end{proof}

\subsection{Brief review of sample path properties}

There is an extensive literature on the sample path properties of SBM $\{X_t\}_{t\geq 0}$.  We do not consider this in detail but briefly mention some basic properties.

\begin{itemize}

 \item In dimensions $d=1$,    $\
X_{t}(dx)=\tilde{X}_t(x)dx$
{has a jointly continuous density} which is given by a weak solution of stochastic partial differential equation (SPDE)
 \[
 {d\tilde X_{t}(x ) =\frac{1}{2}\Delta \tilde X_{t}(x )dt+\sqrt{\tilde X_{t}(x )}%
W(dt,dx)} \] {\[ W(dt,dx) =\text{white noise}\]} { (Konno-Shiga
(1988) \cite{KS-88}, Reimers (1989) \cite{R-89})}

{ {Strong uniqueness is an open problem.}}

\item In dimensions $d\geq 2$ $X_t$ is
a.s. a singular measure (D-Hochberg (1979) \cite{DH-79}). The Hausdorff measure properties of
the support of SBM were established by
Perkins (1989) \cite{P-89} for dimensions  $d\geq 3$ , and by  Le Gall and Perkins (1995) \cite{LP-95} in dimension  $d=2$:
There is a universal constant $c_0\in (0,\infty)$ such that
 \[  X_{t}(A)=c\cdot \phi_d -m(A\cap S(X_{t})),\,\,a.s.
\] where
{ \[ S(X_{t})=\text{closed support of }X_{t}
\]
} and {
\[
\phi_{\geq 3} (r)=r^{2}\log \log \frac{1}{r},\quad \phi_2(r)=
r^2\log^+(1/r)\log^+\log^+(1/r).\] }

\end{itemize}

For a comprehensive development of these and other sample path properties of super-Brownian
motion refer to Perkins (2002) \cite{P-02}.

\section{Spatial dynamical invariance principles}

SBM arose in a natural context as the limit of rescaled branching random walks, branching Brownian motions, interacting Feller CSBP.  There is a natural invariance principle here, namely, for a large class of migration mechanism in the domain of attraction of Brownian motion and branching mechanisms having finite second moments the limit is SBM.  However more surprising, it has been discovered that super-Brownian motion also arises in the scaling limit of systems that have no immediate branching interpretation.   Examples related to finite excursions of SBM arose in lattice trees where David Aldous (1993)
conjectured that if $d>8$ the rescaled tree of size $N$, $X_N$ should  converge to  (ISE)  and Derbez and Slade (1997) \cite{DS-97} proved this in sufficiently high dimensions. Examples also arise in interacting particle systems.  Durrett and Perkins (1999) \cite{DP-99} proved that a class of rescaled contact processes converge to SBM and Cox-Durrett-Perkins (2000) \cite {CDP-00} proved that a class of rescaled voter models converge to SBM.
Other examples are {\em voter model clusters} Bramson-Cox-Le Gall (2001) \cite{BCG-01}, and interacting diffusions Cox-Klenke (2003) \cite{CK-03}.   The proofs of many of these results
involved weak convergence of solutions of martingale problems in the same spirit as the proof of the convergence of branching Brownian motion to SBM.  The proofs in Derbez-Slade \cite{DS-97} and the proof due to van der Hofstad and Slade (2003) \cite{HS-03} that critical oriented percolation above 4+1 dimensions converges to SBM uses lace expansion methods to prove the convergence of the moment measures of the finite dimensional distributions.

\subsubsection{Scaling limits of long-range voter models in $d=1$}

We begin by describing a class of {\em long range voter models } on $\mathbb{Z}^1$.  Let
\be{}  S_N:=\left\{\frac{x}{\wt M_N\cdot N}:x\in\mathbb{Z}^1\right\}\ee and probability distribution on $S_N$ defined by
\bea{}  &&\\p_N(x)&&:=\frac{1}{2[\wt M_N\sqrt{N}]}\nonumber\\&&
\text{if  } x \text{  is one of the }2[\wt{M_N}\sqrt{N}]\text{ equally spaced points  in  } (-\frac{1}{\sqrt{N}},\frac{1}{\sqrt{N}})\nonumber\\&&
:=0,\; \text{otherwise}.\nonumber
\eea

Let   $\xi^N_t(x)=\xi_{Nt}(x\sqrt{N})$   denote the biased voter model with state space  $\{0,1\}^{{S}_N}$ with voting rate for a site with opinion $1$ is
$\gamma_\theta(N)=\wt M_N(N+\theta\sqrt{N}),\;\theta\geq 0$  and is $\wt M_N N$ for a site with opinion $0$ and  voting kernel
$p_N(x,y)=p_N(x-y)$.

Define

\be{}  X^N_t=\frac{1}{N} \sum_{x\in S_N}\xi^N_t(x)\delta_x.\ee

Mueller and Tribe  introduce the space
\be{} \mathcal{C}=\{f:\mathbb{R}\to [0,\infty) \text{  continuous with }|f(x)e^{\lambda |x|}\to 0\text{  as  }|x|\to\infty\; \forall\; \lambda <0\}.\ee

\beT{} (Mueller-Tribe \cite{MT-95} Theorem 2) Let $\wt M_N\equiv 2$ and let the approximate densities of $X^N_t$ be defined by
\be{}  u_N(t,x)=\frac{\sum_{x\sim y}\xi^N_t(y)}{\sum_{x\sim y}1} \quad \text{where   }   x\sim y\text{ iff }|x-y|\leq N^{-1/2}. \ee

Assuming that the the initial approximate densities converge in $\mathcal{C}$ to $u(0)$, then the approximate densities converge in distribution to a continuous $\mathcal{C}$-valued process $u(t)$ which is a solution of the of the stochastic partial differential equation (SPDE) \be{}\frac{\partial u}{\partial
t}=\frac{1}{6}\frac{\partial^2}{\partial x^2}+2\theta_v u(1-u)+\sqrt{4u(1-u)} \dot
W\ee
with initial condition $u(0)$ and where $\dot W$ is space-time white noise.
\end{theorem}
\begin{remark}
This SPDE is a Fisher-KPP equation driven by Fisher-Wright noise.
\end{remark}

In contrast to this, (for the unbiased voter model) the special one-dimensional case of the Cox-Durrett-Perkins long-range voter model in which $\wt M_N\to\infty$ to be presented in the next subsection has as scaling limit SBM. The essential difference is  the CDP scaling is such that the measure is more sparsely distributed on the voters and the  local density of sites occupied by opinion 1 goes to 0.
\beT{} (CDP \cite{CDP-00})  Assume that $\wt M_N\to \infty$ as $N\to\infty$.  Let  $P_N$ denote the law of $X^N_t$ on $D_{M_F(\mathbb{R}^1)}([0,\infty))$.
Assume that  $\sum \xi^N_0(x) <\infty$ and  $X^N_0\to X_0$ in $M_F(\mathbb{R}^1)$ as $N\to\infty$.  Then
\be{}
P_N\Rightarrow P^{1,2,1/3}_{X_0} \ee
as $N\to\infty$ where $P^{d,2,1/3}$ is the law of super-Brownian motion with branching rate $2$ and $\sigma^2=\frac{1}{3}$
where $P^{\gamma,\sigma^2}_{X_0}$ denotes the law of SBM in $\mathbb{R}^d$ with initial measure $X_0\in M_F(\mathbb{R}^d)$ and log-Laplace equation
\be{} \frac{\partial u}{\partial t}= \frac{\sigma^2}{2}\Delta u-\frac{\gamma}{2}u^2.\ee
\end{theorem}

\subsubsection{The rescaled voter model in $\mathbb{Z}^d,\; d\geq 3$}

Assume that the kernel $p(x-y)$ is irreducible and symmetric, $p(0)=0$ and  $\sum_{x\in\mathbb {Z}^d} x^ix^jp(x)=\delta_{ij}\sigma^2$.
Let $\xi_t(x)$ denote the voter model on $\mathbb{Z}^d$ with kernel $p(\cdot)$.

Let  $S_N= \frac{\mathbb{Z}^d}{\sqrt{N}}$,  $\xi^N_t(x) :=\xi_{Nt}(x\sqrt{N}),\;x\in S_N$, and consider the measure-valued process on $\mathbb{R}^d$ defined by
\be{}  X^N_t=\frac{1}{N} \sum_{x\in S_N}\xi^N_t(x)\delta_x.\ee

\beT{FK} (\cite{CDP-00}, Theorem 1.2)
Let $P_N$ denote the law of $X^N_t$ on $D_{M_F(\mathbb{R}^d)}([0,\infty))$ with $d\geq 3$.  Assume that  $\sum \xi^N_0(x) <\infty$ and  $X^N_0\Rightarrow X_0$ in $M_F(\mathbb{R}^1)$ as $N\to\infty$
Then  \be{} P_N\Rightarrow P^{d,2\gamma_e,\sigma^2}_{X_0} \text{   as   } N\to\infty\ee
where $\gamma_e$ is the escape probability  \be{rp}\gamma_e = P_0( p(\cdot)\text{-RW never returns to }0).
\ee
\end{theorem}
\begin{remark}
The case $d=2$ is more subtle.  In this case there is a logarithmic correction and one considers the sequence of measure-valued processes
\be{} X^N_t=\frac{\log N}{N} \sum_{x\in S_N}\xi^N_t(x)\delta_x.\ee
With this scaling Cox-Durrett-Perkins \cite{CDP-00} Theorem 1.2 prove the weak convergence to $P^{2\gamma_e,\sigma^2}_{X_0}$.
\end{remark}

\subsubsection{Rescaled long-range voter models}
To describe the results for the  long range voter model in higher dimensions we let  $M_N$ be a sequence of positive constants with $ M_N\to\infty$ as $N\to\infty$ and consider the voter model on $S_N=\{\frac{x}{M_N\sqrt{N}}:x\in\mathbb{Z}^d\}$  with probability kernels  $p_N(\cdot)$ defined as follows:.

For each $N$ let $W_N$ denote a random variable with values in  $\frac{(\mathbb{Z}^d\backslash \{0\})}{M_N}$ and with  uniform distributed on
$(\frac{(\mathbb{Z}^d\backslash \{0\})}{M_N})\cap I$ where $I=[-1,1]^d$

Then
\begin{itemize}
\item $W_N$ and $-W_N$ have the same distribution,

\item $\lim_{N\to\infty} E(W^i_NW^j_N)=\delta_{ij}\sigma^2$ with $\sigma^2= \frac{1}{3},$ and

\item  the family  $\{|W_N|^2\}_{N\in\N}$ is uniformly integrable.
\end{itemize}

Then define the kernel
 \be{}
p_N(x):=P(\frac{W_N}{\sqrt{N}}=x) \ee

Let  $\xi^N_t(x) =\xi_{Nt}(x\sqrt{N})$ denote the  voter
model on $\{0,1\}^{\mathbb{S}_N}$ with rate $N$ and  voting kernel
$p_N(x,y)=p_N(x-y)$ and let $P_N$ be the law of the measure-valued process \be{}
X^N_t=\frac{1}{N}\sum_{x\in\mathbb{S}_N}\xi^N_t(x)\delta_x\ee

 \beT{} (Cox, Durrett, Perkins (2000) \cite{CDP-00})\\
 Assume  that
\be{} \sum_x \xi^N_0(x)<\infty,\ee
\be{}   X^N_0\Rightarrow X_0 \text{  in  } \mathcal{M}_F(\mathbb{R}^d) \text{  as  }
N\to\infty\ee
and
\be{}\left\{%
\begin{array}{ll}
    M_N/\sqrt{N}\to\infty, & \hbox{ in } d=1 \\
    M^2_N/(\log N)\to\infty , & \hbox{ in }d=2 \\
    M_N\to \infty, & \hbox{ in }d\geq 3, \\
\end{array}.
\right.\ee

Then  $P_N\Rightarrow P^{d,2,1/3}_{X_0}$ as $N\to\infty$.
\end{theorem}

\subsubsection{Interacting diffusions - convergence to SBM}

Assume that the kernel $p(x-y)$ on $S=\mathbb{Z}^d$ is irreducible and symmetric, $p(0)=0$ and  $\sum_{x\in\mathbb {Z}^d} x^ix^jp(x)=\delta_{ij}\sigma^2$.

Now consider the system of interacting diffusions on $\mathbb{Z}^d$:
\begin{align*}
&{d{ X_t(x)}}    = [\sum q_{x,y}({X}_t(y))-{ X}_t(x)]dt  +
\sqrt{{ g(X_t(x))}} \;dW_t(x),
\quad x\in \mathbb{Z}^d\\&\text{where   } q_{x,y}=p(x-y)\text{  and  }\{(W_t(x))_{t\geq 0}\}_{x\in \mathbb{Z}^d} \text{ are  independent Wiener
processes}.
\end{align*}

Let $P_N$ denote the probability laws on $C_{M_F(\mathbb{R}^d)}([0,\infty))$ of the measure-valued processes
\be{}  X^N_t:=\frac{1}{N} \sum_{i\in \mathbb{Z}^d} X_{Nt}(i)\delta_{i/\sqrt{N}}\in M_F(\frac{\mathbb{Z}^d}{\sqrt{N}}).\ee

\beT{} (Cox and Klenke (2003) \cite{CK-03}, Theorem 1)
Assume that $X^N_0\Rightarrow X_0$ in $M_F(\mathbb{R}^d)$ with $d\geq 3$ and that
\be{} g(x)=\kappa x(1-{x})^+,\qquad x\geq 0,\; \kappa >0,\ee

Then \be{} P^N \Rightarrow  P^{d,\gamma,\sigma^2}\ee
 where
\be{}\gamma=\lim_{\theta\to 0}\frac{1}{\theta}\int g(x(0))\nu_\theta(dx)=\gamma_e\cdot \kappa\ee
 $\nu_\theta$ is the invariant measure  for the interacting Wright-Fisher system on  $\mathbb{Z}^d$ with intensity $\theta$ and  $\gamma_e$ is given by (\ref{ncp}).

\end{theorem}

Cox and Klenke \cite{CK-03} conjectured that the same result holds for the more general case with

\begin{itemize}
\item  $g:[0,\infty)\to [0,\infty)\qquad$ is locally Lipschitz continuous
\item $g^{-1}((0,\infty))=(0,b)\qquad $for some $b\in (0,\infty]$,
\item $g(z)\leq C(1+z^2)\qquad$ for some $C<\infty.$

\end{itemize}
In this case the conjectured limit is $P^{d,\gamma,\sigma^2}$ with
\be{}  \gamma =\lim_{\theta\downarrow 0}\frac{1}{\theta} \int g(x(0))\nu_\theta(dx)\ee
where $\nu_\theta$ is the unique stationary measure on $\mathbb{Z}^d$ with intensity $\theta$.

\begin{remark} Note that in Cox-Klenke the limiting branching rate is obtained from the derivative at 0  of the $\mathcal{F}g$. In this thinning out we get super-Brownian excursions. (cf. Bramson, Cox Le Gall).
\end{remark}

\begin{remark}
Cox and Klenke also describe the long range case and note that the invariant measure involved then can be described by the mean-field equation.

\end{remark}

\subsection{Methods of Proof}

The proofs of these results involve three main steps
\begin{itemize}
\item formulation of the sequence of measure-valued processes in terms of martingale problems
\item establishing tightness of the probability laws on $D_{M_F(\mathbb{R}^d)}$
\item verifying any any limit point satisfies the SBM martingale problem.
\end{itemize}

\subsubsection{Step 1: Reformulation as a martingale problem}

Starting from the set of stochastic equations (\ref{vmsie}), we can characterize the rescaled voter model as the solution of the stochastic integral equations
\be{SIE} \xi^N_t(x)=\xi^N_0(x)+\sum_y \int_0^t [\xi^N_{s-}(y)-\xi^N_{s-}(x)]s\Lambda^N_s(x,y),\quad x\in S_N,\ee
 where $\{\Lambda^N_t(x,y):x,y\in S_N\}$ is a system of independent Poisson processes  with rates $Np_N(y-x)$.
 Then  the measure-valued processes $X^N_t$ can be characterized by a martingale problem as follows.

\beT{T2.2}(Cox, Durrett, Perkins (2000) Theorem 2.2)\\ Let
$\phi\in C^{1,3}_b([0,\infty)\times\mathbb{R}^d)$. Then \be{CDP1}
X^N_t(\phi) =X^N_0(\phi)+\int_0^t X^N_s(
\phi_1(s)+\mathcal{A}_N\phi(s))ds+M^N_t(\phi),\ee where
$\phi_1(s)$ denotes the partial derivative $\frac{\partial \phi(s,\cdot)}{\partial s}$ and
\be{} \mathcal{A}_N\phi(s,x)=N\sum_y p_N(x-y)(\phi(s,y)-\phi(s,x)),\ee
and

\be{}M^N_t(\phi)=\frac{1}{N}\sum_x\sum_y\int_0^t\phi(s,x)(\xi_{s-}(y)-\xi_{s-}(x))\hat
{\Lambda}_s(x,y),\quad 0\leq t\leq T, \ee
$\hat{\Lambda}^N_t(x,y)=\Lambda^N_t(x,y)-Np_N(x-y)t$ is a cadlag,
square integrable $(\mathcal{F}_t)$-martingale.

The predictable
square function is given by \be{CDP2}  \langle M^N(\phi)\rangle_t
=\int_0^t[2X^N_s(\phi^2(s)V_{N}(s))+\ve^N_s(\phi)]ds,\ee where
\bea{9.85a} \quad V_N(t,x)&&=\sum_y p_N(y-x)1\{\xi_t(y)=0\} =\text{ density of vacant sites near }x,\\
&&=\quad \lim_{\lambda\to\infty}\sum_{y}p_N(y-x)e^{-\lambda X^N_t(y)}\nonumber\eea
and
$\ve^N_s(\phi)$ satisfies \be{9.85b} E\left(\sup_{s\leq
T}|\ve^N_s(\phi)|^2\right)\to 0,\text{  as  }N\to\infty,\text{ for
any }T>0.\ee

(iii) For any $T>0$, \be{9.86} E\int_0^T
X^N_s(|\mathcal{A}_N\phi(s)-\frac{\sigma^2}{2}\Delta\phi(s)|)ds\to 0\text{
as  }N\to\infty.\ee
\end{theorem}

\subsubsection{Steps 2 and 3: The Cox-Durrett-Perkins Criteria for weak convergence to SBM}

Cox,Durrett and Perkins formulate a general set of additional conditions for a sequence of $M_F(\mathbb{R}^d)$-valued
martingale problems satisfying the conclusions of Theorem \ref{T2.2}  which if satisfied imply that
the solutions to these martingale problems converge weakly to SBM.

Assumptions:

(I1) There is a finite $\gamma >0$ such that, for all $\phi\in
C^\infty_0(\mathbb{R}^d)$ and $T>0$, as $N\to\infty$, \be{}
E\left[\left(\int_0^TX^N_s(\{V_{N}(s)-\gamma\}\phi^2)ds\right)^2\right]\to
0.\ee

(I2)  For all $T>0$ there exists a finite $C_T$ such that
$\lim_{T\downarrow 0}C_T=0$ and for all $N$, \be{}
\int_0^TE^{\xi^N_0}[X^N_s(V_{N}(s))]ds\leq C_T
X^N_0(\mathbf(1)).\ee

(I3) There is a $\theta\in (0,1]$ and  a finite $C(\ve,T,K)$ such
that for all $N\in\mathbb{N}$ all cutoffs $0<\ve$, $K<\infty$ and
all pairs of times $\ve\leq s\leq t\leq T$, we have

\be{}  \sup\left\{E\left[\left(\int_s^t
X^N_r(V_{N}(r))dr\right)^2\right]:X^N_0(\mathbf{1})\leq
K\right\}\leq C(\ve,T,K)|t-s|^{1+\theta}.\ee

\beT{CDP-T3.5} Assume (I1)-(I3).  Then $P_N\Rightarrow
P^{2\gamma,\sigma^2}_{X_0}$.

\end{theorem}
\begin{proof} We refer to \cite{CDP-00} for the details but outline the main ideas here.
The proof involves two parts.  The first is the proof of the
tightness of the processes on $D([0,\infty),M_F(\mathbb{R}^d))$
and all limit points are supported by
$C([0,\infty),M_F(\mathbb{R}^d))$.  This follows the same general lines as in our examples - we omit the details.

The second part shows that every limit point $P$ of $P_N$ satisfies the SBM
martingale problem.
By Skorohod's theorem, given $P_{N_k}\Rightarrow  P$ we may assume there is  $X, X^{N_k}$ defined on a probability space $(\Omega,\mathcal{F},P)$ such that \be{3.10} X^{N_k}\to X,
$P$\text{-a.s.
 in  } D([0,\infty),M_F(\mathbb{R}^d)).\ee

Then a standard argument  shows that for
 $T>0$,

 \be{3.11} \lim_{k\to\infty}\sup_{t\leq
 T}\left|\int_0^tX^{N_k}_s(\mathcal{A}_{N_k}\phi)ds-\int_0^tX_s(\sigma^2\Delta\phi/2)ds\right|=0,\;P\text{
 -a.s.}\ee

Let  \be{}M_t(\phi):=X_t(\phi)-X_0(\phi)-\int_0^t
X_s(\sigma^2\Delta\phi/2)ds.\ee

  Then  using (\ref{9.86}), some calculus and the
fact that  $M_{\cdot}(\phi)$ is continuous to derive uniform
convergence on compacts from convergence on $D$ one can verify
   \be{3.12}
\lim_{k\to\infty}\sup_{t\leq
T}|M^{N_k}_t(\phi)-M_t(\phi)|=0, \;P-\text{  a.s.}.\ee

A stochastic calculus computation starting with the system of stochastic integral equations (\ref{SIE}) yields
\be{} E\langle
M^N(\phi)\rangle^2_t=E\left(\int_0^t\left[2X^N_s(\phi^2V_{N,s})+\ve^N_s(\phi)\right]ds\right)^2\ee
where $|\ve^N_s(\phi)|\leq C_{\phi}X^N_s(\mathbf{1})/\sqrt{N}$.
Then from (I1), Theorem \ref{T2.2}(ii) and the elementary estimate
\be{} E(\sup_{s\leq T} (X^N_s(1))^2)\leq C(X^N_0(1)+(X^N_0(1))^2)\ee
 it follows that that
for $T>0$, $\sup_N E(\langle M^N(\phi)\rangle^2_T)<\infty$.

Then using Burkholder's inequality (see appendix Theorem \ref{Burkholder}) and noting that  $|\Delta
M^N(\phi)(t)|\leq \| \phi\|_\infty(N)^{-1}$, we have

\be{3.13}  \sup_N E\left (\sup_{t\leq
T}|M^N_t(\phi)|^4\right)<\infty.\ee

Fix $0\leq t_1<t_2<\dots <t_n\leq s<t$ and test functions
$h_i:M_F(\mathbb{R}^d)\to\mathbb{R}$ that are bounded and
continuous for $1\leq i\leq n$. Now Theorem \ref{T2.2},
(\ref{3.10}),(\ref{3.12}),(\ref{3.13}) and dominated convergence
imply that
\bean{} E\left((M_t(\phi)-M_s(\phi))\prod_1^n
h_i(X^{N_k}_{t_i})\right)=0.\eean

 Therefore under $P$, $M_t(\phi)$ is a continuous
$\mathcal{F}^X_t$-martingale where $\mathcal{F}^X_t$ is a
canonical right-continuous filtration generated by $X$. Also
(\ref{3.10}),(\ref{3.12}),(\ref{3.13})  imply that
\bean{}
&&E\left(\left(M_t(\phi)^2-M_s(\phi)^2-\int_s^tX_r(2\gamma\phi^2)dr\right)\prod_1^n
h_i(X_{t_i})\right)\\&&
\lim_{k\to\infty}E\left(\left(M^{N_k}_t(\phi)^2-M^{N_k}_s(\phi)^2-\int_s^t
X^{N_k}_r(2\gamma\phi^2)dr\right)\prod_1^n
h_i(X^{N_k}_{t_i})\right) \eean (\ref{9.85a}), (\ref{9.85b})  and (I1) show that
the above equals \be{} \lim_{k\to\infty}
E\left(\left(M^{N_k}_t(\phi)^2-M^{N_k}_s(\phi)^2-(\langle
M^{N_k}(\phi)\rangle_t - \langle
M^{N_k}(\phi)\rangle_s)\right)\prod_1^n
h_i(X^{N_k}_{t_i})\right),\ee which  is 0 by Theorem \ref{T2.2}.
This shows that $\langle M(\phi)\rangle_t=\int_0^t
X_s(2\gamma\phi^2)ds$ for all $t\geq 0$, P-a.s.  Therefore $P$
satisfies the SBM martingale problem
$(MP)^{2\gamma,\sigma^2}_{X_0}$.  Therefore the law of $X$ equals
$P^{2\gamma,\sigma^2}_{X_0}$.

\end{proof}

\subsection{Proof of the invariance principle for interacting Wright-Fisher diffusions}

In this case we have
\be{} M^N_t(\varphi):= X^N_t(\varphi)-X^N_0(\varphi)-\int_0^t X^N_s(\mathcal{A}_N\varphi)ds\ee
is a continuous square integrable martingale with quadratic variation process
\be{} \langle M^N(\varphi\rangle_t=\int_0^t \Gamma^N_s(\varphi^2)ds,\ee
where
\bea{}\qquad  \Gamma^N_s:=&& \frac{1}{N}\sum_{x\in\frac{\mathbb{Z}^d}{\sqrt{N}}}g(NX^N_s(\{x\}))\delta_x\\
&& =\frac{1}{N}\sum_{i\in\mathbb{Z}^d} g(NX_{sN}(i))\delta{i/\sqrt{N}}.\nonumber\eea

Therefore the  key step  Cox-Klenke \cite{CK-03} in verifying the CDG conditions  is to prove that

\be{} \ve^{N,\gamma}_{K,\phi}(t):=\sup \{ E[|(\Gamma^N_t-\gamma X^N_t)(\phi^2)|:X^N_0(1)\leq K]\} \to 0 \text{  as  } N\to\infty.\ee

\subsubsection{Verification via  Duality calculations}

 The verification is achieved using  a duality calculation which involves the coalescing random walk with two particles. The difference between the voter and Wright-Fisher cases in that in the former colliding particles coalesce instantaneously whereas for Wright-Fisher they coalesce {\em with delay}. We now sketch the proof in the Wright-Fisher case.

In the case when $g$ has Wright-Fisher form $g(x)=\kappa x(1-x)^+$, $\gamma_e$ given by (\ref{ncp}) this becomes
\be{} \ve^{N,\gamma}_{K,\varphi}(t)=\kappa\sup   \left\{\left| E\left[\sum_{x\in\mathbb{Z}^d/\sqrt{N}}\left((1-{\gamma_e})X^N_t(\{x\})-NX^N_t(\{x\})^2\right)\varphi^2(x)\right]\right|:X^N_0(1)\leq K\right\}.\ee
Let  $p_t$ be the transition function for the random walk, $Z$.

The dual involves the coalescing random walk when the particles coalesce at rate $\kappa$  when they are at the same site, namely,
\be{} E[X_t(z^1)X_t(z^2)]=E^{(z^1,z^2)}[X_0(Z^1_t)X_0(Z^2_t)].\ee
 Consider two random walks $Z^1$ and $Z^2$ and let  $A_t$ be the event that they have coalesced by time $t$ and $A:=\cup_{t\geq 0}A_t$.
Therefore for $\varphi\in C^2_c(\mathbb{R}^d)$,
\bea{} &&\\ NE[X^N_t(\{i/\sqrt{N}\})^2]&&=\frac{1}{N} E[X_{Nt}(i)^2]\nonumber\\&&
=\frac{1}{N}E^i[X_0(Z^1_{Nt});A_{Nt}]+\frac{1}{N}E^i[X_0(Z^1_{Nt})X_0(Z^2_{Nt}):A_{Nt}^c].\nonumber\eea
By the CLT, there exists a constant $C<\infty$ such that
\be{clt0} p_t(0,j)\leq(\frac{1}{t^{d/2}\wedge 1})\cdot C,\quad j\in\mathbb{Z}^d,\;t>0.\ee

Then since $\frac{1}{N}\sum_j X^N_0(j)\leq K$, $X^N_0(j)\leq NK$,
\bea{} &&\sum_i\frac{1}{N}E^i[X_0(Z^1_{Nt})X_0(Z^2_{Nt})]\varphi(i/\sqrt{N})^2\\&&\leq \frac{1}{N} \sup_{j,k\in\mathbb{Z}^d} p_{Nt}(i,j)p_{Nt}(i,k)X^N_0(j)X^N_0(k)\nonumber\\&&\leq \frac{C}{N}t^{-d}N^{-d}N^2K^2\leq.\nonumber
 CC_{\varphi}K^2t^{-d}N^{1-d/2},\nonumber\eea
where $C_{\varphi}$ depends only on $\varphi$.

Hence by dominated convergence it suffices to show that

\bea{CK2.19} &&\sup_i \wt\ve^{N,\gamma,i}(t)\\&&:=
N^{\frac{d}{2}-1}\sup_i\left\{\left| E^i[X_0(Z^1_{Nt});A^c_{Nt}]-{\gamma}_eE^i[X_0(Z^1_{Nt})]\right|:X_0^N(1)\leq K\right\} \Nto 0\nonumber\eea

The intuitive reason for this is that
\[ \lim_{T\to\infty} P[A^c_T]= {\gamma}_e\]
and the distribution of $Z^1_{Nt}$ and the conditional distribution of $Z^1_{Nt}$ given $A^c_{Nt}$ are close.

To make this precise, let $\delta>0$ and fix $T_0>0$ be such that
\be{CK2.21}|P(A^c_T)-\gamma_e|\leq \frac{(T/2)^{d/2}}{C}\delta, \text{  for all }T\geq T_0\ee
with $C$ as in (\ref{clt0}).

We next obtain an upper bound on the on the probability that$Z^1$ and $Z^2$ coalesce between times $T$ and $tN$ and end at time $tN$ at $j$.

Noting that
\be{} P^i[A_{Nt}\cap A^c_T\cap\{Z^1_{Nt}=j\}]= 1- E^i[e^{-\kappa \int_{T}^{Nt} 1(Z^1_r=Z^2_r)dr}1(Z^1_{Nt}=j)]\ee
and using Jensen's inequality we have
\be{} P^i[A_{Nt}\cap A^c_T\cap\{Z^1_{Nt}=j\}]\leq \kappa E^i[\int_{T}^{Nt} 1(Z^1_r=Z^2_r)dr\cdot1(Z^1_{Nt}=j)].\ee
Therefore
\bea{CK2.22}
&&P^i[A_{Nt}\cap A^c_T\cap\{Z^1_{Nt}=j\}]\\&&
\leq \kappa\int_T^{Nt}dr\sum_{k\in\mathbb{Z}^d}p_r(i,k)p_r(i,k)p_{Nt-r}(k,j)\nonumber\\&&
\leq \kappa p_{Nt}(i,j)\int_T^{Nt}dr \sup_k p_r(0,k)\nonumber\\&&
\leq\kappa C^2t^{-d/2}N^{-d/2}\int_T^{Nt}r^{-d/2}dr\nonumber\\&&
\leq\frac{2\kappa C^2t^{-d/2}}{d-2}T^{1-d/2}N^{-d/2}.\eea
Choosing $T_0$ large enough we can also assume that
\be{CK2.23} \sup_{j\in\mathbb{Z}^d} P^i[A_{Nt}\cap A^c_{T}\cap \{Z^1_{Nt}=j\}]\leq \delta N^{-d/2},\quad T\geq T_0.\ee

Let $R>0$ be such that
\be{CK2.24} P^i[|Z^1_{T_0}|>R]< \frac{\delta}{(1+\gamma_e)(1+(2/t)^{d/2}C)}.\ee
Using (\ref{clt0}) and the Markov property at time $T_0$, we get for $N\geq 2T_0/t$
\be{CK2.25}  (1+\gamma_e)P^i[|Z^1_{T_0}|>R;\, Z^1_{Nt}=j]\leq \delta N^{-d/2}.\ee

Using the CLT again there exists $N_0\geq 2T_0/t$ such that for all $N_0\geq N_)$ and $|k|<R$
\be{CK2.26} |p_{Nt-T_0}(k,j)-p_{Nt-T_0}(0,j)|< \frac{\delta}{1+\gamma_e}N^{d/2}.\ee
Combining (\ref{CK2.23}), (\ref{CK2.25}), (\ref{CK2.26}), (\ref{clt0}), (\ref{CK2.24}) and (\ref{CK2.21}) and using the Markov property we get, for $N\geq N_0$.

\bea{CK2.27} &&\left | P^i[Z^1_{Nt}=j;\,A^c_{Nt}]-\gamma_e P^i[Z^1_{Nt}=j]\right |\\
&& \leq \left |P^i[Z^1_{Nt}=j;\, A^c_{T_0}]-\gamma_eP^i[Z^1_{Nt}=j]\right|+\delta N^{-d/2}\nonumber\\&&\leq
\sum_{|k|<R}\left |P^i[Z^1_{Nt}=j;Z^1_{T_0}=k;\, A^c_{T_0}]-\gamma_eP^i[Z^1_{Nt}=j;\, Z^1_{T_0}=k]\right|+2\delta N^{-d/2}
\nonumber
\\&& \leq  \sum_{|k|<R}p_{Nt-T_0}(k,j)\left|\left(P^i[Z^1_{T_0}=k; A^c_{T_0}]-\gamma_e P^i[Z^1_{T_0}=k]\right)\right |+2\delta N^{-d/2}\nonumber\\&&\leq  \left| P^i[|Z^1_{T_0}|<R; A^c_{T_0}]-\gamma_e  P^i[|Z^1_{T_0}|<R]\right |\cdot p_{Nt-T_0}(0,j) +3\delta N^{-d/2}\nonumber\\&&\leq \left|P^i[A^c_{T_0}]-\gamma_e\right|\cdot C(2/t)^{d/2}N^{-{d/2}}+4\delta N^{-d/2}\nonumber\\&&
\leq 5\delta N^{-d/2}.\eea

Recall that

\be{}\wt\ve^{N,\gamma_e,i}_K(t)= N^{d/2-1}\sup \sum_{j\in\mathbb{Z}^d}X_0(j)\left |[P^i(Z^1_{Nt}=j;A^c_{Nt})]-\gamma_e [ P^i(Z^1_{Nt}=j)]\right|\ee

Since the estimate holds for all $j$, and $\sum X_0(j)\leq NK$,
\be{} \limsup_{N\to\infty} \wt\ve^{N,\gamma_e,i}_K(t)\leq 5K\delta.\ee
Since $\delta >0$ was arbitrary, (\ref{CK2.19}) follows.

\subsection{Applications of the  SBM Invariance Principle}

\subsubsection{Critical parameter of the contact process} Durrett and Perkins \cite{DP-99} used this to obtain sharp asymptotics for the critical parameter of the long-range contact process which improve upon the results of Bramson, Durrett and Swindle (1989) \cite{BDS-89}.

\subsubsection{Extinction time for a voter model cluster}

Under the assumptions of Theorem \ref{FK}
Cox-Perkins (2004) \cite{CP-04}  prove that
\be{} \lim_{N\to\infty} P(\tau^N_0>t) = P(\tau_0>t)=1-\exp(-\frac{X_0(1)}{t\gamma_d}),\quad t>0\ee where
\be{} \tau_0^N=\inf\{ t>0:X^N_t(1)=0\},\quad \tau_0=\inf\{ t>0:X_t(1)=0\}.\ee

Note that this does not immediately follow from the invariance principle since $\tau_0$ is not a continuous function on $D_{M_F(\mathbb{R}^d)}([0,\infty)$.

\subsubsection{Lotka-Volterra models}

In a series of papers Cox, Durrett and Perkins obtain deep results on the region of survival and  coexistence for the Lotka-Volterra model of competing species - see section \ref{s.CCS}.  Recently,

\subsubsection{Voter model clusters}

Bramson, Cox and LeGall \cite{BCG-01} showed that the rescaled voter model cluster on $\mathbb{Z}^d$ in dimensions $d\geq 2$ conditioned on non-extinction at time $t$ converges to the canonical measure of SBM and that the rescaled support converges (wrt Hausdorff metric) to the support of the super-Brownian excursion.

\subsection{Spatial Lotka-Volterra models}

Neuhauser and Pacala \cite{NP-99} introduced a  stochastic spatial  model for the competition between two species based on the Lotka-Volterra equations.  Theirmodel takes into account the local competition between species but allows for spatial segregation of the species.
The model is an interacting particle model on $\mathbb{Z}^d$ with state space $\{0,1\}^{\mathbb{Z}^d}$ and dynamics given by a perturbed voter model with rates dependent at a site on the
density of the two types in a neighbourhood of the site. This models
the effects of short range spatial dispersion and demographic
stochasticity as well as the role of
{\em interspecific and intraspecific competition}.

Neuhauser-Pacala (1999) \cite{NP-99} prove  that the local competitive
interactions reduce the parameter region where coexistence occurs
in the classical mean-field (deterministic) model and spatial segregation of the  species in parts of the parameter region where the classical model predicts coexistence.

Cox and Perkins (2007), (2008) and  Cox, Durrett and  Perkins use
a {super-Brownian invariance principle} for the perturbed voter model to obtain more precise information on the coexistence region
for $d\geq 3$ and  new results on the survival region for $d\geq 2$.

\subsection{Application of SBM  to a spatial  epidemic model}

The SIR epidemic model consider a  population with $S$ susceptible, $I$
infected and $R$  removed individuals.
Recall that the deterministic model is given by the system of ODE
\begin{align*} \frac{dS}{dt}&=-\beta SI,\quad
\frac{dI}{dt}=\beta SI-\gamma I, \quad \frac{dR}{dt}=\gamma I
\end{align*}

The Reed-Frost stochastic model is defined as follows:
\be{} S_{t+1}\sim
\text{Bin}(S_t,(1-p)^{I_t}),\quad
 {1-p=e^{-\lambda/n}},\quad S_0=N\ee

Martin-L\"of (1998) proved that this has the {\em critical threshold} $\lambda =1$. Moreover if $I(0)=N^\alpha$, then the scaling limit as $N\to\infty$ has a phase transition at  $\alpha=1/3$ (Martin-L\"of (1998) \cite{ML-98}, Aldous (1997) \cite{A-97}).

We now consider a spatial analogue of the Reed-Frost model due to Lalley (2008) \cite{L-08}.
At each site on $\mathbb{Z}^d$ there are $N$ individuals (types SIR). Infected models remain infected one unit of time and then recover and become immune. At time $t=0,1,2,\dots,$ for each pair $(i_x,y_y)$ (infected at $x$ and susceptible at $y$ the susceptible individual at $y$ becomes infected with probability $p_N(x,y)$.  Assume that $p_N(.,.)$ is nearest neighbour simple random walk.  Scale the village size $N$ so that the expected number of infections by a contagious individual in a healthy population is $1$ so that the epidemic is critical,
that is
\be{}  p_N(x,y)= \frac{1}{(2d+1)N}\qquad \text{ if  } |x-y|=1,\qquad = 0\text{  otherwise}.\ee

\beT{} (Lalley (2008) \cite{L-08})
Let  $Y^N_t(x)$ be the number infected at time $t$ at site $x$ in a critical SIR epidemic with village size $N$ and initial configuration $Y^N_0(x)$.  Fix $\alpha>0$ and let  $X^N(t,x)$ be the renormalized particle density function process obtained by linear interpolation in $x$ from the valued
\be{}  X^N(t,x)=\frac{Y^N_{[N^\alpha t]}(N^{\alpha/2}x)}{N^{\alpha/2}}\qquad\text{for   } x\in \mathcal{Z}/N^{\alpha/2}.\ee
Assume that there is a compact interval $J$ such that the initial particle density  functions $X^N(0,x)$ have support in $J$ and assume that $X^N(0,x)$ converge in $C_b(\mathbb{R})$ to a function $X(0,x)$.  Then as $N\to\infty$
\be{}X^N(t,x)\Rightarrow X(t,x)\ee
where $X(t,x)$ is the density of a SBM with initial density $X(0,x)$ and killing rate $\theta$ where
\be{}\begin{split}& \theta(x,t)=0\qquad \text{  if  }\alpha <\frac{2}{5}\\
& \theta(x,t)=\int_0^tX(x,s)ds\qquad \text{  if  }\alpha =\frac{2}{5}\end{split}
\ee
\end{theorem}

The idea is that with this scaling up to the extinction time there should be about $O(N^{\alpha/2})$ particles per site and that since the extinction time is $O(N^\alpha)$ the total attrition rate (removed individuals) per generation is $O(N^{5\alpha/2})$. Therefore if $\alpha=2/5$ then a non-trivial proportion of the population is removed.

Lalley and Zheng \cite{LZ-09} have obtained similar results in dimensions 2 and 3 using the absolute continuity of super-Brownian motion  local time in these dimensions. Lalley, Perkins and Zheng (2009) establish the existence of a phase transition in 2 and 3 dimensions.

\chapter{Spatial systems in large space and time scales}

In this chapter we consider critical spatial branching systems and interacting neutral Fleming-Viot processes in large space and time scales. The behaviour of these systems is determined by potential theoretic properties of the migration process such as transience or recurrence. We begin with a brief review of some basic notions.

\section{Migration processes on Abelian groups}

In this section we give a brief review of the basic notions of
random walks and L\'evy processes on groups on abelian groups following \cite{DGW-04}.

Let $S$ be a locally compact (additive) Abelian group with
countable base and with Haar measure $\rho.$  A discrete time
random walk, $\{W_{n}\}_{n\in \mathbb{Z}_+}$, is prescribed by a
transition function
\[
P(x,dy)   :=P(W_{n+1}\in dy|W_{n}=x)=p(d(y-x))
\]
where $p$ is a probability measure on S.  The corresponding
$k$-step transition function is \[ P^{k}(x,dy)  :=P(W_{n+k}\in
dy|W_{n}=x).
\]

A continuous time random walk $\{W_{t}:t\geq0\}$ with jump rate
$1$ is then defined by the transition function
\[
P_{t}(x,dy):=P_{x}(W_{t}\in dy), \;\; t\geq0,
\]
\[
P_{t}(x,dy)=\sum_{k=0}^{\infty}\frac{e^{-t}t^{k}}{k!}P^{k}(x,dy).
\]

A natural generalization of continuous time random walks is the
notion of L\'evy process.
\begin{definition}
A $S$-valued process $\{X_t:t\geq 0\}$ is a \textit{L\'evy
process} if it is stochastically continuous and has stationary and
independent increments.
\end{definition}

 We associate to a  L\'evy process a semigroup $\{T_t:t\geq0\}$
on $\mathcal{B}_c(S)$, the space of bounded measurable functions
on $S$ with compact support, as follows:
\[
T_{t}\varphi(x)=E_{x}(\varphi(X_{t})),
\]
The Green potential of $X$ is the operator
\[  G \varphi = \int_0^\infty
T_t \varphi dt,\quad \varphi\in \mathcal{B}_c(S).\]   The
\textit{fractional operator powers} of $G$ are given by
\[G^{\zeta}\varphi  =\frac{1}{\Gamma(\zeta)}\int_{0}^{\infty}t^{\zeta-1}T_{t}\varphi
dt,\ \ \zeta > 0\,\quad \varphi\in \mathcal{B}_c(S).\quad
\]

\subsection{Transience-Recurrence Properties}
In order to review the definitions of transience and
recurrence, (following \cite{DGW-05}) we consider the \textit{last exit time}, $L_A$,  of
$X$ from a non-empty set $A$ defined by
\[ L_A:=
    \sup\{t > 0:X_t\;\in\; A\} \quad (\text{if}\;\{t> 0:X_t\;\in\; A\}\ne
    \emptyset )
    \]

\begin{definition} The L\'evy process $X_t$ on $S$ is \textit{transient} if for any compact set $K$%
\[
P(L_{K}<\infty)=1.
\]
and \textit{recurrent} if it is not transient.
\end{definition}

The following result in the spirit  of Sato and Watanabe
\cite{SW-01}, \cite{SW1} is the basis for a finer classification of
the transience properties of random walks in terms of the moments
of last exit times.

\begin{proposition} Assume that $X_{t}$ is transient,  for any compact
set
$K\subset S$%
\[
\sup_{x\in K}G1_{K}(x)<\infty
\]
and for any compact set $C$ contained in the interior of $ K$%
\[
\inf_{x\in C}G 1_{K}(x)>0.
\]
Then there exist positive constants \thinspace$c_{1}$ and $c_{2}$ such that
for all $\zeta>0$ and $x\in S$%
\[
c_{1}G^{\zeta+1}1_{C}(x)\leq E_{x}L_{C}^{\zeta}\leq
c_{2}G^{\zeta+1}1_{K}(x).
\]
\end{proposition}
\begin{proof}
See \cite{DGW-05}, Proposition 2.2.1.
\end{proof}

\begin{definition}

The \textit{degree of transience}, $\gamma$, of a transient L\'evy process  $X$   is defined by%
\[
\gamma:=\sup\{\zeta >0:E_{0}L_{K}^{\zeta}<\infty\quad\text{for all
compact}\; K\},
\]
 or equivalently \[ \gamma:=
\sup\{\zeta > 0:G^{\zeta+1}\varphi<\infty \text{ for }\varphi\in
C_{c}^{+}(S)\}
\]
where $C_{c}^{+}(S)$ denotes the space of nonnegative continuous
functions on $S$ with compact support.
\end{definition}
\begin{remark}
 Sato and Watanabe introduced the set
\begin{equation}
\mathcal{T}:=
    \{\zeta > 0:E_0L^\zeta_K<\infty  \;\;\text{for all compact}\;K\}.
\end{equation}
In \cite{DGW-05} we consider the extended set
\[
\mathcal{T}:=
    \{\zeta > -1: \int_1^\infty t^\zeta T_t\varphi dt <\infty \;\;\text{for all}\;\varphi\in
    C^+_c(S)\},
\]
and we call
\[ \gamma:= \sup\{\zeta > -1:\zeta\in \mathcal{T}\}\]
the \textit{degree} of the process. This coincides with the degree
of transience if $\gamma >0$, and if $-1<\gamma<0$, we call
$\gamma$ the \textit{degree of recurrence} of the process.

Given $k\in\mathbb{Z}_+$, the process
is said to be (cf. \cite{DGW-01})
\begin{eqnarray}
    k-\text{strongly transient} & \text{if }  \;k\in\mathcal{T},\quad\text{and}\nonumber\\
    k-\text{weakly-transient} & \text{if }  \;k-1\in\mathcal{T} \;\text{and}\; k\notin
    \mathcal{T}.\nonumber
\end{eqnarray}

\end{remark}
\begin{remark}
The degree of transience can be viewed as a generalization of the
notion of ``critical dimension''. Note  that $G^{\zeta+1}\varphi$
at $\zeta=\gamma$ can be either finite or infinite - we will give
examples of both possibilities below.

\end{remark}

\subsection{Random walks and L\'evy processes in $\mathbb{R}^d$.}
In this section we briefly review the classical results on random
walks and L\'evy processes in $\mathbb{Z}^d$ and $\mathbb{R}^d$.

First recall that symmetric nearest neighbour random walks  in
$\mathbb{Z}^d$ are recurrent in dimensions $d=1,2$ and transient
in dimensions $d\geq3$. Moreover, since the rate of decay of the
transition probabilities for simple symmetric $d$-dimensional
random walk is $p_t(0,0)\sim \textrm{const.} t^{-d/2}$,  its
degree  is $\gamma=d/2-1$.

We next recall the classical  characterization of L\'evy processes
in $\mathbb{R}^d$.
\begin{theorem} (L\'evy-Khintchine representation)  A L\'evy process in $\mathbb{R}^d$ has the
characteristic function (i.e. Fourier transform)
\begin{eqnarray}
&&E[e^{i(z,X_t)}]\nonumber\\
&&=exp\left[t\left(-\frac{1}{2}(z,Az)+\int_{\mathbb{R}^d}
(e^{i(z,x)}-1-i(z,x)1_{\{|x|\leq1\}}(x))\nu(dx)+i(m,z)\right)\right]\nonumber
\end{eqnarray}
where $A$ is a symmetric nonnegative definite $d\times d$ matrix,
$\nu$ is a measure on $\mathbb{R}^d\backslash \{0\}$ satisfying
$\nu(\{0\})=0$ and $\int_{\mathbb{R}^d} (|x|^2\wedge
1)\nu(dx)<\infty$, and $m\in\mathbb{R}^d$.
\end{theorem}
For the proof see \cite{S-99}.

The case $A=Id,\;\nu=0,\;m=0$ is the standard Brownian motion and
the case $A=0,\;m=0$ and $\nu(dx)=|x|^{-\alpha-d}dx$, is the
symmetric $\alpha$-stable process.

\begin{proposition}
For the $\alpha$-stable process on $\mathbb{R}^d$ the degree is
\begin{equation}
\gamma ={\frac{d}{\alpha}}-1
\end{equation}
and in this case
\[\int_0^t s^\gamma T_s\varphi ds \sim \mathrm{const}\cdot \log t
\to\infty\] as $t\to \infty$.
\end{proposition}

The distribution of jumps of the $\alpha$-stable process has
``long tails".

\section{The Persistence-Extinction Dichotomy for Critical Branching Systems}

Consider the super-Brownian motion in ${\mathbb R}^d$ with
initial measure $X_0=m\lambda $, $m>0$ where $\lambda $ is Lebesgue measure
If $\sup |\phi (x)|\cdot (1+|x|^2)^{\frac p2}<\infty ,$
then the solution, $v_t(x)=V[\phi ](t,x),$ to
\be{LLE} \frac{\partial v_t}{\partial t}= Av_t -\frac{\gamma}{2}v_t^2,
 \ee with $A=\frac \Delta 2$
is integrable and integrating both sides with respect to Lebesgue measure
gives
\[
\int v_t(x)dx=\frac \gamma 2\int_0^t\int v_s^2(x)dxds.
\]
Therefore the large time limit of the Laplace functional

\begin{eqnarray*}
\lim_{t\rightarrow \infty }{\mathbb P}_{m\lambda }\left( \exp (-X_t(\phi
))\right) &=&\lim_{t\rightarrow \infty }\exp \left( -m\int v_t(x)dx\right) \\
&=&\lim_{t\rightarrow \infty }\exp \left( -\frac{m\gamma }2\int_0^t\int
v_s^2(x)dxds\right)
\end{eqnarray*}
exists for every $\phi \in {\mathcal B}_{+}$ since the right side is monotone in
$t.$ Therefore $X_t$ converges in distribution as $t\rightarrow \infty $ to
a random measure on ${\mathbb R}^d$ with probability law which we denote by $%
{\mathbb P}_m^{eq}$. Replacing $\phi $ by $\theta \phi $, $\theta >0$, and
evaluating the first and second derivatives with respect to $\theta $ at $%
\theta =0$, we can verify that the first and second moments are given by
\[
{\mathbb P}_{m\lambda }(X_t(\phi ))=m\int \phi (x)dx
\]
and
\begin{eqnarray*}
{\mathbb P}_{m\lambda }\left( X_t(\phi )^2\right) &=&m^2\left( \int \phi
(x)dx\right) ^2+\gamma m\int_0^t\int \left( \int p_s(y-z)\phi (z)dz\right)
^2dyds \\
&=&m^2\left( \int \phi (x)dx\right) ^2+\gamma m\int_0^t\left( \int
p_{2s}(z_1-z_2)\phi (z_1)\phi (z_2)dz_1dz_2\right) ds.
\end{eqnarray*}
Recalling that for the Brownian motion transition kernel $\int_0^\infty
p_s(z)ds$ diverges if $d=1,2$ and is given by $\frac{2c_d}{|z|^{d-2}}$ if $%
d\geq 3,$ we obtain
\begin{eqnarray*}
{\mathbb P}_{m\lambda }\left( X_t(\phi )^2\right) &\uparrow &\infty \text{ if }
d=1,2 \\
&\uparrow & m^2\left( \int \phi (x)dx\right) ^2+\gamma mc_d\int\int \frac{\phi
(z_1)\phi (z_2)}{|z_1-z_2|^{d-2}}dz_1dz_2\text{ if }d\geq 3.
\end{eqnarray*}

If $d\geq 3$, the above imply that $\{X_t(\phi )\}_{t\geq 0}$ are uniformly
integrable and ${\mathbb P}_{m\lambda }(X_\infty (\phi ))=m\lambda (\phi )$,
that is, the limiting equilibrium random measure ${\mathbb P}_m^{eq}$ has the
same intensity, $m$, as the initial intensity - this behaviour is called
{\it persistence}. Bramson, Cox and Greven (1997)  \cite{BCG-97} proved
that $\{{\mathbb P}_m^{eq}:m\in [0,\infty )\}$ is in fact the set of all
extremal invariant measures.

\beT{} \cite{D-77} Let $X_{\infty}$ denote the equilibrium random measure for super-Brownian motion in
$\mathbb{R}^d$ with mean measure $E(X_{\infty}(A))=\lambda(A)$.  Let
\be{} \langle X^K_\infty,\phi\rangle =\int \phi(\frac{x}{K})X_\infty(dx),\ee
and \be{} V(\phi):= \gamma \left( \int\int
|z_1-z_2|^{-(d-2)}\phi (z_1)\phi (z_2)dz_1dz_2\right).\ee

Then the rescaled fluctuations
\be{} \frac{\langle X^K_\infty,\phi\rangle -\langle \lambda,\phi\rangle }{K^{\frac{d+2}{2}}V(\phi)}\Rightarrow Z_\infty\ee
where $Z_\infty$ is the Gaussian free field, that is, the Gaussian random field with covariance kernel
\be{} \frac{1}{|x-y|^{d-2}}.\ee
\end{theorem}

The divergence of the second moment in the low dimensional case suggests
that the behaviour is qualitatively different in these dimensions. It was
proved in Dawson (1977) \cite{D-77} that in this case the spatially homogeneous
super Brownian motion with $X_0=m\lambda $ suffers local extinction, that is, $%
X_t(A)\rightarrow 0$ in probability as $t\rightarrow \infty $ for any
bounded set $A$. Iscoe (1986b) \cite{I-86b} has shown that $X_t(A)\stackrel{%
a.s.}{\longrightarrow }0$ for any bounded set if $d=1$ and that this result
is false if $d=2$.  In dimensions $d=1,2$ Bramson, Cox and Greven (
\cite{BCG-93}) have established that $\delta _0$ is the only measure which is
invariant for the process $X_t$ and that for any locally finite initial
measure the system undergoes local extinction or explodes thus ruling out
the possibility of an invariant measure with infinite mean.

\subsection{Clumping in Low Dimensions}

In order to describe the low dimensional behavior of $X_t$ with $X_0=\lambda
$ (Lebesgue) in more detail we introduce the space-time-mass rescaling
\begin{eqnarray*}
X_t^{K,\xi }(A) &:=&K^{-\xi }X_{Kt}(K^{\frac \xi d}A) \\
X_0^{K,\xi }(A) &=&|A|.
\end{eqnarray*}

Then
\begin{eqnarray*}
{\mathbb P}_\lambda (\exp (-X_t^{K,\xi }(\phi ))) &=&\exp (-\lambda (V_{Kt}\phi
_K))\text{ \thinspace \thinspace \thinspace with} \\
\phi _K(x) &:=&K^{-\xi }\phi (K^{-\frac \xi d}x).
\end{eqnarray*}
Note that
\[
\tilde{v}(t,x):=K^\xi V_{Kt}\phi _K(K^{\frac \xi d}x)
\]
satisfies
\begin{eqnarray*}
\frac{\partial \tilde{v}(t,x)}{\partial t} &=&K^{1-\frac{2\xi }d}\Delta
\tilde{v}(t,x)-\frac \gamma {2}K^{1-\xi  }\tilde{v}%
(t,x)^{2} \\
\tilde{v}(0,x) &=&\phi (x)
\end{eqnarray*}
and therefore $X^{K,\xi }$ is equivalent to a super Brownian motion with
``diffusion coefficient'' $K^{1-\frac{2\xi }d}$ and ``branching
coefficient'' $\frac {\gamma} {2}K^{1-\xi  }$. The branching term
dominates in the $K\rightarrow \infty $ limit and the diffusion term
dominates in the $K\rightarrow 0$ limit if $d< 2$ and the
opposite occurs if $d> 2$.

\begin{theorem}\label{DF1}
(Dawson and Fleischmann 1988) \cite{DF-88} (a) Let $d< 2$. Then $X^{K,\xi
}\stackrel{K\rightarrow \infty }{\Longrightarrow }0$ if $\xi <1$
and $X^{K,\xi }\stackrel{K\rightarrow \infty }{\Longrightarrow }\lambda $ if
$\xi >1$\newline
(b) If $d=1$ and $\xi =1$, then $X_K$ converges in distribution
as $K\rightarrow \infty $ to the pure atomic process $\{X_t^0\}_{t\geq 0}$
in which $X_t^0$ is Poisson with intensity $(\frac {\gamma} {2}
t)X(0)$ and the mass of each atom evolves according to a Feller continuous state branching.\newline
(c) If $d=2,$ then $X^{K,1}\stackrel{\mathcal D}{=}$ $X$, that is $X$
is self-similar.
\end{theorem}

\noindent
\begin{remark}
(b) suggests that for $d=1$ at time $K$ there are clumps of size $K$ with interclump distance $K.$
\end{remark}

In the case $d=2$, the phenomenon of {\it diffusive clustering }arises. The
is made precise in the following result of Klenke.

\begin{theorem}
(Klenke (1997) [\cite{K-97}, Theorem 2])  Let $d=2,\,\,$and $I=(-\infty
,1].$ For $\alpha \in I$, let
\[
X_t^\alpha (B):=t^{-\alpha }X_t(t^{\alpha /2}B).
\]
Then in the sense of finite dimensional distributions
\[
{\mathcal L}^{\frac{(\log t)\lambda }{8\pi }}[\{X_t^\alpha (B)\}_{\alpha \in I}]%
\stackrel{t\rightarrow \infty }{\Longrightarrow }{\mathcal L}^1[\{Z_{1-\alpha
}\}_{\alpha \in I}\cdot \lambda (B)]
\]
where $Z$ is a FB process with $Z_0=1$.
\end{theorem}

In the case $d=2$ Theorem \ref{DF1} (c) provides a link between the small scale
and large scale behaviours. In particular it implies that
\[
X_{Kt}(B(0,1))\stackrel{\mathcal D}{=}\frac{X_t(B(0,K^{-1/2}))}{K^{-1}}.
\]
For $t>0$ the left side goes to zero in probability as $K\rightarrow \infty $
because of the local extinction result which then shows that the local
density at time $t$ is $0$ which implies that it does not have a non-trivial
absolutely continuous component.

\subsection{Ergodic Behaviour}

The extinction-persistence result implies that if $\phi $ has compact
support, then $X_t(\phi )$ converges to zero in probability if $d\leq
2 $ and converges in distribution to a non-degenerate limit if $d>
2$. This can be extended to an ergodic theorem in the latter case.

\begin{theorem}
(Iscoe (1986b) (\cite{I-86b}) , Fleischmann and G\"{a}rtner (1986) (\cite
{FG-86})).\newline
(a) For $d> 2 $ with probability one, $\lim_{t\rightarrow \infty
}\frac 1t\int_0^tX_sds=\lambda $ (in the vague topology). \newline
(b) For $d= 2 $, as ${t\rightarrow \infty }$ $\frac 1t\int_0^tX_sds
$ converges a.s. in the vague topology to $\eta \lambda $ where $\eta $ is a
non-degenerate infinitely divisible random variable with mean one.
\end{theorem}

\noindent
\begin{remark}
(b) implies that in critical dimension $d=2$, $X_t(\phi )$ goes
to zero in probability but there are arbitrarily large values of $t$ for
which $X_t(\phi )$ is large, therefore the large clumps revisit bounded sets
at arbitrarily large times, that is, they have a recurrence property.
\end{remark}

\begin{remark} See  the series of papers of Bojdecki,  Gorostiza and  Talarczyk \cite{BGT-06a}, \cite{BGT-06b}, \cite{BGT-07a}, \cite{BGT-07b} for recent advances in the classification of the occupation time fluctuation structure of branching systems.

\end{remark}

\section{The Equilibrium Clan Decomposition}

In this subsection we consider the structure of equilibrium states of SBM following the development in Dawson-Perkins \cite{DP-91}, \cite{DP-99}.

\subsection{An extension of Historical Brownian motion}

We will consider the historical process with time parameter in $(-\infty
,\infty )$. Let $C_{s,t}^d=C([s,t),{\mathbb R}^d)$, $C^d=C_{-\infty ,\infty }^d$
and $(C^d)^s=\{y\in C^d: y(\cdot)=y(\cdot \wedge s)$.
Also let $M_p(C_{s,t}^d)$ denote the space of measures on $C_{s,t}^d$ with
marginal distributions in $M_p({\mathbb R}^d)$. Let $\mu _\lambda^t $ denote
the law (on $C^d$) of $\bar{\xi}_t:=\{\xi (s\wedge t)\}_{s\in (-\infty
,\infty )}$ when $\{\xi (t):t\in {\mathbb R}\}$ is a Brownian motion with
0-marginal given by Lebesgue measure. The transition Laplace functional is
given by: for $\phi \in p{\mathcal B}(C^d)$
\[
{\mathbb Q}_{s,m}\left( \exp (-H_t(\phi ))\right) =\exp \left( -\int
V_{s,t}\phi (y)m(dy)\right) ,\,\,m\in M_p(C^d)^s
\]
where
\[
V_{s,t}\phi (y)=S_{s,t}\phi (y)-\frac \gamma 2\int_s^tS_{s,r}((V_{r,t}\phi
)^2)dr.
\]

If $H_s=\mu _\lambda ^s,$, then $H_t$ is an infinitely divisible random
measure with canonical measure $\int R_{s,t}(y,A)\mu _\lambda ^s(dy)$ where $%
R$ is characterized by
\be{}V_{s,t}\phi (y)=\int (1-e^{-\nu (\phi
)})R_{s,t}(y,d\nu ).\ee

Let $s<u<t.$ By the Historical Cluster Representation $H_t(\{y:y^u\in \cdot
\})$ (under ${\mathbb Q}_{s,m}$) can be represented as a Poisson random field
of clan measures in $M_p(C_{s,t}^d)$ with Poisson intensity $\frac 2{\gamma
(t-u)}H_u.$ The typical clan measure $\Xi _t(y^{\prime })$ can be
interpreted as the descendent population from an individual alive at time $u$
with history $y^{\prime }$ and $r_{s,u}\Xi _t=\Xi _t(C^d)\delta _{y^{\prime
}}$ where for $\mu \in M_F(C_{s,t}^d),$ $s\leq u\leq t,$ $r_{s,u}\mu
(A):=\mu (\{y:y(\cdot \wedge u)\in A\}).$

\subsection{The historical process conditioned to live forever}

For $t>s$ let $H_t^I$ denote a realization of the historical process $H_t$
conditioned to stay alive forever starting from a finite initial measure $%
H_s=\eta _s\in M_p(C_{-\infty ,s}^d)$. The law ${\mathbb Q}^I(H_t^I\in \cdot )$
is defined rigorously as the weak limit as $T\rightarrow \infty $ of ${\mathbb Q%
}_{s,\eta _s}(H_t\in \cdot |H_T\ne 0)$. Using Remark \ref{Remark2.1} we have
\begin{eqnarray*}
{\mathbb Q}_{s,\eta _s}^I(H_t^I\in A) &=&\lim_{T\rightarrow \infty }{\mathbb Q}%
_{s,\eta _s}(H_t\in A|H_T\ne 0) \\
&=&\lim_{T\rightarrow \infty }\frac{{\mathbb Q}_{s,\eta _s}(1_A(H_t){\mathbb Q}%
_{t,H_t}(H_T\neq 0))}{{\mathbb Q}_{s,\eta _s}(H_T\neq 0)} \\
&=&\lim_{T\rightarrow \infty }\frac{{\mathbb Q}_{s,\eta _s}(1_A(H_t)(1-e^{-%
\frac{2H_t(1)}{\gamma (T-t)}}))}{1-e^{-\frac{2\eta _s(1)}{\gamma (T-s)}}} \\
&=&\eta _s({\bf 1)}^{-1}{\mathbb Q}_{s,\eta _s}({\bf 1}_A(H_t)H_t({\bf 1}))
\end{eqnarray*}
The corresponding {\it normalized Campbell measure} of $H_t$ is
\[
{\ \tilde{{\mathbb Q}}}_{s,\eta _s}(H_t^I\in A,\,\bar{\xi_t} \in B):=(\eta _s(%
{\bf 1)})^{-1}{\mathbb Q}_{s,\eta _s}[{\bf 1}_A(H_t)H_t(B)].
\]
Let
\[
U_{s,t}(\phi ,\psi )(y):=\int \nu (\phi )e^{-\nu (\psi )}R_{s,t}(y,d\nu ).
\]

\beL{Lemma5.1}
Let $\phi ,\psi \in {\mathcal B}_{+}(C).$ Then
\begin{eqnarray*}
U_{s,t}(\phi ,\psi )(y) &=&\frac \partial {\partial \theta }V_{s,t}(\theta
\phi +\psi )|_{\theta =0} \\
&=&P_{s,y}\left( \phi (\bar{\xi}_t)\exp \left( -\int_s^t\gamma V_{u,t}\psi (%
\bar{\xi}_u)du\right) \right)
\end{eqnarray*}
\end{lemma}

\proof%

Since (recall the connections between $R_{s,t}$ and $V_{s,t}$)
\[
\int e^{-\nu (\psi )}\nu (\phi )R_{s,t}(y,d\nu )=\frac \partial {\partial
\theta }V_{s,t}(\theta \phi +\psi )|_{\theta =0},
\]
the first equality is clear. Recall
\[
\text{ }V_{s,t}\phi (y)=S_{s,t}\phi (y)-\frac \gamma
2\int_s^tS_{s,r}((V_{r,t}\phi )^2)dr,
\]
and so, replacing $\phi$ with $\theta \phi+\psi$ and differentiating with
respect to $\theta$ we get
\[
U_{s,t}(\phi ,\psi )=S_{s,t}\phi -\gamma \int_s^tS_{s,r}((V_{r,t}\psi
)U_{r,t})dr.
\]
Then by a Feynman-Kac argument for the Brownian path process we get
\[
U_{s,t}(\phi ,\psi )(y)=P_{s,y}\left( \phi (\bar{\xi}_t)\exp \left(
-\int_s^t\gamma V_{u,t}\psi (\bar{\xi}_u)du\right) \right) .\text{ }%
\]

\noindent
\begin{remark}
Let $R_{s,t}^{\bar{\xi}}(y,d\nu)$ denote the Palm measure associated with $%
R_{s,t}(y,d\nu)$, that is, it is the regular conditional probability on\\ $%
M_p(C_{s,t}^d)$ for $A\times B \rightarrow \int\int 1_A(\xi)d\nu
1_B(\nu)R_{s,t}(y,d\nu)$  given $\bar\xi$. The above Lemma implies that
\be{PLF}
\int \exp (-\nu (\psi ))R_{s,t}^{\bar{\xi}}(y,d\nu )= \exp (-\int_s^t\gamma
V_{u,t}\psi(\bar{\xi}_u)du).
\ee

Intuitively, $R_{s,t}^{\bar{\xi}}(y,d\nu )$ is the law of a random measure
obtained if $\xi $ throws off historical excursions at a constant rate $%
\gamma $ on $[s,t]$. We let $\Xi _{s,t}^I(\xi )$ denote such an immortal
clan, i.e. a random measure with law $R_{s,t}^{\bar{\xi}}(y,\cdot )$.
\end{remark}

\beP{Prop5.1}
Under ${\tilde{{\mathbb Q}}}_{s,\eta _s}$,\newline
(a) $\{\bar{\xi}_t\}_{t\geq s}$ is a path-valued Brownian motion with
initial law $(\eta _s({\bf{1}) }^{-1}\eta _s(dy),and$\newline
(b) the regular conditional law, $Q_{s,t}^{\bar{\xi}} $ of $H_t^I$, given $%
\bar{\xi} $ has Laplace functional
\[
Q_{s,t}^{\bar{\xi}} [\exp ^{-<H_t^I,\psi >}]=\int\left( \int \exp (-\nu
(\psi ))R_{s,t}^{\bar{\xi}}(y,d\nu ) \right)\eta _s(dy)\eta
_s(1)^{-1}Q_{s,\eta _s}(\exp (-<H_t,\psi >)).
\]
\end{proposition}

\proof%

\begin{eqnarray*}
&&{\tilde{{\mathbb Q}}}_{s,\eta _s}\left( \phi (\bar{\xi})\exp ^{-<H_t^I,\psi
>}\right) =(\eta _s({\bf 1)})^{-1}{\mathbb Q}_{s,\eta _s}[\exp ^{-<H_t,\psi
>}H_t(\phi )] \\
&=&-(\eta _s({\bf 1)})^{-1}\frac \partial {\partial \theta }{\mathbb Q}_{s,\eta
_s}\left( \exp ^{-<H_t,\theta \phi +\psi >}\right) \big| _{\theta =0} \\
&=&-(\eta _s({\bf 1)})^{-1}\frac \partial {\partial \theta }\left( \exp
^{-\eta _s(V_{s,t}(\theta \phi +\psi )}\right) \big| _{\theta =0} \\
&=&\exp({-\eta_s(V_{s,t}(\psi))}) P_{s,\eta _s/\eta _s(1)}\left( \phi (\bar{%
\xi}_t)\exp (-\int_s^t\gamma V_{u,t}\psi (\bar{\xi}_u)du)\right) \text{ } \\
&&\,\,\,\,\,\,\,\,\,\,\,\text{( use Lemma \ref{Lemma5.1}) } \\
&=& \int P_{(s,y)}\left( \phi(\bar \xi_t) \int \exp (-\nu (\psi ))R_{s,t}^{%
\bar{\xi}}(y,d\nu )\right) \eta _s(dy)/\eta _s(1) Q_{s,\eta _s}(\exp
(-<H_t,\psi >))
\end{eqnarray*}
which proves the result.

\noindent
\begin{remark}
From the perspective of a typical immortal particle chosen according to $%
H^I_t$, $H^I$ is the sum of $\Xi^I_{s,t}(\xi)$ and an independent copy of $%
H_t$. The former is the contribution to $H^I_t$ of cousins of $\xi$.
\end{remark}

\subsection{Convergence to equilibrium for the historical process}

Recall that if $H_s=\mu _\lambda ^s$, then $H_t$ is an infinitely divisible
random measure with canonical measure $R_{s,t}=\int R_{s,t}(y,A)\mu _\lambda
^s(dy)$ where $R$ is characterized by $V_{s,t}\phi (y)=\int (1-e^{-\nu (\phi
)})R_{s,t}(y,d\nu ).$

Integrate both sides of (\ref{PLF}) with respect to $\mu _\lambda ^s$ to see that $%
R^{\bar \xi}_{s,t}$ is still the Palm measure of the canonical measure $\int
R_{s,t}(y,d\nu) \mu^s_\lambda (dy)$ of $H_t$.

\smallskip

As $s\rightarrow -\infty $, the total mass of $\Xi _{s,t}^I(\xi )$ goes to
infinity (set $\psi =1$ in Proposition \ref{Prop5.1} and let $s\rightarrow -\infty
) $. However

\begin{proposition}
(a) As $s\rightarrow -\infty$, $\Xi _{s,t}^I(\bar{\xi})$ converges to a
locally finite measure iff $d\geq 3.$ \newline
(b) As $s\rightarrow -\infty ,$ $R_{s,t}$ converges vaguely to a measure $%
R_{-\infty ,t}$ such that $\int \mu R_{-\infty ,t}(d\mu )$ is locally finite.%
\newline
(c) $R_{-\infty ,t}$ is ``an entrance law'' for $H$ in the sense that if $%
H_s $ has the infinite law $R_{-\infty ,s}$ and $t>s$, then $H_t$ has law $%
R_{-\infty ,t}.$
\end{proposition}
\proof%

(a) The Laplace functional in (\ref{PLF}) converges as $s\rightarrow -\infty $ by
monotonicity. In the case $d\geq 3,$ we can obtain the mean measure of $\Xi
_{s,0}^I(\bar{\xi})$ for a trajectory $\xi $ as follows. Differentiate the
expression for the Laplace functional in (\ref{PLF})  to see that for $\phi $
continuous with compact support
\[
\lim_{s\rightarrow -\infty }E[\int_{{\mathbb R}^d}\phi (y_t)\Xi _{s,t}^I(\bar{%
\xi},dy)]=\gamma \lim_{s\rightarrow -\infty }\int_{s-t}^0\int p(|r|,y-\xi
(t+r))\phi (y)dydr<\infty
\]
and hence $\Xi _{-\infty ,t}^I$ is a.s. locally finite.

(b) Without loss of generality we can take $t=0$ and $\psi (y)=\tilde{\psi}%
(y_0),\phi (y)=\tilde{\phi}(y_0)$, when $\tilde{\phi},\tilde{\psi}$ are
continuous with compact support. Then for $d\geq 3$, the above gives
\begin{eqnarray*}
\lim_{s\rightarrow -\infty }&&\int \int \mu (\psi )\mu (\phi )R_{s,0}(y,d\mu
)\mu _\lambda ^s(dy) \\
&=&\lim_{s\rightarrow -\infty }\int \int \psi (\xi )\mu (\phi )R_{s,t}^{\bar{%
\xi}} (d\mu )P_{s,\mu _\lambda ^s}(d\xi ) \\
&=&\gamma \lim_{s\rightarrow -\infty }\int \int \int_s^0p(2|r|,z-y))\phi
(z)\psi (y)drdzdy \\
&<&\infty .
\end{eqnarray*}

(c) See \cite{DP-91} Theorem 6.4.%

\begin{theorem} (Dawson-Perkins (1991), \cite{DP-91}, Theorem 6.3) The stationary random measure $X_\infty $ of the super-Brownian
motion in ${\mathbb R}^d$ is an infinitely divisible random measure with
intensity $\lambda $ and with canonical measure given by $\hat{R}%
(B)=R_{-\infty ,0}(\{\pi _0\mu \in B\})$ where $\pi _0\mu (A):=\mu (\{\xi
:\xi _0\in A\}),\ A\in {\mathcal B}({\mathbb R}^d),\ \mu \in M_p(C^d).$
\end{theorem}

\noindent
\begin{remark}
This leads to a description of the super-Brownian motion in equilibrium as a
countable collection of infinite clan measures.
\end{remark}

An infinite clan containing an individual located at $0$ at time $0,$
denoted by $\Xi _{-\infty ,0}^I,$ has law $\int R_{-\infty,0 }^{\bar{\xi}%
}\Pi _0(d\xi )$. It is constructed by running a Brownian trajectory, $\xi $,
backwards to time $-\infty $ and then collecting the mass at time 0
corresponding to $R_{-\infty,0 }^{\bar{\xi}}$.

\subsection{Clan dynamics}

The clan $\{\Xi _{-\infty ,t}^I\}_{t\geq 0}$ then evolves as a historical
Brownian motion with initial condition $\Xi _{-\infty ,0}^I$. It is an easy
consequence of this description that the infinite clan is self-similar in
the sense that
\be{SS}
K^{-2}\Xi _{-\infty ,0}^I(\{y:y_0\in KA\})\stackrel{\mathcal D}{=}\Xi _{-\infty
,0}^I(\{y:y_0\in A\})\text{\thinspace \thinspace }\forall \,\,A\in {\mathcal B}(%
{\mathbb R}^d).
\ee
There has been considerable recent interest in the dynamics of these
infinite clans and this has led to a second dimensional dichotomy - namely
clan recurrence or transience. To describe this consider the total weighted
occupation time that a given clan spends in the unit ball at the origin $%
B(0,1)$. In dimensions $d>2,$ (\ref{SS}) and spherical symmetry implies that the
infinite clan of an individual at the origin and has mean density
proportional to
\be{meanclandensity}
\frac 1{|x|^{d-2}}\text{.}
\ee
Therefore if we exclude the clan mass initially in a ball of radius 2 (which
by the extinction property of Feller branching has a.s. finite lifetime),
then
\begin{eqnarray*}
E[\int_0^\infty && \Xi _{-\infty ,t}^I(\{y :y_t\in B(0,1),\,\,y_0\notin
B(0,2)\})dt] \\
&=&\int_0^\infty \int_{B(0,2)^c}\frac 1{|y|^{d-2}}\left[
\int_{B(0,1)}p_t(x-y)dx\right] dydt \\
&=&\int_{B(0,2)^c}\left[ \int_{B(0,1)}\frac 1{|x-y|^{d-2}}dx\right] \frac
1{|y|^{d-2}}dy \\
&<&\infty \text{ iff }d\geq 5.
\end{eqnarray*}
Hence in dimensions $d\geq 5$, \be{CD1}\int_0^\infty \Xi _{-\infty ,t}^I(\{y:y_t\in
B(0,1)\})dt<\infty ,\,\, \text{a.s.}\ee - this behaviour is called {\it clan
transience.} Analogous results of St\"{o}ckl and Wakolbinger \cite{SW-94} show that the
corresponding infinite clan for a branching particle system in dimensions $%
d\geq 5$ gives positive mass to the unit ball over only a finite time
horizon.

\section{Neutral Stepping Stone Models}

\label{s.NSSM}

\subsection{The two type stepping stone model}
The neutral two type stepping stone model on a countable abelian group $S$ with migration kernel $p(\cdot)$ is given by the system
\begin{align*}
dX_t(x)  &  =\sum_{y\in S_1}p_{y-x}(X_t(y)
-X_t(x)dt\\& +\sqrt{2X_t(x)(1-X_t(x))}dW_t(x)\\
x_0(x)\in [0,1],\; x\in S
\end{align*}

This process can be embedded in the infinitely many types stepping stone model which we now consider.

\subsection{The infinitely many types stepping stone model}

Consider a collection (finite or countable) of subpopulations
(demes),
 indexed by $S$.
 The subpopulation at $\xi \in S$ at time $t$   is described by a
probability distribution $X_\xi(t)$ over a space $E = [0,1]$ of
possible types (alleles). In other words, $X_\xi(t)\in \mathcal{P}(E)$,
the set of probability measures on E so that the state space is \be{}
(\mathcal{P}([0,1]))^{S}.\ee

 Within  each subpopulation
there is mutation, selection and finite population sampling.
Mutation is assumed produce a new type chosen by sampling from a
fixed source distribution $\theta \in \mathcal{P}(E)$. Selection is
prescribed by a fitness function  $V(x)$ in the haploid case or by
$V(x,y)=V(y,x)$ in the diploid case. Migration from site $\xi$ to
site $\xi'$ is assumed to occur via a symmetric random walk with  rates  $q_{\xi,\xi'}= p(\xi -\xi')$.
 Finally
Fleming-Viot continuous sampling is assumed to take place within
each subpopulation. It is a basic property of this model that for
any $t>0$, $X_\xi(t)$ is a purely atomic random measure (with
countably many atoms) and therefore can be represented in the form
$$ X_\xi(t) = \sum_{k\in I} m_{\xi,k}(t) \delta_{y_k}$$
where $m_{\xi,k}(t)\geq 0$  denotes the proportion of the population
in subpopulation $\xi$ of type $y_k\in E$ at time $t$.  Note that in this
model two individuals are related if and only if they are of the
same type.

We denote the vector $\{\mu_\xi\}_{\xi\in S}$ by $\bar \mu$.
The generator is then given by
\bea{}&&\\
GF(\bar \mu) &=&
c\cdot \sum_{\xi \in S}
\int^{}_{[0,1]} {\partial F(\bar\mu)\over \partial \mu _{\xi}(u)}
(\theta(du)-\mu _{\xi}(du))\nonumber\\&& +
\sum
q_{\xi,\xi'}\int_{[0,1]}\frac{\partial F(\bar\mu)}{\partial \mu
_\xi(u)}(\mu _{\xi'}(du)-\mu _\xi(du))\nonumber \\
&&+\frac \gamma 2\sum \int_{[0,1]}\int_{[0,1]}\frac{\partial
^2F(\bar \mu)}{\partial \mu _\xi(u)\partial \mu _\xi(v)}Q_{\mu _\xi}(du,dv)\nonumber
\eea
\[
X_{0,\xi}=\nu\;\forall\;\xi,\quad Q_\mu (du,dv)=\mu (du)\delta
_u(dv)-\mu (du)\mu (dv).
\]
 The first term
corresponds to mutation with source distribution $\theta$, the
second to spatial migration and the last to continuous resampling.
The resampling rate coefficient $\gamma$ is inversely proportional to the
effective population size of a deme.

This existence and uniqueness of this system of interacting Fleming-Viot processes  was established by Vaillancourt \cite{V-90} and Handa \cite{H-90}.

The questions which we wish to investigate are
\begin{itemize}
\item the distribution in a given subpopulation, that is what is the joint distribution of the $\{m_{\xi,k}\}$
\item  the spatial distribution of relatives
\item how are these affected by the migration geometry.
\end{itemize}

\medskip
\subsection{The Dual Process Representation}
Given $n\in\N$ consider the collection
\bea{}&&\\ \Pi_n&&= \{\bar\eta:=(\eta,\pi)\}:\text{ where }\nonumber\\&&\pi\text{ is a partition of }\{1,\dots,n\},\text{ that is, }\nonumber\\&&\pi:\{1,\dots,n\}\to\{1,\dots,|\pi|\}\text{ with }|\pi|\leq n,\nonumber\\&&
\eta:\{1,\dots,|\pi|\}\to S.\nonumber\eea
Now consider the family of functions in $C((\mathcal{P}([0,1]))^S\times\Pi)$ of the form
\be{dfx}F_f(\bar\mu,\bar\eta)
:=\int_{[0,1]}\dots\int_{[0,1]}f(u_{\pi(1)},\dots,u_{\pi(n)})\mu_{\eta_1}(du_1)\dots\mu_{\eta_{|\pi|}}(du_{|\pi|})
\ee
with $f\in C([0,1]^n)$.

We now consider a continuous time Markov chain, $\bar\eta_t=(\eta_t,\pi_t)$,

 with state space $\Pi_n$ and jump rates:

\begin{itemize} \item the partition elements perform continuous time
symmetric random walks on $ S$
 with rates $q_{\xi,\xi'}$ and in addition a
partition element can jump to $\{\infty\} $ with rate $c$ (once a
partition element reaches $\infty$ it remains there without change
of further coalescence). \item each pair of partition elements
during the period they reside at an element of S (but not
$\{\infty\}$) coalesce at rate $\gamma$ to the partition element equal to the union
of the two partition elements.
\end{itemize}

Let $H$ denote the generator of $\bar\eta$. Then for a function of the form (\ref{dfx})
\be{} HF_f(\bar\mu,\bar\eta)=GF_f(\bar\mu,\bar\eta).\ee

We then obtain the dual relationship

\be{DR1} E(F_f(X_t,(\eta ,\pi)))=E(F_f(X_0,(\eta _t,\pi _t))) \ee
  and this proves that the infinitely many types stepping stone martingale problem is well-posed.


\begin{remark} Given the dual we can construct a {\em spatially structured coalescent} that describes the ancestral structure of a sample of a finite number of individuals located at the same or different sites.

  Note
that this is essentially equivalent to  the coalescent
geographically structured populations introduced by developed by
Notohara (1990) \cite{N-90} and Takahata (1991) \cite{T-91}.

\end{remark}

\begin{remark} Note that as $\gamma\to \infty$ the dual converges to the dual of the voter model and we can regard the voter model as the limit as $\gamma\to\infty$ of the interacting Fisher-Wright diffusions.

\end{remark}

\subsection{Spatial homogeneity and the local-fixation coexistence dichotomy}

In this subsection we consider the neutral stepping stone model \underline{without mutation}.

 \beT{DGV} (Dawson-Greven-Vaillancourt (1995) \cite{DGV-95}, Theorem 0.1)

Let $S$ be a countable abelian group and consider the infinitely many types stepping stone model with no mutation ($c=0$). Assume that the initial
random field $\{X_\xi(0)\}_{\xi\in S}$ is spatially stationary, ergodic, weakly mixing and has single site mean
measure satisfying  \be{}E(\int g(u) X_\xi(0,du))=\int g(u)\theta(du),\quad \theta \in \mathcal{P}[0,1].\ee

(a)  If $q_{\xi,\xi'}$ is a symmetric transient
random walk  on $S$, then the stepping stone process $\{X_\xi(t)\}_{\xi\in S}$ converges in
distribution to a nontrivial  invariant $\mathcal{P}([0,1])$-valued  random field $\{X_\xi(\infty)\}_{\xi\in S}$
which also has single site mean measure $\theta $. $\{X_\xi(\infty)\}_{\xi\in S}$   is spatially homogeneous (that is, the law is invariant under translations on $S$), ergodic and weakly mixing, in particular
\be{} E(\langle x_\xi,f\rangle\;\langle x_\zeta,f\rangle)\to \langle \mu,f\rangle^2\text{  as  }d(\xi,\zeta)\to\infty,\quad \forall\; f\in L_\infty([0,1]).\ee

(b)  In (a) the equilibrium state decomposes into countably many coexisting
infinite families, namely,
\be{}  X_\xi(\infty)= \sum_{k=1}^\infty a_{\xi,k}\delta_{y_k}\ee
with $\sum_xi a_{\xi,k}=\infty$ for each $k$.

(c) If $p_\xi$ is
recurrent, then the set of invariant measures is a convex set with extremal invariant measures are $\delta _a$,
$a\in [0,1],$ that is, there is local fixation, and
\be{} \mathcal{L}(\{X_\xi(\infty)\}_{\xi\in S} =\int (\delta_{y})^S)\theta(dy).\ee

\end{theorem}
\begin{proof} We sketch the main steps of the proof.

The proof uses the dual representation (\ref{DR1}), (\ref{DR2}).  Note that $|\pi_t|$ is monotone decreasing so that we can define
\be{3.13} \pi_{\infty}=\lim_{t\to\infty} \pi_t,\quad \pi_\infty=\{\pi_\infty(1),\dots,\pi_\infty(n)\}. \ee
Then we note that $\wh\eta$ is prescribed by a {\em coalescing random walk with delay}.  We let $Z(t)$ be a random walk on $S$ with transition kernel $\{q_{\xi,\xi'}\}$.  Since we have assumed that the random walk is symmetric, then the difference process $Z_1(t)-Z_2(t)$, where $Z_1,Z_2$ are independent copies of the random walk, is a random walk with jump rates $2q_{\xi,\xi'}$. We can assume that the system of coalescing random walks with delay is constructed on a probability space on which the sequence $\{Z_i(t)\}_{i\in\N}$ of independent  random walks and an independent collection of exponentially distributed random variables are defined.

\beL{AL1} If the $q$-random walk is recurrent, then

(a)
\be{3.14} \mathcal{L}(\wh\eta_t)-\mathcal{L}((Z(t);\{1,\dots,n\}))\Rightarrow 0\text{  as  }t\to\infty\ee.

Given two initial sites $0$ and $\xi\ne 0$ and $(\eta,(\{1\},\{2\}))$, $\eta_1=0,\eta_2=\xi$,

\be{} P(\pi_t=\{1,2\})\leq const\cdot \int_0^t P(Z(s)=0)ds\leq \frac{const}{|\xi|^{d-2}}\ee
where $Z(s)$ is a random walk starting at $\xi$ and with jump rate 2.

(b) If the $q$-random walk is transient, then
\be{3.16} \mathcal{L}(\eta_t|\,|\pi_\infty|=k)-\mathcal{L}(Z_1(t),\dots,Z_k(t))\Rightarrow 0\text{  as } t\to\infty\ee
and $P(|\pi_\infty|=1)<1$ provided that $|\pi_0|\ne 1$.

\end{lemma}
\begin{proof} If the random walk is recurrent, then $Z_1(t)-Z_2(t)$ visits $0$ infinitely often and therefore they must coalesce with probability one.

If the random walk is transient, then there exists a random time $\sigma<\infty$ a.s. such that $\sigma$ is the last coalescence time in the system $\wh\eta_t$.  Denote by $\xi^1(t),\dots,\xi^{|\pi_\infty|}(t)$ the position of the partition elements at time $\sigma+t$.
The system $\eta_u$ for times $u=s+t+\sigma$   behaves like a system of $|\pi_\infty|$ random walks in $s$ starting at $\xi^1(t)\dots\xi^{|\pi_\infty|}(t)$ and conditioned on never meeting. Since for every pair $i\ne j$ $\xi^i(t)-\xi^j(t)\to \infty$ as $t\to \infty$, the event that $\xi^i$ and $\xi^j$ never meet after time $t$ tends to one as $t\to\infty$.  It remains to show that the distance between the distributions of the system of $|\pi_\infty|$ independent random walks starting at $(\xi^1(t)\dots\xi^{|\pi_\infty|}(t))$  and starting at $(0,\dots,0)$ tends to $0$ as $s\to\infty$.  This is verified using a coupling by randomized stopping times due to Greven (1987) \cite{Gr-87} and a result of Choquet and Deny on transient random walks (see Spitzer \cite{S-64}, Ch 6. T1) - see (\cite{DGV-95} for details).

\end{proof}
\bigskip

We also note the following elementary result on random probability measures.
\beL{AL2} Let $X_1,X_2$ be a random probability measures on $[0,1]$, having the same mean measures $E(X_i)=\theta\in\mathcal{P}([0,1])$, that is, a measurable map from a probability space $(\Omega,\mathcal{F},P)$ to $\mathcal{P}([0,1])$.

(a) If
\be{} E[\left(\int g(y)X_i(dy)\right)^2]= E[\int g^2(y)X(dy)]\quad \forall g\in C([0,1]),\ee then
\be{} X_i(\omega)=\delta_{y(\omega)}\quad\text{for } a.e. \;\;\omega\in\Omega\text{  and  } \omega\to y\text{  is measurable}.\ee

(b) If in addition,
\be{} E[\left(\int g(y)X_1(dy)\int g(y)X_2(dy)\right)]= E[\int g^2(y)X_1(dy)]\quad \forall g\in C([0,1]),\ee then
\be{} X_1(\omega)=X_2(\omega)=\delta_{y(\omega)}\quad\text{for } a.e. \;\;\omega\in\Omega.\ee
\end{lemma}

We return to the proof of the theorem.

(a) Recurrent Case.

Step 1. Let $m=2$ and take $f(u_1,u_2)=g(u_1)g(u_2), \eta=(\xi,\xi)$. Then by the dual representation and Lemma \ref{AL1}
\be{3.17}
\begin{split}
&E\left( \int_0^1 g(u)X_\xi(t,du)\right)^2
\\& \quad =E\left(\int g^2(u)X_{\eta^1_t}(0,du)1(\pi_t=\{1,2\})\right.\\&\quad \left. +\int g(u)X_{\eta^1}(t,du)\cdot\int g(u)X_{\eta^2}(t,du)1(\pi_t=\{\{1\},\{2\}\}\right)\\& =  E\left(\int g^2(u)X_{\eta^1_t}(0, du)\right)+o(t)
\end{split}
\ee
Therefore in the limit by Lemma \ref{AL2}(a) we have \be{3.19}X_\xi(\infty,du)=\delta_y,\; a.s.\ee


Since $\mathcal{L}(X_\xi(t))\in\mathcal{P}(\mathcal{P}([0,1])$, the set $\{\mathcal{L}(X_\xi(t))\}_{t\geq0}$ is weakly relatively compact.  By (\ref{3.19}) a weak limit point must be concentrated on
\be{} M=\{\delta_u:u\in[0,1]\}\ee
that is, $\mathcal{L}(X_\xi(\infty))=\int_0^1\delta_{\delta_u}H_\xi(du)$ with $H_\xi\in\mathcal{P}([0,1])$. But we have
\be{} E\langle X_\xi(t),f\rangle =\langle \theta,f\rangle \ee
so that for a limit point $\mathcal{L}(\{X_\xi(\infty)\}_{\xi\in S})$
\be{} E\langle X_\xi(\infty),f\rangle=\langle\theta,f\rangle\quad \forall\; f\in C([0,1]). \ee
Therefore $H_\xi=\theta $
\be{3.24} \mathcal{L}(X_\xi(\infty))=\int_0^1\delta_{\delta_u}\theta(du).\ee

Step 2. In order to show consensus  of the components occurs for $t\to\infty$ take $m=2$, $f(u_1,u_2)=g(u_1)g(u_2)$ but use $\eta=(\xi^1,\xi^2)$ with $\xi^1\ne\xi^2$. Then again using Lemma \ref{AL1}
\be{3.25}
\begin{split}
&E\left( \int_0^1 g(u)X_{\xi^1}(t,du)\int_0^1 g(u)X_{\xi^2}(t,du)\right)
\\& \quad =E\left(\int g^2(u)X_{\eta^1_t}(0,du)1(\pi_t=\{1,2\})\right.\\&\quad \left. +\int g(u)X_{\eta^1_t}(0,du)\cdot\int g(u)X_{\eta^2_t}(0,du)1(\pi_t=\{\{1\},\{2\}\})\right)\\&= E\left(\int g^2(u)X_{\eta^1_t}(0, du)\right)+o(t)
\end{split}
\ee
The result then follows from Lemma \ref{AL2}(b), that is
\be{3.19} (X_{\xi^1}(\infty),X_{\xi^2}(\infty))=(\delta_y,\delta_y)\text{  for some  random }y,\;\; a.s.\ee
where
\be{}P(y\in(a,b))=\theta((a,b)).\ee

\medskip

%


Step 3. We can obtain the analogue of  (\ref{3.19}) foe any finite $\xi^1,\dots,\xi^k$.  Therefore  we obtain
\be{} \mathcal{L}((x_\xi(t))_{\xi\in S})\Rightarrow \int\delta_{(\delta_u)^S}\theta(du)\ee
and the proof of (a) is complete.

(b) Transient case. To prove convergence of $\mathcal{L}(t)$ as $t\to\infty$ we first recall that
\be{3.31} \pi_t\to\pi_\infty,  \text(cf. (\ref{3.13})).\ee

Let $n\in\N$ and $f(x_1,\dots,x_n)=\prod_{i=1}^n f_i(x_i)$.  Then by the dual representation

\bean{} && E_{X(0)}(F(X(t),(\eta,\pi)))= E_{(\eta,\pi)}(F(X(0),(\eta_t,\pi_t)))
\\&&= \sum_{m=1}^n E_{(\eta,\pi)}\left(\langle X_{\eta^1_t}(0),\prod_{i\in \pi_t(1)}f_i\rangle,\dots \langle X_{\eta^m_t}(0),\prod_{i\in \pi_t(m)}f_i\rangle 1(|\pi_t|=m)\right)  \\&&
\to E_{(\eta,\pi)}\left(\langle \theta,\prod_{i\in \pi_\infty(1)}f_i\rangle,\dots \langle \theta ,\prod_{i\in \pi_\infty(|\pi_\infty|)}f_i\rangle )\right)
\eean
where we have used the fact that $|Z_i(t)-Z_j(t)|\to \infty$ in probability  as $t\to\infty$ and the weak mixing property of the initial random field so that for $i\ne j$ \bean{}&&\lim_{t\to\infty} E(\int f_1(x)X_{Z_i(t)}(0,dx)\int f_2(y)X_{Z_j(t)}(0,dy))\\&&=\lim_{t\to\infty} E(\int f_1(x)X_{Z_i(t)}(0,dx))E( \int f_2(y)X_{Z_j(t)}(0,dy))\\&&
= \int f_1(x)\theta(dx)\int f_2(y)\theta(dy).\eean
This implies the convergence
of the laws $\mathcal{L}_t$.

The proof of the weak mixing property is obtained  by noting that if  $|\eta_1-\eta_2|\to \infty$, then
\be{}P(\pi_\infty=(\{1\},\{2\}))\to 1\quad\text{as }t\to\infty.\ee




The proof that the limiting law $\mathcal{L}_\infty$ s an invariant measure for the dynamics is standard.
\end{proof}

\begin{remark} (Population structure in 2 dimensions)

The phenomenon of {\em diffusive clustering} in dimension
 {$d=2$} was discovered by Cox and Griffeath  (1986) \cite{CG-86}.

More recently, coalescing random walks used to study the
{coalescence time and identity by descent} between 2
randomly chosen individuals on a 2-d torus ( Cox and Durrett (2002) \cite{CD-02}, Cox, Durrett,
Z\"ahle (2005) \cite{CDZ-05})

\end{remark}

\subsubsection{Homozygosity in large time scales}

Given a probability measure $\mu$ on $[0,1]$ the {\em homozygosity} is defined by
\be{} \int_0^1\int_0^1 1_{x=y} \mu(dx)\mu(dy) =\sum_{i=1}^\infty a_i^2\ee
where $\{a_i\}$ are the masses of the atoms (if any) in $\mu$, that is $\mu=\sum a_i\delta_{y_i}+\mu_{diff}$
and $\mu_{diff}$ is the non-atomic component of the measure.

It follows from Theorem \ref{DGV} that in the recurrent case for any $L\in\N$

\be{}
\lim_{t\to\infty}E\left[\frac 1{N(L)}\sum_{|j|\leq
L}<X_{\xi}(t)\otimes X_0(t),I_\Delta >\right]=1.
\ee
where $I_{\Delta}=\{(x,y):x=y\}$ and $N(L)$ denotes the number of sites in a ball of radius $L$
and for the transient case
\be{}
\lim_{t\to\infty}E\left[
<X_0(t)\otimes X_0(t),I_\Delta >\right]<1
\ee

\begin{theorem} Consider the stepping stone model on $\Bbb{Z}^d$ and random walk kernel given by a nearest neighbour random walk.
Let $d\geq 3$ and $X_0=\nu ,$ with $\nu $
nonatomic. Then\\
(a)
\be{}
\lim_{L\rightarrow \infty }\frac 1{L^d}\sum_{|j|\leq
L}<X_\xi(\infty)\otimes X_0(\infty),I_\Delta >=0,
\ee
(b) Each allelic type present at equilibrium has infinite total mass in $\mathbb{Z}^d$ but has zero spatial density.
\\
In addition, if ${X}(0)$ is given the stationary measure, then\\
(c)
\be{}
\int_0^\infty <X_0(t)\otimes X_0(0),I_\Delta >dt<\infty
\ee
if and only if $d\geq 5$.
\end{theorem}

\begin{proof}
(a) We briefly sketch the argument. We note that if $\pi=(\{1\},\{2\}),\;\eta_1=0,\eta_2=\xi$ then
\be{} \lim_{t\to\infty} E[<X_\xi({t})\otimes X_0(t),I_\Delta >]\leq P_{(\eta,\pi)}(\pi_t=\{1,2\})\ee
since if coalescence does not occur, the expected homozygosity is 0. But the probability that two random walks  $Z_0$ and $Z_\xi$ starting at $0$  and $\xi$  coalesce by time $t$ satisfies
\bea{}&& \lim_{t\to\infty}P(\text{coalesce by time }t)\\&&= \lim_{t\to\infty} E(1-e^{-\gamma\int_0^t(1(Z_0(s)=Z_\xi(s)))ds})\nonumber
\\&&\lim_{t\to\infty}\leq (1- e^{-\gamma\int_0^tP(1(Z_0(s)=Z_\xi(s)))ds}) \sim \frac{1}{|\xi|^{d-2}}.\nonumber\eea
The result follows by summing and dividing by $L^d$.

(c) is the analogue of (\ref{CD1}).
\end{proof}

\subsubsection{Family decomposition and renormalization of the fluctuation field}

The decomposition of the infinitely many types stepping stone model and the related voter model provides a tool for the study of the renormalized fluctuation field. (Recall that the difference between the stepping stone model and the voter model is that coalescence of the random walks  occurs with delay for the stepping stone model but is instantaneous for the voter model. Otherwise the structure of the infinite clusters is similar.)  The following special case of a theorem of I. Z\"ahle illustrates this.

\beT{} \cite{Z-01} Consider the equilibrium voter model $\{X_\xi(t)\}_{\xi\in\mathbb{Z}^d}\in \{0,1\}^{\mathbb{Z}^d} $  with nearest neighbour simple random walk kernel.  For a bounded function $\phi$ with bounded support let
\be{} Z_r(\phi):=\frac{\sum_{\xi\in\mathbb{Z}^d} [X_{\xi}(\infty)-E(X_{\xi}(\infty))]\phi(\frac{\xi}{r})}{r^{\frac{d+2}{2}}}\ee
If $d\geq 3$, then as $r\to\infty$, $Z_r$ converges weakly to the Gaussian free field on $\mathbb{R}^d$, that is, the Gaussian field on $\mathbb{R}^d$ with covariance kernel $\frac{\rm{1}}{r^{\frac{d-2}{2}}}$.
\end{theorem}

\begin{remark}
Recall that the dual of the voter model and the dual for the 2 type Wright-Fisher diffusion differ only in that for the voter model the coalescence is instantaneous and for the Wright-Fisher model coalescence occurs with delay. Using this observation the basic strategy of the proof of this theorem which involves the ``infinite colour'' decomposition can be applied to the case of the Wright-Fisher diffusion.
\end{remark}

\subsection{Historical and genealogical structure of the neutral stepping stone model}

In order to describe the historical and genealogical structure of the neutral stepping stone model we need the analogues of the historical and ancestral processes described in earlier chapters.
In this direction,
Greven, Limic and Winter (2005) \cite{GLW-05} have introduced the historical interacting Moran model and the historical interacting Fisher-Wright diffusions on a countable abelian group $S$.
We now briefly describe their formulation.

Consider the type space $E=\{0,1\}$.  Then the Moran model is given by a locally finite population with the following dynamics
\begin{itemize}

\item each individual moves in $S$ independently according to a random walk with transition function
    \be{}  p_t(x,y)=\sum_{n\geq 0} p^{(n)}(x,y)\frac{t^ne^{-t}}{n!},\quad  x,y\in S\ee

\item resampling: each pair of indivudals dies at rate $\gamma$ and is replaved by a new pair of individual where each new individual adopts a type by choosing the parent independently from the ``dying pair''.

\end{itemize}

This system of particles defines  a measure $\eta\in M(S\times E)$.  As state space we take the Liggett-Spitzer space
\be{} \mathcal{E}_S:= \{\eta\in M(E\times E):\sum_S \eta(\{x\}\times K)\alpha(\{x\})<\infty\ee
where $\alpha$ is a finite measure on $S$ such that
\be{} \sum_y p(x,y)\alpha(\{y\})\leq \Gamma\alpha(\{x\}).\ee

Consider the two type Fisher-Wright system given by
\be{} d\zeta_t(x)= \sum_y (p(x,y)-\delta(x,y))\zeta_t(x)+\sqrt{g(\zeta_t(x)}dw^x_t,\quad x\in S\ee
where $g(z)=\gamma z(1-z)$.

This can be obtained as a limit as follows: let $\theta\in [0,1]$, $(\eta^\rho)_{\rho\in\mathbb{R}_+}$ have law concentrated on $\mathcal{E}_S$ such that $\mathcal{L}[\eta(\cdot\times E)]\in \mathcal{E}_S$ and be translation invariant and ergodic with total mass per site intensity $\rho >0$ and $\mathcal{L}[\eta(\cdot\times \{1\})]\in \mathcal{E}_S$ and translation invariant and ergodic with total mass per site intensity $\theta\rho$.

\be{}  \hat \eta^\rho_t(x):= \frac{\eta^\rho_t(\{x\})\times\{1\})}{\eta_t^\rho(\{x\}\times E)}1_{\{\eta^\rho_t(\{x\}\times E)\ne 0\}}\ee
that is the relative frequency of type 1 individuals.

\beT{}
\be{} \mathcal{L}^\eta[\hat \eta^\rho]\Rightarrow  \mathcal{L}^{\underline{\theta}}[\zeta]\ee
where $\mathcal{L}^{\underline{\theta}}$ denotes that $\zeta_0(x)=\theta$ for all $x\in S$.

\end{theorem}

At each time we can also consider for each individual alive at time $t$t its {\em line of descent} given by a path $y\in D_{S\times E}([0,\infty))$. This path follows the random walk in reversed time from the time $t$ until the birth time  of the individual. At   that time the parent particle from whom the type has been inherited provides the continuation of the path back to its birth place. This is continued until we reach time $0$. The path at times $<0$ and $>t$ are set to be constant equal to their values at times 0 and $t$.

We then get a measure
\be{} \eta^*_t\in\mathcal{N}(D_{S\times E}(([0,\infty)).\ee
Letting $t$ vary we get the historical interacting Moran model.

Consider

\be{}  \hat \eta^\rho_t(x):= \frac{\eta^{*,\rho}_t(A_{\{x\},t}))}{\eta_t^{*,\rho}(E_{\{x\},t})}1_{\{\eta^{\rho,*}_t(E_{\{x\}})\ne 0\}}\ee

where
\be{} E_{A,t}:=\{y\in D_{S\times E}([0,\infty)):y_t\in A \times E\},\ee

We obtain the historical interacting Fisher-Wright process, $\zeta^*=(\zeta^*_t)_{t\geq 0}$, as the limit of of $\eta^{*,\rho}_t$ as $\rho\to\infty$.  It is a $M(D_{S\times E}(\mathbb{R})$-valued process.

Generator: First consider the generator $\wt A$ of the path process. Its generator is defined on a class of functions
\be{}  \Phi(s,y),\quad s\in\mathbb{R},y\in D_{S\times E}(\mathbb{R})\ee as follows.
For $j=1,\dots,n$ let
\be{} g_j:\mathbb{R}\times S\times E\to \mathbb{R}\ee
where $g_j$ are bounded and $C^1$ in the time variable.  For $0<t_2<t_2<\dots <t_n,\;n\in\N$ define
\be{} \Phi(t,y)=\prod_{j=1}^n g_j(t,y_{t\wedge t_j}).\ee

Let $\mathcal{A}$ denote the algebra of functions generated by these functions.
Denote by $A$ the generator of the random walk on $S\times E$ $Af(x,k)=\sum_{z\in S} (a(x,z)-\delta(x,z))f(z,k)$.  Then we define
\be{} \begin{split}\wt A\Phi(t,y)=&&\prod_{j=1}^k[(\frac{\partial}{\partial t}+A)\prod_{j=k+1}g_j(t,y_{t\wedge t_j})]\\&& +[\frac{\partial}{\partial t}\prod_{j=1}^hg_j(t.y_{t\wedge t_j})][\prod_{j=k+1}^n g_j(t,y_{t\wedge t_j})]\end{split}\ee
We also use the notation $y^r_{\cdot}=y_{\cdot\wedge r}$ for $y\in D_{S\times E}([0,\infty))$ and $\eta^{*,r}$ for a measure concentrated on paths stopped at time $r$. Also $pi^*_S,\pi^*_E$ denote the obvious projections on $D_S(\mathbb{R})$, $D_E(\mathbb{R})$.

 The martingale problem for the historical interacting Fisher-Wright system is analogous to the  historical branching martinale problem (Theorem \ref{HMP1}).  It is formulated as follows.

For $\Phi\in \mathcal{A}$ and $(t,s),\; t\geq s$
\be{GLW-M1} \left\{\langle \zeta^*_y,\Phi(t,\cdot)\rangle-\langle \zeta^*_y,\Phi(s,\cdot)\rangle -\int_s^t\langle \zeta^{*,r}_r
,(\wt A\Phi)(r,\cdot)\rangle dr\right\}_{t\geq s}\ee
is a martingale with increasing process
\be{} \left(\int_s^t\int_{(D_{S\times E}(\mathbb{R}))^2} I^r(y,y')\Phi(r,y)\Phi(r,y')\zeta^{*,r}(dy)(\zeta^{*,r}(dy')-\delta_y(dy'))dr\right)_{t\geq s}.
\ee
where
\be{GLW-M2} I^t(y,y')=\left\{\begin{split} &1\quad\text{if} (\pi^*_S y)_t=(\pi^*_S y')_t
\\& 0,\quad \text{otherwise}\end{split}\right.
\ee

\beT{GWL-2}  (Greven-Winter-Limic \cite{GLW-05}, Theorem 2.)
The martingale problem (\ref{GLW-M1}),(\ref{GLW-M2}) is well-posed and arises as the diffusion limit of $\eta^*$.
\end{theorem}

\begin{remark}  Greven-Limic-Winter develop a particle representation starting with the corresponding look-down process in which particles are assigned labels in a countable subset of $[0,\infty)$. They construct this on a probability space with an additional randomization of labels at lookdowns and in this way construct the interacting Moran models and also the interacting Fisher-Wright processes on a common probability space following the program of Donnelly and Kurtz (\cite{DK-96}).
\end{remark}

\chapter{Mutation-Selection Systems}

The basic mechanisms of population biology are mutation, selection, recombination and genetic drift.  In the previous chapter we concentrated on mutation and genetic drift.  In this chapter we introduce mathematical models of recombination and selection.  However it should be emphasized that these are idealization of highly complex biological processes and there is an immense biological literature including empirical investigation, theoretical models of varying degrees of complexity and simulation studies. For example the concept of {\em fitness} is an abstract  notion that in the biological context  can involve fitness at the level of a single gene, genome or phenotype.  At the level of the genome this can involve the interaction between genes ({\em epistasis}) and various models of such interactions have been proposed (see e.g. Gavrilets \cite{G-04}). One of the continuing issues is the
question of the {\em levels of selection} (see e.g. Brandon and Burian (1984) \cite{BB-84}, Lloyd (2005) \cite{L-05}, Okasha (2006), \cite{O-06} ) which include notions of group selection, kin selection, inclusive fitness (see Hamilton (1964) \cite{Ha-64}) and so on. For example, {\em inclusive fitness} represents to effective overall contribution of an individual including its  own  reproductive success as well as its contribution (due to its behavior) to  the fitness of its genetic kin.

Our aim in this chapter  is to introduce some mathematical aspects of the interplay of mutation, selection and genetic drift.

\section{The infinite population dynamics of  mutation, selection and recombination.}

\subsection{Selection}
The investigation of infinite population models with mutation, recombination and selection leads to an interesting class of dynamical systems (see Hofbauer and Sigmund (1988) \cite{HS-88} and B\"urger \cite{B-89}, \cite{B-00}).  These are obtained as special cases of the general FV process by setting $\gamma =0$ and serve as approximations to systems in which the number of individuals $N$ is very large.

One of the objectives of this chapter is to investigate in one setting  the extent to which the behavior of the finite system differs from that of the infinite system.

Consider an infinite diploid population without mutation or recombination (i.e.
$\gamma =0,\; A=0,\; \rho=0)$ with $K$ types of gametes. \ The unordered pair
$\{i,j\}$ represents the genotype determined by the gametes $i$
and $j$. Let $x_{i}(t)$ be the amount of copies of gamete $i$ in
the population at time $t$ and $p_{i}$ denote the frequency
$p_{i}=\frac{x_{i}}{\sum x_{i}}.$

Let $V(i,j)=V(j,i)$ $=b_{i,j}-d_{i,j}$ where $b_{ij}$ and $d_{ij}$
are the birth and death rates \ of the genotype. The fitness,
$V(i)$ of the ith gamete
is defined by%
\[
V(i)=\sum_{j}p_{j}V(i,j)
\]
and the mean fitness is defined by%
\[
\bar V(p)=\bar{V}=\sum_{i}V(i)p_{i}=\sum_{ij}p_{i}p_{j}V(i,j)\text{.}%
\]
Then the population sizes $x_{i}$ satisfy the equations%
\[
\dot{x}_{i}=x_{i}\sum_{j}V(i,j)\frac{x_{j}}{|x|},\;i=1,\dots,K
\]

\begin{proposition}
The proportions $\{p_{i}\}$ satisfy the equations:%
\[
\dot{p}_{i}=p_{i}(V(i)-\bar{V}),\;i=1,\dots,K
\]
\bigskip
\end{proposition}

\begin{proof}
This can be derived from the $\dot{x}$ equations by the
substitution $x_{i}=|x|p_{i}$ giving
\[
\dot{p}_{i}|x|+p_{i}(\sum\dot{x}_{j})=|x|p_{i}V(i)
\]
which yields%
\[
\dot{p}_{i}+p_{i}(\sum_{j}p_{j}V(j))=p_{i}V(i)
\]
and the result immediately follows.
\end{proof}

\subsection{Riemannian structure on $\Delta_{K-1}$}

The deterministic differential equations of selection have played an important role in the development
of population genetics.  A useful tool in their analysis was a geometrical approach developed by Shahshahani and Akin. We
next give a brief introduction to this idea.

Let $M$ be a smooth manifold. The tangent space at $x$, $T_{x}M$
can be identified with the space of tangents at $x$ to all smooth
curves through $x$. The tangent bundle $TM=\{(p,v):p\in M,v\in
T_{p}M\}.$

\begin{definition}
A Riemannian metric on $M$ is a smooth tensor field%
\[
g:C^{\infty}(TM)\otimes C^{\infty}(TM)\rightarrow
C_{0}^{\infty}(M)
\]
such that for each \ $p\in M,$%
\[
g(p)|_{T_{p}M\otimes T_{p}M}:T_{p}M\otimes T_{p}M\rightarrow\mathbb{R}%
\]
with%
\[
g(p):(X,Y)\rightarrow\left\langle X,Y\right\rangle _{g(p)}%
\]
where $\left\langle X,Y\right\rangle _{g(p)}\;$is an inner product
on $T_{p}M$.
\end{definition}

\begin{definition}
The directional derivative in direction $v$ is defined by%
\begin{align*}
\partial_{v}f(x)  &  =\lim_{t\rightarrow0}\frac{f(x+tv)-f(x)}{t}\\
&  =\sum v_{i}\frac{\partial f(x)}{\partial x_{i}}%
\end{align*}
The gradient $\nabla_{g}f(x)$ is defined by%
\[
\left\langle \nabla_{g}f(x),v\right\rangle
_{g}=\partial_{v}f(x)\;\;\forall v\in T_{x}M.
\]
\end{definition}

\begin{example}
Consider the $d$-dimensional manifold $M=\mathbb{R}^{d}$ and $\mathbf{a}%
(\cdot)$ be a smooth map from $M$ to
$\mathbb{R}^{d}\otimes\mathbb{R}^{d}$
($(d\times d)$-matrices). We will write%
\begin{align*}
\mathbf{a}(x)  &  =(a_{ij}(x))\\
\mathbf{a}^{-1}(x)  &  =(a^{ij}(x))
\end{align*}
Assume that%
\[
\sum a^{ij}(x)u_{i}u_{j}\geq\gamma\sum u_{j}^{2},\;\gamma>0.
\]
The tangent space $T_{\mathbf{x}}M\approxeq\mathbb{R}^{d}$ and we
define a Riemannian metric on $M$ by
\[
g_{\mathbf{a}(x)}(\mathbf{u,v)}:=\sum_{i,j=1}^{d}a_{ij}(x)u^{i}v^{j}.
\]
The associated Riemannian gradient and norm are%
\begin{align*}
(\nabla_{\mathbf{a}}f)^{i}  &  =\sum_{j}a^{ij}\frac{\partial
f}{\partial
x_{j}}\\
\Vert u\Vert_{\mathbf{a}(x)}^{2}  &
=\sum_{ij}a_{ij}(x)u^{i}u^{j}.
\end{align*}
\end{example}

\subsubsection{The Shahshahani metric and gradient on $\Delta_{K-1}$}

Let $M_{K}=\mathbb{R}_{+}^{K}:=\{x\in\mathbb{R}^{K},x=(x_{1},\dots
,x_{K}),\;x_{i}>0$ for all $i\}$ is a smooth K-dimensional
manifold.

Shahshahani introduced the following Riemannian metric on $M_{K}$%
\begin{align*}
\left\langle u,v\right\rangle _{g}  &  =g_{x}(u,v):=\sum_{i=1}^{K}%
|x|\frac{u_{i}v_{i}}{x_{i}}\\
|x|  &  =\sum x_{i}%
\end{align*}
\thinspace$\Vert\;\Vert_{g}$ and $\nabla_{g}F$ \ will denote the
corresponding norm and gradient. We have
\[
(\nabla_{g}F)^{i}=\sum_{i}\frac{x^{i}}{|x|}\frac{\partial F}{\partial x^{i}%
}\frac{\partial}{\partial x_{{}}^{i}}%
\]

Recall that the simplex  $\Delta_{K-1}:=\{(p_1,\dots,p_K):p_i\geq 0,\; \sum_{i=1}^K p_i=1\}$.
The interior of the simplex
$\Delta_{K-1}^{0}=\mathbb{R}_{+}^{K}\cap \Delta_{K-1}$ is a
$(K-1)$-dimensional submanifold of $M_{K}$. \ We denote by
$T_{p}\Delta_{K-1}^{0}$ the tangent space to \ $\Delta_{K-1}^{0}$
at $p$. Then $g$ induces a Riemannian metric on
$T_{p}\Delta_{K-1}^{0}$.

\underline{Basic Facts}

We have the Shahshahani inner product on $\Delta{K-1}$ at a point $p \in \Delta{K-1}$:
\be{}  \langle u,v\rangle_p= \sum_{i=1}^K \frac{u_iv_i}{p_i}.\ee

1.  $T_{p}\Delta_{K-1}^{0}$ can be viewed as the subspace of
$T_{p}M_{K}$ of vectors, $v$,  satisfying \ $\left\langle
p,v\right\rangle _{g}=0$  if we identify $p$ with an element
of $T_{p}M_{K}.$

Proof. Recall that $T_{p}\Delta_{K-1}^{0}$ is given by tangents to
all smooth
curves lying in $\Delta_{K-1}^{0}$. \ Therefore if $v\in T_{p}\Delta_{K-1}%
^{0},\;$then $v=q-p$ where $p,q\in\Delta_{K-1}^{0}$ and therefore
$\sum _{i=1}^{K}v_{i}=0.$ Therefore, \
\[
\sum_{i}p_{i}\frac{1}{p_{i}}v_{i}=0.%
\]

2. \ If $F:\Delta_{K-1}^{0}\rightarrow\mathbb{R}$ is smooth, then the Shahshahani gradient is
\[
(\nabla_{g}F)_{i}=p_{i}\left(  \frac{\partial F}{\partial p_{i}}-\sum_{j}%
p_{j}\frac{\partial F}{\partial p_{j}}\right)  .
\]

Proof. From the definition, $\nabla_{g}F$ is the orthogonal
projection on the subspace $T_{p}\Delta_{K-1}^{0}$ of
\[
(\nabla_{g}F)_{i}=p_{i}\frac{\partial F}{\partial p_{i}}%
\]

and therefore we must have $\sum_i(\nabla_{g}F)_{i}=0.$ This then
gives the result.

\begin{remark}
This (Shahshahani) gradient coincides with the gradient on
$\Delta_{K-1}$ associated with the $K$-alleles Wright-Fisher model
and appears in the description of the rate function for large
deviations from the infinite population limit (see below). \
\end{remark}

\begin{theorem}
The dynamical system $\{\mathbf{p}(t):t\geq0\}\;$is given by%
\[
\mathbf{\dot{p}}(t)=\frac{1}{2}(\nabla_{g(\mathbf{p}(t))}\bar{V}%
)(\mathbf{p}(t)).
\]
\bigskip
\end{theorem}

\begin{proof}
From the above, applying the Shahshahani gradient to $\bar{V}$, we get%
\begin{align*}
(\nabla_{g}\bar{V})_{i}  &  =2\left(  p_{i}V(i)-p_{i}\sum p_{j}V(j)\right) \\
&  =2p_{i}(V(i)-\bar{V}).
\end{align*}
\end{proof}

\begin{theorem}
(Fisher's Fundamental Theorem)\newline (a) Mean fitness increases
on the trajectories of $p(t)$.\newline (b) The rate of change of
the mean $\bar {V}(t)$ along orbits is proportional to the
variance.\newline (c) At an equilibrium point the eigenvalues of
the Hessian must be real.\newline
\end{theorem}

\begin{proof}
(a) follows immediately from (b).

(b)
\begin{align*}
d\bar{V}(t)  &  =\left\langle
\nabla_{g}\bar{V}(\mathbf{p}(t)),\mathbf{\dot
{p}}(t)\right\rangle _{g(\mathbf{p}(t))}\\
&  =2\left\langle \dot{p}(t),\dot{p}(t)\right\rangle _{g(\mathbf{p}%
(t))}=2\left(  \sum_{i}p_{i}(t)(V(i)-\bar{V}(t))^{2}\right) \\
&  =2\left(  \sum_{i}p_{i}(t)V(i)^{2}-\bar{V}(t)^{2}\right) \\
&  =2Var_{\mathbf{p}(t)}(\mathbf{V})\geq0.
\end{align*}
(b) It also follows from the gradient form that the Hessian is
symmetric (matrix of mixed second partials of $\bar{V}$).
\end{proof}

\begin{theorem}
(Kimura's Maximum Principle) ``Natural selection acts so as to
maximize the rate of increase in the average fitness of the
population.''
\end{theorem}

\begin{proof}
This simply follows from the property that the directional
derivative $\partial_{v}\bar{V}$ is maximal in the direction of
the gradient.
\end{proof}

\begin{example}
Consider a two type ($\{1,2\}$) population with frequencies $(p_1,p_2)=(p,1-p)$.
\[
V(i,j)=av(i)+av(j)+c\delta_{ij}%
\]
(When $c=0$ we have the additive (or haploid) model. When $a=0$
and $c>0$ we
have the heterozygote advantage model.)\newline In this case%
\begin{align*}
\bar{V}(p_{1},p_{2})  &  =ap_{1}v(1)+ap_{2}v(2)+cp_{1}p_{2}\\
&  =V(p,1-p)=ap(v(1)-v(2))+av(2)+cp(1-p)
\end{align*}
Then depending on the choice of $a,c,v(1),v(2),$the optimum value
of $p$ can range between $0$ and $1.$
\end{example}

\begin{remark} For the multilocus situation there is the Fisher-Price-Ewens version (e.g. Frank (1997) \cite{F-97}, Ewens \cite{E-04}). This is also related to {\em the secondary theorem of natural selection} of Robertson (1966) \cite{R-66}  which relates the rate of change of a quantitative character under selection in terms of the covariance of the character and fitness.

\end{remark}

The above equations are special cases of the class of {\em replicator equations} of the form

\be{} \frac{dp_i(t)}{dt}= p_i(t)(f_i(\mathbf{p}(t))-\sum p_if_i(\mathbf{p}(t)),\quad i=1,\dots, K\ee
where $\{f_i(\mathbf{p}\}_{i=1m\dots, K}$ is a vector field on $\Delta_{K-1}$.
In the linear case  $f_i(\mathbf{p})=\sum_j a_ij p_j$ these are equivalent to the
the Lotka-Volterra equations
\be{} \frac{dx_i(t)}{dt}= x_i(t)\left( r_i +\sum_{j=1}^n K_{ij}x_j(t)\right),\quad i=1,\dots, K-1 \ee
by setting  $p_i(t)= \frac{x_i(t)}{\sum_i x_i(t)}$.

\subsection{Mutation-Selection}

The replicator equations that include both mutation and selection are given by

\be{} \frac{dp_i(t)}{dt}= p_i(t)(V(i)-\bar V)+m(\sum_{j\ne i}q_{ji}p_j -p_i)\ee
where $m$ is the mutation rate and for each $j$, $q_{ji},\;i\ne j$ is the probability that type $j$ mutates to type $i$
and $\sum_{i\ne j} q_{ji}=1$.

\beT{}
The mutation-selection  dynamical system is a Shahshahani gradient system  if and only if
 \be{} q_{ji}=q_i\;\forall \;j,\ee
  (that is type-independent mutation as in the infinitely many alleles model). In the latter case the
potential is

\be{} W(p)=  \bar V(p)-H(q|p),\quad H(q|p)=-\sum_{i=1}^n q_i\log p_i.\ee
\end{theorem}
\begin{proof}
See Hofbauer and Sigmund \cite {HS-88}, Chapt. VI, Theorem 1.
\end{proof}

We will see below that there is a far-reaching analogue of this for the stochastic (finite population) generalizations.

\begin{remark}
In general the deterministic mutation-selection equations  are  not a gradient system and can
exhibit complex dynamics - for example, a stable limit cycle (Hofbauer and Sigmund \cite{HS-88}, 25.4).
An interesting special case is the diploid case with three types - two favourable and mutation.  Baake \cite{B-95} showed that these can exhibit stable limit cycles. Hofbauer (1985) \cite{H-85}  also showed this for selection mutation models with cyclic mutation.

Smale \cite{S-76} pointed out that for n types, $n\geq 5$,   dynamical systems on the simplex can have complex behaviour.  He gave an example that  ``may not be
approximated by a structurally stable, dynamical system, or it may have
strange attractors with an infinite number of periodic solutions''. Some further basic results on competitive  systems are covered by Hirsch (1982), (1985), (1988) \cite{Hir-82} and Liang and Jiang (2003) \cite{LJ-03}.

\end{remark}

\subsection{Multiple loci and recombination}

Multiloci models give rise to dynamical systems that have been extensively studied. They give rise to a large class of dynamical systems that can have complex behaviour.
Akin \cite{A-83} analyzed the simplest two loci model with selection and recombination and proved that in general this is not a gradient system  and that periodic orbits can exist.  We briefly sketch the simplest example.

Consider a two-loci model with two alleles at each loci. We denote the types by  $1=AB,2=Ab,3=aB,4=ab$ and with gamete frequencies \be{} p_{AB},p_{Ab},p_{aB},p_{ab},\quad p_A=p_{AB}+p_{Ab},\; p_B=p_{AB}+p_{aB},\; p_a=p_{aB}+p_{ab},\; p_b=p_{Ab}+p_{ab}.\ee
 Then the {\em measure of linkage disequilibrium} is defined as
\be{}  d:=p_{AB}p_{ab}-p_{Ab}p_{bA}\ee
so that $d=0$ if $p_{AB}=p_Ap_B$, etc. The diploid fitness function is denoted by $V(i,j)$. Some natural assumptions are that
\be{} m_{ij}=m_{ji},\quad m_{14}=m_{23}=0.\ee
There are 10 zygotic types $AB/AB,Ab/AB,\dots,ab/ab$ and the corresponding fitness table
\be{}
\begin{array}{ccccc}
     & AB & Ab & aB & ab \\
  AB & w_{11} & w_{12} & w_{13} & w_{14} \\
  Ab & w_{21} & w_{22} & w_{23} & w_{24} \\
  aB & w_{31} & w_{32} & w_{33} & w_{34} \\
  ab & w_{41} & w_{42} & w_{43} & w_{44}
\end{array}
\ee
The recombination vectorfield
\be{} R=rbd\xi_i, \;i=1,2,3,4\ee
where $r$ is the recombination rate, $b$ is the birth rate for double heterozygotes, $d$ is the linkage disequilibrium and \be{}\xi=(1,-1,-1,1)\ee so that
\be{} d\xi= p-\pi(p)\ee
where $\pi(p)$ has the same marginals as $p$ but in linkage equilibrium (independent loci).

The system of differential equations for the frequencies of types $1,2,3,4$ with selection and recombination are
\be{}  \frac{dp_i}{dt}=p_i(V(i)-\bar V)-rbd\xi_i\,\quad i=1,2,3,4 \ee
where
\be{}V(i)=\sum_{j=1}^4p_jV(i,j),\quad \bar V=\sum_{i=1}^4 p_iV(i),\quad d=p_1p_4-p_2p_3.\ee

In the case $V\equiv 0$ the system approaches linkage equilibrium.  However Akin \cite{A-83} showed that there exist fitness functions $V$ and parameters $b,r$ such that the  system  exhibits a Hopf bifurcation leading to cyclic behaviour. More generally, multilocus systems can exhibit many types of complex behaviour (see for example, Kirzhner, Korol and Nevo (1996) \cite{KKN-96} and Lyubich and Kirzhner (2003) \cite{LK-03}).

\section{Infinitely many types Fleming-Viot}

We now consider the Fleming-Viot process with selection and recombination and establish uniqueness
using a dual representation of Ethier and Kurtz.

In Chapter 6 we showed that the martingale problem for the Fleming-Viot process with mutation selection and recombination is well-posed and defines a $\mathcal{P}(E)$-valued Markov diffusion process. In this chapter we focus on mutation and selection but also give a brief introduction to some aspects of recombination. In evolutionary theory mutation  plays an important role in producing novelty and maintaining diversity while selection eliminates deleterious mutations and makes possible the emergence and fixation of rare advantageous mutations. From a more abstract viewpoint this  can be viewed as a search process which generates new  information.

\subsection{Dual representation with mutation, selection and recombination}

As above we consider the mutation generator $A$ and   the bounded diploid fitness function
For $V\in \mathcal{B}_{\rm{sym}}(E\times E)$, set $\bar
V=\sup_{x,y,z}|V(x,y)-V(y,z)|$.  Without loss of generality we can assume that $\bar V =1$ and
define the selection coefficient $s>0$ and
selection operators
\be{}V_{im}f(x_1,\dots,x_{m+2})= ({V(x_i,x_{m+1})-V(x_{m+1},x_{m+2})})f(x_1,\dots,x_m).\ee

 For $f\in
\mathcal{D}(A^{(n)})\cap\mathcal{B}(E^n)$, define $F(f,\mu)=\int f
d\mu^n$ and

\be{}\begin{split} GF(f,\mu)&&= F(A^{(n)}f,\mu)+\gamma_{1\leq i<j\leq
n}\left(F(\Theta_{ij}f,\mu)-F(f,\mu)\right)+s\sum_{i=1}^nF(V_{in}f,\mu).\end{split} \ee

For $f\in C_{\rm{sim}}(E^\N)$, with $\mathfrak{n}(f)=n$, and
$f\in\mathcal{D}(A^n)\cap\mathcal{B}(E^n)$, let

\be{}  Kf:=\sum_{i=1}^n A_if+\gamma\sum_{j=1}^n\sum_{k\ne
j}[\Theta_{jk}f-f]+s\sum_{i=1}^n[V_{in}f-f].\ee

where $\wt{\Theta}_{jk},\; \mathfrak{n}(f)$ are defined as in section \ref{sdrfv}.

If  $\beta(f):= s\mathfrak{n}(f)$, then \be{}
GF(f,\mu)=F(Kf,\mu)+\bar V(\mathfrak{n}(f))F(f,\mu),\ee and
$\sup_{\mu\in M_1(E)}|F(Kf,\mu)|\leq \text{const}\cdot \mathfrak{n}(f).$

Let $\rho\geq 0$ and
$\eta(x_1,x_2,\Gamma)$ be a transition function from $E\times E \to
E$. For $i=1,\dots,m$ define $R_{im}:\mathcal{B}(E^m)\to
\mathcal{B}(E^{m+1})$ by

\be{} R_{im}f(x_1,\dots,x_{m+1})=\int
f(x_1,\dots,x_{i-1},z,x_{i+1},\dots,x_m)\eta(x_i,x_{m+1},dz)\ee
and assume that $R_{im}:C_b(E^m)\to C_b(E^{m+1})$.  The $R_{im}$
are called the recombination operators for the process and $\rho$
is called the recombination rate.

Given  $V\in \mathcal{B}_{\rm{sym}}(E\times E)$, with $\bar
V:=\sup_{x,y,z}|V(x,y)-V(y,z)|<\infty$,  define the
selection operators
\be{}V_{im}f(x_1,\dots,x_{m+2})=\frac{V(x_i,x_{m+1})-V(x_{m+1},x_{m+2})}{\bar
V}f(x_1,\dots,x_m)\quad \text{for  } i=1,\dots,m.\ee

 For $f\in
\mathcal{D}(A^{(n)})\cap\mathcal{B}(E^n)$, define $F(f,\mu)=\int f
d\mu^n$ and

\be{}\begin{split} GF(f,\mu)&&= F(A^{(n)}f,\mu)+\gamma \sum_{1\leq i<j\leq
n}\left(F(\wt\Theta_{ij}f,\mu)-F(f,\mu)\right)\\&&
+\rho\sum_{i=1}^n\left(F(R_{in}f,\mu)-F(f,\mu)\right)+\bar
V\sum_{i=1}^nF(V_{in}f,\mu).\end{split} \ee

For $f\in C_{\rm{sim}}(E^\N)$, with $\mathfrak{n}(f)=n$, and
$f\in\mathcal{D}(A^n)\cap\mathcal{B}(E^n)$, let

\be{}  Hf:=\sum_{i=1}^n A_if+\gamma\sum_{j=1}^n\sum_{k\ne
j}[\wt\Theta_{jk}f-f]+\rho\sum_{i=1}^n[R_{in}f-f]+\bar
V\sum_{i=1}^n[V_{in}f-f].\ee

If  $\beta(f):=\bar V\mathfrak{n}(f)$, then \be{}
GF(f,\mu)=F(Hf,\mu)+\beta(f))F(f,\mu),\ee and
$\sup_{\mu\in M_1(E)}|F(Hf,\mu)|\leq \text{const}\cdot \mathfrak{n}(f).$

\beT{gfvdual} Let $G$ satisfy the above conditions and assume that the
mutation process with generator $A$ has a version with sample
paths in $D_E[0,\infty)$. Then for each
$\mu\in  \mathcal{P}(E)$ there exists a unique solution $P_\mu$ of the martingale
problem for $G$.
\end{theorem}

\begin{proof} (Ethier-Kurtz (1987) \cite{EK-87}) The uniqueness will be proved by constructing a dual representation.

 Let $N$ be a jump Markov process taking
non-negative integer values with transition intensities
\be{}q_{m,m-1}=\gamma m(m-1),\; q_{n,m+2}=\bar V m,\;
q_{m,m+1}=\rho m,\; q_{i,j}=0\text{   otherwise}. \ee For $1\leq
i\leq m$, let  $\{\tau_k\}$ be the jump times of $N$, $\tau_0=0$,
and let $\{\Gamma_k\}$ be a sequence of random operators which are
conditionally independent given $M$ and satisfy \be{}
P(\Gamma_k=\Theta_{ij}|N)=\frac{2}{N(\tau_k-)N(\tau_k)}\mathbf{1}_{N(\tau_k-)-N(\tau_k)=1},\quad
1\leq i<j\leq N(\tau_k-)\ee \be{}
P(\Gamma_k=R_{im}|N)=\frac{1}{m}\mathbf{1}_{\{N(\tau_k-)=m,N(\tau_k=m+1)\}}\ee

\be{}
P(\Gamma_k=V_{im}|N)=\frac{1}{m}\mathbf{1}_{\{N(\tau_k-)=m,N(\tau_k)=m+2\}}.\ee

For $f\in C_{\rm{sim}}(E^\N)$, define the
$C_{\rm{sim}}(E^\N)$-valued process $Y$ with $Y(0)=f$ by

\be{}
Y(t)=S_{t-\tau_k}\Gamma_kS_{\tau_k-\tau_{k-1}}\Gamma_{k-1}\dots
\Gamma_1 S_{\tau_1}f,\quad \tau_k\leq \tau_{k+1}.\ee Then for any
solution $P_\mu$ to the martingale problem for $G$ and $f\in
C_{\rm{sim}}(E^\N)$ we get the FK-dual representation

\be{} P_\mu[F(f,X(t))]=Q_f\left[F(Y(t),\mu)\exp\left(\bar
V\int_0^t\mathfrak{n}(Y(u))du\right)\right]\ee which establishes that
the martingale problem for $G$ is well-posed. Since the function
$\beta(f)=\bar V\mathfrak{n}(f)$ is not bounded we must verify condition
(\ref{DUI}).  This follows from the following lemma due to Ethier and Kurtz
(1998) \cite{EK-98}, Lemma 2.1.

\end{proof}
\beL{ekl} Let $N(t)=\mathfrak{n}(Y(t))$ be as above,
$\tau_K:=\inf\{t:N(t)\geq K\}$ and $\theta>0$. Then there exists a
function $R(n)\geq \text{const}\cdot n^2$ and a constant $L>0$
such that \be{} E \left[R(N(t\wedge
\tau_K))\exp\left(\theta\int_0^{t\wedge\tau_K}N(s)ds\right){
|}N(0)=n\right]\leq F(n)e^{Lt},\quad\forall\;K\geq1,\ee and given
$N(0)=n$, $\left\{N(t\wedge\tau_K)\exp\left(\bar
V\int_o^{t\wedge\tau_K} N(s)ds\right):\;K\geq 1\right\}$ are
uniformly integrable.

\end{lemma}
\begin{proof} The integer-valued process $N(t)$ is a birth and death process with jump rates
\be{} m\to (m+1)\text{  rate  }\rho m,\quad m\to m+2 \text{  rate  }\bar Vm,\quad m\to m-1 \text{  rate  }\gamma m(m-1).\ee
Let $Q$ denote the corresponding generator.
Take $R(m):=(m!)^{\beta}$, with $\beta<\frac{1}{2}$. Then

\bean{} &&QR(m)+\theta mR(m)\\&&=\gamma m(m-1)(R(m-1)-R(m))+\rho
m(R(m+1)-R(m))\\&&\qquad\qquad +\theta m(R(m+2)-R(m))+\theta m
R(m)\\&& = -\gamma O(m^{2})(m!)^\beta +\rho O(m^2)(m!)^\beta+\theta O(m^{1+2\beta})(m!)^\beta
\eean Since the negative term dominates for large $m$ if $0<\beta<\frac{1}{2}$ and $\gamma>0$, we can
choose $L>0$ such that \be{}  QR(m)+\theta mR(m)\leq L.\ee The
optional sampling theorem implies that for
$\tau_K:=\inf\{t:N(t)\geq K\}$ and $N(0)=m$ \bean{}
&&E\left[\exp\left(\theta\int_0^{t\wedge
\tau_k}N(s)ds\right)|N(0)=m\right]\\&& \leq
E\left[R(N(t\wedge\tau_k)\exp\left(\theta\int_0^{t\wedge
\tau_k}N(s)ds\right)|N(0)=m\right]\\&&\leq R(m)+
E\left[\int_0^{t\wedge\tau_k}\exp\left(\theta\int_0^uN(s)ds\right)(QR(N(u))+\theta
N(u)R(N(u)))du|N(0)=m\right]
\\&&\leq R(m)+L E\left[\int_0^{t\wedge\tau_k}\exp\left(\theta\int_0^uN(s)ds\right)du|N(0)=m\right]\eean
and the lemma follows by Gronwall's inequality.

\end{proof}







\subsection{Girsanov formula for Fleming-Viot with Mutation and Selection}

Recall that the Fleming-Viot martingale problem
$\mathbb{MP}_{(A,\gamma Q,0)}$
corresponds to the case%
\[
\left\langle M(A),M(A)\right\rangle
_{t}=\gamma\int_{0}^{t}Q(X_{s};A,A)ds
\]
where%
\[
Q(\mu;dx,dy)=\mu(dx)\delta_{x}(dy)-\mu(dx)\mu(dy).
\]
and that $M$ is a worthy martingale measure.

Now consider a time-dependent diploid fitness function \
$V:[0,\infty)\times E\times E\rightarrow\mathbb{R}$ with $\Vert
V\Vert_{\infty}<\infty$. . Then the FV martingale problem
$\mathbb{MP}_{(A,Q,V)}$is
\begin{align*}
&  M^{V}(\phi)_{t}\\
&  :=\left\langle X_{t},\phi\right\rangle
-\int_{0}^{t}\left\langle
X_{s},A\phi\right\rangle ds\\
&  \;\;\;\;\;\;\;\;-\int_{0}^{t}\int\left[  \int
V(s,x,y)X_{s}(dy)-\int\int
V(s,y,z)X_{s}(dy)X_{s}(dz)\right]  \phi(x)X_{s}(dx)ds\\
&  =\left\langle X_{t},\phi\right\rangle -\int_{0}^{t}\left\langle X_{s}%
,A\phi\right\rangle ds\\
&  \;\;\;\;\;\;-\int_{0}^{t}\int\int\left[  \left(  \int\frac{V(s,y,z)}%
{\gamma}X_{s}(dz)\right)  \gamma Q(X_{s},dx,dy)\right]  \phi(x)ds
\end{align*}%
\[
\left\langle M^{V}(\phi)\right\rangle
_{t}=\gamma\int_{0}^{t}\int\int \phi(x)\phi(y)Q(X_{s},dx,dy)ds.
\]
We then apply Theorem \ref{DGIR} to conclude that this martingale problem
has a unique solution $\mathbb{P}^{V}$ and that the Radon-Nikodym
derivative
\[
Z_{t}^{V}:=\frac{d\mathbb{P}^{V}}{d\mathbb{P}^{0}}|_{\mathcal{F}_{t}}%
\]
where $\mathbb{P}^{0}$ is the unique solution to
$\mathbb{MP}_{(A,\gamma
Q,0)}$ is given by%
\begin{align*}
Z_{t}^{V}  &  :=\exp\left(  \frac{1}{\gamma}\int_{0}^{t}\int V(s,X_{s}%
,y)M^{0}(ds,dy)\right. \\
&  \left.  -\frac{1}{2\gamma^{2}}\int_{0}^{t}\int\int V(s,X_{s},x)V(s,X_{s}%
,y)\gamma Q(X_{s};dx,dy)ds\right)  .
\end{align*}
where we write%
\[
V(s,X_{s},x)=\int V(s,z,x)X_{s}(dz).
\]

\section[Long-time behaviour]{Long-time behaviour of systems with finite population resampling, mutation and selection}

Systems with finite population resampling can have rather different long-time behaviour than the corresponding
infinite population systems.  One essential difference is that even high fitness types can be lost due to resampling and in the absence of
mutation the system can eventually become unitype.  On the other hand if the mutation process can regenerate all types, then the system can reach
equilibrium in which all types are present.  We now consider these two situations.

\subsection{Fixation in finite population systems without mutation}

In the previous section we have considered the infinite population
system with selection but no mutation. \ In this case Fisher's
fundamental theorem states that such a system evolves to one of
maximal population fitness. \ But what happens in the finite
population case, $\gamma>0$? \ We first observe that if
$V\equiv0,$ then $\{X_{t}(A):t\geq0\}$ is a bounded martingale and
\be{}
X_{t}(A) \ttoo  \left\{ \begin{split} &  1\text{ with probability }%
X_{0}(A)\\
 &  0\text{ with probability }(1-X_{0}(A)).\end{split}\right.
\ee
Therefore%
\[
X_{t}\overset{t\rightarrow\infty}{\Longrightarrow}\delta_{x}\text{
with }x\in A\text{ with probability }X_{0}(A)
\]
that is, the system experiences ultimate ``fixation''. If we add
selection to this, ultimate fixation still occurs. \ However if
$\gamma$ is small then the tendency is for the limiting types to
be those of higher fitness.

\subsection{The Equilibrium Infinitely Many Alleles Model with
Selection}

In order to have a non-degenerate equilibrium a source of new
types through mutation is required. In this section we consider the type independent
infinitely many alleles mutation together with selection.   \ If $\nu_{0}$
is a non-atomic measure, then mutation always leads to a new type
and thus provides a mechanism to guarantee sufficient diversity on
which selection can act.

We first obtain the ergodic theorem in this case.

\beP{} (Ergodicity of IMA Fleming-Viot with selection)\\
Consider the infinitely many allele type Fleming-Viot process on $[0,1]$, mutation source $\nu_0\in\mathcal{P}([0,1])$ and $X_0=\mu$. Then
\be{} X_t\Rightarrow X_\infty\ee
where $X_\infty$ is a random probability measure on $[0,1]$.  The process is ergodic and $\mathcal{L}(X_\infty)$ is the unique stationary measure.
\end{proposition}
\begin{proof}
This follows immediately from the dual representation.  Note that for the type-independent mutation
\be{} Af(x)=\int_0^1 f(y)\nu_0(dy)-f(x)\ee
and then
\be{} q_{1,0}>0.\ee
Due to the quadratic death rate the process $\mathfrak{n}(Y(t)$ is recurrent and starting at $n\geq 1$ returns to $1$ with probability 1. But then the probability that the dual reaches the trap $C^0([0,1])$ is one and the limiting dual is a constant.This implies the result.
\end{proof}

We now identify the resulting stationary measure.

Let \ $\mathbb{P}_{\infty}^{0}$ denote the probability measure on
$C_{\mathcal{P}([0,1]))}(-\infty,\infty)$ corresponding to the reversible
stationary measure, with one dimensional marginal distribution
$\Pi_{\gamma}^{0}(d\mu),$ for the neutral infinitely many alleles
model (recall the representation in terms of the Moran subordinator). Assume that $V$ is symmetric and $V(s,x,y)=V(x,y)=V(y,x).$

The following results is the infinitely many types analogue of a result of Wright \cite{W-49}.

\beT{}
The infinitely many alleles model with selection has a reversible
stationary
measure given by%
\[
\Pi_{\gamma}^{V}(d\mu)=\frac{1}{Z}e^{\frac{V(\mu)}{\gamma}}\Pi_{\gamma}%
^{0}(d\mu)
\]
where $Z$ is a normalizing constant.
\end{theorem}

\begin{proof}
Let \ $X_{0}$ have distribution
\[
\frac{1}{Z}e^{\frac{V(X_{0})}{\gamma}}\Pi_{\gamma}^{0}(dX_{0})
\]
Recall that to verify that this is a reversible equilibrium
measure it suffices to show that for any two continuous functions,
$f$ and $g,$ on
$[0,1]$%
\[
\mathbb{P}_{\infty}(f(X_{0})g(X_{t}))=\mathbb{P}_{\infty}(g(X_{0})f(X_{t})).
\]

But%

\begin{align*}
&  \mathbb{P}_{\infty}(f(X_{0})g(X_{t}))\\
&  =\frac{1}{Z}\int f(X_{0})e^{\frac{V(X_{0})}{\gamma}}g(X_{t})\mathbb{P}%
_{X_{0}}^{V}(d\{X_{s}:0\leq s\leq t\})\Pi_{\gamma}^{0}(dX_{0})\\
&  =\frac{1}{Z}\int f(X_{0})e^{\frac{V(X_{0})}{\gamma}}Z_{t}^{V}%
\mathbb{P}_{\infty}^{0}(dX_{\cdot})g(X_{t}).
\end{align*}

By Theorem \ref{DGIR}
\begin{align*}
Z_{t}^{V}  &  :=\exp\left(  \frac{1}{\gamma}\int_{0}^{t}\int V(X_{s}%
,y)M^{0}(ds,dy)\right. \\
&  \left.  -\frac{1}{2\gamma^{2}}\int_{0}^{t}\int\int V(X_{s},x)V(X_{s}%
,y)\gamma Q(X_{s};dx,dy)ds\right)
\end{align*}
where%
\begin{align*}
M_{s}^{0}(dy)  &  =X_{s}-\int_{0}^{s}A^{\ast}X_{u}du\\
&  =X_{s}-\int_{0}^{s}c[\nu_{0}-X_{u}]du
\end{align*}
As a preparation, note that by Ito's lemma,%
\begin{align*}
&  d_{t}(\int\int V(x,y)X_{t}(dx)X_{t}(dy))\\
&  =\int\int V(x,y)X_{t}(dx)d_{t}X_{t}(dy)+\int\int V(x,y)X_{t}(dy)d_{t}%
X_{t}(dx)\\
&  +\int\int V(X_{s},x)V(X_{s},y)\gamma Q(X_{s};dx,dy)
\end{align*}
Hence by symmetry in $x$ and $y$ and Ito's Lemma,%
\begin{align*}
\frac{1}{\gamma}\int_{0}^{t}\int\int V(x,y)X_{s}(dx)d_{s}X_{s}(dy)
&
=\frac{1}{2\gamma}[\int\int V(x,y)X_{t}(dx)X_{t}(dy)\\
&  \;\;\;-\int\int V(x,y)X_{0}(dx)X_{0}(dy)]\\
&  \;\;\;-\frac{1}{2\gamma^{2}}\int_{0}^{t}\int\int
V(X_{s},x)V(X_{s},y)\gamma Q(X_{s};dx,dy)
\end{align*}
Therefore
\begin{align*}
&  \log(e^{\frac{1}{\gamma}V(X_{0})}Z_{t}^{V})\\
&
=\int_{0}^{t}\int\frac{1}{\gamma}V(X_{s},y)M^{0}(ds,dy)+\frac{1}{\gamma
}\int\int V(x,y)X_{0}(dx)X_{0}(dy)\\
&  -\frac{1}{2\gamma^{2}}\int_{0}^{t}\int\int
V(X_{s},x)V(X_{s},y)\gamma
Q(X_{s};dx,dy)ds\\
&  =\frac{1}{\gamma}\int_{0}^{t}\int\int V(x,y)X_{s}(dx)d_{s}X_{s}%
(dy)-\frac{c}{\gamma}\int_{0}^{t}\int\int V(x,y)X_{s}(dx)(\nu_{0}%
(dy)-X_{s}(dy))ds\\
&  \left.  -\frac{1}{2\gamma^{2}}\int_{0}^{t}\int\int V(X_{s},x)V(X_{s}%
,y)\gamma Q(X_{s};dx,dy)ds\right)  +\frac{1}{\gamma}\int\int V(x,y)X_{0}%
(dx)X_{0}(dy)\\
&  =\frac{1}{2\gamma}[\int\int V(x,y)X_{t}(dx)X_{t}(dy)+\int\int
V(x,y)X_{0}(dx)X_{0}(dy)]\\
&  -\frac{1}{\gamma^{2}}\int_{0}^{t}\int\int
V(X_{s},x)V(X_{s},y)\gamma
Q(X_{s};dx,dy)ds\\
&  -\frac{c}{\gamma}\int_{0}^{t}\int\int V(x,y)X_{s}(dx)(\nu_{0}%
(dy)-X_{s}(dy))ds
\end{align*}
This is symmetric with respect to the direction of time. Also
under $\mathbb{P}^0_{\infty}$, $\{X_{t}:t\in\mathbb{R}\}$ is
stationary and
reversible. \ Therefore we conclude that%

\[
E(f(X_{0})g(X_{t}))=E(f(X_{t})g(X_{0}))
\]
Therefore
$\frac{1}{Z}e^{\frac{V(\mu)}{2\gamma}}\Pi_{\gamma}^{0}(d\mu)$ is a
reversible invariant measure.
\end{proof}

\begin{corollary}
Consider the $K$-allele case with $c=\gamma$ and $\nu_{0}(dx)=dx.$
Assume that $V(p)$ is continuous and has a unique global maximum
$p_{0}\in\Delta_{K-1}.$ Then as $\gamma\rightarrow0$,
$\Pi_{\gamma}^{V}\Longrightarrow\delta_{p_{0}}.$
\end{corollary}

\begin{proof}
In this case $\Pi_{\gamma}^{0}(dp)$ is the Dirichlet (1)
distribution on $\Delta_{K-1}$. \ Let
\[
N_{p}^{\varepsilon}:=\{p:V(p_{0})-V(p)\leq\varepsilon\}
\]
Then for any
$\varepsilon>0,\;\Pi_{\gamma}^{0}(N_{p_{0}}^{\varepsilon})>0.$ It
is then easy to check that%
\[
\Pi_{\gamma}^{V}((N_{p_{0}}^{\varepsilon})^{c})\rightarrow0\text{ as }%
\gamma\rightarrow0\text{.}%
\]
\end{proof}

\subsection{Remarks on further developments}

The Gibbs form of the invariant measure is suggestive.  One can ask if there is a reversible equilibrium for other mutation processes.  However the fact that the type-independent  mutation is the only mutation process for which the equilibrium is reversible was proved by Li, Shiga and Ya (1999) \cite{LSY-99}. This is the analogue of the result of Hofbauer and Sigmund mentioned above.
The reason for the reversibility of the IMA mutation is that this mutation mechanism ``erases'' all historical information.  For the related infinitely many sites mutation historical information is preserved and the stationary state is not reversible.

A natural setting for the study of the macroscopic development of population systems that incorporates finite local capacity, finite population resampling, spatial migration, mutation and selection is the {\em  stepping stone model with infinitely many sites or hierarchical mutation and  general state-dependent fitness}, $V(x,\mu)$  of types.  These systems do not have reversible stationary measures but other methods including mean-field methods can yield partial information on the large scale behaviour of these systems.  This is currently an active field of research.

\section{Remarks on the Bibliography}
The topics covered in these notes have focussed on basic tools from stochastic analysis and the basic mechanisms of reproduction, mutation and selection as well as spatial  migration.  These methods and processes serve as building blocks for the study of more complex systems.  The  Bibliography below contains numerous references to the vast and growing literature on developments along these lines both in stochastic analysis related to population models and references from the biological literature which deal with questions to which stochastic population models have contributed or have the potential to contribute to.



\chapter{Appendix I: Martingales and large deviations}

\section{Martingales and uniform integrability}

\beD{} A family of integrable random variables $\{X_\alpha\}$ is {\em uniformly integrable} if
\be{} \sup_\alpha E[|X_\alpha 1_{|X_\alpha|>L}]\to 0\text{   as  }L\to\infty.\ee

\end{definition}

A sufficient condition for uniform integrability is that there exists a non-negative increasing function such that $\lim_{s\to\infty}\frac{g(s)}{s}=\infty$ and
\be{}\sup_\alpha E(g(|X_\alpha|))<\infty.\ee

\beL{WCM} (a)If $X_n\Rightarrow X$, then $E(|X|)\leq \liminf_{n\to\infty}E(|X_n|)$.  If $p_2>p_1>0$ and
$\sup_n E(|X_n|^{p_2})<\infty$, then $E(|X|^{p_1})=\lim_{n\to\infty} E(|X_n|^{p_1})$.

(b) If a sequence of martingales $X_n(\cdot)$ converge weakly to $X(\cdot)$ in $D_{\mathbb{R}}([0,\infty))$ and for some $\delta >0$
\be{} \sup_n \sup_{0\leq t\leq T} E(|X_n(t)|^{1+\delta})<\infty\ee
for each $T<\infty$, then $X(\cdot)$ is a martingale.
\end{lemma}

\beT{}  Let $\{X_n\}$ be a martingale and $\tau$ a stopping time with $E(\tau)<\infty$.  Assume that
\be{} E(|X_{n+1}-X_n|\,|X_1,\dots,X_n)\leq \alpha<\infty, \; n\leq \tau.\ee
Then
\be{} E X_\tau=E(X_1).\ee

\end{theorem}
\begin{proof} See Breiman \cite{B-68} Prop. 5.33.

\end{proof}

\subsubsection{Application to the Radon-Nikodym Theorem}

Let $P$ and $Q$ be two probability measures on $(\Omega,\mathcal{F})$ and $\mathcal{F}_n$ an increasing sequence of sub-$\sigma$-algebras which generate $\mathcal{F}$.  Assume that $Q<<P$ on $\mathcal{F}_n$ for all $n$ and
\be{} R_n=\left. \frac{dQ}{dP}\right|_{\mathcal{F}_n}.\ee

\beT{RNT}  $((R_n)_{n\in\N},\mathcal{F}_n) $ is a martingale.  $R_n$ is uniformly integrable if and only if  $Q<<P$ on $\mathcal{F}$
and if $Q<<P$, then $R_\infty=\lim_{n\to\infty} R_n$ is the Radon-Nikodym derivative.  If $Q\bot P$, then $R_n\to 0$ a.s. with respect to $P$ and $R_n\to \infty$ with respect to $Q$.

\end{theorem}
Reference: See, for example,  Durrett \cite{Du-96}.

\subsection{Burkholder-Davis-Gundy inequalities}
For every continuous local martingale $M$ and $p>0$ there exists $0<c<C<\infty$ such that

\be{} cE([M]_\infty^{p/2})\leq E((M^*_\infty)^p)\leq CE([M]_\infty^{p/2})\ee
where $M^*_t:=\sup_{s\leq t} |M_s|$.

\beT{Burkholder} (Burkholder's inequality (Burkholder (1973))\\
For $p>1$ and martingale $M$ with jumps $\Delta M$ (Burkholder's inequality (Burkholder (1973)) Theorem 21.1) is
\be{} E[(M^*_t)^p]\leq C E[ (\langle M\rangle^p)+ E(\sup_{s\leq t} |\Delta M_s|^p)].\ee
\end{theorem}

\section{Large Deviations}

\beL{cramer}
(a) Let $X_{1},\dots,X_{N}$ be bounded iid r.v.'s with mean $\mu$ and
$|X_{i}|\leq K$ a.s. \ and $S_{N}=\sum_{i=1}^{N}X_{i}$. Then there exists
constants $C,\eta$ depending on$K$ and $\varepsilon$ so that%
\[
P(|\frac{S_{N}}{N}-\mu|\geq\varepsilon)<Ce^{-N\eta}.
\]
(b) \ Suppose that $\{X_{n}\}$ are random variables such that $E[X_{n+1}%
|S_{n}]=0$ and $E[X_{n+1}^{2}|S_{n}]\leq v.$ Then (Azuma's inequality): for $x>0$ %
\[
P(\frac{|S_{n}|}{n}\geq x)\leq\exp\left(  -nH\left(  \frac{x+v}{1+v}|\frac{v}%
{1+v}\right)  \right)
\]
where
\be{} H(p|p_0)=p\log\frac{p}{p_0}+(1-p)\log\frac{1-p}{1-p_0}.\ee
\end{lemma}

\begin{proof}
(a)(cf. Dembo and Zeitouni \cite{DZ-93}, Theorem 2.2.3) Without loss of generality we can assume that
$K=1$ and $\mu=0$. Then Cram\'{e}r's Theorem states that%
\[
\lim\sup\frac{1}{N}\log P(|\frac{S_{N}}{N}|\geq \varepsilon)\leq\inf
_{|x|>\varepsilon}\Lambda^{\ast}(x)
\]
where $\Lambda^*(x)$ is the Legendre transform
\be{} \Lambda^*(x):= \sup_{x}\{\lambda x-\Lambda(\lambda)\}\ee
\be{} \Lambda(\lambda):=\log E[e^{\lambda X_1}].
\ee
Then a simple calculation shows that%
\begin{align*}
\Lambda^{\ast}(x)  &  \geq H(\frac{x+1}{2}|\frac{1}{2})\\
H(p|\frac{1}{2})  &  =-(p\log2p+(1-p)\log2(1-p))\\
&  >0.
\end{align*}
(b) See Dembo and Zeitouni \cite{DZ-93}  exercise 2.2.29 for hints on the proof.
\end{proof}

\chapter{Appendix II: Measures and Topologies}

\section[Measures and weak convergence]{Measures on Polish spaces and weak convergence}

Throughout these lectures measures on a number of Polish topological spaces will play an essential role.  In addition, the construction of various limiting limiting objects involves weak  convergence of probability measures on these spaces. In this section we collect some basic facts on certain basic classes of Polish spaces and  the topology of
weak convergence of probability measures on Polish spaces.

\subsection{Borel measures on Polish space}

\beD{} A Polish space is a separable completely metrizable topological space;
that is, a space homeomorphic to a complete metric space $(E,d)$
that has a countable dense subset.
\end{definition}

A basic property of a Polish space is that a subset $A$ is
relatively compact iff it is totally bounded with respect to $d$.

\begin{theorem}
Every finite measure $\mu$ on the Borel $\sigma$-algebra $\mathcal{E}=\mathcal{B}(E)$  of a complete
separable metric space $(E,d)$ is Radon, that is, regular relative to the family of compact subsets.
\end{theorem}

\begin{proof} We first prove that $\mu(E)=\sup_{K\subset E,; K\text{ compact}} \mu(K)$.
Let $\{x_{k}\}$ be a dense subset and $B(x,\ve)$ the ball with center $x$ and radius $\ve$. \ For $\varepsilon>0,$ choose integers
$N_{1},N_{k},\dots$ such that%
\[
P(\cup_{k=1}^{N_{n}}B(x_{k},\frac{1}{n}))\geq1-\frac{\varepsilon}{2^{n}}
\]
Let $K$ be the closure of $\cap_{n\geq1}\cup_{k=1}^{N_{n}}B(x_{k},\frac{1}%
{n})$. Then $K$ \ is totally bounded and hence compact and%
\[
P(K^{c})<\varepsilon\text{.}
\]

For the completion, see D.L. Cohn, Measure Theory, Birkh\" auser.
\end{proof}

Sometimes we will consider a special class of metric spaces, the {\em ultrametric   spaces} on which the
strong triangle inequality is satisfied
\be{} d(x,y)\leq \max(d(x,z),d(y,z)).\ee

\begin{example}  Consider a rooted real tree $\mathcal{T}$ with root $\emptyset$.  Then for any $t$ the level set
\be{}X^t:=\{x\in\mathcal{T}:d(\emptyset,x)=t\}\ee
is an ultrametric space.
\end{example}

\subsection{Weak convergence of measures}

A basic reference for weak convergence of probability measures on a complete separable metric space $(E,d)$ is Ethier-Kurtz \cite{EK-86}, Chap. 3.  Let $\mathcal{P}(E)$ denote the set of Borel probability measures on $E$.  The following are basic properties

1.  $(\mathcal{P}(E),\rho)$ is a complete separable metric space where $\rho$ is the Prohorov metric.

2. A subset  $\mathcal{M}\subset \mathcal{P}(E)$ is relatively compact if and only if it is (uniformly) {\em tight}, that is, for each $\ve>0$ there exists a compact set $K\subset E$ such that

\be{} \inf_{P\in\mathcal{M}}P(K)\geq 1-\ve.\ee

\section{Projective Limits of Measures and Processes}

(Reference: P.A. Meyer (1966) \cite{M-66})

Let $(E,\mathcal{E})$ be a measurable space, $T$ an index set and $U$ the
collection of finite subsets of $T$. We denote by $\pi_{u}$ the canonical
projection of $E^{T}$ on $E^{u}$, $u\in U$. The mappings $\pi_{u},\pi_{uv}$
are measurable over the natural product $\sigma$-fields and
\begin{align*}
\pi_{uv}\pi_{v}  &  =\pi_{u}\text{ , }u\subset v\\
\pi_{uv}\pi_{vw}  &  =\pi_{uw},\,\ u\subset v\subset w
\end{align*}
Suppose we are given a probability law $P_{u}$ on every measurable space
$(E^{u},\mathcal{E}^{u})$, $u\in U$. \ The family $(P_{u})_{u\in U}$
constitutes a projective systems of probability laws if for every pair of
elements, $u\subset v,$%
\begin{align*}
\pi_{uv}(P_{v})  &  =P_{u}\\
P_{u}(A)  &  =P_{v}(\pi_{uv}^{-1}(A))
\end{align*}
The projective system admits a \textit{projective limit} if there exists a
probability law $P$ on $(E^{T},\mathcal{E}^{T})$ such that%
\[
\pi_{u}(P)=P_{u}\text{ for every }u\in U
\]
The collection of subsets of $E^{T}$ of the form
$\pi_{u}^{-1}(A_{u})(u\in U,A_{u}\in\mathcal{E}^{u})$, denoted by
$\mathcal{E}_{0}^{T}$, forms a \textit{paving} (i.e. a collection
of subsets of $E^{T}$ containing $\emptyset,E^{T}
$) is closed under finite unions and taking complements.  If we put for $A=\pi_{u}^{-1}(A_{u}%
)\in\mathcal{E}_{0}^{T}$,
\[
P(A)=P_{u}(A_{u})
\]
we obtain \ a function of $A$ independent of the representation of $A$. Given
that the paving $\mathcal{E}_{0}^{T}$ generates the $\sigma$-field
$\mathcal{E}^{T}$ the uniqueness of the projective limit follows.

\begin{definition}
A measure $P$ on $(\Omega,\mathcal{F})$ is said to be \textit{regular} with
respect to a semicompact paving $\mathcal{T}$ (every countable family having
the finite interesection property has the countable intersection property)
whose elements are measurable if%
\[
P(B)=\sup_{\substack{A\in\mathcal{T} \\A\subset B}}P(A)\text{ for every }%
B\in\mathcal{F}
\]
\end{definition}

\begin{theorem}
It suffices to show this for $B$ in a generating algebra $\mathcal{F}%
_{0}\subset\mathcal{F}$ and the extension to a regular probability law is unique.
\end{theorem}

\begin{proof}
Use the monotone class theorem. See Meyer (1966).
\end{proof}

\begin{theorem}\label{PLT}
(Neveu) Suppose that there exists, for each $t\in T$, a semicompact paving
$\mathcal{K}_{t}\subset\mathcal{E}$ such that the law $P_{t}$ is regular wrt
$\mathcal{K}_{t}$. Then the projective system $(P_{u})_{u\in U}$ admits a
projective limit.
\end{theorem}

\begin{proof}
We can suppose that $E$ belongs to each of the pavings $\mathcal{K}_{t}$.
\ Let $\mathcal{K}_{u}$, $\mathcal{K}_{T}$ denote the closure under $(\cup
f,\cap c)$ of the product paving $\prod_{t\in u}\mathcal{K}_{t},\prod_{t\in
T}\mathcal{K}_{t}$ . Each of these pavings is semicompact. \ Next denote by
$\mathcal{K}_{T}^{0}$ the paving on $E^{T}$ \ consisting of subsets of the
form $\pi_{u}^{-1}(A_{u})(u\in U,A_{u}\in\mathcal{K}_{u}),$ $\mathcal{K}%
_{T}^{0}\subset\mathcal{K}_{T}$ and hence is semicompact. \ Let us first show
that the law $P_{u}$ is regular relative to $\mathcal{K}_{u}$. To do this we
apply the previous theorem using $\mathcal{F}_{0}$ the collection of finite
unions of sets of the form $\prod_{t\in u}A_{t}(A_{t}\in\mathcal{E)}$ and for
$\mathcal{K}_{0}$ the collection of finite unions of sets of the form
$\prod_{t\in u}K_{t}(K_{t}\in\mathcal{K}_{t})$. We must verify$P_{u}%
(\prod_{t\in u}A_{y})=\sup_{\substack{K_{t}\in\mathcal{K}_{t} \\K_{t}\subset
A_{t}}}P(\prod_{t\in u}K_{t}).$ For $\varepsilon>0$ take $\varepsilon
_{t}>0(t\in u)$ such that $\sum\varepsilon_{t}=\varepsilon$. Choose for each
$t\in u$ a set $K_{t}\in\mathcal{K}_{t}$ such that $K_{t}\subset A_{t}$ and
$P(A_{t}\backslash K_{t})\leq\varepsilon_{t}.$ The denoting by $B_{t}$ the
inverse image in $E^{u}$ of $A_{t}\backslash K_{t}$ under the projection of
$E^{u}$ onto $E^{\{t\}}$%
\[
P_{u}\left[  (\prod_{t\in u}A_{t})\backslash(\prod_{t\in u}K_{t})\right]  \leq
P_{u}(\cup_{t\in u}B_{t})\leq\sum_{t\in u}P_{t}(A_{t}\backslash K_{t}%
)\leq\varepsilon
\]
Hence $P_{u}$ is regular wrt $\mathcal{K}_{u}$ for each $u\in U$.
Therefore $P$ is defined on   $\mathcal{E}_{0}^{T}$ and can be
extended uniquely to a regular probability law on $\mathcal{E}$.
\end{proof}

We state the following theorem without proof. This implies, for example, that
every probability law on the Baire $\sigma$-field of a compact space is regular.

\begin{theorem}
Let $\mathcal{K}$ be a semicompact paving on the set $\Omega,$ closed under
$(\cup f,\cap c)$ such that the complement of every element of $\mathcal{K}$
belongs to $\mathcal{K}_{\sigma}.$ Every probability law $P$ on the $\sigma
$-fields generated by $\mathcal{K}$ is then regular relative to $\mathcal{K}$.
\end{theorem}

\begin{proof}
See Meyer (1966).
\end{proof}

Remark: Note that we can generalise the above result of Neveu to let $U$ be a
\textit{right-filtering partially ordered set}, i.e. given $u,v\in
U,\;\exists$ $\ \;w\succ u,\;w\succ v$, and a family of probabilities $P_{u} $
on $(E_{u},\mathcal{E}_{u}),\;$\ and measurable maps $\pi_{uv}:(E_{v}%
,\mathcal{E}_{v})\rightarrow(E_{u},\mathcal{E}_{u}),\;u\subset v $.%

\begin{align*}
\pi_{uv}\pi_{v}  &  =\pi_{u}\text{ , }u\subset v\\
\pi_{uv}\pi_{vw}  &  =\pi_{uw},\,\ u\subset v\subset w
\end{align*}

\section{Random measures}

Let $E$ be Polish space with  a countable base $\mathcal{D}$ for the topology, and $\mathcal{A}$ the algebra generated by $\mathcal{D}$.  For $A_1,\dots,A_n\in\mathcal{A}$ let $P_{\{A_1,\dots,A_n\}}$ be a distribution $( R_+)^n$. Assume that
\begin{enumerate}
\item consistency $P_{\{A_1,\dots,A_n\}}= P_{\{A_1,\dots,A_n,E\}}$
\item if $A,B\in\mathcal{A}$ and $A\cap B=\varnothing$, then
\be{} P_{A,B,A\cup B}((x,y,z)\in\mathbb{R}^3: x+y=z)=1,\ee
\item if $\{A_n\}\subset \mathcal{A}$ and $A_n\downarrow \varnothing$, then for all $\ve>0$
\be{} \lim_{n\to\infty} P_{A_n}([\ve,\infty))=0.\ee
\item for $\ve>0$ there exists a compact subset $K_\ve$ such that $P(K_\ve^c([\ve,\infty))<\ve$
\end{enumerate}

\beT{PRPOL}
Assume (1-4) above. Then there is a unique probability, $P$,  on\\ $(M_F(E),\mathcal{B}(M_F(E)))$
 such that
 \be{} P(\{\mu\in M_F(E): (\mu(A_1),\dots,\mu(A_n))\in C\})= P_{A_1,\dots,A_n}(C),\quad C\in\mathcal{B}(\mathbb{R}_+^n).\ee

\end{theorem}

Proof.  See for example, Jagers (1974) \cite{J-74}, Harris (1968) \cite{H-68}.

\section{Topologies on path spaces}

\beD{} Let $\mu_i,\;\mu \in \mathcal{M}_f$.  Then
$(\mu_n)_{n\in\mathbb{N}}$ converges weakly to $\mu$ as
$n\to\infty$, denoted $\mu_n\Rightarrow \mu$ iff and only is \be{}
\int fd\mu_n \nto \int fd\mu\quad\forall \;\; f\in C_b(E)\ee
\end{definition}

Given a Polish space $(E,d)$ we consider the space
$C_E([0,\infty))$ with the metric \be{} \wt
d(f,g)=\sum_{n=1}^\infty 2^{-n} \sup_{0\leq t\leq
n}|f(t)-g(t)|.\ee Then $(C_E([0,\infty),\wt d)$ is also a Polish
space. To prove weak convergence in  $\mathcal{P}((C_E([0,\infty),\wt d))$
it suffices to prove tightness and the convergence of the finite dimensional distributions.

Similarly the space $D_E([0,\infty)$ of c\`adl\`ag functions from
$[0,\infty)$ to $E$ with the Skorohod metric $\wt d$  is a Polish
space where

\be{}  \wt d(f,g)=\inf_{\lambda\in
\Lambda}\left(\gamma(\lambda)+\int_0^\infty e^{-u}\left(1\wedge
\sup_t d(f(t\wedge u),g(t\wedge u))\right)\right)\ee where
$\Lambda$ is the set of continuous, strictly increasing functions
on $[0,\infty)$ and for $\lambda\in\Lambda$, \be{}
\gamma(\lambda)=1+\left(\sup_t |t-\lambda(t)|\vee \sup_{s\ne t}
|\frac{\log (\lambda(s)-\lambda(t))}{s-t}|\right).\ee

\beT{} (Ethier-Kurtz) (Ch. 3, Theorem 10.2) Let $X_n$ and $X$ be processes with sample paths in $D_E([0,\infty)$ and $X_n\Rightarrow X$. Then $X$ is a.s. continuous if and only if $J(X_n)\Rightarrow 0$ where
\be{} J(x)=\int_0^\infty e^{-u} [\sup_{0\leq t\leq u}d(X(t),x(t-))].\ee
\end{theorem}

\subsection{Sufficient conditions for tightness}

\beT{Aldous} (Aldous (1978)) Let $\{P_n\}$ be a sequence of
probability measures on $D([0,\infty),\mathbb{R})$ such that
\begin{itemize}
\item  for each fixed $t$, $P_n\circ X_t^{-1}$ is tight in
$\mathbb{R}$, \item given stopping times $\tau_n$ bounded by
$T<\infty$ and $\delta_n\downarrow 0$ as $n\to\infty$ \be{}
\lim_{n\to\infty} P_n(|X_{\tau_n +\delta_n}-X_{\tau_n}|>\ve)=0,\ee
or \item $\forall\; \eta>0\;\exists \delta,n_0$ such that \be{}
\sup_{n\geq n_0} \sup_{\theta\in [0,\delta]
}P_n(|X_{\tau_n+\theta} -X_{\tau_n}|>\ve)\leq \eta .\ee
\end{itemize}
Then $\{P_n\}$ are tight.
\end{theorem}

\subsection{The Joffe-M\'etivier criteria for tightness of D-semimartingales}
\label{JMC}

We recall the Joffe M\'etivier criterion (\cite{JM-86}) for tightness of locally
square integrable processes.

A c\`adl\`ag adapted process $X$, defined on $(\Omega,\mathcal{F},\mathcal{F}_t,P)$ with values in $\mathbb{R}$ is called a {\em D-semimartingale} if there exists a c\`adl\`ag function $A(t)$, a linear subspace $D(L)\subset C(\mathbb{R})$ and a mapping $L:(D(L)\times \mathbb{R}\times [0,\infty)\times \Omega)\to \mathbb{R}$ with the following properties:
\begin{enumerate}
\item  for every $(x,t,\omega)\in \mathbb{R}\times [0,\infty)\times \Omega$ the mapping $\phi\to L(\phi,x,t,\omega)$ is a linear functional on $D(L)$ and $L(\phi,\cdot,t,\omega)\in D(L)$,
\item for every $\phi\in D(L)$, $(x,t,\omega)\to L(\phi,x,t,\omega)$ is $\mathcal{B}(\mathbb{R})\times \mathcal{P}$-measurable, where $\mathcal{P}$ is the predictable $\sigma$-algebra on $[0,\infty)\times \Omega$, ($\mathcal{P}$ is generated by sets of the form $(s,t]\times F$ where $F\in\mathcal{F}_s$ and $s,t$ are arbitrary)
\item for every $\phi\in D(L)$ the process $M^\phi$ defined by
\be{} M^{\phi}(t,\omega):= \phi(X_t(\omega)-\phi(X_0(\omega))-\int_0^t L(\phi,X_{s-}(\omega),s,\omega)dA_s,\ee
is a locally square integrable martingale on $(\Omega,\mathcal{F},\mathcal{F}_t,P)$,
\item the functions $\psi(x):= x$ and $\psi^2$ belong to $D(L)$.
\end{enumerate}
The functions
\be{} \beta(x,t,\omega):=L(\psi,x,t,\omega)\ee
\be{} \alpha(x,t,\omega):= L((\psi)^2,x,t,\omega)-2x\beta(x,t,\omega)\ee
are called the {\em local characteristics of the first and second order}.

\beT{Joffe-Metivier} Let $X^m=(\Omega^m,\mathcal{F}^m,\mathcal{F}^M_t,P^m)$ be a sequence of D-semimartingales with common $D(L)$ and associated operators $L^m$, functions $A^m,\alpha^m,\beta^m$.  Then the sequence $\{X^m:m\in\mathbb{N}\}$ is tight in $D_{\mathbb{R}}([0,\infty)$ provided the following conditions are satisfied:
\begin{enumerate}
\item $\sup_m E|X_0^m|^2<\infty$,
\item there is a $K>0$ and a sequence of positive adapted processes $\left\{\{C^m_t:t\geq 0
\}\text{ on }\Omega^m\right\}_{m\in\N}$ such that for every $m\in\N, x\in\mathbb{R},\omega\in\Omega^m$,
\be{}|\beta_m(x,t,\omega)|^2+\alpha_m(x,t,\omega)\leq K(C^m_t(\omega)+x^2)\ee
and for every $T>0$,\be{} \sup_m \sup_{t\in[0,T]} E|C^m_t|<\infty,\text{  and  } \lim_{k\to\infty}\sup_m P^m(\sup_{t\in [0,T]} C^m_t\geq k)=0,\ee

\item there exists a positive function $\gamma$ on $[0,\infty)$ and a decreasing sequence of numbers $(\delta_m)$ such that $\lim_{t\to 0}\gamma(t)=0$, $\lim_{m\to\infty} \delta_m=0$ and for all $0<s<t$ and all $m$,
    \be{} (A^m(t)-A^m(s))\leq \gamma(t-s)+\delta_m.\ee
\item  if we set
    \be{} M^m_t:= X^m_t-X^m_0-\int_0^t\beta_m(X^m_{s-},s,\cdot)dA^m_s,\ee
    then for each $T>0$ there is a constant $K_T$ and $m_0$ such that for all $m\geq m_0$,
  then
    \be{} E(\sup_{t\in [0,T]}|X^m_t|^2)\leq K_T(1+E|X^m_0|^2),\ee
    and
    \be{}  E(\sup_{t\in [0,T]}|M^m_t|^2)\leq K_T(1+E|X^m_0|^2),\ee
\end{enumerate}
\end{theorem}

\beC{}

\end{corollary}  Assume that for $T>0$ there is a constant $K_T$ such that
\be{} \sup_m\sup_{t\leq T,x\in\mathbb{R}} (|\alpha_m(t,x)|+|\beta_m(t,x)|)\leq K_T,\;\text{a.s.}\ee
\be{} \sum_m (A^m(t)-A^m(s))\leq K_T(t-s)\text{  if  } 0\leq s\leq t\leq T,\ee
and
\be{}  \sup_m E|X^m_0|^2<\infty,\ee
and $M^m_t$ is a square integrable martingale with $sup_m E(|M^m_T|^2)\leq K_T$. The the $\{X^m:m\in\N\}$ are tight in $D_{\mathbb{R}}([0,\infty)$.

\subsubsection{Criteria for continuous processes}

Now consider  the special case of probability measures on  $C([0,\infty), \R^d)$. This criterion is concerned
with a collection $(X^{(n)}(t))_{t \geq 0}$ of semimartingales
with values in $\R^d$ with continuous  paths. First observe that
by forming \be{jm1a} (<X^{(n)}(t), \lambda>)_{t \geq 0}\quad ,
\quad \lambda \in \R^d \ee we obtain $\R$-valued semi-martingales.
If for every $\lambda \in \R^d$ the laws of these projections are
tight on $C([0,\infty),\R)$ then this is true for $\{[\CL
[(X^{(n)}(t))_{t \geq 0}], n \in \N\}$. The tightness criterion
for $\R$-valued semimartingales is in terms of the so-called local
characteristics of the semimartingales.

For It\^o processes the local characteristics can be
calculated directly from the coefficients.  For example,  if we
have a sequence of semimartingales $X^n$ that are also a Markov
processes with generators:   \be{jm5} L^{(n)} f =
\big(\suml^d_{i=1} a^n_i (x) \frac{\partial}{\partial x_i} +
\suml^d_{i=1} \suml^d_{j=1} b^n_{i,j} (x)
\frac{\partial^2}{\partial x_i
\partial x_j}\big)f \ee  then the local characteristics
are given by  \be{jm6} a^n=(a^n_i)_{i=1,\cdots,d},\;\;
b^n=(b^n_{i,j})_{i,j,=1,\cdots,d}. \ee

The Joffe-M\'etivier criterion implies that if \bea{jm3} &&\sup_n
\supl_{0 \leq t \leq T} E[(|a^n(X^{(n)}(t) | + | b^n (X^{(n)}(t)
|)^2] < \infty,\\&& \lim_{k\to\infty} \sup_n P(\sup_{0\leq t\leq
T}(|a^n(X^{(n)})(t) | + | b^n (X^{(n)})(t) |)\geq k) =0 \eea then
$\{\CL [(X^{(n)}(t))_{t \geq 0}], n \in \N\}$ are tight in
$C([0,\infty),\R)$. See \cite{JM-86} for details.

\beT{DTC} (Ethier-Kurtz \cite{EK-86} Chapt. 3, Theorem 10.2) Let
\be{}  J(x)=\int_0^\infty e^{-u} [J(x,u)\wedge 1]du,\quad J(x,u)=\sup_{0\leq t\leq u}d(x(t),x(t-)).\ee
Assume that a sequence of processes  $X_n\Rightarrow  X$ in $D_E([0,\infty))$.  Then  $X$ is a.s. continuous if and only if $J(X_n)\Rightarrow 0$.

\end{theorem}

\subsection{Tightness of measure-valued processes}

\begin{lemma}\label{L.tightness}
(Tightness \ Lemma).\newline (a) Let $E$ be a compact metric space and
$\{P_{n}\}$ a sequence of probability measures on $D([0,\infty),M_{1}(E))$.
Then $\{P_{n}\}$ is compact if and only if there exists a linear separating
set $D\subset C(E)$ such that $t\rightarrow\int f(x)X_{t}(\omega,dx)$ is
relatively compact in $D([0,\infty),\mathbb{[-\Vert}f\mathbb{\Vert},\Vert
f\Vert]\mathbb{)}$ for each $f\in D$.\newline (b) Assume that$\;$$\{P_{n}\}$
is a family of probability measures on $D([0,\infty),\mathbb{[-}%
K,K]\mathbb{)}$ such that for $0\leq t\leq T$, there are bounded predictable
processes $\{v_{i}(\cdot):i=1,2\}$ such that for each $n$
\[
M_{i,n}(t):=x(\omega,t)^{i}-\int_{0}^{t}v_{i,n}(\omega,s)ds,\;i=1,2
\]
are $P_{n}$-square integrable martingales with
\[
\sup_{n}E_{n}(\sup_{s}(|v_{2,n}(s)|+|v_{1,n}(s)|))<\infty.
\]
Then the family $\{P_{n}\}$ is tight.\newline (c) In (b) we can replace the
$i=2$ condition with: for any $\varepsilon>0$ there exists $f$ and $v_{f,n}$
such that%
\[
\sup_{\lbrack-K,K]}|f_{\varepsilon}(x)-x^{2}|<\varepsilon
\]
and%
\begin{align*}
M_{f,n}(t)  &  :=f_{\varepsilon}(x(\omega,t))-\int_{0}^{t}v_{f_{\varepsilon
},n}(\omega,s)ds\\
\sup_{n}E_{n}(\sup_{s}(|v_{f_{\varepsilon},n}(s)|)  &  <\infty.
\end{align*}
\end{lemma}

\begin{proof}
(a) See e.g. Dawson, \cite{D-93} Section 3.6. \newline (b) Given stopping times $\tau_{n}$
and $\delta_{n}\downarrow0$ as $n\rightarrow\infty.$%
\begin{align*}
&  E_{n}\left[  (x(\tau_{n}+\delta_{n})-x(\tau_{n}))^{2}\right] \\
&  =\{E_{n}[x^{2}(\tau_{n}+\delta_{n})-x^{2}(\tau_{n})]-2E_{n}[x(\tau
_{n})(x(\tau_{n}+\delta_{n})-x_{n}(\tau_{n}))]\}\\
&  \leq E_{n}[\int_{\tau_{n}}^{\tau_{n}+\delta_{n}}|v_{2,n}(s)|ds+2K\int
_{\tau_{n}}^{\tau_{n}+\delta_{n}}|v_{1,n}(s)|ds]\\
&  \leq\delta_{n}\sup_{n}E_{n}(\sup_{s}(|v_{2,n}(s)|+|v_{1,n}(s)|))\\
&  \rightarrow0\text{ \ as }\delta_{n}\rightarrow0.
\end{align*}

The result then follows by Aldous' condition.

(c)%
\begin{align*}
&  E_{n}\left[  (x(\tau_{n}+\delta_{n})-x(\tau_{n}))^{2}\right] \\
&  =\{E_{n}[x^{2}(\tau_{n}+\delta_{n})-x^{2}(\tau_{n})]-2E_{n}[x(\tau
_{n})(x(\tau_{n}+\delta_{n})-x_{n}(\tau_{n}))]\}\\
&  \leq E_{n}(f(x(\tau_{n}+\delta_{n}))-f(x(\tau_{n})))+2K\int_{\tau_{n}%
}^{\tau_{n}+\delta_{n}}|v_{1,n}(s)|ds]+2\varepsilon\\
&  \leq E_{n}[\int_{\tau_{n}}^{\tau_{n}+\delta_{n}}|v_{f_{\varepsilon}%
,n}(s)|ds+2K\int_{\tau_{n}}^{\tau_{n}+\delta_{n}}|v_{1,n}(s)|ds]+2\varepsilon
\\
&  \leq\delta_{n}\sup_{n}E_{n}(\sup_{s}(|v_{f_{\varepsilon},n}(s)|+|v_{1,n}%
(s)|))+2\varepsilon
\end{align*}
Hence for any $\varepsilon>0$%
\begin{align*}
&  \lim_{\delta_{n}\rightarrow0}\sup_{n}E_{n}\left[  (x(\tau_{n}+\delta
_{n})-x(\tau_{n}))^{2}\right] \\
&  \leq\lim_{n\rightarrow\infty}\delta_{n}\sup_{n}E_{n}(\sup_{s}%
(|v_{f_{\varepsilon},n}(s)|+|v_{1,n}(s)|))+2\varepsilon\\
&  =2\varepsilon.
\end{align*}
and the result again follows from Aldous criterion.
\end{proof}

\begin{remark}
These results can be also used to prove tightness in the case of non-compact
$E$. \ However in this case an additional step is required, namely to show
that for $\varepsilon>0$ and $T>0$ there exists a compact subset
$K_{T,\varepsilon}\subset E$ such that%
\[
P_{n}[D([0,T],K_{T,\varepsilon})]>1-\varepsilon\;\;\forall\text{ }n.
\]
\end{remark}

\begin{remark}
Note that if $P_{n}$ is a tight sequence of probability measures on
$D([0,T],\mathbb{R)}$ such that $P_{n}(\sup_{0\leq s\leq T}|x(s)-x(s-)|\leq
\delta_{n})=1$ and $\delta_{n}\rightarrow0$ as $n$ $\rightarrow\infty$, then
for any limit point $P_{\infty}$, $P_{\infty}(C([0,T],\mathbb{R}))=1$.
\end{remark}

\section{The Gromov-Hausdorff metric on the space of compact metric spaces}\label{s.AGH}

Let  $E$ be a metric space and $B_1,\;B_2$ two  subsets.  Then the Hausdorff distance is defined by

\be{}  d_{\rm{H}}(K_1,K_2)= \inf\{ \ve\geq 0: K_1\subset V_\ve(K_2),\;K_2\subset V_\ve(K_1)\}\ee where
$V_\ve(K)$ denotes the $\ve$-neighbourhood of $K$.
This defines a pseudometric, $d_H(B_1,B_2)=0$ iff they have the same closures.

 If $X$ and $Y$ are two compact metric spaces. The {\em Gromov-Hausdorff metric}  $d_{GH}(X,Y)$
is defined to be the infimum of all numbers $d_{H}(f(X),g(Y))$ for
all metric spaces $M$ and all isometric embeddings $f:X\rightarrow
M$ and $g:Y\rightarrow M$ and where $d_{\rm{Haus}}$ denotes Hausdorff
distance between subsets in M. $d_{GH}$ is a pseudometric with $d_{GH}(K_1,K_2)=0$ iff they are isometric.

Now let $(\mathbb{K},d_{GH})$ denote the class of compact metric spaces (modulo isometry) with the Gromov-Hausdorff metric.  Then $(\mathbb{K},d_{GH})$ is complete.

See Gromov \cite{Gr-99} and Evans \cite{E-07}  for detailed expositions on this topic.

\subsection{Metric measure spaces}
The notion of {\em metric measure space} was developed by Gromov \cite{Gr-99} (called mm spaces there).  It is given by a triple $(X,r,\mu)$ where $(X,r)$ is a metric space such that $(\rm{supp}(\mu),r)$ is complete and separable and $\mu\in \mathcal{P}(X)$ is a probability measure on $(X,r)$. Let $\mathbb{M}$ be the space of equivalence classes of metric measure spaces (whose elements are not themselves metric spaces - see remark (2.2(ii)) in\cite{GPW-08}) with equivalence in the sense of measure-preserving isometries.
The distance matrix map is defined for $n\leq \infty$

\be{}  X^n\to \mathbb{R}_+^{\small{(\begin{array}{c} n \\2\end{array})}},\qquad  ((x_i)_{i=1,\dots,n})\to (r(x_i,x_j))_{1\leq i<j\leq n}\ee
and we denote by $R^(X,r)$ the map that sends a sequence of points to its infinite distance matrix.

Then the {\em distance matrix distribution of } $(X,r,\mu)$ (representative of equivalence class) is defined by

\be{}  \nu^{(X,r,\mu)}:= R^{(X,r)}-\text{pushforward of } \mu^{\otimes \mathbb{N}}\quad \in\mathcal{P}(\mathbb{R}_+^{{{\left(\begin{array}{c} \N\\2\end{array}\right)}}}).\ee

Since this depends only on the equivalence class it defined the mapping $\kappa\to \nu^{\kappa} $ for $\kappa\in\mathbb{M}$.
Gromov \cite{Gr-99} (Section 3$\frac{1}{2}$.5) proved that a metric measure space is characterized by its distance matrix distribution.

Greven, Pfaffelhuber and Winter (2008) \cite{GPW-08}  introduced the {\em Gromov-weak topology}.  In this topology  a  sequence $\{\chi_n\}$ converges Gromov-weakly to $\chi$ in $\mathbb{M}$ if and only if $\Phi(\chi_n)$ converges to $\Phi(\chi)$ in $\mathbb{R}$ for all polynomial in $\Pi$.

In \cite{GPW-08}, Theorem 1, they proved that $\mathbb{M}$ equipped with the Gromov-weak topology is Polish.

An important subclass is the set of ultrametric measure spaces given by the closed subset of $\mathbb{M}$

\be{}  \mathbb{U}:=\{\ {\textit{u}} {\in}\mathbb{M}: \textit{u} \text{   is ultra-metric}\}.\ee

\section{Riemannian metrics and gradient}

Let $M$ be a smooth manifold. The tangent space at $x$, $T_{x}M$
can be identified with the space of tangents at $x$ to all smooth
curves through $x$. The tangent bundle $TM=\{(p,v):p\in M,v\in
T_{p}M\}.$

\begin{definition}
A Riemannian metric on $M$ is a smooth tensor field%
\[
g:C^{\infty}(TM)\otimes C^{\infty}(TM)\rightarrow
C_{0}^{\infty}(M)
\]
such that for each \ $p\in M,$%
\[
g(p)|_{T_{p}M\otimes T_{p}M}:T_{p}M\otimes T_{p}M\rightarrow\mathbb{R}%
\]
with%
\[
g(p):(X,Y)\rightarrow\left\langle X,Y\right\rangle _{g(p)}%
\]
where $\left\langle X,Y\right\rangle _{g(p)}\;$is an inner product
on $T_{p}M$.
\end{definition}

\begin{definition}
The directional derivative in direction $v$ is defined by%
\begin{align*}
\partial_{v}f(x)  &  =\lim_{t\rightarrow0}\frac{f(x+tv)-f(x)}{t}\\
&  =\sum v_{i}\frac{\partial f(x)}{\partial x_{i}}%
\end{align*}
The gradient $\nabla_{g}f(x)$ is defined by%
\[
\left\langle \nabla_{g}f(x),v\right\rangle
_{g}=\partial_{v}f(x)\;\;\forall v\in T_{x}M.
\]
\end{definition}

\begin{example}
Consider the $d$-dimensional manifold $M=\mathbb{R}^{d}$ and $\mathbf{a}%
(\cdot)$ be a smooth map from $M$ to
$\mathbb{R}^{d}\otimes\mathbb{R}^{d}$
($(d\times d)$-matrices). We will write%
\begin{align*}
\mathbf{a}(x)  &  =(a_{ij}(x))\\
\mathbf{a}^{-1}(x)  &  =(a^{ij}(x))
\end{align*}
Assume that%
\[
\sum a^{ij}(x)u_{i}u_{j}\geq\gamma\sum u_{j}^{2},\;\gamma>0.
\]
The tangent space $T_{\mathbf{x}}M\approxeq\mathbb{R}^{d}$ and we
define a Riemannian metric on $M$ by
\[
g_{\mathbf{a}(x)}(\mathbf{u,v)}:=\sum_{i,j=1}^{d}a_{ij}(x)u^{i}v^{j}.
\]
The associated Riemannian gradient and norm are%
\begin{align*}
(\nabla_{\mathbf{a}}f)^{i}  &  =\sum_{j}a^{ij}\frac{\partial
f}{\partial
x_{j}}\\
\Vert u\Vert_{\mathbf{a}(x)}^{2}  &
=\sum_{ij}a_{ij}(x)u^{i}u^{j}.
\end{align*}
\end{example}

\subsubsection{The Shahshahani metric and gradient on $\Delta_{K-1}$}

Let $M_{K}=\mathbb{R}_{+}^{K}:=\{x\in\mathbb{R}^{K},x=(x_{1},\dots
,x_{K}),\;x_{i}>0$ for all $i\}$ is a smooth K-dimensional
manifold.

Shahshahani introduced the following Riemannian metric on $M_{K}$%
\begin{align*}
\left\langle u,v\right\rangle _{g}  &  =g_{x}(u,v):=\sum_{i=1}^{K}%
|x|\frac{u_{i}v_{i}}{x_{i}}\\
|x|  &  =\sum x_{i}%
\end{align*}
\thinspace$\Vert\;\Vert_{g}$ and $\nabla_{g}F$ \ will denote the
corresponding norm and gradient. We have
\[
(\nabla_{g}F)^{i}=\sum_{i}\frac{x^{i}}{|x|}\frac{\partial F}{\partial x^{i}%
}\frac{\partial}{\partial x_{{}}^{i}}%
\]

Recall that the simplex  $\Delta_{K-1}:=\{(p_1,\dots,p_K):p_i\geq 0,\; \sum_{i=1}^K p_i=1\}$.
The interior of the simplex
$\Delta_{K-1}^{0}=\mathbb{R}_{+}^{K}\cap \Delta_{K-1}$ is a
$(K-1)$-dimensional submanifold of $M_{K}$. \ We denote by
$T_{p}\Delta_{K-1}^{0}$ the tangent space to \ $\Delta_{K-1}^{0}$
at $p$. Then $g$ induces a Riemannian metric on
$T_{p}\Delta_{K-1}^{0}$.

\underline{Basic Facts}

We have the Shahshahani inner product on $\Delta{K-1}$ at a point $p \in \Delta{K-1}$:
\be{}  \langle u,v\rangle_p= \sum_{i=1}^K \frac{u_iv_i}{p_i}.\ee

1.  $T_{p}\Delta_{K-1}^{0}$ can be viewed as the subspace of
$T_{p}M_{K}$ of vectors, $v$,  satisfying \ $\left\langle
p,v\right\rangle _{g}=0$  if we identify $p$ with an element
of $T_{p}M_{K}.$

Proof. Recall that $T_{p}\Delta_{K-1}^{0}$ is given by tangents to
all smooth
curves lying in $\Delta_{K-1}^{0}$. \ Therefore if $v\in T_{p}\Delta_{K-1}%
^{0},\;$then $v=q-p$ where $p,q\in\Delta_{K-1}^{0}$ and therefore
$\sum _{i=1}^{K}v_{i}=0.$ Therefore, \
\[
\sum_{i}p_{i}\frac{1}{p_{i}}v_{i}=0.%
\]

2. \ If $F:\Delta_{K-1}^{0}\rightarrow\mathbb{R}$ is smooth, then the Shahshahani gradient is
\[
(\nabla_{g}F)_{i}=p_{i}\left(  \frac{\partial F}{\partial p_{i}}-\sum_{j}%
p_{j}\frac{\partial F}{\partial p_{j}}\right)  .
\]

Proof. From the definition, $\nabla_{g}F$ is the orthogonal
projection on the subspace $T_{p}\Delta_{K-1}^{0}$ of
\[
(\nabla_{g}F)_{i}=p_{i}\frac{\partial F}{\partial p_{i}}%
\]

and therefore we must have $\sum_i(\nabla_{g}F)_{i}=0.$ This then
gives the result.

\chapter{Appendix III: Markov Processes}

\section{Operator semigroups}

See Ethier-Kurtz, \cite{EK-86} Chap.1.

Consider a strongly continuous semigroup $\{T_t\}$ with generator $G$ and domain $D(G)$.  A subset $D_0\subset D(G)$ is a {\em core} if the closure of $G|_{D_0}$ equals $G$.  If $D_0$ is dense and $T_t:D_0\to D_0$ for all $t$, then it is a core.

\beT{SGCT} (Kurtz semigroup convergence Theorem \cite{EK-86}, Chap. 1, Theorem 6.5)  Let $L,L_n$ be Banach spaces and $\pi_n:L\to L_n$ is a bounded linear mapping and $\sup_n \|\pi_n\|<\infty$. We say $f_n\in L_n \to f\in L $ if
$\lim_{n\to \infty}\| f_n-\pi_n f \| =0$.

For $n\in\N$ let $T_n$ be a contraction on a Banach space $L_n$, let $\ve_n >0$, $\lim_{n\to\infty} \ve_n =0$. Let $\{T(t)\}$ be a strongly continuous contraction semigroup on $L$ with generator $A$ and let $D$ be a core for $A$. Then the following are equivalent:

(a) For each $f\in L$, $T^{\lfloor t/\ve_n\rfloor}_n \pi_nf \to T(t)f, \text{ for all }t\geq 0,$ uniformly on bounded intervals.

(b)  For each $f\in D$ there exists$f_n\in L_n$ such that $f_n\to F$ and $A_nf_n\to Af$.

\end{theorem}

\beT{SGCT2} \cite{EK-86} Chap. 4, Theorem 2.5.\\
Let $E$ be locally  compact and separable. For $n=1,2,\dots$ let $\{T_n(t)\}$ be a Feller semigroup on $C_0(E)$ and suppose that $X_n$ is a Markov process with semigroup $\{T_n(t)\}$ and sample paths in $D_E([0,\infty))$. Suppose that $\{T(t)\}$ is a Feller semigroup on $C_0(E)$ such that for each $f\in C_0(E)$
\be{} \lim_{n\to\infty} T_n(t)f=T(t)f,\;t\geq 0.\ee
If $\{X_n(0)\}$  has limiting distribution $\nu\in\mathcal{P}(E)$, then there is a Markov process $X$ correspondng to $\{T(t)\}$ with initial distribution $\nu$ and sample paths in $D_E([0,\infty))$ with initial distribution $\nu$ and sample paths in $D_E([0,\infty))$ and $X_n\Rightarrow X$.

\end{theorem}
\section{Some basic result for one dimensional diffusions}

Basic References: It\^o-McKean \cite{IM-65}, Karlin and Taylor \cite{KT-81}, Revuz and Yor \cite{RY-91}.

\subsection{Boundary behaviour classification}
Consider a diffusion process on $[0,L]$ or $[0,\infty)$ with drift and diffusion coefficients $b(x) $ and $\sigma^2(x)$.  The scale function $S(x)$ is defined by
\be{} S(x)=\int_{x_0}^x s(u)du,\qquad  s(x)=\exp\left(-\int_{x_0}^x\frac{2b(u)}{\sigma^2(u)}du\right),\ee
and the speed measure
\be{}  M([x_1,x_2])=\int_{x_1}^{x_2}m(x)dx,\qquad  m(x)=\frac{1}{\sigma^2(x)s(x)}.\ee

Feller introduced a classification of boundary points  as follows:  Applied to the boundary point $0$ this becomes
(following  Karlin and Taylor \cite{KT-81}):

\be{}  \Sigma(0)=\int_0^x S(0,u]m(u)du,\qquad  N(0)=\int_0^x (S(x)-S(u))m(u)du.\ee

Then
\begin{itemize}
\item $0$ is an {\em entrance } boundary if  $S(0,x])=\infty$ and $N(0)<\infty$.
\item $0$ is an {\em exit} boundary if $\Sigma(0)<\infty$ but  $M(0,x]=\infty$
\item $0$ is a regular boundary if $S(0,x]<\infty$ and $M(0,x]<\infty$.

\end{itemize}

For the Feller  CSBP process, $0$  is an exit boundary.  For the Feller CSBP process with immigration (\ref{CSBI})
\be{} s(x)=\exp(-\int_{x_0}^x\frac{2\beta}{4x})= (\frac{x}{x_0})^{-\beta/2},\qquad  m(x)=\frac{1}{4x (\frac{x}{x_0})^{-\beta/2}}=\frac{4}{x_0}x^{1-\frac{\beta}{2}}\ee

\subsection{Excursions}

\subsubsection{The Brownian excursion}

Let  $\{B(t)\}_{t\geq 0}$ be a standard Brownian motion and
\be{}  \tau_1:= \sup \{t\in [0,1]:B(t)=0\},\quad \tau_2:= \inf \{t\geq 1:B(t)=0\}.\ee
Then the Brownian excursion is a nonhomogeneous Markov process defined as follows:
\be{} B^e(t):= \frac{1}{\sqrt{(\tau_2 -\tau_1)}}B(\tau_1 +t(\tau_2 -\tau_1)),\quad 0\leq t\leq 1.\ee
It can be shown (see It\^o-McKean) \cite{IM-65} that the marginal PDF is given by
\be{} f(t,x)= \frac{2x^2}{\sqrt{2\pi t^3(1-t)^3}}e^{-\frac{x^2}{2t(1-t)}}.\ee

 It\^o's excursion measure $n(de)$ is the $\sigma$-finite measure on $C(\mathbb{R}_+,\mathbb{R}_+)$ obtained by
\be{} n(de)=\lim_{\ve\to 0}\frac{1}{2\ve} P_e(de)\ee
where $P_\ve(de)$ is the distribution of Brownian motion started at $\ve$ and stopped at the first time it hits $\zeta(e)=\inf\{s>0:e(s)=0\}$.  Then the {\em normalized Brownian excursion} is given by  $n(de|\zeta(e)=1)$.

\bigskip

\subsubsection{Excursions of non-negative diffusions}
Now consider a general diffusion on $[0,\infty)$ which is regular on $(0,\infty)$ and with $0$ as an absorbing boundary and with laws $\{P_x:0\leq x<\infty]\}$.  Let $T_y$ be the hitting time of $y$ and that $0$ is an exit point for the diffusion.

Revuz and Yor introduced the {\em excursion law of the diffusion} in terms of a $\sigma$-finite measure $\Lambda$ on $C$.  Under
$\Lambda$ the trajectories come in from zero according to an entrance law, then move according to the diffusion.

Let $s(x)$ be a scale function for the diffusion. Since we assume the absorbing point $0$ can be reached with positive probability from $x>0$
we can take
\be{}  s(0)=0,\; s(x)>0\quad\text{for  } x>0\ee

and $s$ is defined uniquely up to a constant factor by

\be{}  P_x(T_y<\infty)=\frac{s(x)}{s(y)},\quad 0<x<y<\infty.\ee

Then there exists  a $\sigma$-finite excursion law ${\mathbb{Q}}$ on
\bea{}  &&W_0:=\{w\in C([0,\infty),\mathbb{R}^+),\;w(0)=0,\;
w(t)>0\text{ for  } 0<t<\zeta\\&& \qquad\qquad\qquad\text{ for some }
\zeta\in(0,\infty)\}\nonumber\eea
obtained as follows.
Denoting by $P^\ve$ the law of the process started
with $w(0)=\ve$ and $\ve> 0, \Q$ is given by:
\be{}
 \mathbb{Q}(\cdot) =\lim_{\ve\to 0}\frac{P^\ve(\cdot)}{S(\ve)},
 \ee
where $S(\cdot)$ is the scale function of the diffusion defined by the relation,
\be{}
P_\ve(T_\eta<\infty)=\frac{S(\ve)}{S(\eta)},\qquad
0<\ve<\eta<\infty.
\ee







\end{document}